\documentclass[11pt,leqno]{amsart}
\usepackage{blindtext}
\usepackage{fancyhdr}
\usepackage{ifpdf}
\usepackage{curves}
\usepackage{tikz}
\usetikzlibrary{arrows}
\input{arrowsnew}
\usetikzlibrary{decorations.markings}
\usepackage{subcaption}
\usepackage{placeins}
\usepackage{tabularx}
\usepackage{url}

\usepackage{amsfonts}
\usepackage{amsthm}
\usepackage{amsmath}
\usepackage{amscd}
\usepackage{amssymb}
\usepackage{graphicx}
\usepackage{rotating}
\usepackage[small,nohug,heads=LaTeX]{diagrams}
\diagramstyle[labelstyle=\scriptstyle]
\usepackage{ulem}
\usepackage{textcomp}
\usepackage{appendix}
\usepackage{chngcntr}
\usepackage{etoolbox}
\usepackage{lipsum}
\usepackage{tikz}
\usepackage{syntonly}
\usepackage[pagebackref]{hyperref}

\usepackage[utf8]{inputenc}
\usepackage[english, russian]{babel}

\usepackage{leftidx}
\usepackage{slashbox}

\usetikzlibrary{matrix,arrows,decorations.pathmorphing}

\AtBeginEnvironment{subappendices}{%
\chapter*{Appendix}
\addcontentsline{toc}{chapter}{Appendices}
\counterwithin{figure}{section}
\counterwithin{table}{section}
}

\def\AA {{\mathbb A}}     
\def\CC {{\mathbb C}}     
\def\KK {{\mathbb K}}     
\def\NN {{\mathbb N}}     
\def\PP {{\mathbb P}}     
\def\RR {{\mathbb R}}     
\def\ZZ {{\mathbb Z}}     

\def\ring#1{\ifmmode \mathaccent'027 #1\else \rm\accent'027 #1\fi}

\def\ol  {\overline}
\def\ul  {\underline}

\def\mc {\mathcal}
\def\mk {\mathfrak}

\def\Hom {\mathrm{Hom}}

\def \Aut {\mathrm{Aut}}

\def \bd {\begin{diagram}}
\def \ed {\end{diagram}}
\def\be  {\begin{eqnarray}}
\def\ee  {\end{eqnarray}}
\def\ben {\begin{eqnarray*}}
\def\een {\end{eqnarray*}}

\def\bpr {\begin{proof}[Proof]}
\def\epr {\end{proof}}
\def\bsp {\begin{split}}
\def\esp {\end{split}}
\def\bprr {\begin{proof}[solution]}
\def\bpru {\begin{proof}[hint]}
\def\bpro {\begin{proof}[answer]}
\def\bcd {\begin{CD}}
\def\ecd {\end{CD}}

\newcommand{\abs}[1]{\left\vert#1\right\vert}
\newcommand{\scal}[1]{\left\langle#1\right\rangle}

\newtheorem{theorem}{Theorem}[section]
\newtheorem{lemma}[theorem]{Lemma}
\newtheorem{prop}[theorem]{Proposition}
\newtheorem{coro}[theorem]{Corollary}
\newtheorem{remark}[theorem]{Remark}
\newtheorem{df}[theorem]{Definition}

\newtheorem{ex}[theorem]{Example}

\newcommand\mapsfrom{\mathrel{\reflectbox{\ensuremath{\mapsto}}}}

\addcontentsline{toc}{chapter}{Appendices}
\begin{document}

\renewcommand{\contentsname}{Contents}
\renewcommand{\refname}{References}

\title[Non-commutative counting invariants and curve complexes]%
{Non-commutative counting invariants and curve complexes}
 
\author{George Dimitrov}
\address[Dimitrov]{Universit\"at Wien\\
	Oskar-Morgenstern-Platz 1, 1090 Wien\\
	\"Osterreich \\
}
\email{george.dimitrov@univie.ac.at}

\author{Ludmil Katzarkov}
\address[Katzarkov]{University of Miami , Coral Gables, USA  \\
 National Research University Higher School of Economics, Russian Federation, University of Miami  }
\email{lkatzarkov@gmail.com}

 \dedicatory{Dedicated to the memory of our friends Vasil Tsanov and Vladimir  Voevodskij}

\renewcommand{\abstractname}{Abstract}

\renewcommand{\figurename}{Figure}

\begin{abstract} 

In our previous paper, viewing $D^b(K(l))$ as a non-commutative curve, where $K(l)$ is the Kronecker quiver with $l$-arrows, we introduced categorical invariants via counting of non-commutative curves. Roughly, these invariants are  sets of subcategories   in a given category and their quotients. The non-commutative curve-counting invariants are obtained by restricting the subcategories to be equivalent to $D^b(K(l))$.  The general definition however defines  much larger class of invariants and many of them behave properly  with respect to fully faithful functors.  Here, after recalling the definition, we focus on examples and extend our studies beyond counting. 
We enrich our invariants  with  structures: the  inclusion of subcategories makes them partially ordered sets, and on the other hand  considering semi-orthogonal pairs of subcategories as edges amounts to directed graphs.

In addition to computing the non-commutative curve-counting invariants in $D^b(Q)$ for two affine quivers, for $ A_n$ and for  $D_4$ we derive  formulas for counting of the subcategories of type $D^b(A_k)$ in $D^b(A_n)$,  whereas  for the two affine quivers and for $D_4$ we determine and count all generated by an exceptional collection subcategories in $D^b(Q)$. It turns out that the problem for counting $D^b(A_k)$ in $D^b(A_n)$ has a geometric combinatorial parallel -  counting of maps between polygons.

Estimating the numbers counting non-commutative curves in $D^b({\mathbb P}^2)$ modulo the group of auto-equivalences  we prove finiteness and that the exact determining of these numbers  leads to solution of Markov problem. Via homological mirror symmetry this gives a new approach to this problem.

Regarding  the  structure of a partially ordered  set mentioned above we initiate intersection theory of non-commutative curves focusing on the  case of non-commutative genus zero. 
 
In the last section  we build an analogue of the classical curve complex, introduced by Harvey and Harrer and used in a spectacular way in Teichmueller theory and Thurston theory by Farb, Minsky, Mazur and many others.   The above-mentioned structure of a directed graph (and related simplicial complex)   is a categorical analogue of the curve complex.
The paper contains pictures of  the graphs in many examples  and presents also an approach to Markov conjecture via counting of subgraphs in a graph associated with $D^b(\PP^2)$. 

Some of the results proved here were announced in the previous work.
 \end{abstract}

\maketitle
\setcounter{tocdepth}{1}
\tableofcontents

\section{Introduction}

 In \cite{DK4} we introduced some  categorical invariants via counting of fully faithful exact  functors. In this paper we focus on examples and  extend  our studies   beyond counting. Here  these invariants produce many numbers and on the other hand they evolve into  categorical versions of  classical geometric structures.  
 
  The examples establish  connections with number theory and classical geometry.   
 In one of the  examples  we show that the  computation  of the  counting invariants is equivalent to computing the  number of equivalence classes  (up to rotation) of $k+1$ subgons in a regular $n+1$-gon (see subsection \ref{Polygons}), for example there are 4 equivalence classes of $3$-subgons in a regular 6-gon and they are: 
     	
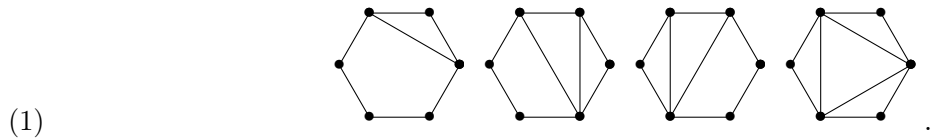
\begin{center}
   		\begin{equation} \label{4 eq classes}
   		\begin{tikzpicture}[scale=1]
   		
   		\begin{scope}[shift={(-4,0)}]
   		\newdimen\R
   		\R=0.8cm
   		\coordinate (center) at (0,0);
   	
   		\draw (0:\R) node[scale=0.8] {{$\bullet$}} --(120:\R) node[scale=0.8]{{$\bullet$}};
   		\draw (0:\R) node[scale=0.8] {{$\bullet$}} \foreach \x in {0,60,...,360} { -- (\x:\R)node [scale=0.8] {{$\bullet$}} };
   		\end{scope} 
   		\begin{scope}[shift={(-2,0)}]
   		\newdimen\R
   		\R=0.8cm
   		\coordinate (center) at (0,0);
   	
   		\draw (60:\R) node[scale=0.8] {{$\bullet$}}--(-60:\R) node[scale=0.8] {{$\bullet$}} --(120:\R) node[scale=0.8]{{$\bullet$}};

   		\draw (0:\R) node[scale=0.8] {{$\bullet$}} \foreach \x in {0,60,...,360} { -- (\x:\R)node [scale=0.8] {{$\bullet$}} };
   		\end{scope}
   		\begin{scope}[shift={(0,0)}]
   		\newdimen\R
   		\R=0.8cm
   		\coordinate (center) at (0,0);
   	
   		\draw (60:\R) node[scale=0.8] {{$\bullet$}}--(-120:\R) node[scale=0.8] {{$\bullet$}} --(120:\R) node[scale=0.8]{{$\bullet$}};

   		\draw (0:\R) node[scale=0.8] {{$\bullet$}} \foreach \x in {0,60,...,360} { -- (\x:\R)node [scale=0.8] {{$\bullet$}} };
   		\end{scope}
   		\begin{scope}[shift={(2,0)}]
   		\newdimen\R
   		\R=0.8cm
   		\coordinate (center) at (0,0);
   	
   		\draw (0:\R) node[scale=0.8] {{$\bullet$}}--(-120:\R) node[scale=0.8] {{$\bullet$}} --(120:\R) node[scale=0.8]{{$\bullet$}}--cycle;

   		\draw (0:\R) node[scale=0.8] {{$\bullet$}} \foreach \x in {0,60,...,360} { -- (\x:\R)node [scale=0.8] {{$\bullet$}} };
   		\end{scope}
   		\end{tikzpicture}.
   		\end{equation}	
   	\end{center}

 In Section \ref{PP2 section} we show a new approach to an old conjecture in number theory - the so called Markov conjecture.  
 
We  equip   our invariants   with a structure of a partially ordered set (Section \ref{Interesection theory intro})  and of a directed graph (Section \ref{section for digraph}).  In particular, Section \ref{Interesection theory intro} concerns  intersection of the objects which   we count in the previous sections.  The ``geometry'' of intersections   in  $D^b(D_4)$ surprisingly relates to a theorem for concurrent  lines\footnote{by concurrent lines we mean three or more lines in a plane which intersect at a single point.}  in the Euclidean plane.    The  structure of a directed graph (and related simplicial complex)   is a categorical analogue of the curve complex used in Teichmueller  and Thurston theory.

 The  idea of  non-commutative counting, introduced in  \cite{DK4}, is unfolded here  in Section \ref{section main Definition}:   for a triangulated  category $\mc T$, a subgroup $\Gamma \subset {\rm Aut}(\mc T)$ of  exact equivalences, and a choice of some restrictions $P$ on fully faithful functors   we  define  the set of subcategories of $\mc T$, which are equivalent to another chosen triangulated  category $\mc A$ via an equivalence  which satisfies   $P$, and modulo $\Gamma$.   The defined in Section \ref{section main Definition}  set of equivalence classes  of  subcategories in $\mc T$    we denote  by $C_{\mc A, P}^{\Gamma}(\mc T)$. 
 
We   discuss    $C_{\mc A, P}^{\Gamma}(\mc T)$ for concrete choices of $\mc A, \Gamma,\mc T, P$.

 When the categories $\mc T$, $\mc A$ are $\KK$-linear, choosing the property $P$ imposed on the functors as follows:\textit{ to be $\KK$-linear}, and restricting  $\Gamma$ to be subgroup of the group of $\KK$-linear auto-equivalences  results in  the set   $C_{\mc A, P}^{\Gamma}(\mc T)$, which we denote by  $C_{\mc A, \KK}^{\Gamma}(\mc T)$ (see Definition \ref{C_A KK(T)} and Lemmas \ref{bijection}, \ref{group action}).

 We prefer to choose first some $\mc A$, which is non-trivial but  well studied.  
 Following Kontsevich-Rosenberg  \cite{KR} we denote $D^b(Rep_\KK(K(l+1)))$ by  $N\PP^l$  for $l\geq -1$.
 For $l\geq 0$ we refer to $ N\PP^l$ as to a non-commutative curve and to  $\dim_{nc}( N\PP^l)=l$ as to its \textit{non-commutative  genus} (see \cite[Remark 12.2]{DK4}  for further motivation and \cite[Section 12.1]{DK4} for the definition of $\dim_{nc}$).  
 For a $\KK$-linear category $\mc T$ and any $l\geq -1$ we denote $C_{N\PP^l, \KK}^\Gamma(\mc T)$ by  $C_{l}^{\Gamma}(\mc T)$ and refer to this set as the set of non-commutative curves, sometimes nc curves for short, of non-commutative genus $l$ in $\mc T$ modulo  $\Gamma$. Note that   $N\PP^{-1}$ has homological dimension $0$ and one  should think of  the ``curves'' of non-commutative  genus $-1$ as a kind of  degenerate curves.  Sometimes, when it is clear from the context that non-commutative genus is meant, we omit writing ``non-commutative'' before genus.

In  Section \ref{Interesection theory intro} we start by  extending the definition of the set  $C_{\mc A, P}^{\{\rm Id\}}(\mc T)$  to define a  set $C_{\mc J, P}(\mc T)$  which consists of  subcategories of $\mc T$ which are allowed to be equivalent not only to a single category $\mc A$ but to any category $\mc A \in \mc J$, where $\mc J$ is a family of pairwise non-equivalent categories. The set $C_{\mc J, P}(\mc T)$ carries a structure of a partially ordered set (see Definition \ref{def of C_J(T)}). 
We refer to the elements in $C_{D^b(A_1),\KK}^{\{\rm Id\}}(\mc T)$ as to derived points (these are subcategories in $\mc T$ generated by a single exceptional object as explained in  \eqref{bijection for eo} and  Remark \ref{no enhancement needed for points}).  Any genus $0$ non-commutative curve contains exactly three derived points and any two of them span the curve (Lemma \ref{derived points in a genus 0}). If the  genus of one of two different nc curves is  less than or equal to $0$, then they intersect (see Definition \ref{interesection} for intersection) either trivially or in  a single derived  point  in $\mc T$ (Lemma  \ref{intersection}).
 Therefore in the cases   $\mc J = \{\mbox{trivial category}, D^b(A_1), D^b(A_2)\}$ and  $\mc J = \{\mbox{trivial category}, D^b(A_1), D^b(A_1) \oplus  D^b(A_1)  , D^b(A_2)\}$   any  two  $x, y \in C_{\mc J, \KK}(\mc T)$ have a  greatest lower bound in $ C_{\mc J, \KK}(\mc T)$  (Corollary \ref{greatest lower bound}). 

In   Section \ref{section for digraph} we enrich the partially ordered set $C_{\mc J, P}(\mc T)$ defined in Section \ref{Interesection theory intro} with a structure of a directed graph, denoted by  $G_{\mc J, P}(\mc T)$ (what we mean by directed graph is explained in the beginning of Section \ref{section for digraph}). The vertices of  $G_{\mc J, P}(\mc T)$  are  $C_{\mc J, P}(\mc T)$ whereas the edges are pairs $(\mc B_1, \mc B_2)$ of subcategories  such that $\Hom(\mc B_2, \mc B_1)=0$. We consider also a simplicial complex  $SC_{\mc J, P}(\mc T)$ whose simplexes are finite subsets in $C_{\mc J, P}(\mc T)$  which can be ordered in  a semi-orthogonal sequences (the information for this simplicial complex is contained in the digraph).

 A good feature of $C_{\mc J, \KK}(\mc T)$  proved here is  (see Proposition \ref{functoriality} and Corollaries  \ref{fully faithful functor on digraphs}, \ref{action on digraphs}), and Subsection \ref{the case J eq DbA1} (b):
 
  \begin{theorem} Any fully faithful  exact  $\KK$-linear  functor  $\mc T_1 \rightarrow \mc T_2$ induces an embedding of  $C_{\mc J, \KK}(\mc T_1)$ into   $C_{\mc J, \KK}(\mc T_2)$ preserving the additional structures introduced in Sections \ref{Interesection theory intro}, \ref{section for digraph}. \end{theorem}
  
  It follows that, if   $\mc T = \langle \mc T_1, \mc T_2 \rangle$ is a semi-orthogonal decomposition, then $C_{\mc J, \KK}(\mc T_i)$, $i=1,2$ are embedded in $C_{\mc J, \KK}(\mc T)$ and they  intersect trivially in $C_{\mc J, \KK}(\mc T)$.

  The  graphs  $G_{D^b(A_1), \KK}\left (N\PP^{l}\right )$ are shown in Figure \ref{graphs for NP} (see subsection \ref{some digraphs}). Further examples  in this paper are as follows.
    \begin{figure} \centering
  		\hspace{-30mm}  \begin{subfigure}[b]{0.2\textwidth}
  			\begin{tikzpicture}  [scale=0.5]
  			\draw[<->]  (-1,0)--(+1,0);
  			\node [left] at (-1 ,0) {$\langle E_1 \rangle $};
  			\node [right] at (+1 ,0) {$\langle E_2 \rangle $};
  			\end{tikzpicture}
  			\caption{$l=-1$}\label{Figure for digraph of NP-1}
  	\end{subfigure}
  	\hspace{5mm}
  	\begin{subfigure}[b]{0.2\textwidth}
  		   		 	\begin{tikzpicture} [scale=0.5] 
  		 	\draw[->]  (-1,0)--(+1,0);
  		 	\node [left] at (-1 ,0) {$\langle s_{0,0}\rangle $};
  		 	\node [right] at (+1 ,0) {$\langle s_{1,1}\rangle $};
  		 	\draw[->]   (1.5,0.6)-- (0.7,1.6);
  		 	\node at (0 ,2) {$\langle s_{0,1}\rangle $};
  		 	\draw[->]  (-0.5,1.6) --  (-1.5,0.5);
  		 	\node [above] at (0 ,0) {$1 $}; 	\node [above] at (0 ,0) {$1 $};
  		 	\node [above] at (1.2 ,1) {$1 $}; 	\node [above] at (-1.2 ,1) {$1 $};
  		 	\end{tikzpicture}
  		 	\caption{$l=0$}\label{Figure for digraph of A2}
  	\end{subfigure} 	\hspace{5mm}
  		\begin{subfigure}[b]{0.2\textwidth} 	\begin{tikzpicture}[scale =0.6]

  			\draw[->]  (-3.8,0)--(-2.3,0);
  			\node [above] at (-3,0) {$l+1 $};	
  			\node [left] at (-3.8,0) {$\cdots$};
  			\draw[->]  (-1,0)--(1,0);
  			\node [left] at (-1 ,0) {$\langle s_i \rangle $};
  			\node [right] at (+1 ,0) {$\langle s_{i+1} \rangle $};
  			\draw[->]  (3,0)--(4.2,0);
  			\node [above] at (3.6,0) {$l+1 $};	
  			\node [right] at (4.2,0) {$\cdots $};	
  			\node [above] at (0,0) {$l+1 $};	
  			\end{tikzpicture}
  	
    		\caption{ $l\geq 1$}\label{Figure for digraph of NPl}	
    		\end{subfigure} 
  	  	\caption{$G_{D^b(A_1), \KK}\left (N\PP^{l}\right )$}\label{graphs for NP}
  \end{figure}
 Let  $\KK$ be algebraically closed field. Till the end of introduction the categories will be $\KK$-linear and when writing  $C_{\mc J}^{\Gamma}(\mc T)$ we mean  $C_{\mc J,\KK}^{\Gamma}(\mc T)$.   Let $\mc T_i=D^b(Rep_\KK(Q_i))$, $i=1,2$,   where:	\be \label{Q1} Q_1=  \begin{diagram}[1em]
 	&       &  2  &       &    \\
 	& \ruTo &    & \rdTo &       \\
 	1  & \rTo  &    &       &  3
 \end{diagram}   \ \qquad \qquad \ Q_2= \begin{diagram}[1em]
 2 &  \rTo  &  3    \\
 \uTo &        & \uTo     \\
 1   & \rTo  &    4\end{diagram}.  \ee 
In  Section \ref{non-trivial examples}    we determine $C_{\mc A, \KK}^{\{\rm Id\}}(\mc T_i)$  for any $\mc A$ which is generated by an exceptional collection and  describe the directed graphs  $G_{D^b(A_1)}(\mc T_i)$ for $i=1,2$. 
 In the case when $\mc A$ is generated by an exceptional pair any non-trivial  $C_{\mc A, \KK}^{\{\rm Id \}}(\mc T_i)$  must be of the form $C_{l}^{\{\rm Id \}}(\mc T_i)$ for  $-1 \leq l\leq 1$ (Remarks \ref{reminder},  \ref{subcats generated by pairs}) and the cardinalities are as follows (here  $S$ denotes the Serre Functor): 
  \begin{gather} \label{table with numbers Q1}
  \abs{C_{l}^{\Gamma}(\mc T_1) }= \begin{tabular}{|c|c|c|c|}\hline
  \backslashbox{$l$}{$\Gamma$}& $\{\rm Id\}$ & $ \langle S \rangle$ & $ {\rm Aut}_{\KK}(\mc T_1)$\\ \hline
  $-1$ & $0$ &   $0$ &  $0$\\ \hline
  $0$ & $\infty$ &   $3$ &  $1$\\ \hline
  $+1$ & $2$ &   $1$ &  $1$\\ \hline
  \end{tabular}
   \abs{C_{l}^{\Gamma}(\mc T_2) }= \begin{tabular}{|c|c|c|c|}\hline
  \backslashbox{$l$}{$\Gamma$}& $\{\rm Id\}$ & $ \langle S \rangle$ & $ {\rm Aut}_{\KK}(\mc T_2)$\\ \hline
  $-1$ & $\infty$ &   $4$ &  $2$\\ \hline
  $0$ & $\infty$ &   $8$ &  $1$\\ \hline
  $+1$ & $4$ &   $2$ &  $1$\\ \hline
  \end{tabular}
  \end{gather}
    \begin{gather}
  \label{vanishings}
  l\geq 2 \ \Rightarrow \abs{C_{l}^{\{\rm Id\}}(\mc T_1) }=\abs{C_{l}^{\{\rm Id\}}(\mc T_2) }=0. \end{gather} 
 
  Let $\mc A$ be a category generated be an exceptional triple.  Since each exceptional triple in $\mc T_1$ is full (proved in \cite{WCB1}),  $C_{\mc A, \KK}^{\{\rm Id\}}(\mc T_1)$ is non-trivial only if $\mc A\cong \mc T_1$.  We show that   non-trivial $C_{\mc A, \KK}^{\{\rm Id\}}(\mc T_2)$ is obtained  only for $\mc A \cong D^b(A_3)$ and $\mc A \cong  D^b(Q_1)$ (see Corollary \ref{triples in Q2}  ) and the cardinalities are:   
  	\begin{gather} \label{triples in Q_2} 
  	\abs{C_{\mc A, \KK}^{\Gamma}\left (D^b(Q_2) \right ) }=   \begin{tabular}{|c|c|c|c|}\hline
  	\backslashbox{$\mc A$}{$\Gamma$}& $\{\rm Id\}$ & $ \langle S \rangle$  & $ {\rm Aut}_\KK(\mc T_1)$\\ \hline
  	$D^b(A_3)$ & $\infty$ &   $4$ & $1$\\ \hline
  	$D^b(Q_1)$ & $4$ &   $2$ &   $1$\\ \hline
  	\end{tabular} .
  	\end{gather}

   The numbers in the column for $\{\rm Id\}$ in tables \eqref{table with numbers Q1},  \eqref{triples in Q_2}  imply that   	$\abs{C_{D^b(Q), \KK}^{\Gamma}\left (D^b(Q_i) \right ) }<\infty$, for any affine acyclic  quiver $Q$ (i.e. with underlying graph an extended Dynkin diagram)  and for $i=1,2$. In     another note, \cite{DK5},   using results by Geigle, Lenzning, Meltzer, Hübner for weighted projective lines we show that for any two affine acyclic quivers $Q$, $Q'$   there are only finitely many triangulated subctegories in $D^b(Rep_\KK(Q))$, which are equivalent to  $D^b(Rep_\KK(Q'))$,  i.e. $\abs{C_{D^b(Q'), \KK}^{\{{\rm Id}\}}\left (D^b(Q) \right ) }<\infty$, furthermore $\abs{C_{D^b(Q'), \KK}^{{\rm Aut}\left (D^b(Q) \right )}\left (D^b(Q) \right ) }=1$. 
    
   In Section \ref{PP2 section} we prove that if we  denote $\mc T = D^b(\PP^2)$ (here we fix $\KK=\CC$) and  $\langle S \rangle \subset  {\rm Aut}(\mc T)$  the subgroup generated by the Serre functor, then $C_{-1}^{\{\rm Id \}}(\mc T) =\emptyset$,   $\forall l\geq 0$ the set $C_{l}^{\langle S \rangle}(\mc T)$ is finite, and: 
    	\begin{gather} \nonumber \left \{l\geq 0:C_{l}^{\langle S \rangle}(\mc T) \neq \emptyset \right \}= \left \{l\geq 0:C_{l}^{\{\rm Id \}}(\mc T) \neq \emptyset \right \}= \{3 m - 1: m \ \mbox{is a Markov number} \}.
    	\end{gather}
    	Furthermore, for any  Markov number\footnote{Recall that a Markov number $x$ is a number $x \in \NN_{\geq 1}$ such that there exist integers $y,z$ with $x^2+y^2+z^2=3 x y z$.}  $m$ we have \begin{gather}  \label{infinite mumbers for PP2} \abs{C_{3 m -1}^{\{\rm Id\}}(\mc T) }=\infty  \\ 
    	\label{finite mumbers for PP2 full group} 1\leq \abs{C_{3 m -1}^{\Aut_{\CC}(\mc T)}(\mc T) }= \abs{\left  \{ y:\begin{array}{l} 0\leq y<m, y \in \ZZ \ \mbox{and there exists an}\\   \mbox{ exceptional vector bundle} \  E \ \mbox{on} \ \PP^2, \\\mbox{with} \  r(E)= m, \  y=c_1(E)\end{array} \right \} } \leq  m,
    	\\ 	 	 	\label{finite mumbers for PP2} 3\leq \abs{C_{3 m -1}^{\langle S \rangle}(\mc T) }=3 \abs{C_{3 m -1}^{\Aut_{\CC}(\mc T)}(\mc T)} \leq 3 m,
    	\end{gather} 
    	where   $c_1(E)$, $r(E)$ are   the first Chern class (which we consider as an integer)  and the rank of $E$.
   
    In Section \ref{PP2 section}  are  computed the first several non-trivial $\abs{C_{3m-1}^{\Aut_\CC(\mc T)}(\mc T)}$, they are (recall that $m$ is a Markov number and on the first row are listed the first 9 Markov numbers): 
    \begin{gather}\label{table for PP2}
    \begin{array}{| c | c | c | c | c | c | c | c | c | c |}
    m                                                & 1 & 2 & 5 & 13  & 29 & 34 & 89 & 169 & 194  \\ \hline
    \abs{C_{3 m -1}^{\Aut_{\CC}(\mc T)}(\mc T) } & 1 & 1 & 2 & 2   & 2  & 2  & 2  & 2   & 2
    \end{array}
    \end{gather}
    	We also prove that  the so called Tyurin conjecture, which is equivalent to the Markov conjecture (\cite[p. 100]{Rudakov1} or \cite[Section 7.2.3 ]{GorKul}), is equivalent to  the following statement:   for all  Markov numbers $m\neq1, m\neq 2$ we have $\abs{C_{3 m -1}^{\Aut_\CC(\mc T)}(\mc T)}=2$, and this is equivalent to the statement:  for all  Markov numbers $m\neq 1, m\neq 2$ we have $\abs{C_{3 m -1}^{\Aut_{\CC}(\mc T)}(\mc T)}\leq 2$.
    Via the latter statement  one could approach  Markov conjecture    using   homological mirror symmetry  and applying  A side techniques  for computing the non-commutative curve-counting invariants for $D^b(\PP^2)$.

     Under some conditions   we can embed  $C_{\mc A}^{\Gamma}(\mc T)$  in a certain set of equivalence classes of subgraphs in the directed graph  $G_{\mc J}(\mc T)$, where $\mc A$ is  fixed   and $\mc J$ is a family of pairwise non-equivalent categories  (Proposition  \ref{from invariants to subgraphs}) and this observation   gives an approach to Markov conjecture by counting subgraphs in $G_{D^b(A_1)}(D^b(\PP^2))$ up to action of ${\rm Aut}_{\CC}\left (D^b(\PP^2) \right )$  as explained in Proposition \ref{markov with graphs}.

    When writing $D^b(A_n)$ or $D^b(D_4)$  we mean categories obtained by any choice of orientations (changing the orientation gives equivalent categories \cite[Section 5]{Miyachi}).   Due to the simple structure of ${\rm Aut}_\KK(D^b(A_n))$ (determined in \cite{Miyachi})  holds   $C_{\mc A, \KK}^{\langle S \rangle}\left (D^b(A_n)\right ) \cong  C^{{\rm Aut}_\KK(D^b(A_n))}_{\mc A,\KK}(D^b(A_n))$  for the formulas that follow, where $S$ is the Serre functor of $D^b(A_n)$.

     In Section \ref{ncc in A_n+1} we obtain formulas  for  $\abs{ C^{\{\rm Id\}}_{D^b(A_k)}(D^b(A_n)) }$ and   $\abs{ C^{{\rm Aut}_\KK(D^b(A_n))}_{D^b(A_k),\KK}(D^b(A_n)) }$ as follows:   for $n\geq1$, $k\geq 1$ the set $C^{\{\rm Id\}}_{D^b(A_k)}(D^b(A_n))$  is non-empty iff $k\leq n$ and    $\abs{C^{\{\rm Id\}}_{D^b(A_k)}(D^b(A_n))}=\binom{n+1}{k+1}$;
      	for $1\leq k < n$ let us denote $D=g.c.d(k+1,n+1)$, then:
    	\begin{gather}  \label{formula third form intro} \abs{C^{{\rm Aut}_{\KK}(D^b(A_n))}_{D^b(A_k),\KK}(D^b(A_n)) }= \sum_{\left \{(x,y) :\begin{array}{c} x,y,x/y \in \ZZ_{\geq 1}  \\ x, y \ \mbox{are  divisors of} \ D \end{array} \right \}} \frac{ D \mu(x/y)}{(k+1)x}   \left ( \begin{array}{c} y \frac{n+1}{D}-1 \\ y \frac{k+1}{D}-1 \end{array} \right), \end{gather}
    	where the summands are binomial coefficients and $\mu(x)$ is the M\"obius function.
    
   For short in the following formulas we write $c_{k,n}$ instead of $ \abs{C^{{\rm Aut}_{\KK}(D^b(A_n))}_{D^b(A_k),\KK}(D^b(A_n)) }$. Some special cases, in which \eqref{formula third form intro} has simple form are explained in Corollaries  \ref{n+2 prime coro}, \ref{n=k+delta}. For example:
   	
   	{ \rm (a) } if $k+1, n+1$ are coprime,	$c_{k,n}= \frac{1}{k+1}\binom{n}{k} =\frac{1}{n+1} \binom{n+1}{k+1};$
   	
   	{\rm (b)} if $g.c.d.(n+1,k+1)=p$ is prime, then 
   	$c_{k,n}=  \frac{1}{k+1}\binom{n}{k}+ \frac{p-1}{k+1} \binom{(n+1)/p -1}{(k+1)/p-1}; 
   	$
   	
   	{\rm (c)} in particular, if  $k+1$ is prime, then 
   	$ c_{k,n}=  \left \{ \begin{array}{c c}  \frac{k}{k+1} + \frac{1}{k+1}\binom{n}{k}&  \mbox{if} \ \frac{n+1}{k+1} \in \ZZ  \\ \frac{1}{k+1}\binom{n}{k}   & \mbox{if} \ \frac{n+1}{k+1} \not \in \ZZ  \end{array} \right. ;
   	$
   	
   		{\rm (d)} if $n=k+\Delta$, where $\Delta$ is $1$ or  prime, then 
   		$ c_{k,k+\Delta}=  \left \{ \begin{array}{c c}  \frac{1}{\Delta} \binom{k+\Delta}{\Delta-1}  &  \mbox{if $\frac{k+1}{\Delta} \not \in \ZZ$}   \\[3mm]  \frac{1}{\Delta} \binom{k+\Delta}{\Delta-1} + \frac{\Delta-1}{\Delta}  &\mbox{if $\frac{k+1}{\Delta}  \in \ZZ$}   \end{array} \right. .
   		$

  Corollary \ref{Polygons coro} shows that the numbers  $c_{k,n}$ have  a  geometric  combinatorial meaning explained before \eqref{4 eq classes}.   
    For example, using (b) above we compute $c_{2,5}=4$, which  by Corollary \ref{Polygons coro}  is equal the number of  equivalence classes (up to rotation) of $3$-subgons in a regular 6-gon  depicted in  \eqref{4 eq classes}.   
   	Another example  is demonstrated in Example \ref{another example}.  To compute the  number of equivalence classes  (up to rotation) of $k+1$ subgons in a regular $n+1$-gon one could use  combinatorial  technique developed by Redfield–Pólya \cite{Stanley}, however we  compute them without referring to it.

   	 In Subsection \ref{subsection ncc in A_n} are determined the numbers $\abs{C^{\Gamma}_{l}(D^b(A_n))}$ for any $l\geq -1$. The genus $0$ case follows from \eqref{formula third form intro}, see {\rm (c)} with $k=2$, whereas the genus $-1$ case is
   	 \begin{gather}     
   	 \label{coro for numC-1} 	\abs{C^{\{\rm Id\}}_{-1}(D^b(A_n))}=2 \binom{n+1}{4}  \quad 
    	 \abs{ C^{{\rm Aut}_\KK(D^b(A_n))}_{-1}(D^b(A_n)) } =\left \{ \begin{array}{c c}\frac{n(n-1)(n-2)}{12}  & \mbox{if} \ n \ \mbox{is even}\\   \frac{(n-1)((n-1)^2+2)}{12} & \mbox{otherwise} \end{array} \right. . \end{gather}
   	 In the end of Section \ref{ncc in A_n+1} we show that the graph $G_{D^b(A_1), \KK}\left (D^b(A_3)\right )$ is as in Figure \ref{Figure for digraph of A3}. 
   	 
   	 \begin{figure}
   	 	\begin{tikzpicture}[scale=1]

   	 	\draw[<-] (0,0)-- (2,0);
   	 	\node [left] at (0,0) {$\langle s_{0,2} \rangle $};
   	 	\node [right] at (2,0) {$\langle s_{2,2} \rangle $};
   	 	\draw[<->] (-0.2,-0.2)-- (1.1,-1.4);	
   	 	\node  at (1.6,-1.5) {$\langle s_{1,1} \rangle $};
   	 	\draw[->] (1.8,-1.3)-- (2.2,-0.2);
   	 	\node  at (1,0.5) {$\langle s_{1,2} \rangle $};	
   	 	\draw[->] (0.5,0.4)-- (0,0.2);
   	 	\draw[<-] (1.5,0.4)-- (2,0.2);
   	 	\draw[->] (1,0.3)-- (1.45,-1.3);
   	 	
   	 	\node  at (1,1.5) {$\langle s_{0,0} \rangle $};	
   	 	\draw[<-] (0.7,1.25)-- (-0.1,0.3);
   	 	\draw[<->] (1.4,1.25)-- (2.1,0.3);
   	 	\draw[->] (1.05,1.3)-- (1.55,-1.3);		
   	 	
   	 	\node  at (1,2.5) {$\langle s_{0,1} \rangle $};	
   	 	\draw[<-] (0.7,2.25)-- (-0.2,0.5);
   	 	\draw[->] (1.3,2.25)-- (2.2,0.5);
   	 	\draw[<-] (1.13,2.3)-- (1.65,-1.3);	
   	 	
   	 	\draw[->] (1,1.3)-- (1,0.7);
   	 	\draw[->] (1,2.35)-- (1,1.7);	
   	 	\end{tikzpicture} 				
   	 	\caption{$G_{D^b(A_1), \KK}\left (D^b(A_3)\right )$}\label{Figure for digraph of A3}
   	 	
   	 \end{figure}
   	In  Section \ref{the numbers for D_4}     we determine $C_{\mc A, \KK}^{\{\rm Id\}}\left (D^b(D_4)\right )$  for any $\mc A$ which is generated by an exceptional collection and describe the graph $G_{D^b(A_1),\KK}(D^b(D_4))$.  In the case when $\mc A\cong D^b(A_1)$ we have:  
  \begin{gather} \label{table with numbers A1 in D4}
    \begin{tabular}{|c|c|c|c|c|}\hline
   $\Gamma$ & $\{\rm Id\}$ & $\langle \kappa \rangle$ & $ \langle S \rangle$ & $ {\rm Aut}_{\KK}(\mc T)$\\ \hline
  $\abs{C_{D^b(A_1), \KK}^{\Gamma}(D^b(D_4)) }$ & $12$ &   $6$ &  $4$ &  $2$\\ \hline
     \end{tabular}.
   \end{gather}
     where $\kappa$ is the auto-equivalence of $\mc T$ due to the rotational  symmetry of the diagram  $D_4$. 
  In the case when $\mc A$ is generated by an exceptional pair, then any non-trivial  $C_{\mc A, \KK}^{\{\rm Id \}}(\mc T)$  must be of the form $C_{l}^{\{\rm Id \}}(\mc T)$ for  $-1 \leq l\leq 0$ (Remark \ref{subcats generated by pairs Dynkin}) and the cardinalities are as follows: 
   \begin{gather} \label{table with numbers D4}
   \abs{C_{l}^{\Gamma}(D^b(D_4)) }= \begin{tabular}{|c|c|c|c|c|}\hline
   \backslashbox{$l$}{$\Gamma$}& $\{\rm Id\}$ & $\langle \kappa \rangle$ & $ \langle S \rangle$ & $ {\rm Aut}_{\KK}(\mc T)$\\ \hline
   $-1$ & $9$ &   $3$ &  $3$ &  $1$\\ \hline
   $0$ & $15$ &   $5$ &  $5$ &  $3$\\ \hline
   \end{tabular}.
   \end{gather} 
    In the case when $\mc A$ is generated by an exceptional triple, then  non-trivial  $C_{\mc A, \KK}^{\{\rm Id \}}(\mc T)$  is obtained only for $\mc A \cong D^b(A_3)$ and  $\mc A \cong D^b(A_1)\oplus D^b(A_1)\oplus D^b(A_1)$ and the cardinalities are (Proposition \ref{triples in D4 prop}): 
   \begin{gather} \label{triples in D4} 
    \abs{C_{\mc A, \KK}^{\Gamma}\left (D^b(D_4) \right ) }=   \begin{tabular}{|c|c|c|c|c|}\hline
   \backslashbox{$\mc A$}{$\Gamma$}& $\{\rm Id\}$ & $\langle \kappa \rangle$  &  $ \langle S \rangle$ & $ {\rm Aut}_\KK(\mc T_1)$\\ \hline
   $D^b(A_3)$ & $9$ &   $3$ &   $3$ & $1$\\ \hline
   $D^b(A_1)\oplus D^b(A_1)\oplus D^b(A_1)$ & $3$ &   $3$ &  $1$&  $1$\\ \hline
   \end{tabular}. 
   \end{gather}
   
   In Section \ref{subsection for examples of intersections} we study  the intersection relations of genus zero nc curves in $D^b(A_3)$ and $D^b(D_4)$ and depict them in  Figures \ref{Figure for  intersection in A3} and \ref{Figure for D4} (here   the genus $0$ nc  curves are presented by lines and the derived points are presented by black dots).  The fact that  each of three triples of diagonals in Figure \ref{Figure for D4} is concurrent (shares  a single point)  is a theorem from Euclidean geometry which is in a  surprising coherence with  the   geometry of genus zero nc curves in   $D^b(D_4)$.  For more details we refer to Subsection \ref{subsection for examples of intersections}.

    In the last section \ref{last section} we draw the graphs  $G_{\mc J, \KK}(\mc T)$ for $\mc J= \{ N\PP^l \}_{l\in \ZZ_{\geq -1}}$ and  $\mc T=D^b(Q_2)$, $D^b(D_4)$.
  
  The invariants $G_{\mc J,\KK}(\mc T)$ and  $SC_{\mc J,\KK}(\mc T)$  are a   categorical analogue  of  curve complexes.  We would like to apply these analogues   to the following problems in future works:

  \begin{enumerate}
  	\item \underline{Markov Conjecture:}   $\abs{C_{3 m -1}^{\Aut_\CC\left (D^b(\PP^2)\right )}\left (D^b(\PP^2)\right ) }=2$, $m\neq 1, m \neq 2$
  	\item \underline{Orlov spectrum:} More precisely, we plan to  compute  Orlov spectrum in new examples. Studying the new invariants in some categories
  	deepens our knowledge for their structure, the gained information will be used as a step towards
  	computing the Orlov spectrum.
  	\item \underline{Iskovskikh Conjecture:} Let $X \neq 3 \dim$ cubic and let $X$ be a conic bundle. Then $X$ is not  rational if and only if $|2K_X+C|\neq \emptyset $. 
  	
   \end{enumerate}
   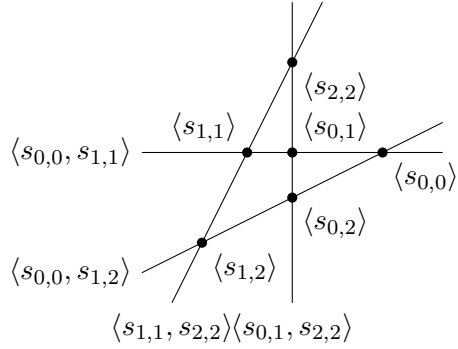
\begin{figure}
   	\begin{tikzpicture}[scale = 0.4] 
   	\draw  (-5,0)--(+5,0);
   	\node [left] at (-5 ,0) {$\langle s_{0,0}, s_{1,1}\rangle $};
   	\draw [fill] (0,0) circle [radius=0.15];
   	\node [above right] at (0 ,0) {$\langle s_{0,1}\rangle $};
   	
   	\draw  (0,-5)--(0,+5);
   	\node [below] at (0 ,-5) {$\langle s_{0,1}, s_{2,2}\rangle $};
   	
   	\draw  (-5,-4)--(+5,+1);
   	\node [left] at (-5 ,-4) {$\langle s_{0,0}, s_{1,2}\rangle $};	
   	\draw [fill] (0,-1.5) circle [radius=0.15];
   	\node [below right] at (0,-1.5) {$\langle s_{0,2}\rangle $};
   	\draw [fill] (3,0) circle [radius=0.15];
   	\node [below right] at (3 ,0) {$\langle s_{0,0}\rangle $};
   	
   	\draw  (-4,-5)--(+1,+5);
   	
   	\draw [fill] (-1.5,0) circle [radius=0.15];
   	\node [above left] at (-1.5,0) {$\langle s_{1,1}\rangle $};
   	\draw [fill] (0,3) circle [radius=0.15];
   	\node [below right] at (0 ,3) {$\langle s_{2,2}\rangle $};	
   	\draw [fill] (-3,-3) circle [radius=0.15];	
   	\node [below right] at (-3 ,-3) {$\langle s_{1,2}\rangle $};
   	\node [below ] at (-4 ,-5) {$\langle s_{1,1}, s_{2,2}\rangle $};
   	
   	\end{tikzpicture}
   	\caption{Intersection relations of the genus zero non-commutative curves in $D^b(A_3)$}
   	\label{Figure for  intersection in A3}
   \end{figure}

   \begin{figure}
   	\begin{tikzpicture} [scale = 0.9] 
   	\draw  (-6,3)--(+6,3);
   	\draw [fill] (0,3) circle [radius=0.07];
   	\node [above right] at (0 ,3) {$\langle s_{23}\rangle $};
   	\draw  (-6,4)--(+6,4);
   	\draw [fill] (0,4) circle [radius=0.07];
   	\node [above right] at (0 ,4) {$\langle s_{2o}\rangle $};		
   	\draw  (-6,-4)--(+6,-4);
   	\draw [fill] (0,-4) circle [radius=0.07];	
   	\node [above left] at (-0.1 ,-4) {$\langle s_{3}\rangle $};
   	
   	\draw  (0,-5)--(0,+5);
   	\node [below] at (0,-5) {$\langle s_{3}, s_{2o} \rangle $};
   	
   	\draw  (-4,-5)--(-4,+5.2);
   	\draw [fill] (-4,-4) circle [radius=0.07];
   	\node [right] at (-4,-5) {$\langle s_{1}, s_{3o} \rangle $};
   	\node [left] at (-6,-4) {$\langle s_{3}, s_{1o} \rangle $};
   	\node [above left] at (-4 ,-4) {$\langle s_{13}\rangle $};
   	\draw [fill] (-4,3) circle [radius=0.07];
   	\node [left] at (-6,3) {$\langle s_{2}, s_{3o} \rangle $};
   	\node [below left] at (-4 ,3) {$\langle s_{3o}\rangle $};
   	\draw [fill] (-4,4) circle [radius=0.07];
   	
   	\node [above right] at (-4 ,4) {$\langle s_{1}\rangle $};
   	\node [left] at (-6,4) {$\langle s_{1}, s_{2o} \rangle $};
   	
   	\draw  (4,-5)--(4,+5);
   	\node [left] at (4,-5) {$\langle s_{2}, s_{1o} \rangle $};
   	\draw [fill] (4,-4) circle [radius=0.07];
   	\node [above right] at (4 ,-4) {$\langle s_{1o}\rangle $};
   	\draw [fill] (4,3) circle [radius=0.07];
   	\node [below right] at (4 ,3) {$\langle s_{2}\rangle $};
   	\draw [fill] (4,4) circle [radius=0.07];
   	\node [above right] at (4 ,4) {$\langle s_{12}\rangle $};

   	\draw [red](-6,4.5)-- (-4,4)--(0,3)--(3.11,2.22)--(6,1.5);	
   	\node [left] at (-6,4.5) {$\langle s_{123}, s_{1} \rangle $};
     	\draw [fill] (3.11,2.22) circle [radius=0.07];	
   	\draw [red] (-0.5,-5)--(0,-4)--(4,4)--(4.5,5);	
   	\node [left] at (-0.5,-5) {$\langle s_{123}, s_3 \rangle $};	
   	\draw [red] (-5,-4.875)--(-4,-4)--(4,3)--(6,4.75);	
   	\node [left] at (-5,-4.875) {$\langle s_{123}, s_{2} \rangle $};		
   	\node [below right] at (2.8,2) {$\langle s_{123} \rangle $};

   	\draw [blue] (-6,2.75)--(-4,3)--(4,4)--(6,4.25);
   	\node [below] at (-6,2.75) {$\langle s_{3o}, \delta \rangle $};
   	\draw [blue](-1.3,5.275)-- (4,-4)--(0,3)--(-0.27,3.47)--(4.5,-4.875);	
   	\node [above] at (-1.3,5.275) {$\langle s_{1o}, \delta \rangle $};
   	\draw [fill] (-0.27,3.47)   circle [radius=0.07];	
   	\draw [blue] (-4.5,-5)-- (-4,-4)--(0,4)--(0.5,5);	
   	\node [below] at (-4.5,-5) {$\langle s_{2o}, \delta \rangle $};
   	\node [above left] at (-0.33,3.4) {$\langle \delta \rangle $};

   	\draw [green] (-5.7,5.7)--(-5.34,5.34)-- (-4,4)--(+4,-4)--(+5,-5);	
   	\node [right] at (5,-5) {$\langle s_{1}, s_{o} \rangle $};
   	\draw [green] (-5.7,5.975)--(-5.34,5.34)--(0,4)--(4,3)--(6,2.5);
   	\node [right] at (6,2.5) {$\langle s_{2}, s_{o} \rangle $};		
   	\draw [green](-5.8,5.45)--(-5.34,5.34)-- (0,-4)--(-4,3)--(0.5,-4.875);	
   	\node [right] at (0.5,-4.875) {$\langle s_3,s_o \rangle $};		
   	\draw [fill] (-5.34,5.34) circle [radius=0.07];
   	\node [above right] at (-5.34,5.34) {$\langle s_o \rangle $};	
   	\end{tikzpicture}
   	\caption{Intersection of the genus zero non-commutative curves in $D^b(D_4)$}\label{Figure for D4}
   \end{figure}

   By their virtues curve complexes behave well in families and they are retracts of Teichmüller spaces. It is reasonable to expect that the categorical analogues of curve complexes are retracts of moduli spaces of stability conditions, and
   then these categorical analogues  and stability conditions should behave well in families.

 \textit{{\bf Acknowledgements:}}
 
 The authors are  thankful to Christian Krattenthaler for bringing to their  attention  Redfield–Pólya theory.
 The authors are thankful to Alexander Efimov, Maxim Kontsevich, So Okada, and Tony Pantev for their encouraging remarks.
 The authors are thankful to Denis Auroux and Paul Seidel for useful conversations.

 \textit{{\bf FUNDING:}} 
 
 The first author  was  supported by  Austrian Science Fund [Grant numbers: P 29178-N35, P 27784]. 
 
 The second author   was  supported by National Science Foundation [Grant numbers: DMS 150908, DMS 1265230, DMS 1201475];
 Simons Foundation [Grant numbers: Investigator Award AWD-003136, Collaboration Award No. AWD-003093];  National Research University Higher School of Economics [Grant number: HMS and automorphic forms];
 Russian Foundation for Fundamental Investigations
 [Grant numbers: Government grant No 14.641.31.0001].

\section{Notations} \label{notations} 

We fix a universe and assume that the  set of objects and  the set of morphisms of any category we consider are elements of this universe. 

 The shift functor  in a triangulated category ${\mathcal T}$ is designated sometimes by $[1]$.
 
 A \textit{trivial triangulated category} is a triangulated category which has only zero objects (for example a category with a single object and a single endomorphism of this object).\footnote{Note that in a trivial triangulated category all triangles are distinguished}
 
 A \textit{triangulated subcategory}  in a triangulated category $\mc T$ is a non-empty full subcategory $\mc B $ in  $\mc T$, s. t. two conditions hold:
 
 {\rm (a) } $Ob(\mc B)[1]=Ob(\mc B)$ ; 
 
 {\rm (b)}  for any $X,Y \in Ob(\mc B) $ and  distinguished triangle $X\rightarrow Z \rightarrow Y \rightarrow X[1]$ in $\mc T$ follws $Z \in  Ob(\mc B)$.

  We write $\langle  S \rangle  \subset \mc T$ for  the triangulated subcategory of $\mc T$ 
 generated by $S$, when $S \subset Ob(\mc T)$. 
  We write $\Hom^i(X,Y)$ for  $\Hom(X,Y[i])$.
  In this paper $\KK$ denotes a field.  
  If $\mc T$ is $\KK$-linear triangulated category  we write  $\hom^i(X,Y)$ for  $\dim_\KK(\Hom(X,Y[i]))$, where $X,Y\in \mc T$. 

A $\KK$-linear  triangulated category $\mc T$ is called  \textit{ proper} if $\sum_{i\in \ZZ} \hom^i(X,Y)<+\infty$ for any two objects $X,Y$ in $\mc T$.

An \textit{exceptional object} in a $\KK$-linear triangulated category  is an object $E\in \mc T$ satisfying $\Hom^i(E,E)=0$ for $i\neq 0$ and  $\Hom(E,E)=\KK $. We denote  by ${\mc T}_{exc}$ the set of all exceptional objects in $\mc T$, 

We will often write \textit{ffe functor} instead of fully faithful exact functor in the sequel.

An \textit{exceptional collection} is a sequence $\mc E = (E_0,E_1,\dots,E_n)\subset \mc T_{exc}$ satisfying $\hom^*(E_i,E_j)=0$ for $i>j$.    If  in addition we have $\langle \mc E \rangle = \mc T$, then $\mc E$ will be called a full exceptional collection.  For a vector $\textbf{p}=(p_0,p_1,\dots,p_n)\in \ZZ^{n+1}$ we denote $\mc E[\textbf{p}]=(E_0[p_0], E_1[p_1],\dots, E_n[p_n])$. Obviously  $\mc E[\textbf{p}]$ is also an exceptional collection. 

If an exceptional collection  $\mc E = (E_0,E_1,\dots,E_n)\subset \mc T_{exc}$ satisfies  $\hom^k(E_i,E_j)=0$ for any $i,j$ and for $k\neq 0$, then it is said to be \textit{strong exceptional collection}.   

For  exceptional collections $\mc E_1$, $\mc E_2$ of equal length we  write $\mc E_1 \sim \mc E_2$ if $\mc E_2 \cong \mc E_1[\textbf{p}]$ for some $\textbf{p} \in \ZZ^{n+1}$.

An abelian category $\mc A$ is said to be hereditary, if ${\rm Ext}^i(X,Y)=0$ for any  $X,Y \in \mc A$ and $i\geq 2$,  it is said to be of finite length, if it is Artinian and Noterian.

By $Q$ we denote an acyclic quiver and   by  $D^b(Rep_\KK(Q))$, or just $D^b(Q)$, -  the derived category of the category of $\KK$-representations of $Q$.

For an integer $l\geq 1$ the $l$-Kronecker quiver  (the quiver with two vertices and  $l$ parallel  arrows)  will be denoted by  $K(l)$.

For а subset $S\subset G$  of a group $G$ we denote by $\langle S \rangle \subset G$ the subgroup   generated by $S$.

Let   $x\in \ZZ$ be a positive integer, let $x = \prod_{i=1}^n p_i^{t_i}$ be the prime factorization of $x$,    denote   $\omega(x)=n$ (in particular $\omega(1)=0$), whereas  denote  $\Omega(x)=\prod_{i=1}^n t_i$, finally recall that the M\"obius function is:
\begin{gather}
	\mu(x) = \left \{ \begin{array}{c c} (-1)^{\omega(x)}  & \mbox{if} \  \omega(x)=\Omega(x) \\ 
		0 & \mbox{otherwise.}\end{array} \right. 
\end{gather}

The number of elements of a finite set $X$ we denote by $\abs{X}$ or by $\#(X)$.

For integers  $a,b \in \ZZ$ we denote by $g.c.d(a,b)$ the greatest common divisor of $a,b$.
\section{Basics}
\subsection{A combinatorial indentity}

The following lemma is  well known and we will use it extensively:
\begin{lemma} \label{ci}  For any $k\geq 1$, $n\geq 0$ holds the equality:
	\begin{gather} \abs{ \left \{ \bd \{1,\dots,k\} & \rTo^{\alpha} &  \{0,1,\dots,n\} \ed : 0\leq \alpha_1 \leq \alpha_2 \leq \dots \leq \alpha_k \leq n \right \} } = \binom{n+k}{k}.\end{gather}
\end{lemma}
\bpr

 There is a bijection from the set we consider to the set  \\ $\left \{ \bd \{1,\dots,k\} & \rTo^{u} &  \{0,1,\dots,n+k-1\} \ed : 0\leq u_0 < u_1 < \dots < u_k \leq n+k-1\right \}$, which assigns  $\alpha \mapsto \{\alpha _j+j-1\}_{j=1}^k$; and the lemma follows. 
\epr

\subsection{Exact functors, graded natural transformations, exact equivalences} \label{exact functors}

We start with some well known facts about exact functors between triangulated categories. In this section  we do not presume that the triangulated categories are  linear over a field. 

 The full definitions of an \textit{exact functor between triangulated categories} and of a \textit{graded natural transformation between two such functors} can be found in the literature and we omit these here. 
We have performed all the  proofs of the statements in this section, however we  omit these proofs, since  they should  be available  in the literature. 

Recall that  an exact functor  $\bd \mc A & \rTo^F & \mc B \ed $, where $\mc A$, $\mc B$ are triangulated categories,  is actually a  pair of an additive functor $Fun$ and  a natural isomorphism between the functors $Fun\circ T_{\mc A}$ and $T_{\mc B}\circ Fun$ subject to some axioms, where $ T_{\mc A}$, $ T_{\mc B}$ are the translation functors of $\mc A$ and $\mc B$, respectively. There is a natural definition  of  composition of exact functors (the underlying functor of the composition is composition of the underlying functors) and one shows that this composition rule is associative. 

 A \textit{graded natural transformation} between two exact functors  $\bd \mc A & \rTo^F & \mc B \ed $,  $\bd \mc A & \rTo^{F'} & \mc B \ed $ is a natural transformation subject to an additional condition (compatibility with the natural isomorphisms, which are part of the data in $F$ and $F'$). The exact functors between $\mc A$ and $\mc B$, viewed as objects, and the graded natural transformations, viewed as morphisms, form a category.  \textit{We write $F \cong F'$ if $F$, $F'$ are isomorphic in this category, $\cong$ is an equivalence relation. }    A graded natural transformation is isomorphism in this category iff it is isomorphism between  the underlying functors. The following statements also hold:
 
 An \textit{exact equivalence between two triangulated categories $\mc A$ and $\mc B$} is an exact functor $\bd \mc A & \rTo^F & \mc B \ed $ such that there exists another exact functor $\bd \mc B & \rTo^G & \mc A \ed $, s. t. $G\circ F \cong {\rm Id}_{\mc A}$,   $F\circ G \cong {\rm Id}_{\mc B}$.
 
 \begin{lemma} \label{comp of functors}
 	Let $\bd \mc A & \rTo^F & \mc B \ed $,  $\bd \mc A & \rTo^{F'} & \mc B \ed $, $\bd \mc B & \rTo^G & \mc C \ed $,  $\bd \mc B & \rTo^{G'} & \mc C \ed $ be exact functors. If $F \cong F'$, $G\cong G'$, then $G\circ F \cong G'\circ F' $.
 \end{lemma}

\begin{lemma} \label{exact equivalences}
	An exact functor  $\bd \mc A & \rTo^F & \mc B \ed $ is an exact equivalence iff the underlying functor in $F$ is an equivalence between categories. In particular composition of exact equivalences is an exact equivalence. 
\end{lemma}

\begin{remark}  \label{KK linear}
	Let  $\mc A$, $\mc B$ be triangulated categories  linear over a field $\KK$.
	
	If  $\bd \mc A & \rTo^F & \mc B \ed $ is an exact equivalence, which is $\KK$-linear, then any exact functor    $\bd \mc B & \rTo^G & \mc A \ed $, s. t. $G\circ F \cong {\rm Id}_{\mc A}$,   $F\circ G \cong {\rm Id}_{\mc B}$ is $\KK$-linear as well. 
	
		Let $\bd \mc A & \rTo^F & \mc B \ed $ be an exact $\KK$-linear functor and $\bd \mc A & \rTo^{F'} & \mc B \ed $  be an exact functor. If $F \cong F'$, then $F'$ is also $\KK$-linear.
	
	Composition of $\KK$-linear exact functors is also $\KK$-linear.

\end{remark}
\begin{df} \label{autoequivalences} For a triangulated category $\mc T$ we denote by $\Aut(\mc T)$ the set of equivalence clases $\{\mbox{exact equivalences from} \ \mc T \mbox{to} \ \mc T\}/\cong$. From Lemma \ref{comp of functors} and the previos discussion it follows that  $\Aut(\mc T)$ is group with the operations: $[F]\circ [G]=[F\circ G]$ and unit $[{\rm Id}_{\mc T}]$.
	
	If $\mc T$ is $\KK$-linear for some field $\KK$, then from Remark \ref{KK linear} it follows that $\{ [F] \in \Aut(\mc T): F \ \mbox{is $\KK$-linear}\}$ is a subgroup of $\Aut(\mc T)$, which we denote by $\Aut_{\KK}(\mc T)$.
\end{df}

\subsection{Fully faithful exact functors (ffe functors)} \label{fff}

When we say for a functor   that it is fully faithful exact functor (\textit{ffe functor} for short), we  mean that the underlying functor is fully faithful.

 By straightforward combining the definitions of  exact functor between triangulated categories, fully faithful functor, and the axioms of a triangulated category  one can show that: 
\begin{lemma} \label{ImF}
	Let $\bd \mc A & \rTo^{F} & \mc B \ed$ be a fully faithful   exact functor. Let $\mc A' \subset \mc A$ be a full triangulated subcategory. \footnote{recall that any full triangulated category is closed under isomorphisms} Then the triangulated subcategory   $\langle \{F(X): X \in Ob(\mc A')\} \rangle$ generated by the set of objects $F(Ob(\mc A))$ is equal to the full  subcategory with objects:  $\{X\in \mc B: \exists A\in Ob(\mc A') \ \mbox{s.t.} \ F(A) \cong X \}$. We will denote this subcategory of $\mc B$ by $F(\mc A')$. Obviously, if  $\mc A'\subset \mc A$ is the trivial subcategory of $\mc A$ (the subcategory whose objects are all zero objects), then $ F(\mc A)$ is the  trivial subcategory of $\mc B$.  We will denote $F(\mc A)$ by  ${\rm Im}(F)$. 
 \end{lemma}	
\begin{remark} \label{composition on subcategories} One can  show that for any triangulated subcategory $\mc A'\subset  \mc A$:
	
	 if  $\bd[1em] \mc A & \rTo^{F} & \mc B \ed$,  $\bd[1em] \mc B & \rTo^{G} & \mc C \ed$ are  ffe functors, then $G\circ F(\mc A')=G(F(\mc A'))$;
	 
	  if  $\bd[1em] \mc A & \rTo^{H} & \mc B \ed$ is  ffe functor with a graded natural isomorphism $F\cong H$, then $F(\mc A') =H(\mc A') $;
	  
	  for any two full triangulated  subcategories $\mc A_1, \mc A_2 \subset  \mc A$ we have  $\mc A_1 \subset  \mc A_2 $  iff  $F(\mc A_1) \subset  F(\mc A_2) $. \end{remark}

\begin{remark} \label{auto eq on perp} It is straightforward to show also  that if $\mc A\subset \mc T$ is a triangulated subcategory and $[\alpha]\in \Aut(\mc T)$, then $\alpha(\mc A^{\perp})=\alpha(\mc A)^{\perp}$.  \end{remark}

\begin{coro} \label{restriction of eq}
	Let $\bd \mc A & \rTo^{F} &\mc B \ed$ be a ffe functor. Then by restricting the underlying functor we obtain a functor $\bd \mc A & \rTo & {\rm Im}(F) \ed$. We claim that this restricted functor together with the same natural isomorphism contained in $F$, defines an exact equivalence, which we denote by $\bd \mc A & \rTo^{F_{\vert}} & {\rm Im}(F) \ed$. 
\end{coro}
\bpr The fact that $\bd \mc A & \rTo^{F_{\vert}} &  {\rm Im}(F)  \ed$ is an exact functor follows directly from the definitions.  Since the underlying functor is fully faithful and by the way we definied ${\rm Im}(F)$ in Lemma \ref{ImF} it follows that the underlying functor is equivalence of categories. The corollary follows from Lemma \ref{exact equivalences}.\epr

\begin{coro} \label{two fff are quiv}
	Let  $\bd \mc A & \rTo^{F} &\mc B \ed$, $\bd \mc A & \rTo^{G} &\mc B \ed$ be as in Lemma \ref{ImF}. If ${\rm Im}(F)={\rm Im}(G)$, then there exists $[\alpha ] \in {\rm Aut}(\mc A)$, s.t. $F\cong G \circ \alpha$. If $\mc A$, $\mc B$, $F$, $G$ are $\KK$-linear, then $[\alpha] \in {\rm Aut}_{\KK}(\mc A)$.
\end{coro} \bpr Let us denote $\mc T={\rm Im}(F)={\rm Im}(G)$. From the previous Corollary we get restrictions $\bd \mc A & \rTo^{F_{\vert}, G_{\vert}} & \mc T \ed$ which are exact equivalences. Hence  there are  exact functors $F_{\vert}^{-1}$,  $G_{\vert}^{-1}$, s.t. $F_{\vert}^{-1}\circ F_{\vert}  \cong{\rm Id}_{\mc A} \cong G_{\vert}^{-1}\circ G_{\vert} $, $F_{\vert}\circ F_{\vert}^{-1}  \cong{\rm Id}_{\mc T} \cong G_{\vert}\circ G_{\vert}^{-1} $. From Lemma \ref{comp of functors} and associativity  we see that $(F_{\vert}^{-1}\circ G_{\vert})\circ (G_{\vert}^{-1}\circ F_{\vert}) \cong {\rm Id}_{\mc A}\cong  (G_{\vert}^{-1}\circ F_{\vert})\circ (F_{\vert}^{-1}\circ G_{\vert})$,  $G_{\vert}\circ (G_{\vert}^{-1}\circ F_{\vert} )\cong  F_{\vert}$, therefore $\alpha = G_{\vert}^{-1}\circ F_{\vert}  $ is an exact equivalence such that $G_{\vert}\circ \alpha \cong  F_{\vert}$. Since $F=j\circ F_{\vert}$, $G=j\circ G_{\vert}$, where $\bd \mc T & \rTo^{j} & B \ed$,  is the embedding functor, it follows that $F\cong G \circ \alpha$. Finally, if  $\mc A$, $\mc B$, $F$, $G$ are $\KK$-linear, then $F_{\vert}, G_{\vert}$ are $\KK$-linear exact equivalences and the last statement follows from Remark \ref{KK linear}.
\epr 
The following lemma is also easy to prove:
\begin{lemma} \label{the generated subcategory by images}
	Let $\bd \mc A & \rTo^{F} & \mc B \ed$ be a fully faithful   exact functor. Let $\mc J$ be a family of  non-empty subsets of $Ob(\mc A)$, then $F(\langle \cup_{S\in \mc J} S \rangle)=\langle \cup_{S\in \mc J} F_{ob}(S) \rangle$, where $F(\langle \cup_{S\in \mc J} S \rangle)$ is as in Lemma \ref{ImF},  $F_{ob}$ is the function determining which object from $\mc B$ is assigned to any object in $\mc A$ via the functor $F$, and  $F_{ob}(S)$ is the set-theoretical  image of the subset $S$ via $F_{ob}$. 
\end{lemma}

\section{Main Definition} \label{section main Definition}
Here  as in the previous section   we do not presume that the triangulated categories are  linear over a field. We note first that 

\begin{lemma} \label{eq relation}
	Let $\mc A$, $\mc T$ be triangulated categories. Let $\Gamma_{\mc A}\subset \Aut(\mc A)$,  $\Gamma_{\mc T}\subset \Aut(\mc T)$ be subgroups. Let $S$ be a set, whose elements are  exact functors from $\mc A$ to $\mc T$. Then the following relation 
	\begin{gather} F,F' \in S \qquad  F \sim F' \iff  F\circ \alpha \cong \beta \circ F'  \ \mbox{for some} \ [\alpha] \in \Gamma_{\mc A}, [\beta] \in \Gamma_{\mc T}. \end{gather}
	is an equivalence relation on $S$. 
\end{lemma}
\bpr Follows by combining  of Lemmas \ref{comp of functors}, the associativity of composition of exact functors and Definition \ref{autoequivalences}. 
\epr
When we apply this lemma, we will have $\Gamma_{\mc A}= \Aut(\mc A)$, as in the following:

\begin{df} \label{C_A(T)}   Let  $\mc A$, $\mc T$ be   triangulated categories. And let 
	 $\Gamma \subset {\rm Aut}(\mc T)$ be  subgroup of the group of autoequivalences.  We denote 
	\begin{gather} \label{C'_l,P} 
	C'_{\mc A,P}(\mc T) = \{\bd \mc A & \rTo^{F} & \mc T \ed : F \ \mbox{is fully faithful exact functor satisfying  property} \ P \}.
	\end{gather} 
	Here we need $P$ to be a property of fully faithful functors, such that \eqref{C'_l,P} is a well defined set,  an example is in  Definition \ref{C_A KK(T)}.  
	We fix an equivalence relation in $C'_{\mc A ,P}(\mc T)$, and  define another set:
	\begin{gather}\label{equivealence}
	C_{\mc A,P}^{\Gamma}(\mc T) = C'_{\mc A,P}  (\mc T)/{\sim} \qquad F \sim F' \iff F \circ \alpha \cong \beta\circ F' \ \mbox{for some} \ \ [\alpha] \in {\rm Aut}(\mc A), [\beta] \in \Gamma  
	\end{gather}  where $ F \circ \alpha \cong \beta \circ F' $ means greaded equivalence of exact functors discussed in Subsection \ref{exact functors}.

\end{df} \begin{df} \label{C_A KK(T)}
		If $\mc A$, $\mc T$ are $\KK$-linear triangulated categories, and the property $P$ of $F$ is: ``$F$ is $\KK$-linear'', then we denote the corresponding sets defined in Definition \ref{C_A(T)}  by $	C'_{\mc A,\KK}(\mc T)$, $	C_{\mc A,\KK}^{\Gamma}(\mc T)$.
\end{df}

\begin{remark} If $\mc A$ is trivial category, and the set $	C'_{\mc A,P}(\mc T)$ in Definition \ref{C_A(T)} is non-empty, then $	C^{\Gamma}_{\mc A,P}(\mc T)$ has a single element for any $\Gamma$. In particular, if $\mc A$ is a trivial category in Definition \ref{C_A KK(T)}, then   $	C^{\Gamma}_{\mc A,\KK}(\mc T)$ has a single element  for any $\Gamma$. 
	
\end{remark}
	
\begin{remark} \label{property}

	 For any ffe functor $\bd[1em] \mc A & \rTo^F & \mc T\ed  $ and  $\Gamma \subset {\rm Aut}_{\KK}(\mc T)$ with    $P$ in  Definition \ref{C_A KK(T)} hold:
	
	(a)  $F$ satisfies $P$ iff its restriction $F_{\vert}:\mc A \rightarrow {\rm Im}(F) $ satisfies $P$:
	
	(b) For any  $[\beta] \in \Gamma$, if $F$ satisfies $P$ then $\beta \circ F$ also satisfies $P$.

\end{remark}

\begin{lemma} \label{bijection} Let $\mc A$, $\mc T$ be as in Definition \ref{C_A(T)}. Let $P$  satisfy (a) in Remark \ref{property}.
	There is a well defined  function 	\begin{gather}  \label{bijection in formula}
	C_{\mc A, P}^{\{\rm Id\}}(\mc T) \rightarrow\left  \{\mc B \subset \mc T:\begin{array}{l} \mc B \ \mbox{is a full tr. subcategory s.t.}  \\  \mbox{there exists an exact equivalence} \\ \bd \mc A & \rTo^{G}& \mc B \ed  \ \mbox{satisfying property} \ P  \end{array}
\right \} \qquad 	C_{\mc A, P}^{\{\rm Id\}}(\mc T) \ni	[F] \mapsto {\rm Im}(F) 	\end{gather}
where ${\rm Im}(F)$ is explained in Lemma \ref{ImF}. The function is bijection. 	
\end{lemma}
\bpr  From Lemmas  \ref{ImF} and \ref{restriction of eq}, and by (a) in Remark \ref{property} it follows that ${\rm Im}(F)$  is an element of the codomain  for any $F\in C_{\mc A, P}'(\mc T)$. If $F, F' \in C_{\mc A, P}'(\mc T)$ and $F \cong F'\circ \alpha$ for some $[\alpha] \in \Aut(\mc A)$, then ${\rm Im}(F)={\rm Im}(F')$, therefore the function is well defined. Lemma \ref{two fff are quiv} implies that the function is injective. If $\mc B \subset \mc T$ is an element of the codomain  with equivalence $G$, then ${\rm Im}(j \circ G)= \mc B$ \footnote{$j$ is the embedding functor of $\mc B$, recall also that each triangulated subcategory is isomorphism closed},  and $(j \circ G)_{\vert} = G$ and again by (a) in Remark \ref{property}    $j\circ G\in C'_{\mc A, P}(\mc T)$, therefore  $[j\circ G]\in C_{\mc A, P}^{\{\rm Id\}}(\mc T)$ is mapped via the function defined here to $\mc B$, hence the function is surjective.  \epr

\begin{prop} \label{functoriality} 	Let $\bd[1em]  \mc T_1 & \rTo^{F_1} & \mc T_2\ed$, $\bd[1em] \mc T_2 & \rTo^{F_2} & \mc T_3 \ed$ be  ffe   functors.   Let  $\mc A$ be a  triangulated category. 	
	
	 Let $P$ be a property of ffe functors such that: 
	 
	 {\rm (a)}  any ffe functor $\bd[1em]\mc A& \rTo^F & \mc T_i \ed $   satisfies $P$ iff its restriction $\bd[1em]\mc A& \rTo^{F_{\vert}} & \mc T_i \ed $  satisfies $P$  for $i=1,2,3$ ; 
	 
	{\rm (b) } for any ffe  functor  $\bd[1em]  \mc A & \rTo^{F} & \mc T_i\ed$ with $P$ it follows that $\bd  \mc A & \rTo^{F_i\circ F} & \mc T_{i+1}\ed$    has  $P$ for   $i=1,2$ .
	
	Then there is   a well defined injective map:
	$
	\bd	C_{\mc A, P}^{\{\rm Id\}}(\mc T_1) & \rTo^{F_{1*}} &	C_{\mc A, P}^{\{\rm Id\}}(\mc T_2) \ed $ defined by $	[F] \mapsto [F_1 \circ F] $,
	 which via the bijection \eqref{bijection in formula} translates to a map  defined by  $ \mc B \mapsto  F_1(\mc B)$.
	 
In particular, we have also injections $F_{2*}$, $(F_{2} \circ F_{1})_*$. 
	Furthermore, we have  $(F_{2} \circ F_{1})_*=F_{2*} \circ F_{1*}$.
\end{prop}
\bpr  The fact that the function is well defined follows from the assumption for $P$. From the associativity of the composition,  and from  Lemma \ref{comp of functors}.  Let $[F_1\circ F] = [F_1\circ F']$ in $C_{\mc A, P}^{\{\rm Id\}}(\mc T_1)$. Then from Lemma \ref{bijection}  it follows that ${\rm Im}(F_1\circ F)={\rm Im}(F_1\circ F')$ in $\mc T_2$, hence, see Remark \ref{composition on subcategories},  $F_1 ({\rm Im}( F))= F_1({\rm Im}( F'))$, and  the fully-faithfulness of $F_1$ implies that ${\rm Im}( F)= {\rm Im}( F')$ in $\mc T_1$, and now Lemma \ref{comp of functors} implies that  $[ F] = [ F']$ in $C_{\mc A, P}^{\{\rm Id\}}(\mc T_1)$. Thus, we proved that the function is injective. 

The last statement follows from the associativity of composition of exact functors. 
\epr

\begin{coro} \label{group action}
	Let $\mc A$, $\mc T$, $P$,  $\Gamma$ be as in Definition \ref{C_A(T)}. Let $P$, $\Gamma$ satisfy (a), (b) in Remark  \ref{property}.
	
	 There is a group action $\Gamma \times C_{\mc A,P}^{\{\rm Id\}}(\mc T) \rightarrow  C_{\mc A, P}^{\{\rm Id\}}(\mc T)$ defined by $([\beta],[F]_{\{\rm Id\}}) \mapsto [\beta \circ F]_{\{\rm Id\}}$, which via the bijection \eqref{bijection in formula} translates to action on the codomain defined by  $([\beta], \mc B) \mapsto \beta(\mc B)$.  
	
	 The natural projection  $\begin{array}{c} C_{\mc A, P}^{\{\rm Id\}}(\mc T) \rightarrow C_{\mc A,P}^{\Gamma}(\mc T) \end{array}$  descends to a bijection $ C_{\mc A, P}^{\{\rm Id\}}\left ( \mc T\right )/\Gamma \rightarrow C_{\mc A, P}^{\Gamma}\left ( \mc T\right ) $. 
	
	In particular, if $\mc A, \mc T$ are $\KK$-linear and $\Gamma \subset {\rm Aut}_{\KK}(\mc T)$, then we get a bijection (see Definition \ref{C_A KK(T)})   $ C_{\mc A, \KK}^{\{\rm Id\}}\left ( \mc T\right )/\Gamma \rightarrow C_{\mc A, \KK}^{\Gamma}\left ( \mc T\right ) $. 
\end{coro} 
\bpr  The first part follows from Proposition \ref{functoriality}. The rest is also easy combining of facts from subsections \ref{exact functors} and \ref{fff} and (a), (b) in Remark  \ref{property}.
\epr

\section{Counting non-commutative curves in a given triangulated category} \label{counting non-commutative curves section}
	Following Kontsevich-Rosenberg  \cite{KR} we denote  $D^b(K(l+1))$ by  $N\PP^l$  for $l\geq -1$,  this is a $\KK$-linear triangulated category. In this Section we specialize the Definition \ref{C_A(T)} by choice $\mc A = N\PP^l$.
	
	\begin{remark} \label{exce pairs in NPl}
		For each $l\geq -1$  $N\PP^l$ has a full strong  exceptional pair $(A,B)$. Each exceptional pair in  $N\PP^l$ is full (follows from \cite{WCB1}). Each strong exceptional   pair $(A,B)$ has $\hom(A,B)=l+1$.
	\end{remark}

\begin{remark} \label{full exact has lr adjoints} Let $l\geq -1$ and $\mc T$ be any  triangulated cateogry linear over $\KK$,  let $\bd N \PP^l & \rTo^{F} & \mc T \ed$ be any fully faithful  exact $\KK$-linear  functor.  Then ${\rm Im}(F)$ is a triangulated subcategory of $\mc T$ generated by two exceptional objects, hence due to \cite[Theorem 3.2]{B} the funcor $F$ has   left and right adjoints, the subcategory $\mc A= {\rm Im}(F)$ is admissible, and  there are SOD $\mc T=\langle \mc A, \mc A^{\perp} \rangle $, $\mc T=\langle  ^{\perp}\mc A, \mc A \rangle $, in particular $\mc A$ is closed under direct summands.
	
\end{remark}

\begin{remark} \label{functor from NP tp NP} From the previous remarks it follows that, if for some integer  $j\geq -1$ there exists a ffe  $\KK$-linear functor $\bd N \PP^j & \rTo^{F} &  N \PP^l \ed$, then $l=j$ and $F$ is an exact  equivalence.
\end{remark}
\begin{df} \label{C_l} Let $l\in \ZZ_{\geq -1}$ and   let $\mc T$ be any  triangulated category linear over $\KK$. And let $\Gamma \subset {\rm Aut}_{\KK}(\mc T)$ (see Definition \ref{autoequivalences} ).   We denote $C_{N\PP^l, \KK }^\Gamma(\mc T)$ (deifined in Definition \ref{C_A KK(T)}) by $C_{l}^{\Gamma}(\mc T)$ and refer to the elements of $C_{l}^{\Gamma}(\mc T)$ as to non-commutative curves  of non-commutative  genus $l$ in $\mc T$ modulo $\Gamma$.  Unless explicitly specified otherwise,  writing genus in this paper we always mean non-commutative genus in the sense explained here.
	
\end{df}

\begin{prop} \label{bijection 123}
	Let $\mc T$ be a triangulated category linear over $\KK$.  Let $\mc T$ admits an enhancement. Let $C_{l}^{\{\rm Id\}}(\mc T)$ be as in Definition \ref{C_l}. Then the bijcetion from Lemma \ref{bijection}  can be described as follows:
	\begin{gather}
\label{bijection 1} 	C_{l}^{\{\rm Id\}}(\mc T) \rightarrow\left  \{\mc A\subset \mc T:\begin{array}{l} \mc A \ \mbox{is a full tr. subcategory s.t.} \ \mc A =\langle E_1,E_2 \rangle \\  \mbox{for some strong exceptional pair} \ (E_1,E_2) \\ \mbox{s.t.} \ \hom(E_1,E_2)=l+1 \end{array} \right \} \ 
	[F] \mapsto  {\rm Im}(F). 
	\end{gather}
	For any  full exceptional pair $(A,B)$ in $N\PP^l$ this function is defined by the assignement:\\ $ C_l^{\{\rm Id\}} \ni [F] \mapsto \langle F(A),F(B) \rangle $. 
	
	 When $\mc A = D^b(A_1)=D^b(point)$, then    the bijection from Lemma \ref{bijection}  can be described as follows:
	 	\begin{gather}
	 	\label{bijection for eo} 	C_{ D^b(A_1), \KK }^{\{\rm Id\}}(\mc T) \rightarrow\left  \{\mc A\subset \mc T:\begin{array}{l} \mc A \ \mbox{is a full tr. subcategory s.t.} \ \mc A =\langle E\rangle \\  \mbox{for some exceptional object} \ E \end{array} \right \} \ 
	 	[F] \mapsto  {\rm Im}(F). 
	 	\end{gather}
\end{prop}
\bpr
 $N\PP^l$ has a full strong exceptional pair $(A,B)$ with $\hom(A,B)=l+1$ (Remark \ref{exce pairs in NPl}). The definition of a fully faithfull exact $\KK$-linear functor implies that $(F(A), F(B))$ is a strong exceptional collection in $\mc T$ with $\hom(F(A), F(B))=l+1$, furthermore also quickly from the definitions one shows that ${\rm Im}(F)=\langle F(A), F(B) \rangle $, so we see that ${\rm Im}(F)$ is indeed in the codomain in \eqref{bijection 1} and that the last statement in the Proposition holds.

Let $\mc A \subset \mc T$ be an element of  the codomain of \eqref{bijection 1}, therefore it has a full strong exceptional pair $(E_1,E_2)$. Since $\mc T$ admits an enhancement, it follows that $\mc A$ also admits an enhancement. Now we can apply \cite[Corollary 1.9]{Orlov} to $\mc A$, and we conclude that there is an exact $\KK$-linear equivalence $\bd N\PP^l & \rTo^F & \mc A\ed$. Therefore $\mc A$ is in the codomain of the bijection of Lemma \ref{bijection} and the proposition is proved. 
The proof of \eqref{bijection for eo} is analogous. 
\epr

\begin{remark} \label{no enhancement needed for points}
	The bijection \eqref{bijection for eo}  can be proved for any proper $\KK$-linear category $\mc T$, without restricting it  to have  an enhancement.
\end{remark}

\section{Introduction to intersection theory} \label{Interesection theory intro}

\subsection{Definition of $C_{{\mc J}, P}(\mc T)$ for a family of categories ${\mc J}$ and partial order in  $C_{{\mc J}, P}(\mc T)$}  \label{IT}

Taking into account Lemma \ref{bijection}  we can extend Definition \ref{C_A(T)} to the case of family of source categories as follows:

\begin{df} \label{def of C_J(T)} Let  $\mc J $ be a small set of  of  pairwise non-equivalent  triangulated categories, let $\mc T$ be a triangulated category.  Let $P$ be a property of fully faithful functors from $\mc A$  to subcategories in $\mc T$ satisfying  (a) in Remark \ref{property} for any $\mc A \in \mc J$.  From this data we  define a partially ordered set $C_{\mc J, P}(\mc T)$  as follows. The set is: 
	$$ C_{\mc J, P}(\mc T)= \left  \{\mc B \subset \mc T:\begin{array}{l} \mc B \ \mbox{is a full tr. subcategory s.t. for some} \ \mc A \in \mc J  \\  \mbox{there exists an exact equivalence} \\ \bd \mc A & \rTo^{G}& \mc B \ed  \ \mbox{satisfying property} \ P  \end{array}
	\right \}.$$
	  We define the partial order in this set as follows: for  $x,y \in  C_{\mc J, P}(\mc T)$ we define  $x\leq y$ iff the set of objects of the subcategory $x$ is a subset of the set of objects of the subcategory $y$. 
	
\end{df}

\begin{df} \label{C_J KK(T)}
	Let   $\mc J $, $\mc T$ be as in Definition \ref{def of C_J(T)}. If  any $\mc A \in \mc J$ and $\mc T$ are $\KK$-linear triangulated categories, and the property $P$ of $F$ is: ``$F$ is $\KK$-linear'', then we denote the corresponding partially ordered set defined in Definition \ref{def of C_J(T)}   by $		C_{\mc J,\KK}(\mc T)$. 
\end{df}
\begin{remark}
If $\mc J$ contains a trivial category, then  $C_{\mc J, \KK}(\mc T)$ contains the trivial subcategory in $\mc T$ and denoting this subcategory by $T$ we have $T\leq \mc B$ for any $\mc B \in C_{\mc J, \KK}(\mc T)$.
\end{remark}

\subsection{The case $\mc J = \{\mbox{trivial category}, D^b(A_1), D^b(A_2)\}$}  We will show here that in the partially ordered set $C_{\mc J, \KK}(\mc T)$ any two elements have a greatest lower bound. 

\begin{df} \label{interesection} In this section $\mc T$ is a $\KK$-linear proper triangulated category, and $C_l^{\{\rm Id\}}(\mc T)$ is the set of non-commutative curves in $\mc T$ of genus $l\geq -1$ as defined in Definition \ref{C_l}. Let $[F_i]\in C_{l_i}^{\{\rm Id\}}(\mc T)$, let $\mc A_i \subset \mc T $ be the subcategories to which $[F_i]$ are mapped via the bijection  
	\eqref{bijection in formula}, $i=1,2$. 	 From Remark \ref{full exact has lr adjoints} the full subcategory with objects  $Ob(\mc A_1) \cap Ob(\mc A_2)$ is a triangulated  subcategory, which is closed under direct summands. We  call it the \underline{intersection of the non-commutative curves $[F_1]$, $[F_2]$ in $\mc T$}.
	
	A subcategory in $\mc T$ generated by an exceptional object will be called a \underline{derived  point} of $\mc T$ . Via the bijection \eqref{bijection for eo}
	the derived points in $\mc T$ are exactly the elements in $C_{D^b(A_1),\KK}(\mc T)$.
\end{df}

This intersection is easy to describe when one of the curves is of genus $0$ or genus $-1$.
\begin{lemma} \label{derived points in a genus 0} Let $\mc T$ be any $\KK$-linear proper category. 
	Any non-commutative curve of genus $0$, i.e. element in the set  $C_{0}^{\{\rm Id\}}(\mc T)$, contains  exactly three elements of the set $C_{D^b(A_1),\KK}^{\{\rm Id\}}(\mc T)$ (in terms of the terminology from the definition contains exactly three derived points in $\mc T$). Any two of these three  derived points  span the curve (the subcategory generated by them is the initial curve).  
\end{lemma}
\bpr  Let  $[F]\in C_{0}^{\{\rm Id\}}(\mc T)$ be a  non-commutative curve of genus zero and let    $\mc A$ be the subcategory of $\mc T$ to which $[F]$ is mapped via the bijection  
\eqref{bijection in formula}. $\mc A$ is  a full triangulated subcategory of $\mc T$ and   $\bd D^b(A_2)=D^b(K(1)) & \rTo^{F_{\vert}} & \mc A \ed$, is a $\KK$-linear  equivalence of triangulated categories.  By Lemma \ref{C_eoA_n+1} it follows that $C_{D^b(A_1),\KK}^{\{\rm Id\}}(\mc A)$ has three elements (using the notations of Lemma \ref{C_eoA_n+1} these are $\langle F(s_{0,0}) \rangle$,$ \langle F(s_{1,1}) \rangle$, $ \langle F(s_{0,1}) \rangle$).   On the other hand  $C_{D^b(A_1),\KK}^{\{\rm Id\}}(\mc A)$ is exactly the set \\ $\left \{\alpha \in C_{D^b(A_1),\KK}^{\{\rm Id\}}(\mc  T): \alpha \subset \mc A\right \}$, hence the first part of the lemma is proved.  

Since $(s_{0,0},s_{1,1})$, $(s_{1,1},s_{0,1})$, $(s_{0,1},s_{0,0})$ are full exceptional pairs in $D^b(K(1))= D^b(A_2)$ (see Lemma \ref{help lemma 1} and  Remark \ref{exce pairs in NPl}) it follows that any two of the three derived points  $\langle F(s_{0,0}) \rangle$,$ \langle F(s_{1,1}) \rangle$, $ \langle F(s_{0,1}) \rangle$ span the non-commutative curve.
\epr

\begin{lemma}  \label{intersection}  Let $\mc T$ be as in Definition \ref{interesection}. If the genus of one of two different non-commutative curves is  less or equal to $0$, then they intersect either trivially or in  a single derived  point  in $\mc T$.	
\end{lemma} 
\bpr Let  $[F_1]\in C_{l_1}^{\{\rm Id\}}(\mc T)$, $ [F_2] \in C_{l_2}^{\{\rm Id\}}(\mc T)$, let $\mc A, \mc B \subset \mc T $ be the subcategories to which $[F_1], [F_2]$ are mapped via the bijection  
\eqref{bijection in formula}. Now $\bd D^b(K(l_1+1)) & \rTo^{F_{1\vert}} & \mc A \ed$, $\bd D^b(K(l_2+1)) & \rTo^{F_{2\vert}} & \mc B \ed$ are equivalences of triangulated categories.  Assume $l_1=0$ or $l_1=-1$. 

Note that $\mc A\cap \mc B$ is a proper subcategory of $\mc A$ (otherwise $Ob(\mc A) \subset Ob(\mc B)$ and then the exceptional pair of $\mc A$ should generate the category $\mc B$, since by Remark \ref{exce pairs in NPl} each exceptional pair in $\mc B$ is full, and it would follow that $\mc A = \mc B$). 

It remains to show that if $\mc A \cap \mc B$ is not trivial (i.e. contains a nonzero object ), then it $\mc A \cap \mc B = \langle E \rangle$ for some exceptional object $E$ in $\mc T$. 
Now $\mc A$ is equivalent to  $ D^b(K(1))$ or $ D^b(K(0))$ and it follows that:
\begin{gather}\label{property 1 for A}  \mbox{\textit{The indecomposable objects  in $\mc A$ are exceptional objects.}}    \\ \label{property 2 for A}  \mbox{\textit{ Each object is direct sum of exceptional objects.}} \\  \mbox{\textit{For any two exceptional objects $E_1$, $E_2$ in $\mc A$, if $\langle E_1 \rangle \not = \langle E_2 \rangle$, then}} \nonumber \\[-3mm] \label{property 3 for A}  \\[-3mm]
\mbox{ 	\textit{ either $(E_1, E_2)$ or $(E_2, E_1)$ is a full exceptional collection in $\mc A$.}} \nonumber
\end{gather}

Since $\mc A \cap \mc B$ is closed under direct summands and non-trivial it follows from \eqref{property 2 for A} that it must contain at least one exceptional object $E$.  Since $\mc A \cap \mc B$ is proper subcategory in $\mc A$ it follows from \eqref{property 3 for A} that for any other exceptional object $E'\in \mc A\cap \mc B$  we have $\langle E' \rangle = \langle E \rangle $, now using again \eqref{property 2 for A} and the fact that $\mc A \cap \mc B$ is closed under direct summands it follows that $\mc A \cap \mc B=\langle E \rangle$. \epr
\begin{coro} \label{interesection genus zero} Let $\mc T$ be any proper  $\KK$-linear  category.

	If $\mc A$, $\mc B$ are two different non-commutative curves of genus zero, let $\alpha_1, \alpha_2, \alpha_3$,  $\beta_1, \beta_2, \beta_3$ be the derived points in $\mc A$ and $\mc B$ respectively as in Lemma \ref{derived points in a genus 0}.
	
	Then  $\mc A$ intersects  $\mc B$ non-trivially iff $\{\alpha_1, \alpha_2, \alpha_3\}  \cap \{\beta_1, \beta_2, \beta_3\}$ contains a single element and the intersection of   $\mc A$ and  $\mc B$ is exactly this element.
\end{coro}
\bpr 
Follows from   Lemmas \ref{intersection} and \ref{derived points in a genus 0}. 
\epr
\begin{coro} \label{greatest lower bound} Let $\mc T$ be any proper  $\KK$-linear  category, and let $\mc J = \{\mbox{trivial category}, D^b(A_1), D^b(A_2)\}$ or $\mc J = \{\mbox{trivial category}, D^b(A_1), D^b(A_1) \oplus  D^b(A_1)  , D^b(A_2)\}$. For any two $x,y \in C_{\mc J, \KK}(\mc T)$ the full subcategory whose objects are the set theoretical intersection of the sets of objects of $x$ and $y$ is also an element of  $ C_{\mc J, \KK}(\mc T)$ and it is a greatest lower bound of  $x,y$
\end{coro}
\bpr 
The fact that $x\cap y \in C_{\mc J,\KK}(\mc T)$ follows from Lemma \ref{intersection}. The rest is trivial.  
\epr

\section{Directed graph  of subcategories in a given category } 
\label{section for digraph} 
First we fix some terminology. By a  \textit{directed graph (or for short digraph)} $G$ we mean a pair $G=(A,B)$ of a set $A$, whose elements we call vertices,  and a subset $B\subset A\times A \setminus diag ( A\times A)$, whose elements we call edges. When we draw pictures, we represent the elements in $B$ by points (in Euclidean plane or in Euclidean space), and a pair $(a,b) \in B$, such that $(b,a) \not \in B$, we present by an arrow starting from the point representing $a$ and directed towards the point representing $b$, we call such pairs \textit{directed edges}. If both $(a,b)\in B$ and $(b,a) \in B$, then we draw  a double sided arrow between $a$ and $b$.

\textit{From any directed graph $G = (A,B)$ we define a simplicial complex $SC(G)$  as follows:}
\begin{gather} SC(G) = (A, B') \nonumber \\[-2mm] \label{sc from digraph} \\[-2mm]  \nonumber  B'=\{X\subset A: \exists n\in \ZZ_{\geq 1}  \ \mbox{and a bijection } \{1,\dots,n\} \rightarrow X \ \mbox{s.t.} \  (x_i,x_j) \in B \  \mbox{ for }  \ 1\leq j < j\leq n \ \} \end{gather}

\textit{Throughout this section    $\mc J $ denotes  a set of  of  pairwise non-equivalent  triangulated categories and  $\mc T$ is a triangulated category. We say that $\mc J$ is $\KK$-linear if $\mc A$ is so for any $\mc A \in \mc J$.}
\begin{df} \label{SC_A(T)}   Let $P$ be a property of fully faithful functors from $\mc A$  to subcategories in $\mc T$ satisfying  (a) in Remark \ref{property} for any $\mc A \in \mc J$.  From this data we  define here a directed graph $G_{\mc J, P}(\mc T)=(A,B)$ whose set  of vertices is equipped with a partial order. 
	
	The set of vertices   is $C_{\mc J, P}(\mc T)$ as defined in Definition  \ref{def of C_J(T)} with the partial order explained there. 
	
	The set of  edges is: 	
	$B=\left \{(\mc B_1,\mc B_2) : \mc B_1, \mc B_2 \in C_{\mc J, P}^{\{\rm Id\}}(\mc T) \ \mbox{and} \ \Hom(\mc B_2, \mc B_1)=0 \right  \}. $
	
	The simplicial complex  $SC(G_{\mc J, P}(\mc T))$, defined in \eqref{sc from digraph}, will be denoted for short by $SC_{\mc J, P}(\mc T)$. 
\end{df}
\begin{df} \label{G_A KK(T)}  Let   $\mc T$, $\mc J$ be  $\KK$-linear and the property $P$ of $F$ be: ``$F$ is $\KK$-linear'', then we denote the corresponding graph and simplicial complex defined in Definition \ref{SC_A(T)}  by $	G_{\mc J,\KK}(\mc T)$, $SC_{\mc J,\KK}(\mc T)$.
\end{df}
As a corollary of Proposition \ref{functoriality} and Remark \ref{composition on subcategories} we get  
\begin{coro} \label{fully faithful functor on digraphs}  Let  $\bd[1em]  \mc T_1 & \rTo^{F_1} & \mc T_2\ed$, $\bd[1em] \mc T_2 & \rTo^{F_2} & \mc T_3 \ed$ be  ffe  functors.   Let $P$ be a property of ffe functors  satisfying  {\rm (a), (b)}  in Proposition \ref{functoriality} for any $\mc A \in \mc J$, so Definition  \ref{SC_A(T)} gives graphs $\{G_{\mc J, P}(\mc T_i)=(A_i,B_i)\}_{i=1}^3$.

Let $i=1$ or $i=2$.
 The   assignments $A_i \ni \mc B \mapsto F_i (\mc B)$, $B_i \ni (\mc B, \mc B') \mapsto  (F_i(\mc B), F_i(\mc B')) $\footnote{ $F_i(\mc B)$ is explained in Definition \ref{ImF}} define a full embedding of the graph $G_{\mc J, P}(\mc T_i)$ into  $G_{\mc J, P}(\mc T_{i+1})$  in the sense that it is injective on veritces and for any pair $(a,b) \in A_i\times A_i $ we have $(a,b) \in B_i$ if and only if $(F_i(a),F_i(b)) \in B_{i+1}$. This assignment preserves the partial orders in $A_i$, $A_{i+1}$ as well, more precisely  for any $x,y\in C_{\mc A, P}^{\{\rm Id\}}(\mc T_i)$ we have $x\leq y$ iff $F_{i}(x)\leq F_{i}(y)$.
		Thus, we can view $G_{\mc J, P}(\mc T_i)$ as a full subgraph of $G_{\mc J, P}(\mc T_{i+1})$.
		
		Let us denote by $\bd[1em] G_{\mc J, P}(\mc T_i) & \rTo^{F_{i*}} & G_{\mc J, P}(\mc T_{i+1})\ed$, $i=1,2$ the induced embeddings of digraphs preserving the partial orders.  We have $(F_{2} \circ F_{1})_*=F_{2*} \circ F_{1*}$. If  $F_i$ is equivalence, then $F_{i*}$ is isomorphism.
	\end{coro}

\begin{coro} \label{action on digraphs} We get an action of $\Aut_{\KK}(\mc T)$ on the digraph  $G_{\mc J, \KK}(\mc T)$ for any  $\KK$-linear  $\mc T$ and $\mc J$. 
\end{coro}

Under some conditions we can embed $C_{\mc A,\KK}^{\Gamma}(\mc T)$ (see Definition \ref{C_A KK(T)}) in a certain set of equivalence classes of subgraphs in the digraph  $G_{\mc J, \KK}(\mc T)$. More precisely:

\begin{prop} \label{from invariants to subgraphs} 
	Let $\mc A$, $\mc T$ be $\KK$-linear triangulated categories, let $\mc J$ be $\KK$-linear as well, and $\Gamma \subset {\rm Aut}_{\KK}(\mc T)$ be subgroup. 
	 Let $\mc A$ be generated by the vertices of $G_{\mc J, \KK}(\mc A)$.
	
	There is a well defined injective function 	\begin{gather}  \label{injection graphs in formula}
	C_{\mc A, \KK}^{\rm \Gamma }(\mc T) \rightarrow\left  \{\mc B \subset G_{\mc J, \KK}(\mc T):\begin{array}{l} \mc B \ \mbox{is a full subgarph s.t.}  \\  \mbox{there exists an  isomorphism} \\ \bd  G_{\mc J, \KK}(\mc A) & \rTo^{\mu}& \mc B \ed  \ \mbox{of digraphs} \end{array}
	\right \}/\Gamma \quad 	C_{\mc A, \KK}^{\Gamma}(\mc T) \ni	[F] \mapsto [{\rm Im}(F_{*})] 	\end{gather}	 
	where for a fully faithful functor $F$ by $F_{*}$ is denoted the induced embedding of digraphs as in Corollary  \ref{fully faithful functor on digraphs} and  ${\rm Im}(F_{*})$ is the full subgraph of $  G_{\mc J, \KK}(\mc  T)$ determined by this embedding. The quotient by $\Gamma$ on the right hand side means that we take equivalence classes of subgraphs, where two subgraphs $B_1$, $B_2$ are in the same class iff there exists $[\beta] \in \Gamma$, s.t. $\beta_*(B_1)=B_2$.
\end{prop}
\bpr If $\mc A$ is trivial then both sets have single elements and the function is obviously bijective. So let $\mc A$ be non-trivial. Let us first show that the assignment is well defined. If $[F_1] = [F_2]$, then there exists an auto-equivalence $\alpha$ of $\mc A$ and $[\beta] \in \Gamma$, such that $\beta \circ F_1 \cong F_2 \circ \alpha$ and since $F_1$, $F_2$, $\beta$ are $\KK$-linear, $\alpha$ is also $\KK$-linear (as in the proof of Corollary \ref{two fff are quiv}). Now using the second statement in Remark \ref{composition on subcategories} and the last statement in Corollary \ref{fully faithful functor on digraphs} we see that $\beta_* \circ F_{1*} = F_{2*} \circ \alpha_*$. On  the other hand from Corollary \ref{action on digraphs} we know that $\alpha_*$ is an auto-equivalence, hence $\beta_*({\rm Im}( F_{1*})) = {\rm Im}(\beta_* \circ F_{1*})={\rm Im}(F_{2*} \circ \alpha_*)={\rm Im}(F_{2*})$, hence $[{\rm Im}(F_{1*})]=[{\rm Im}(F_{2*})]$.

To prove that the function is injective, denote  by $A$ the set of vertices of  $G_{\mc J, \KK}(\mc A)$. Using Lemma \ref{the generated subcategory by images} and the fact that $\mc A$ is generated by  the elements in $A$  we can write  for any ffe  functor $\bd  \mc A & \rTo^{F}& \mc T \ed$ 
\begin{gather} \nonumber 
{\rm Im}(F)= F\left (\langle \cup_{\mc B \in A} Ob(\mc B) \rangle \right )=\langle \cup_{\mc B \in A} F_{ob}\left ( Ob(\mc B) \right ) \rangle = \langle \cup_{\mc B \in A} Ob( F\left ( \mc B \right )) \rangle =  \langle \cup_{\mc B \in V({\rm Im}(F_{*}))} Ob( \mc B ) \rangle, 
\end{gather}
where  we  denote  by $V({\rm Im}(F_{*}))$  the vertices of the digraph $G_{\mc J, \KK}({\rm Im}(F_{*}))$.
Hence we see that for any two  fully faithful exact  functors $\bd  \mc A & \rTo^{F_1, F_2}& \mc T \ed$ we have 
$
{\rm Im}(F_{1*})= {\rm Im}(F_{2*}) \ \  \Rightarrow \ \ {\rm Im}(F_{1})= {\rm Im}(F_{2})
$. Assume now that $[{\rm Im}(F_{1*})]=[{\rm Im}(F_{2*})]$, therefore $\beta_*({\rm Im}(F_{1*}))={\rm Im}(F_{2*})$ for some  $[\beta] \in \Gamma$, hence ${\rm Im}((\beta \circ F_{1})_{*})={\rm Im}(\beta_* \circ F_{1*}))=\beta_*({\rm Im}(F_{1*}))={\rm Im}(F_{2*})$ and we proved that this implies  ${\rm Im}(\beta \circ F_{1})={\rm Im}( F_{2})$, and now from Corollary \ref{two fff are quiv} and the definition of  $C_{\mc A, \KK}^{\rm \Gamma }(\mc T)$ it follows that $[F_1]=[F_2]$.
\epr

\subsection{The case $\mc J=\{D^b(A_1)\}$}  \label{the case J eq DbA1}
By a \textit{valued directed graph}  we mean a directed graph  $G = (A,B)$  together with a positive  integer attached to each  edge $(a,b)\in B$ such  that $(b,a) \not \in B$.

Let the family $\mc J$ contains only the category $D^b(A_1)$. 
	For a $\KK$-linear  $\mc T$	the vertices of $G_{D^b(A_1), \KK}(\mc T) =(A,B)$ are  $A=C_{D^b(A_1), \KK}(\mc T)$, which is the same as  equivalence classes of exceptional objects (the equivalence is $\sim$ as explained in Section \ref{notations} ).
	
	In this case one can make  $G_{D^b(A_1), \KK}(\mc T) =(A,B)$  a valued directed graph  as follows. For any such pair $(a,b)$ we  take any exceptional objects $X$, $Y$ with $\langle X \rangle = a$, $\langle Y \rangle = b$ and  attach  the number  $\sum_{i\in \ZZ} \hom(X,Y[i])$ to it.  We have:
	
	{\rm (a)} The  embedding of digraphs $F_*$ induced by a ffe functor is such that  the number attached to an edge $(a_1,a_2)$ in $G_{D^b(A_1),\KK}(\mc T_1)$ is the same as to   $(F(a_1),F(a_2))$.

	{\rm (b)}  Proposition  \ref{from invariants to subgraphs} can be strengthened as follows:  If  $\mc A$ is generated by its exceptional objects, then in \eqref{injection graphs in formula} we can modify the right hand side by imposing on $\mu$ to be isomorphism of valued digraphs and then yet the function is well defined and it is injective.

\section{First  examples: two affine  quivers} \label{non-trivial examples} In \cite{DK4}  we postponed  writing  down the details of the  proof of Proposition  \ref{Cl(PP2)}  and of part of Proposition  \ref{the numbers for two quivers}   for a future work. 

In this Section $\KK$ is algebraically closed. 
Using  results of \cite{DK1} and \cite{DK3}  we 
will prove:

\begin{prop} \label{the numbers for two quivers}  Let $\mc T_i=D^b(Q_i)$, $i=1,2$,   where:
 $Q_1$, $Q_2$ are in \eqref{Q1}.
Then the numbers $ \abs{C_{l}^{\Gamma}(\mc T_i) }$  for $l\in \{-1,0,+1\}$, $i=1,2$,  $\Gamma \in \{ \{\rm Id\}, \langle S \rangle, {\rm Aut}_{\KK}(\mc T_i) \}$ are in \eqref{table with numbers Q1},   \eqref{vanishings}. Furthermore holds  \eqref{triples in Q_2}.  
\end{prop} 
	The vanishings \ref{vanishings} were already explained in \cite[Proposition 12.9]{DK4}.

To compute the numbers in these tables we first collect some  material from   \cite{Miyachi}.

\subsection{Preliminaries from \cite{Miyachi} } \label{material from Mi}
Here $Q$ is one of the quivers in \eqref{Q1}, hence it is  acyclic and representation infinite. In \cite{Miyachi} is described   the group of auto-equivalences of $D^b({\rm mod} A)$, where ${\rm mod} A$ is the category of finite dimensional left $A$-modules for $A$, which are  path algebras of some acyclic quivers. With the convention  for defining the path algebra of a quiver in  \cite{Miyachi}, \cite{Ha} (see \cite[p. 347]{Miyachi} and \cite[p. 47]{Ha})  $mod A$ is equivalent to the category of contravariant representations of $Q$ (see \cite[p.47,48]{Ha}). The category of representations $Rep_{\KK}(Q)$, which we consider, is the same as the category of  contravariant representations of the quivers $Q'$ which are obtained from $Q$ by inverting the arrows. So in fact, our category $\mc T = D^b(Q)$ is equivalent to $D^b({\rm mod} A)$ from  \cite{Miyachi}, where $A$ is the path algebra as defined in \cite{Miyachi} for the quiver $Q'$ dual to $Q$. However $Q_1'$, $Q_2'$ are isomorphic to $Q_1$, $Q_2$. 

Before  \cite[Definitions 2.2]{Miyachi} is explained that  the Auslander Reiten quiver of $\mc T=D^b(Q)$ contains the quiver $Q'$ (and hence the quiver $Q$ which is isomorphic to $Q'$) as the full  subquiver with vertices corresponding to the indecomposable projective representations,
under the inclusion $ Rep_\KK(Q) \subset D^b(Q)$.

In \cite[Definitions 2.2, 2.3]{Miyachi}  they   define a quiver, denoted by ${\bf \Gamma}^{irr}$,  which is the disjoint union of all connected components of the  Auslander-Reiten quiver of $D^b(Q)$, which are  isomorphic to the connected component  containing the indecomposable projective modules.     \cite[Theorem 2.4 ]{Miyachi}  gives an isomorphism of quivers $\rho: \ZZ\times (\ZZ Q)\rightarrow {\bf \Gamma}^{irr}$\footnote{this Theorem  gives an isomorphism from $\ZZ\times (\ZZ Q'))$ to  ${\bf \Gamma}^{irr}$. however $\ZZ\times (\ZZ Q')$ is isomorphic to  $\ZZ\times (\ZZ Q)$}, where:  

(a) $\ZZ Q$ is  isomorphic to quiver with set of vertices $\ZZ\times V(Q)$ and for  every arrow $x \rightarrow y$ in $Q$ there are arrows $(n,x) \rightarrow (n,y)$ and $(n,y) \rightarrow (n+1,x)$ in $\ZZ Q$.  In particular,  the subquiver with vertices with zero first coordinate is isomorphic to $Q$.

(b) $\ZZ \times (\ZZ Q) $ is quiver whose connected components are $\{i\}\times \left ( \ZZ Q \right )$ and each of them is a labeled copy of $\ZZ Q$. 

(c) The identification $\rho: \ZZ\times (\ZZ Q)\rightarrow {\bf \Gamma}^{irr} $ is such that:

(c.1) Its restriction  to the full  subquiver of $\ZZ\times (\ZZ Q)$ with vertices whose first and second coordinates vanish give isomorphism with the full subquiver  of $ {\bf \Gamma}^{irr}$ with vertices the indecomposable projective representations.

(c.2)  The translation functor of $\mc T$ induces action on ${\bf \Gamma}^{irr}$ which by $\rho$ corresponds to an automorphism $\sigma$ of $\ZZ \times (\ZZ Q)$ acting by increasing the first component by $1$. 

(c.3)   The Serre functor of $\mc T$ shifted by $[-1]$ induces an automorphism  $\tau $ of $\ZZ \times (\ZZ Q)$, which on vertices maps $(a,(i,j))$ to $(a,(i-1,j))$. 
\subsection{On the   Serre functor and the exceptional objects}

Recall that throughout the paper we use the notation $A\sim B$, where $A$, $B$ are  exceptional collections, explained in Section \ref{notations}.

Let $\mc T = \mc T_i$ with $i=1$ or $2$. We recall first that   every exceptional collection  in $\mc T$ can be extended to a full exceptional collection in $\mc T$ (this is proved in \cite{WCB1}).

If $\mc E=(E_1,E_2.\dots,E_n)$ is a full exceptional collection in $\mc T$ then after appropriate number of right mutations of $E_1$ we obtain another full exceptional collection    $(E_2.\dots,E_{n},E_{n+1})$ (the first $n-1$ elements are the last $n-1$ elements of  $\mc E$), analogously enough number of left mutations of $E_n$ result a full exceptional collection $(E_0,E_1.\dots,E_{n-1})$, by induction one obtains an infinite sequence in both directions $\{E_i\}_{i\in \ZZ}$ - the helix induced by  $\mc E=(E_1,E_2.\dots,E_n)$, as defined in \cite[p. 222]{BP}. A property of such a helix is that  $\Phi(E_{i+n})\cong E_i$ for any $i$, where $\Phi=S[1-n]$ (see \cite[p. 223]{BP})  and $S$ is the Serre functor. Taking into account that if $\mc E$ and  $\mc E'$ are two full exceptional collections such that $\mc E(i)\sim \mc E'(i)$ for all but one $i$, say $k$, then it follows that  $\mc E(k)\sim \mc E'(k)$ as well, we deduce that: 
\begin{gather} 
\mbox{let $(E_1,E_2.\dots,E_n)$ and $(E_2.\dots,E_{n},E_{n+1})$ be full exceptional collections in $\mc T$,} \nonumber \\[-3mm] \label{Serre}\\[-3mm] \mbox{ then  $S(E_{n+1})\sim  E_1$, where $S$ is the Serre functor  of $\mc T$}. \nonumber
\end{gather}
	
From	\cite[Corollary 2.6]{DK1}  and since the homological dimension of $D^b(Q_i)$ is one, it follows that:
\begin{gather} \mbox{For any two exceptional objects}   \ X, Y \in \mc T_i \quad i=1,2   \nonumber \\[-3mm] \label{at most one non-vanisching}  \\[-3 mm] \mbox{ at most one  element in}  \ \{ \hom^p(X,Y) \}_{p\in \ZZ} \ \mbox{is non-zero.}\nonumber \end{gather}

\begin{remark} \label{reminder}
	 In \cite{DK4} was shown that (see \cite[(7) and (8)]{DK4}):
	 \begin{gather}
	 \label{vanishings general}
	 \mbox{If} \  l\geq 2 \ \mbox{and $Q$ is affine aciclic quiver}  \ \Rightarrow \abs{C_{l}^{\{\rm Id\}}(D^b(Q)) }=0, \\
	 \label{vanishings general 2}  \mbox{If} \  l\geq 1 \ \mbox{and $Q$ is Dynkin  quiver}  \ \Rightarrow \abs{C_{l}^{\{\rm Id\}}(D^b(Q)) }=0.  \end{gather} 
\end{remark}

\begin{remark} \label{subcats generated by pairs} Note that from \eqref{at most one non-vanisching} and Proposition \ref{bijection 123} it follows that if $C_{\mc A, \KK}^{\{\rm Id \}}(D^b(Q_i)) \neq \emptyset$ and $\mc A$ is a category generated by an exceptional pair, then  $C_{\mc A, \KK}^{\{\rm Id \}}(\mc T_i)$ must be  $C_{l}^{\{\rm Id \}}(\mc T_i)$ for some $l\geq -1$.  Furthermore, \eqref{vanishings general} implies that $l\leq 1$.
	\end{remark}

\subsection{	The case $Q_1$.} \label{The case Q1}

Here we derive the numbers in the first table  of \eqref{table with numbers Q1}.

Using the bijection \eqref{bijection 1}  we will identify $C_l^{\{\rm Id\}}(\mc T_1)$ with the codomain of the bijection, and we will denote the codomain by  $C_l^{\{\rm Id\}}(\mc T_1)$ as well ($l\in \ZZ_{\geq -1}$.) 

We will need  information about  exceptional collections in $\mc T_1$ obtained in 
 \cite{DK1} and  \cite{DK3}. In    \cite{DK3}  are listed  the full  exceptional collections in $\mc T_1$ up to isomorphism and shifts, they are:	
		\begin{gather} \label{mc T}  {\mk T} = \left \{  \begin{array}{c  c  c} (M',a^{m},a^{m+1})  & (a^m, b^{m+1},a^{m+1}) & (a^m,a^{m+1},M) \\
			(M,b^{m},b^{m+1}) &(b^m,a^m,b^{m+1})  & (b^m,b^{m+1},M')\\
			(b^{m},M',a^{m}) &(a^m , M, b^{m+1}) &  \end{array}: m\in \NN \right \},
		\end{gather}
where $\{a^m\}_{m\in \ZZ}$,  $\{b^m\}_{m\in \ZZ}$, $\{M,M'\}$ are all the exceptional objects in $\mc T_1$ up to shift organized in a convenient way (see Section 3.2 in  \cite{DK3} for details).

From \cite[formula (70)]{DK3} and \cite[table (4) in Proposiotion 2.4]{DK1}	 we see that: 
\begin{gather} \mbox{\textit{Let}  $E_1,E_2$ \textit{ be an exc. pair in} $\mc T_1$, \textit{then there exists} $p\in \ZZ$  \textit{ s. t.} $\hom^p(E_1,E_2)\neq 0$  \textit{and} (\textit{see} \eqref{at most one non-vanisching}): } \nonumber\\
\label{options} \mbox{\textit{either} $\hom^p(E_1,E_2)=2$ \textit{and}} \ \ (E_1,E_2)\sim X, \ X \in  \left \{ (a^m,a^{m+1}), (b^m,b^{m+1}): m\in \ZZ  \right \}  \\ 
\mbox{\textit{or} $\hom^p(E_1,E_2)=1$ \textit{and}} \ \ (E_1,E_2)\sim X, \ X \in    \left \{ \begin{array}{c}  (M',a^m), (a^m,b^{m+1}), (b^m,a^{m}), \\ (a^m,M), (M,b^m), (b^m,M') \end{array} : m\in \ZZ \right \}. \nonumber  \end{gather}

The zeros in the row with $l=-1$ in the first table of \ref{table with numbers Q1} are due to the lack of orthogonal pairs in $\mc T_1$ (see \eqref{options}). Indeed,   	if $\mc A\subset \mc T$, $\mc A \in C_{-1}^{\{\rm Id\}}(\mc T_1)$, then from bijection \eqref{bijection 1}  there is a strong    exc. pair $(E_1,E_2)$ with  $\langle E_1,E_2 \rangle =\mc A $ and $\hom(E_1,E_2)=0$, thus $(E_1,E_2)$ must be  orthogonal, therefore   $C_{-1}^{\{\rm Id\}}(\mc T_1)=\emptyset$.

	Now we can draw the graph  $G_{D^b(A_1), \KK}(\mc T_1)$. Using the explanations related to \eqref{mc T}  we deduce that  set of vertices is $\{\langle  M \rangle , \langle  M' \rangle \} \cup \{ \langle  a^m \rangle, \langle  b^{m} \rangle : m \in \ZZ   \}$. All the exceptional triples up to shift are listed in \eqref{mc T}, from this list one sees also all the exceptional pairs up to shift. We consider first  the subgraph of  containing $\langle  M \rangle , \langle  M' \rangle $, $ \langle  a^m \rangle, \langle  b^{m} \rangle$, $\langle  a^{m+1} \rangle, \langle  b^{m+1} \rangle$ for some $m\in \ZZ$. From   \eqref{mc T} we see that  all the exceptional pairs consisting of some of these vertices are $(M', a^m)$, $(M', a^{m+1})$, $( b^m, M')$, $(b^{m+1}, M')$, $(b^m,b^{m+1})$, $(a^m,a^{m+1})$,  $(a^m,b^{m+1})$,  $(b^m,a^{m})$, $(b^{m+1},a^{m+1})$, $(M, b^m)$, $(M, b^{m+1})$, $( a^m, M)$, $(a^{m+1}, M)$ and this corresponds to 13 edges, which are the middle part in Figure \ref{Figure for digraph of Q1}. 	The entire digraph $G_{D^b(A_1), \KK}(\mc T_1)$ consists of infinitely many figures as this, indexed by $m$ as demonstrated in Figure \ref{Figure for digraph of Q1}.
		\begin{figure}

		\begin{tikzpicture}[scale=1.5] 
		\node [left] at (0,0) {$\langle a^m \rangle $};
		\draw[->]  (0,0)--(1,0);
		\node [right] at (1,0) {$\langle a^{m+1} \rangle $};
		\draw[<-]  (-0.4,0.2)--(-0.7,1);	
		\node [above] at (-0.8,1) {$\langle b^{m} \rangle $};
		
		\draw[->]  (-0.55,1.2)--(0.45,1.2);	
		\node [right] at (0.45,1.2) {$\langle b^{m+1} \rangle $};
		\draw[<-]  (1.3,0.2)--(1,1);
		\draw[<-]  (0.6,1)--(-0.1,0.2);		
		\node at (0.2,3) {$\langle M \rangle $};	
		\draw[->]  (0.1,2.9)--(-0.8,1.4);	
		\draw[->]  (0.3,2.9)--(0.8,1.4);
		\draw[->]  (-0.2,0.2)--(0.15,2.9);	
		\draw[->]  (1.1,0.2)--(0.2,2.9);
		\node at (0.2,-2.5) {$\langle M' \rangle $};
		\draw[->]  (0.1,-2.3)--(-0.2,-0.2);	
		
		\draw[->]  (0.3,-2.3)--(1.1,-0.2);	
		\draw[->]  (-0.8,1)--(0.05,-2.3);	
		\draw[->]  (0.7,1)--(0.2,-2.3);		
		
		\node [left] at (-1.5,0) {$\langle a^{m-1} \rangle $};
		\draw[->]  (-1.5,0)--(-0.7,0);
	
		\draw[<-]  (-1.9,0.2)--(-2.2,1);	
		\node [above] at (-2.3,1) {$\langle b^{m-1} \rangle $};
		
		\draw[->]  (-1.9,1.2)--(-1.1,1.2);

		\draw[<-]  (-0.9,1)--(-1.6,0.2);				
		\draw[->]  (-0.1,-2.3)--(-1.7,-0.2);
		\draw[->]  (-2.1,1)--(-0.05,-2.3);
		\draw[->]  (-0.1,2.9)--(-2.1,1.4);	
		\draw[->]  (-1.7,0.2)--(0,2.9);

		\draw[->]  (1.9,0)--(2.7,0);
		\node [right] at (2.7,0) {$\langle a^{m+2} \rangle $};

		\draw[->]  (1.3,1.2)--(2.15,1.2);	
		\node [right] at (2.15,1.2) {$\langle b^{m+2} \rangle $};
		\draw[<-]  (3,0.2)--(2.7,1);
		\draw[<-]  (2.3,1)--(1.6,0.2);				
		
		\draw[->]  (0.4,-2.3)--(2.7,-0.2);
		
		\draw[->]  (2.4,1)--(0.35,-2.3);
		\draw[->]  (0.5,2.9)--(2.3,1.4);
		\draw[->]  (2.6,0.2)--(0.4,2.9);	
		\node  at (-2.6,0.6) {$\cdots$};	
		\node  at (-2.6,2.3) {$\cdots$};
		\node  at (-2.6,-1) {$\cdots$};											
		\node  at (-2.4,0) {$\cdots$};		
		\node  at (-2.8,1.2) {$\cdots$};	
		\node  at (3.4,0.6) {$\cdots$};										
		\node  at (3.6,0) {$\cdots$};		
		\node  at (3.2,1.2) {$\cdots$};		
		\node  at (3.4,2.3) {$\cdots$};
		\node  at (3.4,-1) {$\cdots$};	
		\node  [above] at  (0.5,0) {$2$}; 	
		\node  [above] at  (2.2,0) {$2$}; 		
		\node  [above] at  (-1,0) {$2$}; 		
		\node  [above] at  (0.1,1.2) {$2$}; 	
		\node  [above] at  (1.8,1.2) {$2$}; 		
		\node  [above] at  (-1.4,1.2) {$2$}; 																						
		\end{tikzpicture} 				
		\caption{$G_{D^b(A_1), \KK}\left (D^b(Q_1)\right )$}\label{Figure for digraph of Q1}
		
	\end{figure}
		The two dimensional simplexes in  $SC_{\mc J, \KK}(\mc T)$ are the triples of vertices in triangles  as the quiver $Q_1$. 
	Such triangles correspond to exceptional triples in $\mc T$. All of these triangles are seen in  Figure \ref{Figure for digraph of Q1} and they intersect with each other  in the proper sides.   So in this figure is presented also a geometric realization of the  complete simplicial complex $SC_{\mc J, \KK}(\mc T)$.
	
	In Subsection \ref{the case J eq DbA1} is explained how to make the digraph a valued graph. In this case using \eqref{options} we see that to the edges from $a^m$ to $a^{m+1}$ and from $b^m$ to $b^{m+1}$ is attached $2$ for each $m\in \ZZ$, and to all other edges is attached $1$.

Next we derive the numbers  the column of the first table in \eqref{table with numbers Q1} with  $\Gamma = \{\rm Id\}$. 	Let $\mc A\subset \mc T$, $\mc A \in C_1^{\{\rm Id\}}(\mc T_1)$. Therefore we have a  strong exceptional pair $(E_1,E_2)$, s.t. $\hom(E_1,E_2)=2$,  $\langle E_1,E_2 \rangle =\mc A $. From \eqref{at most one non-vanisching}, \eqref{options} we see that $\mc A = \langle a^m,a^{m+1}\rangle$ or $\mc A = \langle b^m,b^{m+1}\rangle$ for some $m\in \ZZ$, now \eqref{mc T} implies that  $\mc A = \langle M\rangle^{\perp}$ or  $\mc A = \langle M'\rangle^{\perp}$, furthermore $ \langle M'\rangle^{\perp}\neq \langle M\rangle^{\perp}$
 since for example $\hom^1(M,a^m)\neq 0$(see \cite[formula (72)]{DK3}), thus we see that for any $m\in \ZZ$:
 \begin{gather}  \label{C_1^{{rm Id}}(mc T_1)} C_1^{\{\rm Id\}}(\mc T_1) =\{\langle M\rangle^{\perp}=\langle a^m,a^{m+1}\rangle, \quad   \langle M'\rangle^{\perp}=\langle b^m,b^{m+1}\rangle\} \qquad \abs{C_1^{\{\rm Id\}}(\mc T_1)}=2. \end{gather} 
  	Let $\mc A\subset \mc T$, $\mc A \in C_0^{\{\rm Id\}}(\mc T_1)$. Therefore we have a  strong exceptional pair $(E_1,E_2)$, s.t. $\hom(E_1,E_2)=1$ and   $\langle E_1,E_2 \rangle =\mc A $. From \eqref{at most one non-vanisching}, \eqref{options} we see that $\mc A = \langle X\rangle$ for some of the exceptional pairs in the third row of \eqref{options}. Sine the triples in   \eqref{mc T} are full exceptional, it follows: 
 	\begin{gather}  \langle M',a^{m}\rangle = \langle a^{m}, b^{m+1}\rangle=\langle b^{m+1}, M'\rangle =\langle a^{m+1}\rangle^{\perp} \\ \langle M,b^{m}\rangle = \langle b^{m}, a^{m}\rangle=\langle a^{m}, M\rangle =\langle b^{m+1}\rangle^{\perp}. 
 	\end{gather}
 	Thus, we see that if  $\mc A\subset \mc T$, $\mc A \in C_0^{\{\rm Id\}}(\mc T_1)$, then $\mc A=\langle a^m\rangle^{\perp}$ or $\mc A=\langle b^m\rangle^{\perp}$ for some $m \in \ZZ$. On the other hand for  $i\neq j$ we have $\langle a^i\rangle^{\perp}\neq \langle a^j\rangle^{\perp}$ , since $(a^i,a^{j+1})$ and $(a^j,a^{i+1})$ cannot be both exceptional pairs (see \cite[Corollary 3.7]{DK3}), and similarly $\langle b^i\rangle^{\perp}\neq \langle b^j\rangle^{\perp}$  for  $i\neq j$. Similar  arguments using \cite[Corollary 3.7]{DK3} show that  for any $i,j \in \ZZ$  we have  $\langle a^i\rangle^{\perp}\neq \langle b^j\rangle^{\perp}$. So we deduce: 
  	\begin{gather}  \label{C_0^{{rm Id}}(mc T_1)}  C_0^{\{\rm Id\}}(\mc T_1) =\{\langle a^m\rangle^{\perp}:m\in \ZZ \}\cup\{\langle b^m\rangle^{\perp}:m\in \ZZ\} \ \ \mbox{disjoint union, injective sequences}. \end{gather} 
  	In particular, $  \ \ \abs{C_0^{\{\rm Id\}}(\mc T_1)}=\infty$. 
  	
  	To compute $C_0^{\Gamma}(\mc T_1)$ and $C_1^{\Gamma}(\mc T_1)$ for non-trivial $\Gamma$  we  use Corollary \ref{group action}.
  	
  	 We   determine first the action of $\langle S \rangle$ on the sets  $C_0^{\{\rm Id\}}(\mc T_1)$ and $C_1^{\{\rm Id\}}(\mc T_1)$.
  	In \eqref{mc T} we see that there are full exceptional collections  $(b^m,M',a^m)$, $(M',a^m,a^{m+1})$ and applying \eqref{Serre} we deduce that $S(a^{m+1})\sim b^m$. Similarly the pairs of full exceptional collections  $(a^{m-1},a^{m}, M)$, $(a^{m}, M, b^{m+1})$ and  $(M',a^m,a^{m+1})$,  $(a^m,a^{m+1}, M)$ and $(M,b^{m},b^{m+1})$, $(b^m,b^{m+1},M')$ imply the rest incidences here: 
  	\begin{gather} \label{Serre in Q_1}
  	S(a^{m+1})\sim b^m \qquad  S(b^{m+1})\sim a^{m-1} \qquad   S(M)\sim M' \qquad S(M')\sim M \qquad m\in \ZZ.
  	\end{gather}

  	Therefore  $S\left (\langle M'\rangle^{\perp}\right )= \langle M\rangle^{\perp}$ and from \eqref{C_1^{{rm Id}}(mc T_1)} we see that the two elements of $C_1^{\{\rm Id\}}(\mc T_1)$ project to  a single element  in $C_1^{\{\rm Id\}}(\mc T_1)/\langle S \rangle \cong C_1^{\langle S \rangle}(\mc T_1)$, hence $\abs{C_1^{\langle S \rangle}(\mc T_1)}=1$. On the other hand from  \eqref{Serre in Q_1} follow the formulas  $S( \langle b^m\rangle^{\perp})\sim  \langle a^{m-2}\rangle^{\perp}$ and $S^2( \langle b^m\rangle^{\perp})\sim  \langle b^{m-3}\rangle^{\perp}$, and it  follows that  the two sequences of $C_0^{\{\rm Id\}}(\mc T_1)$ project to  three equivalence classes  $[\langle b^0\rangle^{\perp}], [\langle b^1\rangle^{\perp}], [\langle b^2\rangle^{\perp}]$  in $C_0^{\{\rm Id\}}(\mc T_1)/\langle S \rangle \cong C_0^{\langle S \rangle}(\mc T_1)$, furthermore using \eqref{C_0^{{rm Id}}(mc T_1)} and these formulas one easily proves that if for \\ $x,y\in \{ \langle b^0\rangle^{\perp}, \langle b^1\rangle^{\perp}, \langle b^2\rangle^{\perp} \}$, $i\in \ZZ$ we have $S^i(x)=y$ then $i=0$, therefore    $\abs{C_0^{\langle S \rangle}(\mc T_1) }=3$. 
  	
  	Since   $C_0^{\{\rm Id\}}(\mc T_1)/{\rm Aut}_{\KK}(\mc T_1) \cong C_0^{{\rm Aut}_{\KK}(\mc T_1)}(\mc T_1)$, in order to completely fill in  the first  table in \eqref{table with numbers Q1}, it is enough to show that there  exists  $\KK$-linear auto-equivalence $u: \mc T_2 \cong \mc T_2$, such that  $u(b^{m})\cong b^{m-1} $ for each $m\in \ZZ$. Indeed, such $u$  will ensure that  the three elements in $C_0^{\{\rm Id\}}(\mc T_1)/\langle S \rangle $  project to a single element in  $C_0^{\{\rm Id\}}(\mc T_1)/{\rm Aut}_{\KK}(\mc T_1) $. We find such an  auto-equivalence  in the following:
  	
    \begin{prop}
  	There exist $\KK$-linear autoequivalences $\zeta$ of $ \mc T_1$, such that for each $m\in \ZZ$:
  	\begin{gather} \label{autoeq zeta formula}
  	\zeta(a^m) \cong b^m \qquad \zeta(b^m) \cong a^{m-1}. 
  	\end{gather} In particular $\zeta^2$ satisfy $\zeta^2(b^m)\cong b^{m-1}$ for each $m\in \ZZ$. 
  \end{prop}
  \bpr
  
  In \cite[Theorem 4.5 (3)]{Miyachi} is proved that ${\rm Aut}_{\KK}(\mc T_1) \cong \ZZ \times \left (\ZZ \ltimes \left [\begin{array}{c c} \KK^{\times} & \KK \\ 0 & 1\end{array} \right ] \right )$, where one $\ZZ$-factor is generated by the translation functor, and we will show that a generator of the other $\ZZ$-factor is the desired auto-equivalence  $\zeta$.

  The Auslander Reiten quiver of $\mc T_1=D^b(Q_1)$ contains the quiver $Q_1$  as explained in the beginning of  Subsection \ref{material from Mi}.  
  We will show that this subquiver is  
  \begin{gather} \label{subquiver of ARQ Q_1}
      \begin{diagram}[1em]
  &       &   a^1[-1] &       &    \\
  & \ruTo &    & \rdTo &       \\
  b^1[-1]  & \rTo  &    &       &  b^2[-1]
  \end{diagram}  
  \end{gather}

  We  show that the indecomposable projective representations are in the exceptional collection $(b^1[-1],a^1[-1],b^2[-1])$ (see \eqref{mc T}). In terms of the notations of  \cite[Proposition 2.2]{DK1}  the simple representations  are $ E_1^0, M,  E_3^0$, which with the notations, which we use here   (following \cite[(70)]{DK3}) are $(a^0,M,b^1[-1])$.  From \eqref{at most one non-vanisching}  and \cite[Corollary 3.7]{DK3} we know that:
  
 $ \forall X\in \{b^1[-1],a^1[-1],b^2[-1]\} \ \ \forall Y \in  \{a^0,M,b^1[-1]\} \quad \hom^1(X,Y)=0   $ therefore: \\
 $ \forall X\in \{b^1[-1],a^1[-1],b^2[-1]\}$ $ \forall Y \in Rep_\KK(Q_1) \quad \hom^1(X,Y)=0. \nonumber  $
  Furthermore, from \eqref{at most one non-vanisching}, \cite[Corollary 3.7]{DK3}, \eqref{options}  we see that $ \hom(b^1,a^1)= \hom(a^1,b^2)=1 $, $ \hom(b^1,b^2)=2$.  These arguments imply that  the full subquiver of the Auslander-Reiten quiver of $D^b(Q_1)$  with vertices  the indecomposable projective representations is \eqref{subquiver of ARQ Q_1}.

  Next, we will describe the connected component of the quiver  ${\bf \Gamma}^{irr}$ discussed in Subsection  \ref{material from Mi} which contains $b^1,a^1,b^2$. An identification $\rho: \ZZ\times (\ZZ Q_1)\rightarrow {\bf \Gamma}^{irr}$is explained in  Subsection \ref{material from Mi} (c).
  
  The full subquiver containing  $b^1,a^1,b^2$  is just shift of \eqref{subquiver of ARQ Q_1} and via the identification $\rho$ it corresponds to a subquiver in $\{1\}\times \ZZ Q_1$ (see (c.2) in  Subsection \ref{material from Mi}), and the connected component containing  $b^1,a^1,b^2$ is  $\rho\left ( \{1\}\times \ZZ Q_1 \right )$.   Using (c.3) in Subsection \ref{material from Mi} we obtain the desired connected component by acting on a shift of \eqref{subquiver of ARQ Q_1} with $\langle [-1] S \rangle$.  We first make the formulas  \eqref{Serre in Q_1} more precise
  
  \begin{gather} \label{Serre in Q_1 prec}
  S(a^{m+1})\cong b^m[1] \qquad  S(b^{m+1})\cong a^{m-1}[1]
  \qquad  m\in \ZZ.
  \end{gather} 
  Indeed,  using the property of the Serre Functor we  can write $1=\hom(a^{m+1}, a^{m+1})=$ \\ $\hom(a^{m+1}, S(a^{m+1}))=\hom(a^{m+1}, b^{m}[p])$, where $S(a^{m+1})\cong b^{m}[p]$ due to  $ S(a^{m+1}) \sim b^{m}$ in \eqref{Serre in Q_1}, and now from $\hom^1(a^{m+1}, b^{m})\neq 0$ in  \cite[Corollary 3.7]{DK3} and  \eqref{at most one non-vanisching} it follows $p=1$.  Similarly one derives the other  incidence in  \eqref{Serre in Q_1 prec}.

\begin{figure}
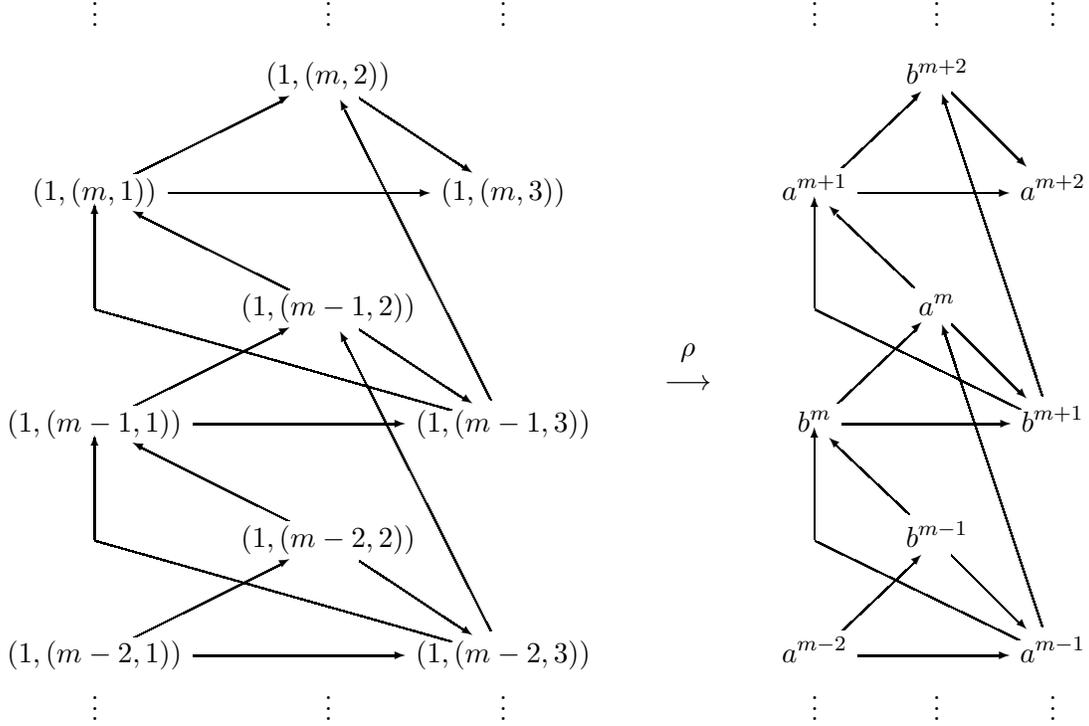


  \begin{gather} \nonumber
  \begin{diagram} \vdots &      & \vdots  &     &  \vdots\\
             &                  &  (1,(m,2))   &                   &    \\
             & \ruTo            &              & \rdTo  \luTo(2,6) &     \\
  (1,(m,1))  & \rTo             &              &                   &  (1,(m,3)) \\
     \uTo[abut]    & \luTo            &              &                   &            \\
   \HonV     &                  &  (1,(m-1,2)) &                   &    \\
             & \ruTo \luLine(4,2)[abut]    &              & \rdTo  \luTo(2,6)            &       \\
(1,(m-1,1))  & \rTo             &              &                   &  (1,(m-1,3)) \\
   \uTo[abut]  &   \luTo          &              &                   &            \\ 
    \HonV     &                  &  (1,(m-2,2)) &                   &    \\
             & \ruTo \luLine(4,2)[abut] &              & \rdTo             &       \\
(1,(m-2,1))  & \rTo             &              &                   &  (1,(m-2,3)) \\
  \vdots &      & \vdots  &     &  \vdots
    \end{diagram} \quad \qquad\qquad \begin{array}{c} \rho \\ \longrightarrow \end{array}  \qquad
   \begin{diagram}[]\vdots &      & \vdots  &     &  \vdots\\
            &                  &   b^{m+2}     &                   &    \\
            & \ruTo            &               & \rdTo  \luTo(2,6) &     \\
  a^{m+1}   & \rTo             &               &                   &  a^{m+2}   \\
    \uTo[abut] & \luTo  &              &                   &            \\
   \HonV   &                  &  a^m &                   &    \\
   & \ruTo   \luLine(4,2)[abut]         &              & \rdTo   \luTo(2,6)           &       \\
   b^m  & \rTo             &              &                   &  b^{m+1} \\
    \uTo[abut]  &   \luTo          &              &                   &            \\ 
   \HonV   &                  &   b^{m-1} &                   &    \\
   & \ruTo  \luLine(4,2)[abut]          &              & \rdTo             &       \\
   a^{m-2}  & \rTo             &              &                   &   a^{m-1} \\
  \vdots &      & \vdots  &     &  \vdots\\
   \end{diagram} 
  \end{gather} \caption{The identification $\rho$ for $Q_2$, $m$ is odd}  \label{extend subquiver of ARQ Q_1}
  \end{figure}

 We shift the subquiver \eqref{subquiver of ARQ Q_1} of ${\bf \Gamma}^{irr}$ and by acting on it with $\langle [-1] S \rangle$ using \eqref{Serre in Q_1 prec}  we see that the identification $\rho$ for odd  $m\in \ZZ$ can be depicted locally as shown in Figure \ref{extend subquiver of ARQ Q_1},   where although in the picture some arrows look intersecting, they actually do not intersect.

  Looking in the proof of \cite[Theorem 4.5 (2)]{Miyachi} and correcting a typo there\footnote{Note noe that in our case $q=1$, $p=2$ and in the proof of \cite[Theorem 4.5 (2)]{Miyachi} is written $\rho(m,1)=(m-1,p)$, which we think must be $\rho(m,1)=(m-1,p+q)$}  we see that  a generator of the other $\ZZ$-factor in ${\rm Aut}_{\KK}(\mc T_1)$  acts on the left hand side in Figure \ref{extend subquiver of ARQ Q_1}  as follows (for any $n\in \ZZ$)
  \begin{gather}  (1,(n,1)) \mapsto (1,(n-1,3)) \quad  (1,(n,2)) \mapsto (1,(n,1))  \quad (1,(n,3)) \mapsto (1,(n,2))
  \end{gather}
   Note that we have observed the rule for labeling of the  the vertices   in  \cite[Figure 6]{Miyachi}.   Via the correspondence in Figure \ref{extend subquiver of ARQ Q_1} we get an auto-equivalence $\zeta$ of $\mc T_1$ such that the first two entities in \eqref{autoeq zeta formula} hold for any $m\in \ZZ$ (not only odd). 
     \epr

  \subsection{	The case $Q_2$.} \label{The case Q2}	Here we derive the numbers in the second table in \eqref{table with numbers Q1}. The  approach  is  as in subsection \ref{The case Q1}. 
  
  The exceptional objects of $\mc T_2$ are classified in \cite[Proposition 2.3]{DK1} and they are organised in $8$ series $\{E_1^m\}_{m\in \ZZ}$,  $\{E_2^m\}_{m\in \ZZ}$, $\dots$,  $\{E_8^m\}_{m\in \ZZ}$, and two Ext-nontrivial couples $\{F_+,G_-\}$, $\{F_-,G_+\}$ (all are objects in $Rep_\KK(Q_2)$). Furthermore in \cite[Proposition 2.5]{DK1} are computed  $\hom(X,Y)$, $\hom^1(X,Y)$ for any pair of exceptional representations. It is convenient to group the eight series in four series as follows:
    \begin{gather}  a^m=\left \{  \begin{array}{c c} E_1^{-m}  & m\leq 0 \\  E_2^{m-1}[1]  & m\geq 1 \end{array} \right. ; \qquad  b^m= \left \{ \begin{array}{c c} E_4^{-m}  & m\leq 0 \\  E_3^{m-1}[1]  & m\geq 1  \end{array} \right . ; \nonumber \\[-2mm] \label{abcd}  \\[-2mm] c^m= \left \{ \begin{array}{c c} E_6^{-m}  & m\leq 0 \\  E_5^{m-1}[1]  & m\geq 1  \end{array} \right . ;  \qquad d^m= \left \{ \begin{array}{c c} E_7^{-m}  & m\leq 0 \\  E_8^{m-1}[1]  & m\geq 1  \end{array} \right  .. \nonumber \end{gather}
   And then using  \cite[Proposition 2.5]{DK1} one proves: 
   \begin{lemma} \label{exceptional pairs in T2} For any $x\in \{a,b,c,d\}$, $m\in \ZZ$ the pair $(x^m,x^{m+1})$ is  strong exceptional  with $\hom(x^m,x^{m+1})=2$
  	
  	Let  $E_1,E_2$ be any exceptional pair in $\mc T_2$,  there is a trichotomy
  	
  	{\rm (a)}  $\hom^p(E_1,E_2)=2$ for some $p\in \ZZ$ which is equivalent to  $(E_1,E_2)\sim (x^m,x^{m+1})$ for some $x\in \{a,b,c,d\}$ and some $m\in \ZZ$;
  	
  		{\rm (b)}  $\hom^p(E_1,E_2)=0$ for all $p\in \ZZ$ which is equivalent to  $(E_1,E_2)\sim X$ or  $(E_1,E_2)\sim \ol{X}$ for some $X\in \{(a^m,b^{m+1}), (c^m,d^m) :m\in \ZZ\} \cup \{ (F_+,F_-), (G_+,G_-), (F_+,G_+), (F_-,G_-) \}$, where $\ol{(A,B)}=(B,A)$;
  		
  			{\rm (c)}  $\hom^p(E_1,E_2)=1$ for some $p\in \ZZ$ which is equivalent to  $(E_1,E_2)\sim X$, where $X$ is some of the pairs in the following table for some $m\in \ZZ$:
  			\begin{gather} \label{1 pairs}\begin{array}{| c | c | c |}
  			\hline
  		(a^{m-1},c^m)	&        (c^m,G_-)        &        (G_-,a^{m-1})           \\ \hline
  			 (c^m,a^m)  	&        (a^m,F_+)        &       (F_+,c^{m})          \\ \hline
  			  (a^{m-1},d^m)  	&        (d^m,G_+)        &       (G_+,a^{m-1})          \\ \hline
  			  (d^{m},a^m)	&        (a^m,F_-)        &        (F_-,d^{m})           \\ \hline
  			   (b^{m},c^m)	&        (c^m,F_-)        &        (F_-,b^{m})           \\ \hline
  			    (c^{m-1},b^m)	&        (b^m,G_+)        &        (G_+,c^{m-1})           \\ \hline
  			     (b^{m},d^m)	&        (d^m,F_+)        &        (F_+,b^{m})           \\ \hline
  			     (d^{m-1},b^m)	&        (b^m,G_-)        &        (G_-,d^{m-1})           \\ \hline
  			\end{array}
  			\end{gather}
  	Furthermore, for any pair $(A,B)$ taken from the first column and for aby $m\in \ZZ$ the only $p$ such that $\hom^p(A,B)\neq 0$ is  $p=0$.
  	\end{lemma}

  Now we can discuss the graph $G_{D^b(A_1), \KK}(\mc T_2)$. 
  		Using the explanations related to \eqref{abcd}  we deduce that the  set of vertices of $G_{D^b(A_1), \KK}(\mc T)$ is $\{\langle  F_+ \rangle , \langle  F_- \rangle, \langle  G_+ \rangle , \langle  G_- \rangle \}$ \\ $\cup \{ \langle  a^m \rangle, \langle  b^{m} \rangle, \langle  c^{m} \rangle, \langle  d^{m} \rangle   : m \in \ZZ   \}$. All the exceptional pairs  up to shift are listed in Lemma  \ref{exceptional pairs in T2}. Let us draw the full subgraph of $G_{D^b(A_1), \KK}(\mc T_2)$ with vertices $\langle  a^m \rangle, \langle  b^{m} \rangle, \langle  c^{m} \rangle, \langle  d^{m} \rangle$,  $\langle  a^{m+1} \rangle, \langle  b^{m+1} \rangle, \langle  c^{m+1} \rangle$,\\ $ \langle  d^{m+1} \rangle$. Using Lemma \ref{exceptional pairs in T2} (a), (b), (c) we draw the  cube in Figure \ref{cube for Q2} presenting the  subgraph, to the  arrows from $\langle x^m \rangle$ to  $\langle x^{m+1} \rangle$ is attached $2$, to the rest one-sided arrows is attached $1$.

  		\begin{figure}
  		\begin{tikzpicture}[scale=1.5]

  		\node [left] at (0,0) {$\langle d^m \rangle $};
  		\draw[->]  (0,0)--(1,0);
  		\node [right] at (1,0) {$\langle d^{m+1} \rangle $};
  		\node [left] at (0.7,1) {$\langle b^m \rangle $};
  		\draw[->]  (0.7,1)--(1.7, 1);
  		\node [right] at (1.7, 1) {$\langle b^{m+1} \rangle $};
  		\draw[->]  (0.2,0.8)--(-0.2, 0.2);	
  		\draw[->]  (2,0.8)--(1.7, 0.2);

  		\node [left] at (0,1.4) {$\langle a^m \rangle $};
  		\draw[->]  (-0.05,1.4)--(1,1.4);
  		\node [right] at (1,1.4) {$\langle a^{m+1} \rangle $};
  		\node [left] at (0.7,2.4) {$\langle c^m \rangle $};
  		\draw[->]  (0.7,2.4)--(1.7, 2.4);
  		\node [right] at (1.7, 2.4) {$\langle c^{m+1} \rangle $};
  		\draw[->]  (0.2,2.2)--(-0.2, 1.6);	
  		\draw[->]  (2,2.2)--(1.6, 1.6);					
  		
  		\draw[->]  (-0.4,0.2)--(-0.4, 1.2);		
  		\draw[->]  (0.4,1.2)--(0.4,2.2);
  		\draw[->]  (1.5,0.2)--(1.5, 1.2);		
  		\draw[->]  (2.2,1.2)--(2.2,2.2);
  		
  		\draw[->]  (0,0.2)--(1.7, 0.8);		
  		\draw[->]  (0,1.6)--(1.7, 2.2);	
  		
  		\draw[->]  (-0.1,1.1)--(1, 0.2);
  		\draw[->]  (0.7,2.2)--(1.7, 1.2);
  		
  		\draw[<->]  (-0.3,0.2)--(0.3, 2.2);	
  		\draw[<->]  (1.6,0.2)--(2.1, 2.2);	
  		\draw[<->]  (-0.1,1.3)--(1.7, 1.1);										
  		\end{tikzpicture} 
  		\caption{Subgraph of $G_{D^b(A_1),\KK}(D^b(Q_2))$} \label{cube for Q2}
  		\end{figure}
  		In this subgraph we have $18$ edges coming from:  4 exceptional collections in Lemma \ref{exceptional pairs in T2} (a), $2$ exceptional pairs in Lemma \ref{exceptional pairs in T2} (a), and $12$ exceptional pairs coming from the first column in table \eqref{1 pairs}.   The 18 edges in the cube above  are pairwise disjoint.
  		
  		On the other hand using columns 2 and 3 in table \eqref{1 pairs} we can represent  four subgraphs in $G_{D^b(A_1), \KK}(\mc T_2)$  (the arrows  from $\langle x^m \rangle$ to  $\langle x^{m+1} \rangle$ are labeled by $2$, the rest by $1$) as in Figure \ref{four pyra 1 for Q2}.

  		\begin{figure}
  		
  		\begin{tikzpicture}[scale=1.1]

  		\node [left] at (0,0) {$\langle a^m \rangle $};
  		\draw[->]  (0,0)--(1,0);
  		\draw[->]  (-0.55,1.2)--(0.45,1.2);
  		\node [right] at (1,0) {$\langle a^{m+1} \rangle $};
  		\draw[<-]  (-0.4,0.2)--(-0.7,1);	
  		\node [above] at (-0.8,1) {$\langle c^{m} \rangle $};

  		\node [right] at (0.45,1.2) {$\langle c^{m+1} \rangle $};
  		\draw[<-]  (1.3,0.2)--(1,1);
  		\draw[<-]  (0.6,1)--(-0.1,0.2);		
  		\node at (0.2,3) {$\langle F_+ \rangle $};	
  		\draw[->]  (0.1,2.9)--(-0.8,1.6);	
  		\draw[->]  (0.3,2.9)--(0.8,1.6);
  		\draw[->]  (-0.2,0.2)--(0.15,2.9);	
  		\draw[->]  (1.1,0.2)--(0.2,2.9);
  		\node at (0.2,-2.5) {$\langle G_- \rangle $};
  		\draw[->]  (0.1,-2.3)--(-0.2,-0.2);	
  		
  		\draw[->]  (0.3,-2.3)--(1.1,-0.2);	
  		\draw[->]  (-0.8,1)--(0.05,-2.3);	
  		\draw[->]  (0.7,1)--(0.2,-2.3);		
  		\end{tikzpicture} 				
  		\begin{tikzpicture}[scale=1.1]

  		\node [left] at (0,0) {$\langle d^m \rangle $};
  		\draw[->]  (0,0)--(1,0);
  		\draw[->]  (-0.55,1.2)--(0.45,1.2);
  		\node [right] at (1,0) {$\langle d^{m+1} \rangle $};
  		\draw[<-]  (-0.4,0.2)--(-0.7,1);	
  		\node [above] at (-0.8,1) {$\langle b^{m} \rangle $};

  		\node [right] at (0.45,1.2) {$\langle b^{m+1} \rangle $};
  		\draw[<-]  (1.3,0.2)--(1,1);
  		\draw[<-]  (0.6,1)--(-0.1,0.2);		
  		\node at (0.2,3) {$\langle F_+ \rangle $};	
  		\draw[->]  (0.1,2.9)--(-0.8,1.6);	
  		\draw[->]  (0.3,2.9)--(0.8,1.6);
  		\draw[->]  (-0.2,0.2)--(0.15,2.9);	
  		\draw[->]  (1.1,0.2)--(0.2,2.9);
  		\node at (0.2,-2.5) {$\langle G_- \rangle $};
  		\draw[->]  (0.1,-2.3)--(-0.2,-0.2);	
  		
  		\draw[->]  (0.3,-2.3)--(1.1,-0.2);	
  		\draw[->]  (-0.8,1)--(0.05,-2.3);	
  		\draw[->]  (0.7,1)--(0.2,-2.3);		
  		\end{tikzpicture} 				
  		\begin{tikzpicture}[scale=1.1] 
  		
  		\node [left] at (0,0) {$\langle c^m \rangle $};
  		\draw[->]  (0,0)--(1,0);
  		\draw[->]  (-0.55,1.2)--(0.45,1.2);
  		\node [right] at (1,0) {$\langle c^{m+1} \rangle $};
  		\draw[<-]  (-0.4,0.2)--(-0.7,1);	
  		\node [above] at (-0.8,1) {$\langle b^{m} \rangle $};

  		\node [right] at (0.45,1.2) {$\langle b^{m+1} \rangle $};
  		\draw[<-]  (1.3,0.2)--(1,1);
  		\draw[<-]  (0.6,1)--(-0.1,0.2);		
  		\node at (0.2,3) {$\langle F_- \rangle $};	
  		\draw[->]  (0.1,2.9)--(-0.8,1.6);	
  		\draw[->]  (0.3,2.9)--(0.8,1.6);
  		\draw[->]  (-0.2,0.2)--(0.15,2.9);	
  		\draw[->]  (1.1,0.2)--(0.2,2.9);
  		\node at (0.2,-2.5) {$\langle G_+ \rangle $};
  		\draw[->]  (0.1,-2.3)--(-0.2,-0.2);	
  		
  		\draw[->]  (0.3,-2.3)--(1.1,-0.2);	
  		\draw[->]  (-0.8,1)--(0.05,-2.3);	
  		\draw[->]  (0.7,1)--(0.2,-2.3);		
  		\end{tikzpicture} 				
  		\begin{tikzpicture}[scale=1.1]

  		\node [left] at (0,0) {$\langle a^m \rangle $};
  		\draw[->]  (0,0)--(1,0);
  		\draw[->]  (-0.55,1.2)--(0.45,1.2);
  		\node [right] at (1,0) {$\langle a^{m+1} \rangle $};
  		\draw[<-]  (-0.4,0.2)--(-0.7,1);	
  		\node [above] at (-0.8,1) {$\langle d^{m} \rangle $};

  		\node [right] at (0.45,1.2) {$\langle d^{m+1} \rangle $};
  		\draw[<-]  (1.3,0.2)--(1,1);
  		\draw[<-]  (0.6,1)--(-0.1,0.2);		
  		\node at (0.2,3) {$\langle F_- \rangle $};	
  		\draw[->]  (0.1,2.9)--(-0.8,1.6);	
  		\draw[->]  (0.3,2.9)--(0.8,1.6);
  		\draw[->]  (-0.2,0.2)--(0.15,2.9);	
  		\draw[->]  (1.1,0.2)--(0.2,2.9);
  		\node at (0.2,-2.5) {$\langle G_+ \rangle $};
  		\draw[->]  (0.1,-2.3)--(-0.2,-0.2);	
  		
  		\draw[->]  (0.3,-2.3)--(1.1,-0.2);	
  		\draw[->]  (-0.8,1)--(0.05,-2.3);	
  		\draw[->]  (0.7,1)--(0.2,-2.3);		
  		\end{tikzpicture} 	
  		\caption{Subgraphs of $G_{D^b(A_1),\KK}(D^b(Q_2))$} \label{four pyra 1 for Q2}			
  		\end{figure}
  		Comparing with Figure \ref{Figure for digraph of Q1} we see that there are four different subgraphs in $G_{D^b(A_1), \KK}(\mc T_2)$, which are isomorphic to  $G_{D^b(A_1), \KK}(D^b(Q_1))$. Taking into account the arrows in the cube above and that any of the  four subgraphs in Figure \ref{four pyra 2 for Q2} is not isomorphic to the middle part  in Figure \ref{Figure for digraph of Q1}
  		\begin{figure}
  		\begin{tikzpicture}[scale=1.1] 
  		  		\node [left] at (0,0) {$\langle a^m \rangle $};
  		\draw[->]  (0,0)--(1,0);
  		\draw[->]  (-0.45,1.2)--(0.45,1.2);
  		\node [right] at (1,0) {$\langle a^{m+1} \rangle $};
  		\draw[<-]  (-0.4,0.2)--(-0.7,1);	
  		\node [above] at (-0.8,1) {$\langle c^{m} \rangle $};

  		\node [right] at (0.45,1.2) {$\langle c^{m+1} \rangle $};
  		\draw[<-]  (1.3,0.2)--(1,1);
  		\draw[<-]  (0.6,1)--(-0.1,0.2);		
  		\node at (0.2,3.1) {$\langle F_- \rangle $};	
  		\draw[<-]  (0.1,2.9)--(-0.8,1.6);	
  		\draw[<-]  (0.3,2.9)--(0.8,1.6);
  		\draw[->]  (-0.2,0.2)--(0.15,2.9);	
  		\draw[->]  (1.1,0.2)--(0.2,2.9);
  		\node at (0.2,-2.6) {$\langle G_+ \rangle $};
  		\draw[->]  (0.1,-2.3)--(-0.2,-0.2);	
  		
  		\draw[->]  (0.3,-2.3)--(1.1,-0.2);	
  		\draw[<-]  (-0.8,1)--(0.05,-2.3);	
  		\draw[<-]  (0.7,1)--(0.2,-2.3);		
  		\end{tikzpicture} 				
  		\begin{tikzpicture}[scale=1.1] 
  	  		\node [left] at (0,0) {$\langle d^m \rangle $};
  		\draw[->]  (0,0)--(1,0);
  		\draw[->]  (-0.45,1.2)--(0.45,1.2);
  		\node [right] at (1,0) {$\langle d^{m+1} \rangle $};
  		\draw[<-]  (-0.4,0.2)--(-0.7,1);	
  		\node [above] at (-0.8,1) {$\langle b^{m} \rangle $};

  		\node [right] at (0.45,1.2) {$\langle b^{m+1} \rangle $};
  		\draw[<-]  (1.3,0.2)--(1,1);
  		\draw[<-]  (0.6,1)--(-0.1,0.2);		
  		\node at (0.2,3.1) {$\langle F_- \rangle $};	
  		\draw[->]  (0.1,2.9)--(-0.8,1.6);	
  		\draw[->]  (0.3,2.9)--(0.8,1.6);
  		\draw[<-]  (-0.2,0.2)--(0.15,2.9);	
  		\draw[<-]  (1.1,0.2)--(0.2,2.9);
  		\node at (0.2,-2.6) {$\langle G_+ \rangle $};
  		\draw[<-]  (0.1,-2.3)--(-0.2,-0.2);	
  		
  		\draw[<-]  (0.3,-2.3)--(1.1,-0.2);	
  		\draw[->]  (-0.8,1)--(0.05,-2.3);	
  		\draw[->]  (0.7,1)--(0.2,-2.3);		
  		\end{tikzpicture} 				
  		\begin{tikzpicture}[scale=1.1]

  		\node [left] at (0,0) {$\langle c^m \rangle $};
  		\draw[->]  (0,0)--(1,0);
  		\draw[->]  (-0.45,1.2)--(0.45,1.2);
  		\node [right] at (1,0) {$\langle c^{m+1} \rangle $};
  		\draw[<-]  (-0.4,0.2)--(-0.7,1);	
  		\node [above] at (-0.8,1) {$\langle b^{m} \rangle $};

  		\node [right] at (0.45,1.2) {$\langle b^{m+1} \rangle $};
  		\draw[<-]  (1.3,0.2)--(1,1);
  		\draw[<-]  (0.6,1)--(-0.1,0.2);		
  		\node at (0.2,3.1) {$\langle F_+ \rangle $};	
  		\draw[->]  (0.1,2.9)--(-0.8,1.6);	
  		\draw[->]  (0.3,2.9)--(0.8,1.6);
  		\draw[<-]  (-0.2,0.2)--(0.15,2.9);	
  		\draw[<-]  (1.1,0.2)--(0.2,2.9);
  		\node at (0.2,-2.6) {$\langle G_- \rangle $};
  		\draw[<-]  (0.1,-2.3)--(-0.2,-0.2);	
  		
  		\draw[<-]  (0.3,-2.3)--(1.1,-0.2);	
  		\draw[->]  (-0.8,1)--(0.05,-2.3);	
  		\draw[->]  (0.7,1)--(0.2,-2.3);		
  		\end{tikzpicture} 				
  		\begin{tikzpicture}[scale=1.1] 
  		  		\node [left] at (0,0) {$\langle a^m \rangle $};
  		\draw[->]  (0,0)--(1,0);
  		\draw[->]  (-0.45,1.2)--(0.45,1.2);
  		\node [right] at (1,0) {$\langle a^{m+1} \rangle $};
  		\draw[<-]  (-0.4,0.2)--(-0.7,1);	
  		\node [above] at (-0.8,1) {$\langle d^{m} \rangle $};

  		\node [right] at (0.45,1.2) {$\langle d^{m+1} \rangle $};
  		\draw[<-]  (1.3,0.2)--(1,1);
  		\draw[<-]  (0.6,1)--(-0.1,0.2);		
  		\node at (0.2,3.1) {$\langle F_+ \rangle $};	
  		\draw[<-]  (0.1,2.9)--(-0.8,1.6);	
  		\draw[<-]  (0.3,2.9)--(0.8,1.6);
  		\draw[->]  (-0.2,0.2)--(0.15,2.9);	
  		\draw[->]  (1.1,0.2)--(0.2,2.9);
  		\node at (0.2,-2.6) {$\langle G_- \rangle $};
  		\draw[->]  (0.1,-2.3)--(-0.2,-0.2);	
  		
  		\draw[->]  (0.3,-2.3)--(1.1,-0.2);	
  		\draw[<-]  (-0.8,1)--(0.05,-2.3);	
  		\draw[<-]  (0.7,1)--(0.2,-2.3);		
  		\end{tikzpicture} 							
  			\caption{Subgraphs of $G_{D^b(A_1),\KK}(D^b(Q_2))$} \label{four pyra 2 for Q2}			
  		\end{figure}
  		we deduce:
  		
  		\begin{lemma} \label{four subgraphs} Let $Q_1$, $Q_2$ be the quivers in \eqref{Q1}. 
  			  			There are exactly four subgraphs in $G_{D^b(A_1), \KK}(D^b(Q_2))$, which are isomorphic to $G_{D^b(A_1), \KK}(D^b(Q_1))$.	
  		\end{lemma}

 \begin{coro} \label{some exceptional quadruples}
 For each $m\in \ZZ$	there exist the following full exceptional  collections in $\mc T_2$
 		\begin{gather} \label{Q2 table}\begin{array}{| c | c | c |}
 		\hline
 	 (G_+,a^{m},a^{m+1},F_+)   &    	(a^{m},a^{m+1},F_+,F_-)	&        (a^{m+1},F_+,F_-,b^{m+2})                \\ \hline
 	  	(F_+,b^{m},b^{m+1},G_+)	 & 	(b^{m},b^{m+1},G_+,G_-)	&        (b^{m+1},G_+,G_-,a^{m})                    \\ \hline
 	(G_+,c^{m},c^{m+1},G_-)	 & 	(c^{m},c^{m+1},G_-,F_-)	&        (c^{m+1},G_-,F_-,d^{m+1})                      \\ \hline
 	(F_-,d^{m},d^{m+1},F_+) & 	(d^{m},d^{m+1},F_+,G_+)	&        (d^{m+1},F_+,G_+,c^{m+1})                      \\ \hline
 	(G_-,d^{m},d^{m+1},G_+) & 	(d^{m},d^{m+1},G_+,F_+)	&                             \\ \hline
 		 		\end{array}
 		\end{gather}
 		Furthermore, we have  for any $m\in \ZZ$
 		\begin{gather} \label{row 1 Serre}  S(b^{m+2}) \cong a^{m}[1] \quad S(a^{m}) \cong b^{m}[1]\quad S(d^{m+1}) \cong c^{m}[1]\quad S(c^{m+1}) \cong d^{m}[1]\\ \label{row 2 Serre}
 		  S(F_-) \sim G_+ \quad S(G_-) \sim F_+\quad S(G_+) \sim F_-\quad S(F_+) \sim G_-\end{gather}
 \end{coro} \bpr These are exceptional collections, since any subsequence of length two in any of them is an exceptional pair (see  Lemma \ref{exceptional pairs in T2} or Figures \ref{cube for Q2}, \ref{four pyra 1 for Q2}, \ref{four pyra 2 for Q2} describing subgraphs of $G_{D^b(A_1),\KK}(\mc T_2)$). Since these sequences have maximal possible length - four, from \cite{WCB1} we know that they are full exceptional collections.
  From \eqref{Serre} it follows that  for any $m\in \ZZ$ we have  $S(b^{m+2}) \sim a^{m}$,  $ S(a^{m}) \sim b^{m}$, $ S(d^{m+1}) \sim c^{m}$, $ S(c^{m+1}) \sim d^{m}$, $
  S(F_-) \sim G_+$, $S(G_-) \sim F_+$, $ S(G_+) \sim F_- $, $S(F_+) \sim G_-$.
 
Finally to prove the more precise incidences in \eqref{row 1 Serre}  we use first  \cite[Proposition 2.5]{DK1} to derive:
\begin{gather}  \hom^p(c^{m+1}, d^m) = \left \{ \begin{array}{c c} 1 & \mbox{if} \ p=1 \\ 0 & \mbox{otherwise}  \end{array} \right. \qquad  \hom^p(d^{m+1}, c^m) = \left \{ \begin{array}{c c} 1 & \mbox{if} \ p=1 \\ 0 & \mbox{otherwise}  \end{array} \right. \nonumber \\[-2mm] \label{help for Serere for Q_2} \\[-2mm] \nonumber 
\hom^p(a^{m}, b^m) = \left \{ \begin{array}{c c} 1 & \mbox{if} \ p=1 \\ 0 & \mbox{otherwise}  \end{array} \right. \qquad  \hom^p(b^{m+2}, a^m) = \left \{ \begin{array}{c c} 1 & \mbox{if} \ p=1 \\ 0 & \mbox{otherwise}  \end{array} \right. \end{gather}
and then using the main property of the Serre Functor we  can write $1=\hom(c^{m+1}, c^{m+1})=\hom(c^{m+1}, S(c^{m+1}))=\hom(c^{m+1}, d^{m}[p])$, where $S(c^{m+1})\cong d^{m}[p]$ due to  $ S(c^{m+1}) \sim d^{m}$, and now from \eqref{help for Serere for Q_2} it follows $p=1$.  Similarly one derives the rest  incidences in  \eqref{row 1 Serre}.
\epr
 
\begin{prop}\label{autoeq sigma prop} 
	There exist $\KK$-linear  autoequivlences $\theta, \zeta: \mc T_2 \cong \mc T_2$, such that for each $m\in \ZZ$: 
	\begin{gather} \label{autoeq sigma} \theta(a^m) \cong a^m \quad \theta(b^m) \cong b^m \quad  \theta(c^m) \cong d^m \quad \theta(d^m) \cong c^m \quad \theta(F_\pm) \cong F_\mp \quad \theta(G_\pm) \cong G_\mp. \\
	\label{autoeq zeta} \zeta(a^m) \cong d^m \quad \zeta(b^m) \cong c^{m-1}  \quad \zeta(c^m) \cong b^m \quad \zeta(d^m) \cong a^{m-1} \ \begin{array}{c}\zeta(F_+) \sim F_+ \ \zeta(G_-) \sim G_-\\ \zeta(G_+) \sim F_- \ \zeta(F_-) \sim G_+.\end{array}\end{gather}
\end{prop}
 \bpr Recall that the exceptional objects of $\mc T_2$ are  $\{E_1^m\}_{m\in \ZZ}$,  $\{E_2^m\}_{m\in \ZZ}$, $\dots$,  $\{E_8^m\}_{m\in \ZZ}$, and  $\{F_+,G_-\}$, $\{F_-,G_+\}$ and that $a^m, b^m, c^m, d^m$ are defined in \eqref{abcd}.  The symmetry of the quiver $Q_2$ along the diagonal gives an exact equivalence $Rep_{\KK}(Q_2) \rightarrow Rep_\KK(Q_2)$, which  extends to  $\theta $.  
 
In \cite[Theorem 4.5 (3)]{Miyachi} is proved that ${\rm Aut}_{\KK}(\mc T_2) \cong \ZZ^2 \times (S_2 \ltimes \KK^{\times})$, where one $\ZZ$-factor is generated by the translation functor, and we will show that a generator of the other $\ZZ$-factor is the desired auto-equivalence  $\zeta$.

 The Auslander Reiten quiver of $\mc T_2=D^b(Q_2)$ contains the quiver $Q_2$  as explained in the beginning of  Subsection \ref{material from Mi}.  
   We will show that this subquiver is  (up to symmetry):
   \begin{gather} \label{subquiver of ARQ}
   \begin{diagram}
   c^1[-1] &  \rTo  &  b^2[-1]   \\
   \uTo &        & \uTo     \\
   b^1[-1]  & \rTo  &     d^1[-1]
   \end{diagram}.
   \end{gather}

  We  show that the indecomposable projective representations are in the exceptional collection $(b^1[-1], d^1[-1], c^1[-1],b^2[-1])$. Indeed, from \eqref{abcd}  we see that this exceptional collection is \\ $(E_3^0,E_8^0,E_5^0,E_3^1)$ and from  \cite[Proposition 2.5]{DK1} we see that  the simple representations  are $F_+,F_-, E_1^0, E_3^0$, and from the table in  \cite[Proposition 2.5]{DK1} we know that: $\hom^1(X,Y)=0 $ for   
  $  X\in \{E_3^0,E_8^0,E_5^0,E_3^1\}$, $ Y \in  \{F_+,F_-, E_1^0, E_3^0\}  $, therefore  $\hom^1(X,Y)=0 $ for 
  $X\in \{E_3^0,E_8^0,E_5^0,E_3^1\} $,  $  Y \in Rep_\KK(Q_2). $
  Furthermore, from table in  \cite[Proposition 2.5]{DK1} we check that $\hom(E_3^0,E_8^0)=\hom(E_3^0,E_5^0)=\hom(E_8^0,E_3^1)=\hom(E_5^0,E_3^1)=1$, These arguments imply that  the full subquiver of the Auslander-Reiten quiver of $\mc T_2$  with vertices  the indecomposable projective representations is \eqref{subquiver of ARQ}.

 Next, we will describe a connected component of the quiver  ${\bf \Gamma}^{irr}$ (this quiver was discussed in Subsection  \ref{material from Mi} ) and we will  depict locally the identification $\rho: \ZZ\times (\ZZ Q_2)\rightarrow {\bf \Gamma}^{irr}$	from (c) in Subsection \ref{material from Mi}.
  
\begin{figure}
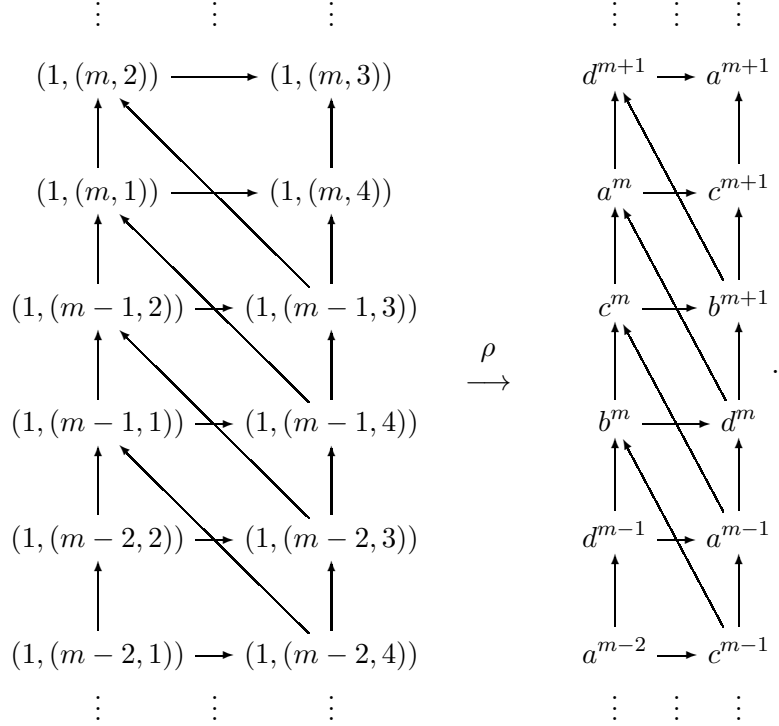

   \begin{gather}  \nonumber 
 \begin{diagram} \vdots & \vdots & \vdots \\
 (1,(m,2)) &  \rTo  &  (1,(m,3))  \\
 \uTo &    \luTo(2,4)     & \uTo     \\
 (1,(m,1) ) & \rTo  &     (1,(m,4))\\
 \uTo &     \luTo(2,4)    & \uTo     \\
 (1,(m-1,2))  &  \rTo  & (1, (m-1,3) )  \\
 \uTo &     \luTo(2,4)    & \uTo     \\
(1, (m-1,1))  & \rTo  &    (1,(m-1,4) )\\
 \uTo &    \luTo(2,4)      & \uTo     \\
(1, (m-2,2)) &  \rTo  &  (1,(m-2,3) )  \\
 \uTo &      & \uTo     \\
(1, (m-2,1) ) & \rTo  &    (1, (m-2,4)) \\
 \vdots & \vdots & \vdots 
 \end{diagram} \qquad \begin{array}{c} \rho \\ \longrightarrow \end{array}  \qquad  \begin{diagram} \vdots & \vdots & \vdots \\
    d^{m+1} &  \rTo  &  a^{m+1}   \\
   \uTo &    \luTo(2,4)     & \uTo     \\
    a^{m}  & \rTo  &     c^{m+1}\\
     \uTo &     \luTo(2,4)    & \uTo     \\
    c^m &  \rTo  &  b^{m+1}   \\
   \uTo &     \luTo(2,4)    & \uTo     \\
   b^m  & \rTo  &     d^m \\
     \uTo &    \luTo(2,4)      & \uTo     \\
     d^{m-1} &  \rTo  &  a^{m-1}   \\
     \uTo &      & \uTo     \\
     a^{m-2}  & \rTo  &     c^{m-1} \\
      \vdots & \vdots & \vdots 
   \end{diagram}.
   \end{gather}
    \caption{The identification $\rho$ for $Q_1$}  \label{extend subquiver of ARQ}
   \end{figure}

 Now we shift the subquiver \eqref{subquiver of ARQ} of ${\bf \Gamma}^{irr}$ and by acting on it with with $\langle  [-1]S \rangle$ using \eqref{row 1 Serre} we see that the identification $\rho$ for odd  $m\in \ZZ$ can be depicted locally as shown in Figure \ref{extend subquiver of ARQ}, where although in the picture some arrows intersect they actually do not intersect. 
  
  Looking in the proof of \cite[Theorem 4.5 (3)]{Miyachi} (note that we have observed the rule for labeling of the  the vertices   in  \cite[Figure 6]{Miyachi}) we see that  a generator of the other $\ZZ$-factor in ${\rm Aut}_{\KK}(\mc T_2)$  acts on the left hand side in Figure  \ref{extend subquiver of ARQ}  as follows (for any $n\in \ZZ$)
  \begin{gather}  (1,(n,1)) \mapsto (1,(n-1,4)) \quad  (1,(n,2)) \mapsto (1,(n,1))  \\  (1,(n,3)) \mapsto (1,(n,2))\quad  (1,(n,4)) \mapsto (1,(n-1,3)).
  \end{gather}
  Via the correspondence in Figure \ref{extend subquiver of ARQ} we get an auto-equivalence $\zeta$ of $\mc T_2$ such that the first four entities in \eqref{autoeq zeta} hold for any $m\in \ZZ$ (not only odd). To show the rest, we use Lemma \ref{exceptional pairs in T2}  to deduce that $(a^m,d^{m+1},a^{m+1},F_+)$, $(d^m,a^{m},d^{m+1},F_+)$ are both full exceptional collections in $\mc T_2$. From the already proved we see that  $(\zeta(a^m),\zeta(d^{m+1}),\zeta(a^{m+1}),\zeta(F_+))\cong (d^m,a^m,d^{m+1},\zeta(F_+))$ is a full exceptional collection, and it follows that $\zeta(F_+)\sim F_+$. From table \ref{Q2 table} we see that  $(a^{m},a^{m+1},F_+,F_-)$ and   $(d^{m},d^{m+1},F_+,G_+)$ are both full exceptional collections and from the already proved \\  $(\zeta(a^{m}),\zeta(a^{m+1}),\zeta(F_+),\zeta(F_-))\sim  (d^m,d^{m+1},F_+,\zeta(F_-))$, which implies $\zeta(F_-)\sim G_+$. In table \ref{Q2 table} we have also exceptional sequences $(c^{m},c^{m+1},F_-,G_-)	$ and $(b^{m},b^{m+1},G_+,G_-)$ and we already proved that $ (\zeta(c^{m}),\zeta(c^{m+1}),\zeta(F_-),\zeta(G_-))\sim (b^{m},b^{m+1},G_+,\zeta(G_-))	$, therefore $\zeta(G_-)\sim G_-$. Finally, from the full exceptional collections $(b^{m},b^{m+1},G_-,G_+)$, $(c^{m-1},c^{m},G_-,F_-)$ in table \ref{Q2 table}  and similar \\ arguments we deduce  $\zeta(G_+)\sim F_-$.
  \epr
\begin{coro}  \label{coro for Q2 genus 1} For any $m,n \in \ZZ$ and any $x\in \{a,b,c,d\}$ we have $\langle x^m,x^{m+1} \rangle = \langle x^n,x^{n+1} \rangle$. 
		In particular the following subcategories do not depend on $m$
		\begin{gather}
		A = \langle a^m,a^{m+1} \rangle \quad B=\langle b^m,b^{m+1} \rangle \quad  C=\langle  c^m,c^{m+1} \rangle \quad D=\langle d^m,d^{m+1} \rangle.
		\end{gather}
		The action of the Serre functor on these subcateories is described in the following table:
		\begin{gather} \label{Serre on the 2 pairs}\begin{array}{| c | c | c | c | c |}
				\hline
				X    & A     & B &  C  & D   \\ \hline
				S(X) & B     & A &  D  & C \\ \hline				\end{array} ,
		\end{gather}
		and the set of non-commutative curves of genus one in $\mc T_2$ is:
		\begin{gather}  \label{C_1^{{rm Id}}(mc T_2)}  C_1^{\{\rm Id\}}(\mc T_2) =\{ A,B,C,D \} \qquad  \abs{C_1^{\{\rm Id\}}(\mc T_2)}= 4.   \end{gather} Furthermore, \eqref{Serre on the 2 pairs},  \eqref{C_1^{{rm Id}}(mc T_2)} imply that $C_1^{\langle S \rangle}(\mc T_2) \cong  C_1^{\{\rm Id\}}(\mc T_2)/\langle S \rangle $ are finite sets with $2$ elements and:
		\begin{gather}   \label{C_1^{S}(mc T_2)}  C_1^{\{\rm Id\}}(\mc T_2)/\langle S \rangle =\{[A],[C]\} \qquad  \abs{C_1^{\langle S \rangle}(\mc T_2)}= 2.  \end{gather}
	 The auto-equivalences $\theta, \zeta$ in Proposition \ref{autoeq sigma prop} satisfy  $\zeta \circ \theta (C)= A$ and hence the 	$ \abs{C_1^{{\rm Aut}_{\KK}(\mc T_2)}(\mc T_2)}= 1$. 
\end{coro}
\bpr From table \eqref{Q2 table} we see that for any $m \in \ZZ$ we have: $$ \langle a^m,a^{m+1} \rangle =  \langle F_+, F_- \rangle^\perp \ \  \langle b^m,b^{m+1} \rangle =  \langle G_+, G_- \rangle^\perp \ \ \langle c^m,c^{m+1} \rangle =  \langle G_-, F_- \rangle^\perp \ \ \langle d^m,d^{m+1} \rangle = \langle F_+, G_+ \rangle^\perp $$ and the first statement follows.  Recalling that $(F_+,F_-), (G_+,G_-), (F_+,G_+), (F_-,G_-)$ are orthogonal pairs and $ F_+,F_-, G_+,G_- $ are pairwise non-isomorphic it follows that $\langle F_+, F_- \rangle, $  $ \langle G_+, G_- \rangle $ ,$ \langle G_-, F_- \rangle $ , $ \langle F_+, G_+ \rangle$  are pairwise different, and hence $A,B,C,D$ are  pairwise different subcategories.

	Let $\mc A\subset \mc T_2$, $\mc A \in C_1^{\{\rm Id\}}(\mc T_2)$. Therefore we have a  strong exceptional pair $(E_1,E_2)$, s.t. $\hom(E_1,E_2)=2$,   $\langle E_1,E_2 \rangle =\mc A $. From Lemma \ref{exceptional pairs in T2}  we see that $\mc A = \langle x^{m}, x^{m+1}\rangle$ for $x\in\{a,b,c,d\}$. i.e. $\mc A \in \{A,B,C,D\}$.  From Lemma \ref{exceptional pairs in T2} (a) and \eqref{at most one non-vanisching}  it follows that $\{A,B,C,D\} \subset  C_1^{\{\rm Id\}}(\mc T_2)$,   and \eqref{C_1^{{rm Id}}(mc T_2)} follows. 
	 
	 From \eqref{row 1 Serre} and the already proven first part of the Lemma follows table \eqref{Serre on the 2 pairs} and \eqref{C_1^{S}(mc T_2)}.
	 
	  Finally, using the formulas in Proposition \ref{autoeq sigma prop} we compute $\zeta \circ \theta (\langle c^m, c^{m+1} \rangle) = \zeta  (\langle d^m, d^{m+1} \rangle) = \langle a^{m-1}, a^{m} \rangle $  and the corollary is  proved. 
\epr
 \begin{coro} \label{coro extensions of 1 pairs}
 	For each $m\in \ZZ$	the exceptional pairs from table \eqref{1 pairs} can be extended to the following full exceptional  collections in $\mc T_2$
 
 	\begin{gather} \label{extensions of 1 pairs} \begin{array}{| c | c | c |}
 	\hline
 	(a^{m-1},c^m,F_-,d^m)	&        (c^m,G_-,F_-,d^m)        &        (G_-,a^{m-1},F_-,d^m)           \\ \hline
 	(c^m,a^m,F_-,b^{m+1})  	&        (a^m,F_+,F_-,b^{m+1})        &       (F_+,c^{m},F_-,b^{m+1})          \\ \hline
 	(a^{m-1},d^m,a^{m},F_+)  	&        (d^m,G_+,a^{m},F_+)        &       (G_+,a^{m-1},a^{m},F_+)          \\ \hline
 	(d^{m},a^m,d^{m+1},F_+)	&        (a^m,F_-,d^{m+1},F_+)        &        (F_-,d^{m},d^{m+1},F_+)           \\ \hline
 	(b^{m},c^m, b^{m+1},G_-)	&        (c^m,F_-, b^{m+1},G_-)        &        (F_-,b^{m}, b^{m+1},G_-)           \\ \hline
 	(c^{m-1},b^m,G_-,a^{m-1})	&        (b^m,G_+,G_-,a^{m-1})        &        (G_+,c^{m-1},G_-,a^{m-1})           \\ \hline
 	(b^{m},d^m,G_+,c^{m})	&        (d^m,F_+,G_+,c^{m})        &        (F_+,b^{m},G_+,c^{m})           \\ \hline
 	(d^{m-1},b^m,d^{m},G_+)	&        (b^m,G_-,d^{m},G_+)        &        (G_-,d^{m-1},d^{m},G_+)           \\ \hline
 	\end{array}
 	\end{gather}
  \end{coro} \bpr The proof that these are  exceptional collections is as in the  Corollary \ref{some exceptional quadruples}.  Since these sequences have maximal possible length - four, from \cite{WCB1} we know that they are full exceptional collections. \epr
 \begin{coro}  \label{coro for Q2 genus 0}
 	If we denote:
 		\begin{gather}
 		cG_-^m=\langle c^m,G_-\rangle \quad 	aF_+^m=\langle a^m,F_+\rangle \quad dG_+^m=\langle d^m,G_+\rangle \quad aF_-^m=\langle a^m,F_-\rangle \\
 		cF_-^m=\langle c^m,F_-\rangle \quad 	bG_+^m=\langle b^m,G_+\rangle \quad dF_+^m=\langle d^m,F_+\rangle \quad bG_-^m=\langle b^m,G_-\rangle
 		\end{gather}
 		Then the action of the Serre functor on these subcategories is described in the following table:
 			\begin{gather} \label{Serre on the 1 pairs}\begin{array}{| c | c | c | c | c | c | c | c | c |}
 					\hline
X    & cG_-^m     & aF_+^m &  dG_+^m     & aF_-^m & cF_-^m	   &  	bG_+^m     & dF_+^m    &bG_-^m  \\ \hline
S(X) & dF_+^{m-1} & bG_-^m &  cF_-^{m-1} & bG_+^m & dG_+^{m-1} &  	aF_-^{m-2} & cG_-^{m-1}&aF_+^{m-2} \\ \hline				\end{array} ,
 			\end{gather}
 		and the set of non-commutative curves of genus zero in $\mc T_2$ is:
 		 			\begin{gather}   C_0^{\{\rm Id\}}(\mc T_2) =\{cG_-^m:m\in \ZZ \}\cup\{	aF_+^m:m\in \ZZ \}\cup\{	dG_+^m:m\in \ZZ \} \cup\{	aF_-^m:m\in \ZZ \}\nonumber \\  \label{C_0^{{rm Id}}(mc T_2)} \cup \{cF_-^m:m\in \ZZ \}\cup\{	bG_+^m:m\in \ZZ \}\cup\{	dF_+^m:m\in \ZZ \} \cup\{	 bG_-^m:m\in \ZZ \}
 		\\ \mbox{disjoint union, injective sequences. In particular, $  \ \ \abs{C_0^{\{\rm Id\}}(\mc T_2) }=\infty$}. \nonumber\end{gather} Furthermore, \eqref{Serre on the 1 pairs},  \eqref{C_0^{{rm Id}}(mc T_2)} imply that $C_0^{\langle S \rangle}(\mc T_2) \cong  C_0^{\{\rm Id\}}(\mc T_2)/\langle S \rangle $ are finite sets with $8$ elements and:
 		\begin{gather}   \label{C_0^{S}(mc T_2)}  C_0^{\{\rm Id\}}(\mc T_2)/\langle S \rangle =\{[cG_-^0],[cG_-^1]\}\cup\{	[aF_+^0],	[aF_+^1] \}\cup\{	[dG_+^0],	[dG_+^1] \} \cup\{	[aF_-^0],[aF_-^1] \}. \end{gather}
 		Furthermore, the   auto-equivalences $\theta, \zeta$ in Proposition \ref{autoeq sigma prop} satisfy:  
 		\begin{gather}
 	\label{help formula for full group 22}	\theta\left (cG_-^0\right )=dG_+^0 \quad 	\theta\left (cG_-^1\right )=dG_+^1 \quad \theta\left (aF_+^0\right )= aF_-^0 \quad \theta\left (aF_+^1\right )= aF_-^1 \\  	\label{help formula for full group 23}
 		\zeta^2\left (cG_-^1\right )=cG_-^0 \quad  \zeta^2\left (aF_+^1\right )=aF_+^0 \quad S\circ \zeta \left (cG_-^0\right ) = aF_+^{-2} =  \zeta^4\left (aF_+^0\right ) 
 		\end{gather}
 		 and hence the 	$ \abs{C_0^{{\rm Aut}_\KK(\mc T_2)}(\mc T_2)}= 1$.
 		 		Thus  the row with $l=0$  in the second table of  \eqref{table with numbers Q1} is proved. 
 \end{coro} 
 \bpr  From \eqref{row 1 Serre} \eqref{row 2  Serre} one computes table \eqref{Serre on the 1 pairs}.
 
 	Let $\mc A\subset \mc T_2$, $\mc A \in C_0^{\{\rm Id\}}(\mc T_2)$. Therefore $\langle E_1,E_2 \rangle =\mc A $ for a  strong exceptional pair $(E_1,E_2)$ with $\hom(E_1,E_2)=1$ . From Lemma \ref{exceptional pairs in T2}  we see that $\mc A = \langle X\rangle$ for some of the exceptional pairs in table \eqref{1 pairs}. However from Corollary \ref{coro extensions of 1 pairs}  it follows that:  
 \begin{gather} \label{equal on each row} \langle X\rangle=\langle X'\rangle \ \mbox{ whenever $X,X'$ are exceptional pairs from the same row of table \eqref{1 pairs}}
 \end{gather}
  hence follows the equality in \eqref{C_0^{{rm Id}}(mc T_2)} (by Lemma \ref{exceptional pairs in T2}, \eqref{at most one non-vanisching}, and Proposition \ref{bijection 123} follows that any of the elements on the RHS of\eqref{C_0^{{rm Id}}(mc T_2)}  is an element of $C_0^{\{\rm Id\}}(\mc T_2)$). Next we show that the sequences in \eqref{C_0^{{rm Id}}(mc T_2)} are injective and that the union is disjoint.
 
 If $cG^m_- = cG^n_-$ for some $m,n \in \ZZ$ then from the first row in table \eqref{extensions of 1 pairs} it follows that  ${}^\perp(cG^m_-) = \langle F_-, d^m\rangle ={}^\perp(cG^n_-) = \langle F_-, d^n\rangle$, which implies that  $(c^m,G_-,F_-,d^n)$, $ (c^n,G_-,F_-,d^m)$ are exceptional collections and hence $(c^m,d^n)$, $(c^n,d^m)$ are exceptional pairs, therefore by Lemma \ref{exceptional pairs in T2} it follows that $m=n$, therefore the first sequence in \eqref{C_0^{{rm Id}}(mc T_2)} is injective. 
 
 If $aF^m_+ = aF^n_+ $ for some $m,n \in \ZZ$ then by similar arguments using the second row in table  \eqref{extensions of 1 pairs} it follows that $(a^{m},b^{n+1})$, $(a^{n},b^{m+1})$  are exceptional  pairs, therefore by Lemma \ref{exceptional pairs in T2} it follows that $m=n$ and the second sequence in \eqref{C_0^{{rm Id}}(mc T_2)} is injective. 
 
 If $dG^m_+ = dG^n_+ $ for some $m,n \in \ZZ$ then by similar arguments using the third row in table \eqref{extensions of 1 pairs} it follows that $(d^{m},a^{n})$,$(d^{n},a^{m})$  are exceptional  pairs, therefore by Lemma \ref{exceptional pairs in T2} it follows that $m=n$. 
 
 If $aF^m_- = aF^n_- $ for some $m,n \in \ZZ$ then by similar arguments using the fourth row in table \eqref{extensions of 1 pairs} it follows that $(a^{m},d^{n+1})$is an exceptional  pair, therefore by Lemma \ref{exceptional pairs in T2}  it follows that $m=n$. 
 
 Any of the equalities  $cG^m_- =  dG^n_+$,  $cG^m_- =  aF^n_-$, $aF^m_+ =  dG^n_+$,  $aF^m_+ =   aF^n_-$,   for some $m,n \in \ZZ$ contradicts that $\{  (F_+,F_-), (G_+,G_-), (F_+,G_+), (F_-,G_-) \}$ are orthogonal pairs - see Lemma \ref{exceptional pairs in T2}. 
 
 If $cG^m_- =  aF^n_+$ or $dG^m_+ =  aF^n_-$ for some $m,n$ then from \eqref{equal on each row} and  \eqref{1 pairs} it follows that  $\langle a^{m-1},c^m \rangle = \langle F_+,c^{n} \rangle$ or  $\langle a^{m-1},d^m \rangle = \langle F_-,d^{n} \rangle$, which implies $m=n$ and $a^{m-1} \sim F_{+/-}$ -  contradiction.

 So far we showed that the first four sequences are injective and the corresponding images are pairwise disjoint.  However from table \eqref{Serre on the 1 pairs} we see that the union of the images of the  last four sequences  is obtained bijectively from the already discussed union via the Serre functor. To prove \eqref{C_0^{{rm Id}}(mc T_2)} it remains to show that no an element on the first row equals an element in the second row. 
  
  Recalling that $\{  (F_+,F_-), (G_+,G_-), (F_+,G_+), (F_-,G_-) \}$  are orthogonal pairs it remains to show that non of the equalities $cG_-^m=dF_+^n$, $aF_+^m=bG_-^n$, $dG_+^m=cF_-^n$,  $aF_+^m=bG_+^n$ for any $m,n \in \ZZ$ can happen.  Using table \eqref{extensions of 1 pairs} we derive the following contradictions to Lemma \ref{exceptional pairs in T2}
  
  $cG_-^m=dF_+^n$ implies that  
  $d^n, F_+,F_-,d^m$ is exc. collection, hence $F_+,d^m$ is exc. pair,
  
  $aF_+^m=bG_-^n$ implies that  
  $b^n, G_-,F_-,b^{m+1}$ is exc. collection, hence $G_-,b^{m+1}$ is exc. pair,
  
  $dG_+^m=cF_-^n$ implies that  
  $c^n, F_-,a^m,F^{+}$ is exc. collection, hence $F_-,a^{m}$ is exc. pair,
  
  $aF_+^m=bG_+^n$  implies that $b^n,G_+, F_-, b^{m+1}$ is an exc. collection, hence  $G_+,b^{m+1}$ is exc. pair.

Table \eqref{Serre on the 1 pairs} shows that the action of $\langle S \rangle$ on $C_0^{\{\rm Id\} }(\mc T_2)$ (discussed in Corollary \ref{group action}) has four invariant subsets, each of them  quotient by  $\langle S \rangle$ is a two-elements set,  and the corresponding splitting  of the quotient   $C_0^{\{\rm Id\}}(\mc T_2)/\langle S \rangle $ into four  quotients is as in \eqref{C_0^{S}(mc T_2)}. 
  
  Finally, the equalities  in \eqref{help formula for full group 22}, \eqref{help formula for full group 23} follow by combining formulas \eqref{row 1 Serre}, \eqref{row 2 Serre} and the formulas in  Proposition \ref{autoeq sigma prop} 
 \epr
 
 \begin{coro}  \label{coro for Q2 genus -1}
 	If we denote:
 	\begin{gather} \nonumber 
 	AB^m=\langle a^m,b^{m+1} \rangle \quad 	CD^m=\langle c^m,d^m\rangle \\ \nonumber  F_\pm=\langle F_+,F_-\rangle \quad G_\pm=\langle G_+,G_-\rangle \quad
 	FG_+=\langle F_+,G_+\rangle \quad 	FG_-=\langle F_-,G_-\rangle.
 	\end{gather}
 	Then the action of the Serre functor on these subcateories is described in the following table:
 	\begin{gather} \label{Serre on the -1 pairs}\begin{array}{| c | c | c | c | c | c | c | c | c |}
 	\hline
 	X    & AB^m     & 	CD^m  &   F_\pm    &  G_\pm & FG_+	   &  	FG_-      \\ \hline
 	S(X) & AB^{m-1} & CD^{m-1} &   G_\pm & F_\pm & 	FG_- &  	FG_+.  \\ \hline				\end{array} ,
 	\end{gather}
 	and the set of non-commutative curves of genus $-1$ in $\mc T_2$ is:
 	\begin{gather}   C_{-1}^{\{\rm Id\}}(\mc T_2) =\{AB^m :m\in \ZZ \}\cup\{	CD^m :m\in \ZZ \}\cup\{ F_\pm, G_\pm, FG_+, 	FG_- \} \nonumber \\[-2mm]  \label{C_{-1}^{{rm Id}}(mc T_2)} \\[-2mm]
 	 \nonumber \mbox{disjoint union, injective sequences. In particular, $  \ \ \abs{C_{-1}^{\{\rm Id\}}(\mc T_2) }=\infty$}. \nonumber\end{gather} Furthermore, \eqref{Serre on the -1 pairs},  \eqref{C_{-1}^{{rm Id}}(mc T_2)} imply that $C_{-1}^{\langle S \rangle}(\mc T_2) \cong  C_{-1}^{\{\rm Id\}}(\mc T_2)/\langle S \rangle $ are finite sets with $4$ elements and for any $m, n \in \ZZ$ we have:
 	\begin{gather}   \label{C_{-1}^{S}(mc T_2)}  C_{-1}^{\{\rm Id\}}(\mc T_2)/\langle S \rangle =\{[AB^m],[CD^n],	[F_\pm],	[FG_+] \}. \end{gather}
 		The   autoequivalences $ \zeta$ in Proposition \ref{autoeq sigma prop} satisfy for each $m \in \ZZ$:  
 		\begin{gather}
 		\label{help formula for full group -12}	\zeta\left (AB^m\right )=CD^m \quad 	\zeta\left (F_\pm\right )=FG_+  \ \Rightarrow \ C_{-1}^{\{\rm Id\}}(\mc T_2)/{\rm Aut}_{\KK}(\mc T_2)  =\{[AB^m],	[F_\pm]\} 
 		\end{gather}
 		and we claim that the two elements in $C_{-1}^{\{\rm Id\}}(\mc T_2)/{\rm Aut}_{\KK}(\mc T_2) $ are differeint, i.e. 	$ \abs{C_{-1}^{{\rm Aut}_{\KK}(\mc T_2)}(\mc T_2) }= 2$.

 		Thus  the row with $l=-1$, in the second  table of  \eqref{table with numbers Q1} is proved. 
  \end{coro} 
 \bpr

   From \eqref{row 1 Serre}, \eqref{row 2  Serre} one computes table \eqref{Serre on the -1 pairs}.
  
 Recalling  Proposition \ref{bijection 123} we see that   $\mc A\subset \mc T_2$, $\mc A \in C_{-1}^{\{\rm Id\}}(\mc T_2)$ is equivalent to  existence of an orthogonal exceptional pair $(E_1,E_2)$, s.t.   $\langle E_1,E_2 \rangle =\mc A $. In Lemma \ref{exceptional pairs in T2} (b) are listed  all the orthogonal exceptional pairs in $\mc T_2$. Hence follows   \eqref{C_{-1}^{{rm Id}}(mc T_2)}.
 
\begin{remark} \label{orthogonal pairs}  It is easy to show that for any two orthogonal exceptional pairs $(\alpha, \beta)$, $(\gamma, \delta)$ from $\langle \alpha, \beta\rangle = \langle \gamma, \delta \rangle$ it follows that either $\alpha \sim \gamma$, $\beta \sim  \delta$ or  $\alpha \sim \delta$, $\beta \sim  \gamma$. 
\end{remark}
Hence follows  that the sequences in \eqref{C_{-1}^{{rm Id}}(mc T_2)} are injective and that the union is disjoint. From the already proved facts follows \eqref{C_{-1}^{S}(mc T_2)}.
To compute the equialities in \eqref{help formula for full group -12} we apply  Proposition \ref{autoeq sigma prop}.

It remains to show that for any $m\in \ZZ$ there is no a $\KK$-linear auto-equivalence $v$ of $\mc T_2$, such that $v(\langle F_+,F_-\rangle) = \langle a^m,b^{m+1}\rangle$. Indeed, if such $m\in \ZZ$ and an auto-equivalence $v$ existed, then we note that there are full exceptional collections  (see  table \eqref{Q2 table} and Corollary \ref{exceptional pairs in T2}) 	$(a^{m},a^{m+1},F_+,F_-)$, $(c^m,d^m,a^m,b^{m+1})$ and we derive $ v\left ( \langle a^{m},a^{m+1} \rangle \right ) =v\left ( \langle F_+,F_- \rangle^\perp \right )=v\left ( \langle F_+,F_- \rangle \right )^\perp=\langle  a^m,b^{m+1}\rangle^\perp =\langle c^m,d^m \rangle $,
which contradicts the fact that $ ( a^m,a^{m+1} )$ is not orthogonal exceptional pair, whereas $ (c^m,d^m) $ is an orthogonal exceptional pair (see Corollary \ref{exceptional pairs in T2}). The corollary is proved.
 \epr
 
Lemma \ref{four subgraphs} suggests the following

\begin{prop} \label{triples in Q_2 part 1}  Let $Q_1$, $Q_2$ be the quivers in \eqref{Q1}, then the numbers $\abs{C^{\Gamma}_{D^b(Q_1), \KK}(D^b(Q_2))}$ for $\Gamma=\{\rm Id\}$,  $\Gamma=\langle S \rangle$, $\Gamma={\rm Aut}_{\KK}(D^b(Q_2))$ are $4,2,1$, respectively. 
   \end{prop} \bpr Lemma \ref{four subgraphs}  and   Proposition \ref{from invariants to subgraphs}  imply  that  $\abs{C_{D^b(Q_1), \KK}(D^b(Q_2))}\leq 4$. On the other hand in Figure  \ref{four pyra 2 for Q2}  we  read the following four exceptional sequences (here $m$ can be  any integer):
\begin{gather}
(b^m,d^m,b^{m+1}, G_+), (b^m,c^m,b^{m+1}, G_-),  (d^m,a^m,d^{m+1}, F_+),  (c^m,a^m,c^{m+1}, F_-) 
\end{gather}
In particular, we have four different subcategories  $\langle G_+ \rangle^{\perp} $,  $\langle G_- \rangle^{\perp} $,  $\langle F_+ \rangle^{\perp} $,  $\langle F_- \rangle^{\perp} $. The auto-equivalences  \eqref{row 2 Serre} and \eqref{autoeq sigma} give us
\begin{gather}
S\left ( \langle G_+ \rangle^{\perp} \right )  =\langle F_- \rangle^{\perp} \quad S\left (\langle F_- \rangle^{\perp}  \right )  =\langle G_+ \rangle^{\perp} \quad  S\left ( \langle G_- \rangle^{\perp} \right )  =\langle F_+ \rangle^{\perp}  \quad \theta\left ( \langle F_+ \rangle^{\perp} \right )  =\langle F_- \rangle^{\perp}
\end{gather} and therefore these four subcategories are pairwise equivalent. We will show that they are equivalent to $D^b(Q_1)$, then the already explained implies \eqref{expl for Q1 in Q2} and the Proposition will be proved:
\begin{gather}
\label{expl for Q1 in Q2} C^{\{\rm Id\}}_{D^b(Q_1), \KK}(D^b(Q_2))=\{\langle G_+ \rangle^{\perp},  \langle G_- \rangle^{\perp},  \langle F_+ \rangle^{\perp} ,  \langle F_- \rangle^{\perp}\}  
\\
C^{\langle S \rangle}_{D^b(Q_1), \KK}(D^b(Q_2))=\{[\langle F_+ \rangle^{\perp}], [ \langle F_- \rangle^{\perp}] \} \nonumber  \qquad \abs{C^{{\rm Aut}_{\KK}(D^b(Q_2))}_{D^b(Q_1), \KK}(D^b(Q_2))}=1
\end{gather}

The subcategory $\langle G_+ \rangle^{\perp}$ has a full strong exceptional collection $(b^1,d^1,b^{2})$ (see   \eqref{subquiver of ARQ} and Figure \ref{extend subquiver of ARQ}, and \eqref{at most one non-vanisching}) and the category $D^b(Q_1)$ has a full strong exceptional collection $(b^1,a^1,b^2)$ (see \eqref{subquiver of ARQ Q_1}). Let  $(A,B,C)$ be any of these two exceptional triples, then we  have  $\hom(A,B)=\hom(B,C)=1$, $\hom(A,C)=2$, and we will show that  the composition of any nonzero morphisms in $\Hom(A,B)$, $\Hom(B, C)$ is a non-zero morphism in $\Hom(A,C)$, and this implies  that both the triples have the same  endomorphism  algebra  ${\rm End}(E_1 \oplus E_2 \oplus E_3)$, then from   \cite[Corollary 1.9]{Orlov} would follow that $\langle G_+ \rangle^{\perp}$ is equivalent to $D^b(Q_1)$.
In \eqref{abcd} we see that the triple $(b^1[-1],d^1[-1],b^{2}[-1])$ is $(E_3^0, E_8^0, E_3^1)$, and in  \cite[Proposition 2.3]{DK1} we see that these are the following representations :
\begin{gather}
E_3^0 = \begin{diagram}[1em]
0 &  \rTo &  \KK   \\
\uTo &        & \uTo     \\
0  & \rTo^{Id}   &    0
\end{diagram}  \ \ E_8^0 = \begin{diagram}[1em]
\KK &  \rTo^{Id}  &  k   \\
\uTo &        & \uTo     \\
0  & \rTo^{Id}   &    0
\end{diagram} \ \  E_3^1 = \begin{diagram}[1em]
\KK &  \rTo^{j_+^1}  &  \KK^{2}   \\
\uTo^{Id} &        & \uTo^{j_-^1}     \\
\KK  & \rTo^{Id}   &    \KK
\end{diagram}
\end{gather} 
and one checks that the composition of   non-zero morphisms from $E_3^0$ to $E_8^0 $ and from  $E_8^0$ to $E_3^1 $ is non-zero.  In Section 3.2 in  \cite{DK3} we see that $(b^1[-1],a^1[-1],b^2[-1])$ is $(E_3^0,E_2^0,E_3^1)$ and from \cite[Proposition 2.2]{DK1}:
\begin{gather}
E_3^0 = \begin{diagram}[1em]
&       &  \KK &       &    \\
& \ruTo &    & \luTo &       \\
0  & \rTo^{Id}  &    &       &  0
\end{diagram}
\ \  \ E_2^0 =  \begin{diagram}[1em]
&       &  \KK &       &    \\
& \ruTo &    & \luTo^{Id} &       \\
0  & \rTo  &    &       &  \KK
\end{diagram}
\ \  \ 
E_3^1 =\begin{diagram}[1em]
&       &  \KK^{2} &       &    \\
& \ruTo^{j_+^1} &    & \luTo^{j_-^1} &       \\
\KK  & \rTo^{Id}  &    &       &  \KK
\end{diagram}
\end{gather} 
and one checks that the composition of   non-zero morphisms from $E_3^0$ to $E_2^0 $ and from  $E_2^0$ to $E_3^1 $ is non-zero. 
\epr
Having observed \eqref{expl for Q1 in Q2} it is natural to determine what is $\langle x^m \rangle^{\perp}$, for $x \in \{a,b,c,d\}$:
\begin{prop}  \label{triples in Q_2 part 2} Let  $Q_2$ be the quiver in \eqref{Q1}, then: 
	\begin{gather} 
	\label{expl for A3 in Q2} C^{\{\rm Id\}}_{D^b(A_3), \KK}(D^b(Q_2))=\left \{\langle x^m \rangle^{\perp} : m \in \ZZ, x\in \{a,b,c,d\}\right \} 
	\end{gather} and  the cardinalities  $\abs{C^{\Gamma}_{D^b(A_3), \KK}(D^b(Q_2))}$ for $\Gamma=\{\rm Id\}, \langle S \rangle,{\rm Aut}_{\KK}(D^b(Q_2))$ are $\infty,4,1$, respectively. \end{prop} \bpr From the Figures \ref{cube for Q2}, \ref{four pyra 1 for Q2}, \ref{four pyra 2 for Q2} (or using Lemma \ref{exceptional pairs in T2}) one concludes  that for each $m$ we have a full exceptional collection $(c^m,G_-,d^m,a^m)$, hence $\langle a^m \rangle^{\perp} = \langle c^m,G_-,d^m \rangle $, also these pictures contain the data, that after shifts of the elements of $ (c^m,G_-,d^m)$ one gets a strong triple $(A,B,C)$ with $\hom(A,B)=\hom(B,C)=1$, $\hom(A,C)=0$, and now  \cite[Corollary 1.9]{Orlov} ensures that $ \langle a^m \rangle^{\perp} \cong D^b(A_{3})$, and hence $\langle a^m \rangle^{\perp} \in C^{\{\rm Id\}}_{D^b(A_3), \KK}(D^b(Q_2))$. We use again some auto-equivalences of $D^b(Q_2)$. More preccisely, from \eqref{row 1 Serre}, \eqref{autoeq sigma}  we have:

\begin{gather} \label{actions help for A3 in Q2 1}
\begin{array}{| c | c | c | c | c | c | c | c | c |}
\hline
X    & \langle a^m \rangle^{\perp}     & \langle b^m \rangle^{\perp} &  \langle d^m \rangle^{\perp}    &  \langle c^{m-1} \rangle^{\perp}	   &  \langle b^m \rangle^{\perp} & \langle c^m \rangle^{\perp}    \\ \hline
S(X) & \langle b^m \rangle^{\perp} & \langle a^{m-2} \rangle^{\perp} &  \langle c^{m-1} \rangle^{\perp} & \langle d^{m-2} \rangle^{\perp} & \langle c^{m-1} \rangle^{\perp} &  	\langle b^m \rangle^{\perp}.  \\ \hline				\end{array} 
\end{gather} 
which implies that $\langle x^m \rangle^{\perp}\cong D^b(A_3)$ and $\langle x^m \rangle^{\perp} \in C^{\{\rm Id\}}_{D^b(A_3), \KK}(D^b(Q_2))$, for $x \in \{a,b,c,d\}$ and $m \in \ZZ$.

On the other hand if we take any $\mc B \in C^{\{\rm Id\}}_{D^b(A_3), \KK}(D^b(Q_2))$, then viewing it as a subcategory it must be generated by an exceptional triple, say $(A,B,C)$. Due to \cite{WCB1} this triple can be extended to a full exceptional collection   $(A,B,C,D)$ and hence $\mc B= \langle D \rangle^{\perp}$ for some $\langle D \rangle \in C_{D^b(A_1), \KK}(D^b(Q_2))$, and \eqref{expl for A3 in Q2} follows from  $C_{D^b(A_1), \KK}(D^b(Q_2) = \{\langle  F_+ \rangle , \langle  F_- \rangle, \langle  G_+ \rangle , \langle  G_- \rangle \}$  $\cup \{ \langle  a^m \rangle, \langle  b^{m} \rangle, \langle  c^{m} \rangle, \langle  d^{m} \rangle   : m \in \ZZ   \}$  and \eqref{expl for Q1 in Q2}. The rest of the Proposition follows from \eqref{actions help for A3 in Q2 1}.
\epr
\begin{coro} \label{triples in Q2}   Let  $Q_1$, $Q_2$ be the quivers in \eqref{Q1}.
	If $C_{\mc A, \KK}^{\{\rm Id\}}(D^b(Q_2))\neq \emptyset$ and $\mc A$ has a full exceptional triple, then either $\mc A \cong D^b(A_3)$ or  $\mc A \cong D^b(Q_1)$ and   we have table \eqref{triples in Q_2}.	
\end{coro} 
\bpr The fact that $\mc A$ must be either equivalent to $D^b(A_3)$ or to $D^b(Q_1)$ follows from the arguments in the last paragraph of the proof of Proposition  \ref{triples in Q_2 part 2}.   \epr

\section{Non-commutative curve counting  in $D^b(\PP^2)$ and Markov conjecture} \label{PP2 section}

Here we choose $\KK=\CC$ and  initiate our study of the non-commutative curve counting in $D^b(\PP^2)$ and show that the completion of this project leads to proving (or disproving) Markov conjecture.

\begin{prop} \label{Cl(PP2)} Denote $\mc T = D^b(\PP^2)$. Let  $\langle S \rangle \subset  {\rm Aut}_{\CC}(\mc T)$  be the subgroup generated by the Serre functor. Then $C_{-1}^{\{\rm Id \}}(\mc T) =\emptyset$ and  $\forall l\geq 0$ the set $C_{l}^{\langle S \rangle}(\mc T)$ is finite. It is non-empty for: 
	\begin{gather} \label{numbers for P2} \left \{l\geq 0:C_{l}^{\langle S \rangle}(\mc T) \neq \emptyset \right \}= \left \{l\geq 0:C_{l}^{\{\rm Id \}}(\mc T) \neq \emptyset \right \}= \{3 m - 1: m \ \mbox{is a Markov number} \}.
	\end{gather}
	Furthermore, for any  Markov number $m$ hold 
	\eqref{infinite mumbers for PP2}, 
	\eqref{finite mumbers for PP2 full group}, 	\eqref{finite mumbers for PP2}, 
	where   $c_1(E)$, $r(E)$ are   the first Chern class (which we consider as an integer)  and the rank of $E$.
\end{prop}
\bpr
	We first collect known facts for exceptional collections in $\mc T$.  Let $\{\mc O(i)\}_{i\in \ZZ}$ be the sequence of line bundles on $\PP^2$.    Beilinson   \cite{Beili} gave a full exceptional collection  $(\mc O(i), \mc O(i+1), \mc O(i+2))$ in $D^b(\PP^2)$  for eqch $i\in \ZZ$.  The full exceptional collections in $\mc T$ have length 3  and ${\rm rank}(K_0(\mc T))=3$.

	\cite[Proposition 2.10]{OrlovKuleshov} ensures that each exceptional object $E \in \mc T$ after shifting  becomes an exceptional sheaf, on the other hand each exceptional sheaf on $\PP^2$ is locally free (see \cite[Subsection 1.1]{GoroRuda}). Thus \textit{each exceptional object in $\mc T$ is a shift of an exceptional vector bundle}.  Now using \cite[Theorem 4.6.1]{GorKul}  we conclude that:
	\begin{gather} \label{extension to full} \begin{array}{l} \mbox{ any exceptional collection  in $\mc T$ can be extended to
	a full collection,} \\ \mbox{ after shifting any full exceptional collection is a collection of 3 vector bundles.} \end{array}
	\end{gather}

		In \cite{BP} a strong exceptional collection  $\mc E$ which  remains strong under all mutations is called \textit{non-degenerate}.   Furthermore  in \cite{BP} are  defined  so called \textit{ geometric}  exceptional collections and \cite[Corollary 2.4]{BP} says that geometricity implies non-degeneracy.  Furthermore, \cite[Proposition 3.3]{BP} 
		claims that a full exceptional collection of length $m$ of coherent sheaves on a smooth projective variety $X$ of dimension $n$ is geometric if and only if $m=n+1$.  Combining with \eqref{extension to full} we conclude:  
		\begin{gather} \label{strong pairs}
		\mbox{each exceptional pair can be shifted to a strong exceptional pair} \  (E_1,E_2) \  \mbox{ of bundles}.
		\end{gather}
	Recall that we have  a bijection \ref{bijection 1}. Take any $\mc A =\langle E_1,E_2 \rangle \subset \mc T$ from the codomain of this bijection, in particular $\hom(E_1,E_2)=l+1$  (from \eqref{strong pairs} we know that $E_1,E_2$ can be chosen so that they form an exceptional pair of bundles). Using  \eqref{extension to full} we can extend the pair to a full exceptional collection  of vector bundles $(E_1,E_2,E_3)$ and so we conclude that there exists unique up to isomorphism exceptional vector bundle $E_3$ such that $^{\perp}\mc A=\langle E_3\rangle$. On the other hand   from \cite[Subsection 7.2.2]{GorKul} we know that  $\hom(E_1,E_2)=3 r(E_3)$, and that \begin{gather} \label{Markov equation} r(E_1)^2+r(E_2)^2+r(E_3)^2= 3 r(E_1) r(E_2) r(E_3)\end{gather} where $r(E_i)$ is the rank of $E_i$ (this holds  for each full exceptional collection $(E_1,E_2,E_3)$). So to  each $\mc A$ we attached a subcategory generated by an exceptional vector bundle $E$ with $r(E)=(l+1)/3$, so that $\mc A=\langle E \rangle^\perp$, in particular this assignment is injective. Using  \eqref{strong pairs}, \eqref{extension to full} one can show that the assignment is also surjective, more precisely, we obtain a bijection:
	  From  \eqref{extension to full}, \eqref{strong pairs} it follows that we have a bijection, we choose to denote the codomain by $EB_l$
	 	\begin{gather}
	 	\label{bijection 2} 	C_{l}^{\{\rm Id\}}(\mc T) \rightarrow\left  \{\langle E \rangle\subset  \mc T:\begin{array}{l} E \ \mbox{is an exc. v. bundle with} \ r(E)=\frac{l+1}{3}  \end{array} \right \}:=EB_l \ \  \ 
	 	[F] \mapsto  {}^\perp{\rm Im}(F). 
	 	\end{gather}
	 	In particular we see that $C_{-1}^{\{\rm Id \}}(\mc T)=\emptyset$. 
 The  group action $\Gamma \times C_{\mc A, \CC}^{\{\rm Id\}}(\mc T) \rightarrow  C_{\mc A, \CC}^{\{\rm Id\}}(\mc T)$ defined by $([\beta],[F]_{\{\rm Id\}}) \mapsto [\beta \circ F]_{\{\rm Id\}}$, from Corollary \ref{group action} translates to action (see Remark \ref{auto eq on perp} ) \begin{gather} \label{gr action} \Gamma \times EB_l \rightarrow  EB_l \ \ ([\beta],\langle E \rangle ) \mapsto \langle \beta (E) \rangle \qquad C_l^\Gamma(\mc T)\cong EB_l/\Gamma.  \end{gather}Having the bijection \eqref{bijection 2}, the equation \eqref{Markov equation} and knowing that every Markov number is a rank of an exceptional vector bundle (\cite[Theorem 3.2]{Rudakov1}), we obtain \eqref{numbers for P2}, and from now on we assume that  \begin{gather} \label{Markov number} m=\frac{l+1}{3} \ \mbox{is a Markov number.}\end{gather} 
 
 For a vector bundle $E$ on $\PP^2$ we denote $E\otimes \mc O(n)$ by $E(n)$, and the rational number $\frac{c_1(E)}{r(E)}$ is denoted by $\mu(E)$ and called slope of $E$. We have $c_1(E(n))=c_1(E)+n r(E)$, $\mu(E(n))=\mu(E)+n$ for $n\in \ZZ$ (see for example \cite[p.100]{Rudakov1}).
For each $x\in EB_l$ there is unique up to isomorphism exceptional vector bundle $E_x$ on $\PP^2$ with $r(E_x)=\frac{l+1}{3}$, s.  t. $x=\langle E_x \rangle$, we will write $c_1(x)$ for   $c_1(E_x)$, and $x(n)$ for $\langle E_x(n) \rangle$. 
  \cite[Proposition 4.6]{GoroRuda} or \cite[Theorem 2.2]{Rudakov1}  say that: if   $E, E'$ are exceptional vector bundles s.t.  $\mu(E)=\mu(E')$  then \ $E\cong E'$, therefore
  \begin{gather} \label{from equal chern equal things}
\mbox{if} \  x, y \in EB_l \ \mbox{then} \ c_1(x)=c_1(y) \iff x=y;  \qquad \bd EB_l & \rTo^{c_1} & \ZZ \ed \ \ \mbox{is injective}.
  \end{gather}
Since we have $c_1(x(n))=c_1(x)+n m$, then denoting:   \begin{gather}  \label{disj union} EB_l^{[0,m)}=\{x\in EB_l: 0\leq c_1(x)<m\} \ \ \mbox{we have a disjoint union}\ \
EB_l=\bigcup_{n\in \ZZ} EB_l^{[0,m)}(n) 
\end{gather} 
we already showed that $EB_l\neq \emptyset$ (under the assumption \ref{Markov number}), and from \eqref{disj union} it follows that $EB_l$ is an infinite set (hence \eqref{infinite mumbers for PP2}). Now let $[S]\in \Aut(\mc T)$ be the Serre functor of $\mc T$, since the canonical bundle on $\PP^2$ is $\mc O(-3)$, then for a vector bundle $E$ we have $S(E)=E\otimes \mc O(-3)[2]$ (see \cite[2) in  3.2  Examples]{BK}), hence the action of $[S]$ on $EB_l$ \eqref{gr action} is:
\begin{gather}
\forall x \in EB_l \qquad ([S],x ) \mapsto x(-3) \ \Rightarrow \ \left ([S], EB_l^{[0,m)}(n) \right ) \mapsto  EB_l^{[0,m)}(n-3),
\end{gather}
therefore  from \eqref{gr action}, \eqref{disj union} we see that: 
$ \abs{C_l^{\langle S \rangle}(\mc T)}=\abs{EB_l/\langle S \rangle} = 3 \abs{EB_l^{[0,m)}}
$ and   part of \eqref{finite mumbers for PP2} follows.  

Finally, let $\Gamma=\Aut_{\CC}(\mc T)$, we want to compute $EB_l/\Gamma$ (see \eqref{gr action}). \cite[Theorem 3.1]{BondalOrlov} says that $\Aut_{\CC}(\mc T)$ is generated by automorphisms of the variety, translations, and by the functor $F=\_\otimes \mc O(1)$.  Since the exceptional vector bundles are homogeneous,\cite[Subsection 1.1]{GoroRuda}, then the automorphisms of the variety and the translations act trivially on $EB_l$, therefore $EB_l/\Gamma=EB_l/\langle F \rangle$. On the other hand, since $F(E(n))=E(n+1)$ for any vector bundle $E$, it follows
\begin{gather}
\forall x \in EB_l \qquad ([F],x ) \mapsto x(1) \ \Rightarrow \ \left ([F], EB_l^{[0,m)}(n) \right ) \mapsto  EB_l^{[0,m)}(n+1),
\end{gather} and therefore  $ \abs{C_l^{\Aut_\CC(\mc T)}(\mc T) }=\abs{EB_l/\langle F \rangle} =  \abs{EB_l^{[0,m)}}
$ and we deduce \eqref{finite mumbers for PP2 full group} and the equality in \eqref{finite mumbers for PP2}.
\epr
\begin{coro} \label{corollary for Markov}  Denote $\mc T = D^b(\PP^2)$. The first several non-trivial $\abs{C_{3m-1}^{\Aut_{\CC}(\mc T)}(\mc T)}$ are (recall that $m$ is a Markov number and on the first row are listed the first 9 Markov numbers) listed in \eqref{table for PP2}. 

Furthermore, the so called Tyurin's conjecture, which is equivalent to the Markov's conjecture  (this equivalence is proved in \cite[p. 100]{Rudakov1}, see also \cite[Section 7.2.3 ]{GorKul}) is equivalent to  the following statement:   for all  Markov numbers $m\neq1, m\neq 2$ we have $\abs{C_{3 m -1}^{\Aut_{\CC}(\mc T)}(\mc T)}\leq 2$ and this is equivalent to the statement:  for all  Markov numbers $m\neq 1, m\neq 2$ we have $\abs{C_{3 m -1}^{\Aut_{\CC}(\mc T)}(\mc T)}= 2$.
	\end{coro}
\bpr We first recall that, \cite[Proposition 4.6]{GoroRuda} or \cite[Theorem 2.2]{Rudakov1}, if   $E, E'$ are exceptional vector bundles s.t.  $\mu(E)=\mu(E')$  then \ $E\cong E'$.
In \cite[p. 100]{Rudakov1} a vector bundle on $\PP^2$ is called representative (``представительньiми''), if $0\leq \mu(E)\leq \frac{1}{2}$. Since $\mu(E(n))=\mu(E)+n$, and $\mu(E^{\star})=-\mu(E)$, and since starting from an exceptional vector bundle $E$ the dual $E^{\star}$ and the twisted bundles $E(j)$ are exceptional bundles as well (see e.g. \cite[Subsection 1.2]{GoroRuda}) it follows that all the exceptional vector bundles up to isomorphism) on $\PP^2$ are obtained by the representative exceptional vector bundles using dualisation and twisting. Taking into account these remarks, we can rewrite  \eqref{finite mumbers for PP2 full group} as follows (here $m$ is a Markov number):
\begin{gather} 
\abs{C_{3 m -1}^{\Aut_{\CC}(\mc T)}(\mc T)}= \abs{\left  \{ \mu(E):\begin{array}{l} E \ \mbox{is representative exceptional} \nonumber \\   \mbox{ bundle on} \ \PP^2, \mbox{with} \  r(E)= m, \end{array} \right \}} \\[-2mm] \label{formula for the numbers of PP2} \\[-2mm]  + \abs{\left  \{ \mu(E):\begin{array}{l} E \ \mbox{is representative exceptional bundle}\\    \mbox{  on} \ \PP^2, \mbox{with} \  r(E)= m, \ \mbox{and} \  0<\mu(E)<\frac{1}{2}\end{array} \right \} }\nonumber
\end{gather}
Furthermore, \cite[Proposition 4.4]{GoroRuda} says that for an exceptional vector bundle $E$ the numbers  $c_1(E)$, $r(E)$ have no common divisors, hence they can be found from $\mu(E)$. On \cite[p. 100]{Rudakov1} is written, that  all the possible slopes of representative exceptional bundles of rank $\leq 200$ are:
\begin{gather}\label{liest of slopes} \frac{0}{1}, \frac{1}{2}, \frac{2}{5},  \frac{5}{13}, \frac{12}{29},  \frac{13}{34},\frac{34}{89},\frac{70}{169},\frac{75}{194},\end{gather} from these arguments it follows table \ref{table for PP2}.
 Furthermore, Tyurin's conjecture (\cite[p. 100]{Rudakov1}) says that a representative exceptional bundle is uniquely determined by its rank. Knowing that for each Markov number $m$ there exists such a vector bundle with rank $m$, and since the slopes $0$, $\frac{1}{2}$ are already in the list  \eqref{liest of slopes}, then  for $m\neq 1$ and $m\neq 2$     such a vector bundle  $E$ has  $0<\mu(E)<\frac{1}{2}$, and  we see that the Turin's conjecture is equivalent to the statement that for each Markov number $m\neq 1, \neq 2$ there exists no more than one (up to isomorphism)  representative exceptional vector bundle $E$ with $r(E)=m$. Now using \eqref{formula for the numbers of PP2} and that   $0<\mu(E)<\frac{1}{2}$ for such $E$  we deduce the last statement of the corollary.
\epr

It follows an approach to Markov conjecture by counting subgraphs in a given valued digraph:

\begin{prop} \label{markov with graphs} If for any Markov number $m\neq 1$ and $m\neq 2$ one can find two subgraphs  of the valued graph $G_{D^b(A_1), \CC}(D^b(\PP^2))$ of the form:
	
	\begin{tikzpicture} 
	\node [left] at (-7,0) {$\cdots$};
	\draw[->]  (-7,0)--(-5,0);	
	\node [above] at (-6,0) {$3 m $};	
	\draw[->]  (-3.8,0)--(-1.8,0);
	\node [above] at (-2.8,0) {$3 m $};	
	\node [left] at (-3.8,0) {$\langle x_{i-1} \rangle $};
	\draw[->]  (-1,0)--(+1,0);
	\node [left] at (-1 ,0) {$\langle x_i \rangle $};
	\node [right] at (+1 ,0) {$\langle x_{i+1} \rangle $};
	\draw[->]  (2.2,0)--(4.2,0);
	\node [above] at (3.2,0) {$3 m $};	
	\node [right] at (4.2,0) {$\cdots $};	
	\node [above] at (0,0) {$3 m $};	
	
	\end{tikzpicture}

	\begin{tikzpicture} 
	\node [left] at (-7,0) {$\cdots$};
	\draw[->]  (-7,0)--(-5,0);	
	\node [above] at (-6,0) {$3 m $};	
	\draw[->]  (-3.8,0)--(-1.8,0);
	\node [above] at (-2.8,0) {$3 m $};	
	\node [left] at (-3.8,0) {$\langle y_{i-1} \rangle $};
	\draw[->]  (-1,0)--(+1,0);
	\node [left] at (-1 ,0) {$\langle y_i \rangle $};
	\node [right] at (+1 ,0) {$\langle y_{i+1} \rangle $};
	\draw[->]  (2.2,0)--(4.2,0);
	\node [above] at (3.2,0) {$3 m $};	
	\node [right] at (4.2,0) {$\cdots $};	
	\node [above] at (0,0) {$3 m $};	
	
	\end{tikzpicture}
	
	such that for any   subgraph of $G_{D^b(A_1), \CC}(D^b(\PP^2))$  of the form: 
	
	\begin{tikzpicture} 
	\node [left] at (-7,0) {$\cdots$};
	\draw[->]  (-7,0)--(-5,0);	
	\node [above] at (-6,0) {$3 m $};	
	\draw[->]  (-3.8,0)--(-1.8,0);
	\node [above] at (-2.8,0) {$3 m $};	
	\node [left] at (-3.8,0) {$\langle z_{i-1} \rangle $};
	\draw[->]  (-1,0)--(+1,0);
	\node [left] at (-1 ,0) {$\langle z_i \rangle $};
	\node [right] at (+1 ,0) {$\langle z_{i+1} \rangle $};
	\draw[->]  (2.2,0)--(4.2,0);
	\node [above] at (3.2,0) {$3 m $};	
	\node [right] at (4.2,0) {$\cdots $};	
	\node [above] at (0,0) {$3 m $};	
	
	\end{tikzpicture}
	
	there exists $\beta \in \Aut_{\CC}(D^b(\PP^2))$ such that either $\beta(\langle z_i \rangle) =  \langle x_{i+N} \rangle$ for each $i\in \ZZ$ or  $\beta(\langle z_i \rangle) =  \langle y_{i+N} \rangle$ for each $i\in \ZZ$ for some $N \in \ZZ$, then Markov conjecture holds. 
\end{prop}
\bpr Let $\mc A=N\PP^{3m-1}$ and  
$\Gamma=\Aut_{\CC}(D^b(\PP^2))$.  If the condition holds, then the right hand side of \eqref{injection graphs in formula}, modified as in Subsection \ref{the case J eq DbA1} (b), has no more than  two elements, therefore$\abs{ C_{3m-1}^{\rm \Gamma }(\mc T)} \leq 2$ and by Corollary \ref{corollary for Markov} Markov conjecture follows.
\epr

\section{Non-commutative  counting  in $D^b(A_{n+1})$}  \label{ncc in A_n+1}
 Here again $\KK$ is an arbitrary algebraically closed field.   Let  $\mc T_n = D^b(Rep_{\KK}(Q_n))$,   where $Q_n$ is the  quiver $Q_n= \bd[1em] 0 & \rTo & 1 & \rTo & \cdots & \rTo & n-1 & \rTo & n \ed$,   $n\geq 0$.

Here we study the sets  $C^{ {\rm Aut}_{\KK}(\mc T_n)}_{\mc A,\KK}({\mc T_n}) $, $C^{ \{\rm Id\}}_{\mc A,\KK}({\mc T_n}) $  for   $\mc A = \mc T_{k-1}\cong D^b(A_k)$, $k\geq 1$  (see Definition \ref{C_A KK(T)}) and the sets  $C_l^{\rm Id}(\mc T_n)$,  $C_l^{ {\rm Aut}_{\KK}(\mc T_n) }(\mc T_n)$ for $l\geq -1$. 
  In \cite{Miyachi} is shown that  ${\rm Aut}_\KK(\mc T_n)$ is generated by the translation  and the Serre functor,  it follows that for the choices of $\mc A$ in this  section we have \begin{gather} \label{Serre as the entire for Dynkin} C^{{\rm Aut}_{\KK}(\mc T_n)}_{\mc A, \KK}(\mc T_n) \cong  C^{\langle S \rangle}_{\mc A,\KK}(\mc T_n). \end{gather} From now on we denote  $\Gamma_n = {\rm Aut}_{\KK}(\mc T_n)$ and
for shorter formulas we will omit $\KK$ in the notations  $C^{\{\rm Id\}}_{\mc T_{k-1}, \KK}(\mc T_n)$ and $C^{\Gamma_n}_{\mc T_{k-1},\KK}(\mc T_n)$ and write just   $C^{\{\rm Id\}}_{\mc T_{k-1}}(\mc T_n)$ and $C^{\Gamma_n}_{\mc T_{k-1}}(\mc T_n)$.

Sometimes we denote the value of a sequence $\alpha: \{0,1,\dots,k\}\rightarrow X$ at $i\in \{0,1,\dots,k\}$ via $\alpha_i$ for the sake of  brevity, otherwise we denote it by $\alpha(i)$. 

 For $0\leq i\leq j \leq n$ we denote by $s_{i,j}$ the representation in $Rep_{\KK}(Q_n)$,  whose dimension vector satisfies $\ul{\dim}_x(s_{i,j})=1$ iff $i\leq x\leq j$   and  $\ul{\dim}_x(s_{i,j})=0$ iff $x<i$ or $x>j$, and which assigns the identity map to the arrows $\bd[1em] x & \rTo & y \ed$ for $i\leq x$, $y\leq j$, and zero to all other arrows. It is known that:  
 \begin{gather}\label{exceptional objects} \mbox{\rm up to shifts, all exceptional objects in} \  \mc T_n \ \mbox{\rm are } \  \{s_{i,j}\}_{0\leq i\leq j \leq n }, \quad s_{i,j} \sim s_{u,v} \iff i=u, j=v. \end{gather}
Using the bijection \eqref{bijection for eo} in Proposition \ref{bijection 123} and Lemma  \ref{ci}  we see that: 

  \begin{lemma} \label{C_eoA_n+1}
  Let $n\geq 0$, then	the set $C^{\{\rm Id\}}_{\mc T_0, \KK}(\mc T_n)$ can be described as follows 
  	\begin{gather}C^{\{\rm Id\}}_{\mc T_0, \KK}(\mc T_n) =\{\langle  s_{\alpha,\beta} \rangle  : 0\leq \alpha \leq \beta\leq n \} \qquad \abs{C^{\{\rm Id\}}_{\mc T_0, \KK}(\mc T_n) } = \frac{(n+1)(n+2)}{2}. \end{gather}
  \end{lemma}
 In particular $\abs{C^{\{\rm Id\}}_{\mc T_0, \KK}(\mc T_0) } =1$ and we can assume from now on $n\geq 1$.
 The following remark will be useful to determine the exceptional pairs:
 
 \begin{remark} \label{hom and hom1} Recall that for any acyclic $Q$ and any   $A$, $B \in Rep_\KK(Q)$ holds the Euler formula: 
\begin{gather} \label{Euler}  \hom(A,B) - \hom^1(A,B) = \sum_{v\in V(Q)} \ul{\dim}_v(A)  \ul{\dim}_v(B) - \sum_{\bd[1em] x & \rTo & y \ed \in Arrows(Q)} \ul{\dim}_x(A)  \ul{\dim}_y(B).  \end{gather}

 Furthermore if $Q$ is Dynkin, then using  the table in \cite[page 5 in Version 3]{DK2} we see that  for any pair of exceptional representations $(A,B)$ we have 
 \begin{gather} \hom(A,B) = \left \{ \begin{array}{c c} \langle \ul{\dim}(A), \ul{\dim}(B) \rangle & \mbox{if} \ \langle \ul{\dim}(A), \ul{\dim}(B) \rangle >0 \\ 0 & \mbox{otherwise} \end{array}\right .
 \\
 \hom^1(A,B) = \left \{ \begin{array}{c c} -\langle \ul{\dim}(A), \ul{\dim}(B) \rangle & \mbox{if} \ \langle \ul{\dim}(A), \ul{\dim}(B) \rangle <0 \\ 0 & \mbox{otherwise} \end{array}\right .
 	\end{gather}
in particular $(A,B)$ is  exceptional iff $\langle \ul{\dim}(B), \ul{\dim}(A) \rangle =0 $.
 \end{remark}
 
  using this remark one proves Lemmas \ref{help lemma 1}, \ref{noEP}: 
 \begin{lemma} \label{help lemma 1}  Let $(s_{\alpha,\beta}, s_{i,j})$ be an exceptional pair in $\mc T_n$.
 	
 	{\rm (a)} If $\alpha < i$ and  $\beta < j$, then $\beta < i$, furthermore, if $\beta = i-1$, then  $\hom^1(s_{\alpha,\beta}, s_{i,j})=1$, $\hom^i(s_{\alpha,\beta}, s_{i,j})=0$ for $i\neq 1$ otherwise $(s_{\alpha,\beta}, s_{i,j})$ is orthogonal pair,
 	
 		{\rm	(b) } If $i< \alpha $ and  $j<\beta$, then $j<  \alpha -1$ and  $(s_{\alpha,\beta}, s_{i,j})$ is orthogonal pair,
 		
 		{\rm	(c) }  If $ \alpha \leq i\leq j \leq \beta $, then either ( $\alpha = i$ and $j<\beta $ and  $\hom(s_{\alpha,\beta}, s_{i,j})=1$, $\hom^i(s_{\alpha,\beta}, s_{i,j})=0$ for $i\neq 0$ ) or ( $\alpha < i$ and $j<\beta $ and  $(s_{\alpha,\beta}, s_{i,j})$   is orthogonal pair ),
 		
 		{\rm	(d) } If $ i \leq \alpha \leq \beta \leq j $, then either ( $i <\alpha $ and $j=\beta $ and  $\hom(s_{\alpha,\beta}, s_{i,j})=1$, $\hom^i(s_{\alpha,\beta}, s_{i,j})=0$ for $i\neq 0$ ) or ( $i<\alpha$ and $\beta < j$ and  $(s_{\alpha,\beta}, s_{i,j})$   is orthogonal pair ),
 
 \end{lemma}
 
 \begin{lemma} \label{noEP}
 	If 	$\alpha < i\leq \beta < j$, then neither $(s_{\alpha, \beta}, s_{i,j})$ nor $( s_{i,j}, s_{\alpha, \beta})$ is an exceptional pair. If 	$i<\alpha \leq j< \beta $, then neither $(s_{\alpha, \beta}, s_{i,j})$ nor $( s_{i,j}, s_{\alpha, \beta})$ is an exceptional pair.
 \end{lemma}

\begin{remark} \label{subcats generated by pairs Dynkin} Note that from Remark \ref{hom and hom1}  it follows that property \eqref{at most one non-vanisching} holds for $D^b(Q)$ with a Dynkin quiver $Q$. This property and Proposition \ref{bijection 123} imply that if $C_{\mc A, \KK}^{\{\rm Id \}}(D^b(Q_i)) \neq \emptyset$ and $\mc A$ is a category generated by an exceptional pair, then  $C_{\mc A, \KK}^{\{\rm Id \}}(\mc T_i)$ must be  $C_{l}^{\{\rm Id \}}(\mc T_i)$ for some $l\geq -1$. Furthermore, \eqref{vanishings general 2} implies that $l\leq 0$.  
\end{remark}
 
 \begin{lemma} \label{help lemma 2} For $0\leq \alpha \leq \beta < j \leq n $  the pairs  $(s_{\beta+1,j},s_{\alpha,j})$, $(s_{\alpha,\beta},s_{\beta+1,j})$, $(s_{\alpha,j}, s_{\alpha,\beta})$ are exceptional and $\langle s_{\beta+1,j},s_{\alpha,j} \rangle = \langle  s_{\alpha,\beta},s_{\beta+1,j} \rangle = \langle s_{\alpha,j}, s_{\alpha,\beta}\rangle $ 
  \end{lemma}
 \bpr  The fact that the pairs are exceptional one proves using \eqref{Euler}. If we denote $(A,B)=(s_{\alpha,\beta},s_{\beta+1,j})$, then one shows that $(B,R_B(A))\sim (s_{\beta+1,j},s_{\alpha,j}) $ and that  $(L_A(B),A)\sim (s_{\alpha,j}, s_{\alpha,\beta})$, hence the lemma follows. 
 \epr

  \subsection{Combinatorial descriptions of $C^{\{\rm Id\}}_{\mc T_{k-1}}(\mc T_n)$ and $C^{\Gamma_n}_{\mc T_{k-1}}(\mc T_n)$}
  \begin{lemma} \label{Cid_A_kinA_n+1}
  	The set $C^{\{\rm Id\}}_{\mc T_{k-1}}(\mc T_n)$ can be described as follows: $C^{\{\rm Id\}}_{\mc T_{k-1}}(\mc T_n) =$ 
  	\begin{gather}\nonumber  \{\langle  s_{\alpha_0,\alpha_1}, s_{\alpha_1+1,\alpha_2+1},\dots, s_{\alpha_j+j,\alpha_{j+1}+j},\dots, s_{\alpha_{k-1}+k-1,\alpha_{k}+k-1}\rangle  : 0\leq \alpha_0 \leq\alpha_1 \dots \leq \alpha_{k}  \leq n+1-k \}\\ \nonumber \cong  \left \{ \bd \{0,1,\dots,k\} & \rTo^{\alpha} &  \{0,1,\dots,n-k+1\} \ed : 0\leq \alpha_0 \leq \alpha_1 \leq \dots \leq \alpha_k \leq n+1-k \right \}. \end{gather}
  \end{lemma}
  \bpr  The case $k=1$ is in Lemma  \ref{C_eoA_n+1}. 
  
  Assume that the statement holds for some  $k\geq 1$. 
  
  Let $\mc B \in C^{\{\rm Id\}}_{\mc T_{k}}(\mc T_{n})$. We will show first  that $\mc B$ must have the desired form.  Now  $\mc B \subset \mc  T_{n}$,  $\mc B \cong \mc T_{k}$, and there is a SOD $\mc B =\langle \mc B', E \rangle$, where $\mc B'\cong \mc T_{k-1}$, and $E$ is an  exceptional object and the decompositions are not orthogonal. By the induction assumption \\ $\mc B' = \langle  s_{\alpha_0,\alpha_1}, s_{\alpha_1+1,\alpha_2+1},\dots, s_{\alpha_j+j,\alpha_{j+1}+j},\dots, s_{\alpha_{k-1}+k-1,\alpha_{k}+k-1}\rangle $ 
  and  $E=s_{u,v}$, and we can write:
  \begin{gather} 0\leq \alpha_0 \leq\alpha_1 \leq \dots \leq \alpha_{k}  \leq n+1-k \qquad 0\leq u\leq v \leq n, \\
  \label{excpairs} ( s_{\alpha_j+j,\alpha_{j+1}+j}, s_{u,v}) \ \ \mbox{are  exc. pairs for all} \  0\leq j \leq k-1, \  \mbox{not all pairs are orthogonal}. \end{gather}
  
  If $u\geq \alpha_{k}+k-1$, then  Lemma \ref{help lemma 1}  and \eqref{excpairs} imply that $u=\alpha_{k}+k\leq v$, and then we obtain the desired form of $\mc B$. If $v\leq \alpha_0$, then using     Lemma \ref{help lemma 1} and the fact that  $(s_{\alpha_0,\alpha_1}, s_{u,v})$ is an exceptional pair  we see that $u$ must be equal to $\alpha_0$ and $v<\alpha_1$ and using Lemma \ref{help lemma 2} we get $\langle s_{\alpha_0,\alpha_1}, s_{u,v} \rangle =\langle s_{\alpha_0,\alpha_1}, s_{\alpha_0,v} \rangle  = \langle s_{\alpha_0,v}, s_{v+1,\alpha_1} \rangle$  and then we obtain the desired form of $\mc B$.
  
  Therefore we can assume that $u<\alpha_{k}+k-1$, $v>\alpha_0$.  Now if $v>  \alpha_{k}+k-1$, then Lemmas \ref{help lemma 1} (d), \ref{noEP}, the inequality  $u<\alpha_{k}+k-1$, and since $(s_{\alpha_j+j,\alpha_{j+1}+j}, s_{u,v})$ is an exceptional pair for all $0\leq j \leq k-1$ imply that $u<\alpha_0$ (one shows that $u\geq \alpha_0$ is impossible), but now  Lemma \ref{help lemma 1} (c), (d) implies that all pairs are orthogonal, which contradicts \eqref{excpairs}.  
  
  Therefore we can assume that $u<\alpha_{k}+k-1, \alpha_0 < v\leq \alpha_{k}+k-1$.

  Now we consider the case $\alpha_0 \leq  u$, then  there exists    $0\leq j\leq k-1$, s.t. $\alpha_j +j \leq u \leq \alpha_{j+1} +j$. Now one of these must hold
  
  (i)   $\alpha_j+j< u$, $v=\alpha_{j+1}+j$, however now   $(s_{\alpha_{j}+j,\alpha_{j+1}+j}, s_{u,v})$ is not an exceptional pair (see Lemma \ref{help lemma 1} (c)), which contradicts \eqref{excpairs}.  
  
  (ii) $\alpha_j+j< u$, $v<\alpha_{j+1}+j$, then all pairs in \eqref{excpairs} are orthogonal, which is a contradiction.
  
  (iii)  $\alpha_j+j< u$, $v>\alpha_{j+1}+j$, however now  $(s_{\alpha_{j}+j,\alpha_{j+1}+j}, s_{u,v})$ is not an exceptional pair (see Lemma \ref{noEP}), which contradicts \eqref{excpairs}.

  (iv)  $u=\alpha_j+j$. Now due to Lemma \ref{help lemma 1} the condition that  $(s_{\alpha_{j}+j,\alpha_{j+1}+j}, s_{u,v})$ is  an exceptional pair imposes that $v<\alpha_{j+1}+j$ and therefore $s_{u,v}$ is orthogonal to $s_{\alpha_{t}+t,\alpha_{t+1}+t}$ for $t>j$ and using mutations (see also Lemma \ref{help lemma 2}) we can derive the desired form of $\mc B$: 
  \begin{gather}
  \mc B = \left  \langle  s_{\alpha_0,\alpha_1}, s_{\alpha_1+1,\alpha_2+1},\dots, s_{\alpha_j+j,\alpha_{j+1}+j},\dots, s_{\alpha_{k-1}+k-1,\alpha_{k}+k-1},s_{u,v}\right \rangle \nonumber =\\
  \left  \langle  s_{\alpha_0,\alpha_1}, s_{\alpha_1+1,\alpha_2+1},\dots, s_{\alpha_j+j,\alpha_{j+1}+j},s_{u,v},s_{\alpha_{j+1}+j+1,\alpha_{j+2}+j+1}\dots, s_{\alpha_{k-1}+k-1,\alpha_{k}+k-1} \right \rangle  \nonumber \\ 
  = \left  \langle  s_{\alpha_0,\alpha_1}, s_{\alpha_1+1,\alpha_2+1},\dots, s_{\alpha_j+j,v},  s_{v+1,\alpha_{j+1}+j},s_{\alpha_{j+1}+j+1,\alpha_{j+2}+j+1}\dots, s_{\alpha_{k-1}+k-1,\alpha_{k}+k-1} \right \rangle.  \nonumber\end{gather}
  It remains to consider the case   $u<\alpha_{0}, \alpha_0 < v\leq \alpha_{k}+k-1$. There exists    $0\leq j\leq k-1$, s.t. $\alpha_j +j \leq v \leq \alpha_{j+1} +j$. Since $u<\alpha_{0}$, then $ v < \alpha_{j+1} +j$ is incompatible with the fact that $(s_{\alpha_{j}+j,\alpha_{j+1}+j}, s_{u,v})$ is  an exceptional pair and Lemma \ref{noEP}. Thus $ v = \alpha_{j+1} +j$.  Now if $j<k-1$, then $(s_{\alpha_{j+1}+j+1,\alpha_{j+2}+j+1}, s_{u,v})$ would not be an exceptional pair, which is a contradiction. Therefore we reduce to the case  $u<\alpha_{0},v= \alpha_{k} +k-1 $.
  
  Now we apply $k$ times Lemma \ref{help lemma 2} to  deduce  the desired form of $\mc B$:
  \begin{gather}
  \mc B = \left  \langle  s_{\alpha_0,\alpha_1}, s_{\alpha_1+1,\alpha_2+1},\dots, s_{\alpha_j+j,\alpha_{j+1}+j},\dots, s_{\alpha_{k-1}+k-1,\alpha_{k}+k-1}, s_{u,v}\right \rangle \nonumber =\\
  \left  \langle  s_{\alpha_0,\alpha_1}, s_{\alpha_1+1,\alpha_2+1},\dots, s_{\alpha_j+j,\alpha_{j+1}+j},s_{\alpha_{j+1}+j+1,\alpha_{j+2}+j+1}\dots,    
  s_{u,\alpha_{k-1}+k-2}, s_{\alpha_{k-1}+k-1,\alpha_{k}+k-1} \right \rangle  \nonumber \\ 
  \cdots = \left  \langle s_{u,\alpha_0-1}, s_{\alpha_0,\alpha_1}, s_{\alpha_1+1,\alpha_2+1},\dots, s_{\alpha_j+j,\alpha_{j+1}+j},s_{\alpha_{j+1}+j+1,\alpha_{j+2}+j+1}\dots, s_{\alpha_{k-1}+k-1,\alpha_{k}+k-1} \right \rangle.  \nonumber\end{gather}

  Conversely, if $\mc B \subset \mc T_n$ has the form: $$\mc B= \left  \langle  s_{\alpha_0,\alpha_1}, s_{\alpha_1+1,\alpha_2+1},\dots, s_{\alpha_j+j,\alpha_{j+1}+j},\dots, s_{\alpha_{k-1}+k-1,\alpha_{k}+k-1}, s_{\alpha_{k}+k,\alpha_{k+1}+k}\right \rangle, $$ then denoting $E_i=s_{\alpha_j+j,\alpha_{j+1}+j}[j]$ for $0\leq j \leq k$ we get a strong full exceptional collection (here we use again Lemma \ref{help lemma 1}) in  $\mc B$ which with the help of \cite[Corollary 1.9]{Orlov} ensures that $\mc B \cong D^b(A_{k+1})\cong \mc T_k$, and hence $\mc B \in C^{\{\rm Id\}}_{\mc T_{k}}(\mc T_{n})$.
  
  To prove the lemma completely, it remains to show that if 
  \begin{gather}
  \left  \langle  s_{\alpha_0,\alpha_1}, s_{\alpha_1+1,\alpha_2+1},\dots, s_{\alpha_j+j,\alpha_{j+1}+j},\dots, s_{\alpha_{k-1}+k-1,\alpha_{k}+k-1}, s_{\alpha_{k}+k,\alpha_{k+1}+k}\right \rangle \nonumber  \\ = \left  \langle  s_{\alpha'_0,\alpha'_1}, s_{\alpha'_1+1,\alpha'_2+1},\dots, s_{\alpha'_j+j,\alpha'_{j+1}+j},\dots, s_{\alpha'_{k-1}+k-1,\alpha'_{k}+k-1}, s_{\alpha'_{k}+k,\alpha'_{k+1}+k}\right \rangle \nonumber
  \end{gather}
  then the sequences $\alpha$, $\alpha'$ coincide, i.e. $\alpha=\alpha'$. Indeed let the equality above holds. Inequality  $\alpha_0<\alpha'_0$ is impossible, since it gives contradiction: $\Hom(\mc B,s_{\alpha_0,\alpha_0})=0$,  $\Hom(\mc B,s_{\alpha_0,\alpha_0})\neq 0$, similarly $\alpha_0>\alpha'_0$, $\alpha_{k+1}>\alpha'_{k+1}$,  $\alpha_{k+1}<\alpha'_{k+1}$ are impossible. Therefore $\alpha_0=\alpha'_0$, $\alpha_{k+1}=\alpha'_{k+1}$. If $\alpha_1<\alpha'_1$, then all the pairs  $ (s_{\alpha_0,\alpha_1},  s_{\alpha'_j+j,\alpha'_{j+1}+j} )  $, $j\geq 1$  are orthogonal, thus we obtain an exceptional collection in $\mc B \cong \mc T_k$ of length $k+1$ in which one element is orthogonal to all the rest, which is known to be impossible. Therefore we see that $\alpha_1=\alpha'_1$, taking the right orthogonal to $s_{\alpha_0,\alpha_1}=s_{\alpha'_0,\alpha'_1}$ in $\mc B$ we see that: \begin{gather}
  \left  \langle  s_{\alpha_1+1,\alpha_2+1},\dots, s_{\alpha_j+j,\alpha_{j+1}+j},\dots, s_{\alpha_{k-1}+k-1,\alpha_{k}+k-1}, s_{\alpha_{k}+k,\alpha_{k+1}+k+1}\right \rangle \nonumber  \\ = \left  \langle  s_{\alpha'_1+1,\alpha'_2+1},\dots, s_{\alpha'_j+j,\alpha'_{j+1}+j},\dots, s_{\alpha'_{k-1}+k-1,\alpha'_{k}+k-1}, s_{\alpha'_{k}+k,\alpha'_{k+1}+k+1}\right \rangle. \nonumber
  \end{gather}
  
  Now by the induction assumption we deduce that $\alpha = \alpha'$, and the lemma is proved.
  \epr
  \begin{coro} \label{identity group case for A_kinA_n+1}
  	For $n\geq1$, $k\geq 1$ the set $C^{\{\rm Id\}}_{\mc T_{k-1}}(\mc T_n)$  is non-empty iff $k\leq n+1$ and then   $\abs{C^{\{\rm Id\}}_{\mc T_{k-1}}(\mc T_n) }=\binom{n+2}{k+1}$. In particular if $k=n+1$ we have $\abs{C^{\{\rm Id\}}_{\mc T_{k-1}}(\mc T_n)}=\abs{C^{\Gamma_n}_{\mc T_{k-1}}(\mc T_n)}=1$. 
  	 
  \end{coro}
  \bpr Follows from the previous lemma and Lemma \ref{ci}.
  \epr
   We  study now  the action of the Serre functor.
  \begin{prop} \label{Serre on exceptional objects}
  	Let $S$ be the Serre functor of $\mc T_n$. Then for \begin{gather} \label{Serre1 for A_n} 0\leq i\leq j < n \quad \Rightarrow \quad  S(s_{i,j})\cong s_{i+1,j+1}[1] \\\label{Serre2 for A_n} 0\leq i\leq n \quad  \Rightarrow \quad   S(s_{i,n}) \cong s_{0,i} \quad \end{gather}
  \end{prop}
  \bpr Using Lemmas  \ref{help lemma 1}, \ref{help lemma 2} we see that the sequences in Figure \ref{exc col in An}  are exceptional collections. 
  \begin{figure} \begin{gather}\nonumber  (s_{0,0},s_{1,1}, \dots, s_{n,n}) \\
  \nonumber  (s_{1,1},s_{2,2} \dots, s_{n,n},s_{0,n}), \\  
  \nonumber  (s_{2,2},s_{3,3} \dots, s_{n,n},s_{0,n}, s_{0,0}) \\ 
  \nonumber \dots  \\ \nonumber  (s_{i+1,i+1},s_{i+2,i+2} \dots, s_{n,n},s_{0,n}, s_{0,0},s_{1,1},s_{2,2},\dots,s_{i-1,i-1}) \\
  \nonumber  (s_{i+2,i+2},s_{i+3,i+3} \dots, s_{n,n},s_{0,n}, s_{0,0},s_{1,1},s_{2,2},\dots,s_{i,i})
  \\ 
  \nonumber \dots \\ \nonumber  ( s_{n,n}, s_{0,n}, s_{0,0},s_{1,1},s_{2,2},\dots,s_{n-2,n-2})
  \\ \nonumber  (s_{0,n}, s_{0,0},s_{1,1},s_{2,2},\dots,s_{n-1,n-1})\end{gather} \caption{Some full exceptional collections in $\mc T_n$} \label{exc col in An}
  \end{figure}
  Since the first exceptional collection is full and all others are obtained from it via mutations, it follows that all of the listed exceptional collections are full (actually in \cite{WCB1} is shown that all exceptional collections of length $n+1$ in $\mc T_n$ are full). Therefore \eqref{Serre} implies that 
  \begin{gather}
  0\leq i < n \Rightarrow  S(s_{i,i}) \sim s_{i+1,i+1} \\  S(s_{0,n}) \sim s_{0,0} \quad  S(s_{n,n}) \sim s_{0,n}
  \end{gather} Let $0\leq i<n$. Using the main property of the Serre functor we see that $1=\hom(s_{i,i}, s_{i,i})=\hom\left (s_{i,i}, S(s_{i,i})\right )$, on the other hand in  Lemma \ref{help lemma 1} (a) we have $\hom^p(s_{i,i}, s_{i+1,i+1})=1$ iff $p=1$, therefore $S(s_{i,i})=s_{i+1,i+1}[1]$. Now it follows \eqref{Serre1 for A_n}, since $S$ is an exact autoequivalence of $\mc T_n$ and since $s_{i,j}$ is built from $s_{i,i}, s_{i+1,i+1},\dots,s_{j,j}$ by a tower of triangles. 
  
  It remains to prove \eqref{Serre2 for A_n}. We have already proved that $ S(s_{0,n}) \sim s_{0,0} \quad  S(s_{n,n}) \sim s_{0,n} $, on the other hand $1=\hom(s_{0,n},s_{0,n})=\hom\left (s_{0,n},S(s_{0,n})\right )$ and $1=\hom(s_{n,n},s_{n,n})=\hom\left (s_{n,n},S(s_{n,n})\right )$, and the vanishings $ \hom^i\left (s_{0,n},s_{0,0}\right )= \hom^i\left (s_{n,n},s_{0,n}\right )=0 $ for $i\neq 0$ from  Lemma \ref{help lemma 1} (c), (d)   ensure that $ S(s_{0,n})\cong  s_{0,0} $ and $ S(s_{n,n}) \cong s_{0,n} $. Choose now $0 <i < n$. From Lemma \ref{help lemma 2} we know that $(s_{i,n-1},s_{n,n})$ and $(s_{0,i},s_{i+1,n})$ are  exceptional pairs and that $\langle s_{i,n-1},s_{n,n} \rangle = \langle s_{n,n},s_{i,n} \rangle$, $ \langle s_{i+1,n}, s_{0,n} \rangle =  \langle s_{0,i},s_{i+1,n} \rangle = \langle s_{0,n}, s_{0,i} \rangle$. Using the already proved equalities we derive:
  \begin{gather}\nonumber  \langle s_{0,n}, S(s_{i,n}) \rangle = S\left ( \langle s_{n,n}, s_{i,n} \rangle \right ) =S\left ( \langle s_{i,n-1}, s_{n,n} \rangle \right ) =  \langle s_{i+1,n}, s_{0,n} \rangle =  \langle  s_{0,n}, s_{0,i} \rangle  \Rightarrow S(s_{i,n}) \sim  s_{0,i}. \end{gather} 
  On the other hand using \eqref{Euler} on shows that:
  \begin{gather} \nonumber \hom(s_{i,n}, s_{0,i})=1 \qquad \hom^j(s_{i,n}, s_{0,i})=0 , j\neq 0 \end{gather}
  which together with the equalities  $1=\hom(s_{i,n}, s_{i,n})=\hom\left (s_{i,n}, S( s_{i,n})\right )$ coming from the fact that $S$ is a Serre functor  imply  $S(s_{i,n}) \cong  s_{0,i}$, hence the lemma is completely proved. \epr
  
  \begin{coro}
  	\label{lemma orbits for DbPoint}
  	Any orbit in $C_{\mc T_0,\KK}^{\{\rm Id\}}(\mc T_n)$ of  $\langle S \rangle$ has the following form (here $0\leq \beta \leq n$):
  	\begin{gather}
  	\nonumber \langle s_{0,\beta}\rangle  \mapsto \langle s_{1,\beta+1} \rangle
  	\mapsto\cdots \mapsto \langle s_{n-\beta,n} \rangle  \nonumber \\[- 2mm] \label{orbits for Dbpoint}   \\[-2mm] \nonumber
  	\mapsto  \langle s_{0,n-\beta} \rangle  \mapsto \langle s_{1,n-\beta+1} \rangle 
  	\mapsto\dots\mapsto   \langle s_{\beta,n} \rangle \mapsto \ \mbox{return to} \ \langle s_{0,\beta}\rangle.
  	\end{gather}
  	The orbit reduces to the first  row iff $n$ is even and $\beta =\frac{n}{2}$  and then the orbit has $\frac{n}{2}+1$ elements. 
  	
  	Otherwise the elements in \eqref{orbits for Dbpoint}  are pairwise distinct and in this case the orbit has $n+2$ elements. 	
  \end{coro}
  \bpr From Lemma \ref{C_eoA_n+1} we know that any $\mc A \in C_{\mc T_0,\KK}^{\{\rm Id\}}(\mc T_n)$ has the form  $\mc A=\langle s_{\alpha,\beta} \rangle$ for some $0\leq \alpha \leq \beta\leq n$, furthermore the numbers $\alpha, \beta$ are uniquely determined from  $\mc A$.   The action of $S$, described in Proposition \ref{Serre on exceptional objects}, allows us to map any $\langle s_{\alpha',\beta'} \rangle$ to an element of the form $\langle s_{0,\beta} \rangle$, and starting from  $\langle s_{0,\beta} \rangle$ applying repeatedly the Serre functor generates the  sequence in \eqref{orbits for Dbpoint}.
  The elements in the given sequences are of the form $ S^j( \langle s_{0,\beta}\rangle )$,  it remains to show the minimal $j$ s.t. $S^j(\langle s_{0,\beta}\rangle )=\langle s_{0,\beta}\rangle $. 
  The first element, $\langle s_{0,\beta}\rangle $,  in the sequence \eqref{orbits for Dbpoint}  can appear again in this sequence, only if $n-\beta=\beta$ and then the orbit reduces to the first row, and has $1+\beta=1+\frac{n}{2}$ elements, otherwise all the elements in \eqref{orbits for Dbpoint}  are pairwise different and their number is $1+n-\beta +1 +\beta = n+2$.
  \epr

  \begin{coro} \label{C_A1inA_n+1aut}  Let $n\geq 1$ and $\Gamma_n = \rm Aut_{\KK}(\mc T_n)$, then 
  	$
  	\abs{ C^{\Gamma_n}_{\mc T_0,\KK}(\mc T_n)  } =\left \{ \begin{array}{c c}\frac{(n+1)}{2}  & \mbox{if} \ n \ \mbox{is odd}\\ \frac{(n+2)}{2} & \mbox{otherwise} \end{array} \right. .
  	$	
  \end{coro}
  \bpr First we recall \eqref{Serre as the entire for Dynkin} and therefore it is enough to count the orbits of $\langle S \rangle $ on  $ C^{\{\rm Id\}}_{\mc T_0, \KK}(\mc T_n) $.  From Corollary \ref{lemma orbits for DbPoint} we know that if $n$ is odd, then all the orbits 
  contain the same number of elements $n+2$, on the other hand from Lemma \ref{C_eoA_n+1}  we know that $\abs{ C^{\{\rm Id\}}_{\mc T_0,\KK}(\mc T_n)}=\frac{(n+1)(n+2)}{2} $, therefore the number of orbits, which is the same as $\abs{  C^{\Gamma_n}_{\mc T_0, \KK}(\mc T_n)  }$, is $\frac{(n+1)}{2}$. 
  
  If $n$ is even, then from Corollary \ref{lemma orbits for DbPoint}  we know that there is  unique  orbit with  $\frac{n}{2}+1$ elements and the rest  have $n+2$ elements. Therefore $\abs{ C^{\Gamma_n}_{\mc T_0, \KK}(\mc T_n) }-1=\left (\frac{(n+1)(n+2)}{2} - \frac{n+2}{2} \right )/(n+2)$.
  \epr

 We have already clarified the case  $k=1$ (see Corollary \ref{C_A1inA_n+1aut}).  From now on , unless else specified, we identify $C_{\mc T_{k-1}}^{\{\rm Id \}}(\mc T_n)$ with $X_n^k$ (see below) 
 using Lemma \ref{Cid_A_kinA_n+1}, we assume $2\leq k < n+1$, and denote  
 \begin{gather} \label{kn} 2\leq k < n+1, \qquad  \delta_n^k=  n+1-k\\ 
 \label{Xnk} X_n^k=\left \{ \bd \{0,1,\dots,k\} & \rTo^{\alpha} &  \{0,1,\dots,\delta_n^k\} \ed : 0\leq \alpha_0 \leq \alpha_1 \leq \dots \leq \alpha_k \leq \delta_n^k \right \}.\end{gather}

 \begin{coro} \label{Serre on C_0A_kinA_n}
 	The action (explained in Corollary \ref{group action}) of $ S $ on the elements in   $C_{\mc T_{k-1}}^{\{\rm Id \}}(\mc T_n)$ is:
 	\begin{gather}
 	\alpha_k<\delta_n^k \ \ \Rightarrow \ \ 	S \left((\alpha_0,\dots,\alpha_{k-1},\alpha_k) \right )=(\alpha_0+1,\dots,\alpha_{k-1}+1,\alpha_k+1) \label{Serre 22} \\
 	S\left ((\alpha_0,\dots,\alpha_{k-1},\delta_n^k)\right ) = (0,\alpha_0,\dots,\alpha_{k-1}) \label{Serre 23}.
 	\end{gather} \end{coro}
 \bpr Iff $\alpha_k<\delta_n^k$ we use \eqref{Serre1 for A_n}. If $\alpha_k=\delta_n^k$ we use  \eqref{Serre1 for A_n}, \eqref{Serre2 for A_n}, and  Lemma \ref{help lemma 2} to  compute \begin{gather}\nonumber 
 S\left (\langle  s_{\alpha_0,\alpha_1}, s_{\alpha_1+1,\alpha_2+1},\dots, s_{\alpha_{k-2}+k-2,\alpha_{k-1}+k-2}, s_{\alpha_{k-1}+k-1,n}\rangle\right )=
 \left \langle  S( s_{\alpha_0,\alpha_1}),\dots, S(s_{\alpha_{k-1}+k-1,n}) \right \rangle = \nonumber \\ 
 \left \langle  s_{\alpha_0+1,\alpha_1+1}, s_{\alpha_1+2,\alpha_2+2},\dots, s_{\alpha_{k-2}+k-1,\alpha_{k-1}+k-1}, s_{0,\alpha_{k-1}+k-1} \right \rangle \nonumber \\ =   \left \langle s_{0,\alpha_0},  s_{\alpha_0+1,\alpha_1+1}, s_{\alpha_1+2,\alpha_2+2},\dots, s_{\alpha_{k-2}+k-1,\alpha_{k-1}+k-1} \right \rangle \nonumber.
 \end{gather} \epr
 \begin{coro} \label{orbits as elements in our set} 
 	The set $C^{\Gamma_n}_{\mc T_{k-1}}(\mc T_n)$ can be identified with the set of orbits of the action of $\langle S \rangle \subset Perm(X_n^k)$ on $X_n^k$ described in Corollary \ref{Serre on C_0A_kinA_n}.
 \end{coro} \bpr This follows from the results in this subsection and  \eqref{Serre as the entire for Dynkin}. \epr
 \subsection{Introducing  $d$-additive sequences}  \label{introducing d-additive}
 We give a name \textit{$d$-additive  sequences} of some elements in $X_n^k$ (see  \eqref{Xnk} for the notation $X_n^k$), and show that we need them to study  the orbits of the action of $\langle S \rangle $ on $X_n^k$ described in Corollary \ref{Serre on C_0A_kinA_n}. The goal is proof of Proposition \ref{main prop for d-additive} and its corollaries.

 Let $\mc O$ be any orbit of $\langle S \rangle$ in $X_n^k$.

 Using Corollary \ref{Serre on C_0A_kinA_n} one shows that there is  an element in $\mc O$ with vanishing first coordinate. Take any such element $\alpha=(0,\alpha_1,\alpha_2,\dots, \alpha_k) \in \mc O$.  Applying $S$ recursively to $\alpha$    produces the  cycle shown in Figure \ref{orbit general}
 \begin{figure}
 \begin{gather}
 \nonumber \alpha=(0,\alpha_1,\alpha_2,\dots, \alpha_k)\mapsto\cdots \mapsto \left (\delta_n^k-\alpha_k,\delta_n^k+\alpha_1-\alpha_k,\dots, \delta_n^k+\alpha_{k-1}-\alpha_k, \delta_n^k \right ) \\ 
 \mapsto \left (0,\delta_n^k-\alpha_k,\delta_n^k+\alpha_1-\alpha_k,\dots, \delta_n^k+\alpha_{k-1}-\alpha_k\right )\mapsto \quad \qquad \qquad \qquad \cdots  \qquad \qquad \qquad \quad  \nonumber  \\ 
 \quad \qquad \qquad  \cdots  \qquad \qquad \qquad \quad  \mapsto \left (\alpha_{k}-\alpha_{k-1},\delta_n^k-\alpha_{k-1},\delta_n^k+\alpha_1-\alpha_{k-1},\dots, \delta_n^k+\alpha_{k-2}-\alpha_{k-1},  \delta_n^k \right ) \nonumber \\ \mapsto \left (0,\alpha_{k}-\alpha_{k-1},\delta_n^k-\alpha_{k-1},\delta_n^k+\alpha_1-\alpha_{k-1},\dots, \delta_n^k+\alpha_{k-2}-\alpha_{k-1} \right )\mapsto \quad \qquad \qquad  \cdots   \qquad \qquad \quad  \nonumber  \\ 
 \quad \qquad   \cdots  \qquad \qquad    \mapsto \left (\alpha_{k-1}-\alpha_{k-2},\alpha_{k}-\alpha_{k-2},\delta_n^k-\alpha_{k-2},\delta_n^k+\alpha_1-\alpha_{k-2},\dots, \delta_n^k +\alpha_{k-3}-\alpha_{k-2} , \delta_n^k \right ) \nonumber \\
  \vdots  \nonumber \\ 
 \mapsto \left (0,\alpha_{i}-\alpha_{i-1},\dots,\alpha_{k}-\alpha_{i-1},\delta_n^k-\alpha_{i-1}, \delta_n^k+\alpha_1-\alpha_{i-1}, \dots, \delta_n^k+\alpha_{i-2}-\alpha_{i-1} \right )\mapsto  \qquad \qquad \cdots \qquad  \qquad  \nonumber  \\ 
 \quad \cdots \quad  \mapsto \left (\alpha_{i-1}-\alpha_{i-2},\alpha_{i}-\alpha_{i-2},\dots,\alpha_{k}-\alpha_{i-2},\delta_n^k-\alpha_{i-2}, \delta_n^k+\alpha_1-\alpha_{i-2}, \dots,\delta_n^k+\alpha_{i-3}-\alpha_{i-2}, \delta_n^k \right ) \nonumber  \\
 \vdots  \nonumber \\ 
 (0,\alpha_2-\alpha_1,\alpha_3-\alpha_1,\dots, \alpha_k-\alpha_1,\delta_n^k-\alpha_1)\mapsto \qquad \qquad \cdots  \qquad \qquad  \mapsto \left (\alpha_1,\alpha_2,\dots, \alpha_k, \delta_n^k \right ) \nonumber \\ 
 \mapsto \ \mbox{return to} \ \alpha=(0,\alpha_1,\alpha_2,\dots, \alpha_k) \nonumber
 \end{gather} \caption{Acting with $S$ on an element in $X_n^k$}  \label{orbit general}
 \end{figure}
 where  in the part  indexed by $i$, the index $i$ varies decreasingly  from $k-1$ to $3$, in particular if $k=3$ or $k=2$ this part does not exists, furthermore, if $k=2$ the cycle looks like this (here we put $\beta=\alpha_1$, $\gamma=\alpha_2$):
 	\begin{gather}
 	\nonumber (0,\beta,\gamma)\mapsto (1,\beta+1,\gamma+1)\mapsto\cdots \mapsto (n-1-\gamma,\beta+n-1-\gamma,n-1) \\ \label{orbit} 
 	\mapsto (0,n-1-\gamma,\beta+n-1-\gamma)\mapsto(1,n-\gamma,\beta+n-\gamma)\mapsto\dots\mapsto(\gamma-\beta,n-1-\beta,n-1)\\ \nonumber 
 	\mapsto (0,\gamma-\beta,n-1-\beta)\mapsto\dots\mapsto(\beta,\gamma,n-1) \mapsto \ \mbox{return to} \ (0,\beta,\gamma).
 	\end{gather}
    Each element in Figure \ref{orbit general} is obtained from the previous via the formulas in Corollary \ref{Serre on C_0A_kinA_n}.

 Note that in the cycle from Figure \ref{orbit general} the elements with a zero leading coordinate are in the following sequence, and they appear in this order: 
 \begin{gather}\label{elements with zero first coordinate}  \alpha,S^{\delta_n^k+1-\alpha_k }\left (\alpha \right ),
 \dots, S^{\delta_n^k+1+u-\alpha_{k-u}}\left (\alpha \right ), \dots, S^{\delta_n^k+k-\alpha_{1}}\left (\alpha \right ), S^{\delta_n^k+k+1}\left (\alpha \right )=\alpha.
 \end{gather}
 where $\delta_n^k+1-\alpha_k< \delta_n^k+2-\alpha_{k-1} <  \dots < \delta_n^k+k-\alpha_{1} < \delta_n^k+k+1=n+2 $.
 Therefore: 
 \begin{lemma} \label{help lemma for general orbits 1}   We have  $S^{n+2}\left (\alpha \right )=\alpha$  (and therefore $ \frac{n+2}{ \abs{\mc O}}\in \NN $).  Furthermore  (recall that $\alpha_0=0$):  
 	\begin{gather} \label{rough formula for number of el in orbit} \abs{\mc O}=\min\left \{\delta_n^k+1+u-\alpha_{k-u}: 0\leq u \leq k  \ \mbox{and}  \  S^{\delta_n^k+1+u-\alpha_{k-u}}\left (\alpha \right )=\alpha \right \}. \end{gather}
 	In particular, $\abs{\mc O}<n+2$ iff $ S^{\delta_n^k+1+u-\alpha_{k-u}}\left (\alpha \right )=\alpha$ for some $0\leq u \leq k-1$.
 \end{lemma}
 Next step is to prove:
 \begin{lemma} \label{lemma equivalences}  Here $v$ denotes an integer.  For any $1\leq v \leq  k$ we have an equivalence: 
 	\begin{gather}  \label{equivalences}  		
 	S^{\delta_n^k+v-\alpha_{k+1-v}}\left (\alpha \right )= \alpha \qquad \iff \qquad \begin{array}{c}  \{\alpha_{j+(k+1-v)}= \alpha_j+\alpha_{k+1-v}\}_{j=1}^{v-1}\\  \mbox{and} \ \ \delta_{n}^k=\alpha_{v}+\alpha_{k+1-v} \\ \mbox{and} \ \ \{\alpha_{j+v} = \alpha_{j} + \alpha_{v} \}_{j=1}^{k-v}.  \end{array}    \end{gather}
 	Furthermore:
 	\begin{gather} \label{better formula for number of el in orbit} \abs{\mc O}=\min \left ( \{n+2\} \cup \left \{\delta_n^k+v-\alpha_{k+1-v}: 1\leq v \leq \frac{k+1}{2}  \ \mbox{and}  \  S^{\delta_n^k+v-\alpha_{k+1-v}}\left (\alpha \right )=\alpha \right \} \right). \end{gather}
 \end{lemma}
 \bpr   
 Now we take the elements with a zero leading term in Figure \ref{orbit general}, remove the first  coordinate (which equals zero), and order the rest coordinates  in the table below according to the order in which the corresponding element appear in the cycle from Figure  \ref{orbit general}:
 { \small	\begin{gather} \label{table with zero term elements} \nonumber
 	\begin{array}{c c c c c c c c}
 	\alpha_1 &\alpha_2 & \dots & \alpha_{k+1-i} & \alpha_{k+2-i} &  \alpha_{k+3-i} &\dots & \alpha_k  \\
 	\delta_n^k-\alpha_k & \delta_n^k+\alpha_1-\alpha_k &\dots &  \delta_n^k+\alpha_{k-i}-\alpha_k & & & & \delta_n^k+\alpha_{k-1}-\alpha_k  \\
 	\{ \ \alpha_{i}-\alpha_{i-1} & \alpha_{i+1}-\alpha_{i-1}&\dots &\alpha_{k}-\alpha_{i-1} &\delta_n^k-\alpha_{i-1} & \delta_n^k+\alpha_1-\alpha_{i-1} & \dots & \delta_n^k+\alpha_{i-2}-\alpha_{i-1}  \ \}_{i=k}^2.
 	\end{array}
 	\end{gather}
 }
 On the one hand,  using the locations of these elements  in Figure \ref{orbit general} we deduce that  the first row equals to the  second, or to the third, ..., or to the $k+1$-th row  if and  only if the following (see also \eqref{elements with zero first coordinate}) hold, respectively :
 \begin{gather}\label{help for orbits} S^{\delta_n^k-\alpha_k + 1  }\left (\alpha \right )=\alpha, S^{\delta_n^k+2-\alpha_{k-1}}\left (\alpha \right )=\alpha, \dots, S^{\delta_n^k+1+u-\alpha_{k-u}}\left (\alpha \right )=\alpha, \dots, S^{\delta_n^k+k-\alpha_{1}}\left (\alpha \right )=\alpha\end{gather}
 on the other hand, equalizing the coordinates in the table above   shows that  the $k$ equalities in \eqref{help for orbits} are equivalent to the following $k$ systems, respectively: 
 \begin{gather}
 \alpha_1= \delta_n^k-\alpha_k,  \quad \{\alpha_j= \delta_n^k+\alpha_{j-1}- \alpha_k\}_{j=2}^k\nonumber \\
 \{\alpha_j=\alpha_{j+(i-1)}-\alpha_{i-1}\}_{j=1}^{k+1-i}, \ \ \delta_{n}^k=\alpha_{k+2-i}+\alpha_{i-1}, \ \{\delta_n^k+\alpha_{j}- \alpha_{i-1}= \alpha_{j+k+2-i}\}_{j=1}^{i-2}  \ \ k\geq i \geq 2 \nonumber 
 \end{gather}
 these systems of equations in turn are equivalent to the following, respectively: 
 
 \begin{gather}
 \alpha_1+\alpha_k= \delta_n^k,  \quad \{\alpha_j= \alpha_1+\alpha_{j-1}\}_{j=2}^k\nonumber \\
 \{\alpha_{j+(i-1)}= \alpha_j+\alpha_{i-1}\}_{j=1}^{k+1-i}, \ \ \delta_{n}^k=\alpha_{k+2-i}+\alpha_{i-1}, \ \{\alpha_{j+k+2-i} = \alpha_{j} + \alpha_{k+2-i} \}_{j=1}^{i-2}  \quad  k\geq i \geq 2. \nonumber 
 \end{gather} The latter sequence of $k$ systems of equations  is the same as a sequence indexed  by $i$, where $i$ starts with $k+1$ and decreases to $2$ and for any such $ k+1\geq i \geq 2$ the system is:
 \begin{gather}\label{sequaence of systems of equations} \{\alpha_{j+(i-1)}= \alpha_j+\alpha_{i-1}\}_{j=1}^{k+1-i} \ \mbox{and} \  \delta_{n}^k=\alpha_{k+2-i}+\alpha_{i-1} \ \mbox{and} \  \{\alpha_{j+k+2-i} = \alpha_{j} + \alpha_{k+2-i} \}_{j=1}^{i-2}.   \end{gather}
 Next we note that 
 for any $0\leq u\leq k-1$ the system for $i=k+1-u$ in \eqref{sequaence of systems of equations} coincides with the system for $i=2+u$, and is the following:
 \begin{gather}  \label{help for orbits 3} \{\alpha_{j+(k-u)}= \alpha_j+\alpha_{k-u}\}_{j=1}^{u} \ \  \mbox{and} \ \ \delta_{n}^k=\alpha_{1+u}+\alpha_{k-u} \ \ \mbox{and} \ \ \{\alpha_{j+1+u} = \alpha_{j} + \alpha_{1+u} \}_{j=1}^{k-u-1}.  \end{gather}
 So far ( see  also \eqref{help for orbits}), we proved that 
 \begin{gather} 		
 S^{\delta_n^k+1+u-\alpha_{k-u}}\left (\alpha \right )= \alpha \qquad \iff \qquad \begin{array}{c}  \{\alpha_{j+(k-u)}= \alpha_j+\alpha_{k-u}\}_{j=1}^{u}\\  \mbox{and} \ \ \delta_{n}^k=\alpha_{1+u}+\alpha_{k-u} \\ \mbox{and} \ \ \{\alpha_{j+1+u} = \alpha_{j} + \alpha_{1+u} \}_{j=1}^{k-u-1}.  \end{array}    \end{gather}
 By replacing $u$ with $v-1$ we obtain \eqref{equivalences}   for $ 1\leq v \leq k$ .
 
 Finally, note that  for any  $ 0\leq u \leq k-1 $ by replacing $u$ with $k-1-u$  the system \eqref{help for orbits 3} remains the same, and that the set of integers in $[0,k-1]$  can be represented as follows  (disjoint unions):
  \begin{gather} \left \{u:0\leq u\leq \frac{k-2}{2}\right \} \cup \left \{k-1-u:0\leq u\leq \frac{k-2}{2}\right \}  \quad \mbox{if} \ k \ \mbox{is even}\\
 \left \{u:0\leq u\leq \frac{k-1}{2}\right \} \cup \left \{k-1-u:0\leq u\leq \frac{k-1}{2}-1\right \}  \quad \mbox{if} \ k \ \mbox{is odd}.\end{gather} 
  From  \eqref{rough formula for number of el in orbit} and \eqref{equivalences}\footnote{taking into account that $\delta_n^k+k+1-\alpha_{0}=n+2$ and that for integers $u, \gamma$ the inequalities $u\leq \gamma-1$ and $u\leq \gamma-\frac{1}{2}$ are equivalent} we see that $\abs{\mc O}$ is  either $n+2$ or the minimal integer among  $  \delta_n^k+1+u-\alpha_{k-u} $ s.t. $ 0\leq u \leq \frac{k-1}{2}$  and  $ S^{\delta_n^k+1+u-\alpha_{k-u}}\left (\alpha \right )=\alpha $, if such $u$ exists, hence we get  \eqref{better formula for number of el in orbit} (we again prefer to use the variable $v=u+1$).  
  \epr
 
 \begin{lemma} \label{lemma for inc periodic} Let $1\leq v \leq k$ be integer, here $k$ can be any positive integer (not restricted as in \eqref{kn}). Let  
 	$\bd  \{ 1,\dots,k \} & \rTo^{\beta} & \ZZ \ed $ be a sequence which satisfies: $\beta(j+v)= \beta(j)+\beta(v)$ for  $1\leq j \leq k-v$. Then there exists unique extension of $\beta$ to a sequence $\bd[1em] \ZZ & \rTo^\beta & \ZZ \ed$ satisfying 
 	\begin{gather} \label{inc periodic}  \beta(j+i v)=\beta(j)+i \beta(v)=\beta(j)+ \beta(i v) \qquad i,j\in \ZZ. \end{gather}
 	
 \end{lemma}
 \bpr  If such extension exists, then $\beta(v)=\beta(0+v)=\beta(0)+\beta(v)$, hence  it must satisfy $\beta(0)=0$.   Furthermore, the values of such extension are uniquely determined by its values in $[1,v]$, hence if the extension  exists, it is unique.   
 
 So let us define $\beta(0)=0$. Then from the assumption for $\beta$ we have $\beta(j+v)=\beta(j)+\beta(v)$ for any $j\geq 0$, s.t. $j+v\leq k$. Therefore $\beta(2 v)= 2 \beta(v)$ if $2 v \leq k$, and by induction one proves, that for any $0\leq j$, s.t. $j v \leq k$ we have $\beta(jv) = j \beta(v)$. 
 Next, let $j \geq 0$ be such that $j+ 2 v \leq k$, then $\beta(j+2v)=\beta(j+v) + \beta(v)=\beta(j) +\beta(v)+ \beta(v)=\beta(j) + 2 \beta(v)$, and now by induction  on proves
 \begin{gather} \label{property for the proof of orbits} i,j\geq 0, \ j+ i v \leq k \ \ \Rightarrow \beta(j+iv)= \beta(j) +\beta(i v) =  \beta(j) +i \beta( v).  
 \end{gather}
 We can define a sequence $\beta':\ZZ \rightarrow \ZZ$  by the rule  $\beta' (j + iv)=i\beta (v) +\beta(j)$ for $i\in \ZZ$, $0\leq j <v$. Then one easily shows that this $\beta'$ satisfies $\beta' (j + iv)=\beta' (j) +\beta'(iv)=\beta' (j ) +i \beta'(v)$ for any $j\in \ZZ$, $ i \in \ZZ$. Since for $x\in [1,k]$ we have  $i,j\geq 0$, $j<v$, s. t. $x=i v + j$, therefore by the property \eqref{property for the proof of orbits} of $\beta$ and of the definition of $\beta'$ one computes $\beta'(x)=\beta(x)$, hence $\beta'$ is an extension of $\beta$.
 \epr
 \begin{lemma} \label{equal extensions} Let $1\leq v \leq k$, $1\leq u \leq k$ be integers, such that $u+v \leq k+1$, here $k$ can be any positive integer. Let  $\beta:[1,\dots,k]\rightarrow \ZZ$ be a sequence which satisfies:
 	\begin{gather}   \begin{array}{c}  \{\beta(j+u)= \beta(j)+\beta(u)\}_{j=1}^{k-u}   \mbox{and} \ \ \{\beta(j+v) = \beta(j) + \beta(v) \}_{j=1}^{k-v}.  \end{array}    \end{gather}
 	Then there exists unique extension of $\beta$ to a sequence $\bd[1em] \ZZ & \rTo^\beta & \ZZ \ed$ satisfying 	for any $i,j\in \ZZ$ both
 	\begin{gather} \label{inc periodic 3}  \beta(j+i v)=\beta(j)+i \beta(v)=\beta(j)+ \beta(i v) \qquad \beta(j+i u)=\beta(j)+i \beta(u)=\beta(j)+ \beta(i u).   \end{gather}
 \end{lemma}
 \bpr  Let us extend $\beta$ using Lemma \ref{lemma for inc periodic} and the part $\{\beta(j+v) = \beta(j) + \beta(v) \}_{j=1}^{k-v}$ of the given system. In particular we have a sequence $\bd[1en] \ZZ & \rTo^{\beta} & \ZZ\ed$ satisfying \eqref{inc periodic}. We have also an extension due to the other part of the given system $\{\beta(j+u)= \beta(j)+\beta(u)\}_{j=1}^{k-u}$, however we will show that this extension coincides with the already constructed $\beta$. 
 Indeed from  $\{\beta(j+u)= \beta(j)+\beta(u)\}_{j=1}^{k-u}$ and $\beta(0)=0$ (see also the proof of  Lemma \ref{lemma for inc periodic} ) it follows that 
 \begin{gather} \label{inc periodic 2}  \beta(j+i u)=\beta(j)+i \beta(u)=\beta(j)+ \beta(i u) \qquad i\geq 0, j\geq 0, \quad j+i u\leq k, \end{gather}
 using this we will show that $\beta(j+u)=\beta(j)+\beta(u)$ for any  $j\in \ZZ$ and then  the lemma  follows by induction.   Indeed, for any  $j\in \ZZ$ we have   $j=j' v + i'$ for unique $j'\in \ZZ$ and $0\leq i' <v$ and using  \eqref{inc periodic} we get $\beta(j+u)=\beta(j' v + i'+u) =\beta(j' v) + \beta( i'+u)  $, on the other hand $i'+u < v+u\leq k+1$, i.e  $i'+u\leq k$, and by \eqref{inc periodic 2} $\beta( i'+u)=\beta(i')+\beta(u)$, therefore we can combine and continue, utilizing once again  \eqref{inc periodic}, the equalities as follows $\beta(j+u)=\beta(j' v) + \beta( i')+\beta(u) = \beta(j' v+i')+\beta(u)=\beta(j)+\beta(u)$.
 \epr 
 \begin{lemma} \label{solution}  Let $1\leq v \leq \frac{k+1}{2}$, $u = k+1-v$ be integers, here $k$ can be any positive integer.  Let us  denote $d=g.c.d.(u,v)$, in particular $1\leq d  \frac{k+1}{2}$.
 	A sequence   $\beta:[0,1,\dots,k]\rightarrow \ZZ$  satisfies:
 	\begin{gather}  \label{system in lemma}   \begin{array}{c} \left  \{\beta(j+u)= \beta(j)+\beta(u) \right \}_{j=1}^{k-u}\\ \mbox{and} \ \ \beta(0)=0 \ \  \mbox{and} \ \ \delta_{n}^k=\beta(v)+\beta(u) \\ \mbox{and} \ \ \left  \{\beta(j+v) = \beta(j) + \beta(v) \right \}_{j=1}^{k-v}  \end{array}    \end{gather}
 	if and only if the following hold:  
 	\begin{gather}\label{first row for good sequences}
 	\frac{k+1}{d} \in \ZZ \qquad \beta(d) = d \frac{\delta_n^k}{k+1} \in \ZZ \\ \label{second row for good sequences}
 	i\geq 0, j \geq 0, j+i d \leq k \Rightarrow  \beta(j + i d) = \beta(j) + \beta(i d)=\beta(j) + i \beta(d),
 	\end{gather}
 	which, in turn, is equivalent to:
 	\begin{gather}  \label{system 2 in lemma} \frac{k+1}{d} \in \ZZ, \ \    \begin{array}{c} \left  \{\beta(j+(k+1-d))= \beta(j)+\beta(k+1-d) \right \}_{j=1}^{d-1}\\ \mbox{and} \ \ \beta(0)=0 \ \  \mbox{and} \ \ \delta_{n}^k=\beta(d)+\beta(k+1-d) \\ \mbox{and} \ \ \left  \{\beta(j+d) = \beta(j) + \beta(d) \right \}_{j=1}^{k-d}.  \end{array}    \end{gather}
 \end{lemma}
 \bpr We show first that the system \eqref{system in lemma} implies \eqref{first row for good sequences} and \eqref{second row for good sequences}.  Using Lemma \ref{equal extensions} we see that we have unique  extension $\bd[1em] \ZZ & \rTo^\beta & \ZZ \ed$ of the given sequence such that \eqref{inc periodic 3} holds. There are integers $t_u, t_v \in \ZZ$, $\gamma_u ,\gamma_v \in \NN$ such that: \ \ $
 t_u u + t_v v = d$, $ u = \gamma_u d $, $ v = \gamma_v v.
 $
 We are given also $u+v=k+1$, therefore  \begin{gather} \label{help for orbits 5} \frac{k+1}{d} = \gamma_u + \gamma_v \in \ZZ.  \end{gather}
 Using \eqref{inc periodic 3} we deduce that the extended function $\bd[1em] \ZZ & \rTo^\beta & \ZZ \ed$  satisfies for any $i,j \in \ZZ$:  $$\beta(j d + i)=\beta\left (j (t_u u + t_v v) + i \right )= \beta \left (j (t_u u + t_v v) \right ) + \beta(i)= j \beta(t_u u + t_v v) + \beta(i) = j \beta(d) + \beta(i).$$
 Thus we obtain \eqref{second row for good sequences}. We are given also  $\beta(u)+\beta(v) = \delta_n^k$, hence using the already proved \eqref{second row for good sequences} we obtain $(\gamma_u + \gamma_v)\beta(d)=\delta_n^k$ and now combining with \eqref{help for orbits 5} we derive \eqref{first row for good sequences}.
 
 The converse (that \eqref{system in lemma} follows from \eqref{first row for good sequences}, \eqref{first row for good sequences}) is even easier to prove.  
 
 Now we apply the already proved equivalence to the case $v = d$, $u=k+1-d$, and provided that $\frac{k+1}{d} \in \ZZ$ it follows that $g.c.d.(u,v)=d$, follows  that \eqref{system 2 in lemma} is equivalent to \eqref{first row for good sequences}, \eqref{second row for good sequences}.
 \epr
 Lemmas \ref{help lemma for general orbits 1}, \ref{lemma equivalences}  and \ref{solution} show that the special kind of sequences $\alpha \in X_n^k$ defined now will play an important role in our discussion.  
 \begin{remark} \label{note for integers}  Note in advance that for integers $k$, $d$  such that $n+1>k\geq 1$ the condition  $\frac{k+1}{d} \in \ZZ$, $1\leq d \leq \frac{k+1}{2}$ is equivalent to:  $\frac{k+1}{d} \in \ZZ$, $1\leq d < k+1$. Also   $d \frac{\delta_n^k}{k+1} \in \ZZ$ is equivalent to $d \frac{n+2}{k+1} \in \ZZ$.
 \end{remark}
 \begin{df}\label{def of d-additive}  Let $1\leq k < n+1$ be integers.  We call an integer $d$ a  \underline{divisor of $(k,n)$}, if   $1\leq d \leq (k+1)$, $(k+1)/d \in \ZZ$, $d (n+2)/(k+1) \in \ZZ$. In particular $k+1$ is always divisor of $(k,n)$.  We say that  $d$ is  \underline{a  proper  divisor of $(k,n)$}, if it is divisor and $d<(k+1)$. 
 	
 	For a proper divisor $d$ of $(k,n)$ we define  a  \underline{$d$-additive sequence}  as  an element  $\beta \in X_n^k$ satisfying \begin{gather} \label{d-additive}
 	\beta(d) = d \frac{\delta_n^k}{k+1}  \qquad   i\geq 0, j \geq 0, j+i d \leq k \Rightarrow  \beta(j + i d) = \beta(j) + \beta(i d)=\beta(j) + i \beta(d).
 	\end{gather}
 	If $d=k+1$, then we say that $\beta \in X_n^k$ is a $d$-additive sequence iff $\beta(0)=0$.
 \end{df} 
 \begin{remark}  \label{from d' to d}
 	Note that any $d$-additive sequence $\beta$ satisfies $\beta(0)=0$. If $d$ and $d'$ are two divisors of $(k,n)$ and $\frac{d}{d'} \in \ZZ$, then any $d'$-additive sequence is also $d$-additive.    	  	
 \end{remark}

Using Lemmas \ref{help lemma for general orbits 1}, \ref{lemma equivalences}  and \ref{solution} we prove the main result of this subsection:
 
 \begin{prop} \label{main prop for d-additive} Let $\mc O$ be any orbit of $\langle S \rangle $ in $X_n^k$. Then there exists $\alpha \in \mc O$ with  $\alpha_0 = 0$. Take any $\alpha \in \mc O$ such that $\alpha_0 = 0$, 
 	then:   	  	$ \abs{\mc O}=\min  \left \{ d \frac{n+2}{k+1}: d   \ \mbox{ is a divisor of $(k,n)$ and $\alpha$ is $d$-additive} \right \}.$
 \end{prop}
 \bpr Lemmas \ref{solution}, \ref{lemma equivalences} (see also Remark \ref{note for integers}) imply that 
 \begin{gather} \nonumber  \abs{\mc O}=\min \left ( \{n+2\} \cup \left \{\delta_n^k+d-\alpha_{k+1-d}:  d   \ \mbox{ is a proper divisor of $(k,n)$ and $\alpha$ is $d$-additive} \right \} \right). \end{gather}
 however, for a $d$-additive sequence $\alpha$ and a proper divisor $d$ of $(k,n)$
 we have 
 \begin{gather}
 \delta_n^k+d-\alpha(k+1-d)= \delta_n^k+d-\alpha\left (\left (\frac{k+1}{d}-1 \right )d \right)=\delta_n^k+d-\left (\frac{k+1}{d}-1 \right )\alpha\left (d \right) \\
 =\delta_n^k+d-\left (\frac{k+1}{d}-1 \right )d \frac{\delta_n^k}{k+1} = d+ d \frac{\delta_n^k}{k+1}=d+ d \frac{n-k+1}{k+1} = d  \frac{n+2}{k+1}.
 \end{gather}  
 Recalling that $k+1$ is a diviosr of $(k,n)$ and that each $\beta \in X_n^k$ with $\beta(0)=0$ is $(k+1)$-additive we deduce the equality.
 \epr
 
 In view of Proposition \ref{main prop for d-additive} we give another definition:
 \begin{df} \label{def of period}
 	Let $\alpha \in X_n^k$ satisfy $\alpha(0)=0$. We call the following  number   \underline{ period of $\alpha$}: $$P(\alpha) = \min  \left \{ d : d   \ \mbox{ is a divisor of $(k,n)$ and $\alpha$ is $d$-additive} \right \}.$$  
 \end{df}
 We remind that we view the  elements of $C^{\rm Aut_{\KK}(\mc T_n)}_{\mc T_{k-1}}(\mc T_n)$ as orbits of $\langle S \rangle$ in $X_n^k$  (Corollary \ref{orbits as elements in our set}  )
 
 \begin{coro} \label{from number to period} Let $\mc O \in C^{\Gamma_n}_{\mc T_{k-1}}(\mc T_n)$, then the following three are equivalent: 	
 	{\rm (a)}  $\abs{\mc O} = d \frac{n+2}{k+1}$; \\ 
 	{\rm (b)}   each  $\alpha \in \mc O$ with $\alpha(0)=0$ has period $d$; \  	
 	{\rm (c)}   some  $\alpha \in \mc O$ with $\alpha(0)=0$ has  period $d$.   
 	
 	Furthermore, if  $d$ is the divisor of $(k,n)$ such that  $\abs{\mc O} = d \frac{n+2}{k+1}$, then
 	$ \abs{\{\alpha \in \mc O : \alpha(0)=0\}}=d. $
 \end{coro}
 \bpr The first part follows straightforward  from Proposition \ref{main prop for d-additive}. 
 To prove the last sentence   we take any $\alpha \in \mc O$ with $\alpha(0)=0$  (such $\alpha$ exists) and by the first part this $\alpha$ is $d$-additive and not $d'$-additive for $d'<d$, hence by Lemmas  \ref{solution}  and \ref{lemma equivalences}
 it follows that in the sequence \eqref{elements with zero first coordinate} containing all the elements with a zero leading coordinate in $\mc O$ we have $S^{\delta_n^k +v-\alpha_{k+1-v}}(\alpha)\neq \alpha$ for $v<d$ and $S^{\delta_n^k +d-\alpha_{k+1-d}}(\alpha) =  \alpha$ (recall that in this formula $\alpha$ was any element in $\mc O$ with $\alpha(0)=0$). Since all other elements in $\mc O$ have leading coefficient which is non-zero, it follows that $\delta_n^k +d-\alpha_{k+1-d}$ is the minimal integer, such that $S^{\delta_n^k +d-\alpha_{k+1-d}}(\alpha) =  \alpha$, and therefore the following is a complete list of elements with zero leading term in $\mc O$ and they are all different:
 $$\alpha,S^{\delta_n^k+1-\alpha_k }\left (\alpha \right ), S^{\delta_n^k+2-\alpha_{k-1}}\left (\alpha \right ), \dots, S^{\delta_n^k+1+u-\alpha_{k-u}}\left (\alpha \right ), \dots, S^{\delta_n^k+d-1-\alpha_{k+2-d}}\left (\alpha \right ) $$ hence the last statement follows.  
 \epr
 \begin{coro} \label{number of orbits with d times}
 	Let  $d$ be a divisor $(k,n)$, then:
 	$$
 	\abs{ \left  \{\mc O \in C^{\Gamma_n}_{\mc T_{k-1}}(\mc T_n): \abs{\mc O}= d \frac{n+2}{k+1} \right  \} } = \frac{1}{d} \abs{ \left  \{\alpha  \in  X_n^k: \alpha(0)= 0,  P(\alpha) = d \right \}}.
 	$$
 \end{coro}
 \bpr From  Corollary \ref{from number to period}  follows that this function is surjective and each fiber has $d$ elements:
 \begin{gather} \nonumber  \left  \{\alpha  \in  X_n^k: \alpha(0)= 0, \quad P(\alpha) = d \right \} \rightarrow \left \{\mc O \in C^{\Gamma_n}_{\mc T_{k-1}}(\mc T_n): \abs{\mc O}= d \frac{n+2}{k+1} \right  \}  \\ 
 \alpha \mapsto \{S^j(\alpha): j\in \ZZ\}.  \nonumber  \end{gather}
 \epr
 In the next subsection   we  study  properties of $d$-additive sequences in order to  determine $\abs{C^{\Gamma_n}_{\mc T_{k-1}}(\mc T_n) }$. 
   \subsection{Propoerties of $d$-additive sequences.}

  \begin{lemma} \label{number of d-additive sequences}  Let $1\leq k < n+1$ be integers.
  	Let $d$ be  any  divisor of $(k,n)$.	 Then
  	
  	{\rm (a)} The set $\{\beta \in X_n^k : \beta \ \mbox{is $d$-additive}  \}$ has $ \left ( \begin{array}{c}  d \frac{n+2}{k+1}-1  \\ d-1 \end{array} \right)$ elements. 
  	
  	{\rm (b)}  There exists  $\beta \in X_n^k $ with $\beta(0)=0$, such that  $P(\beta)=d$ ($P(\beta)$ is in definition \ref{def of period}). 
  \end{lemma}
  \bpr (a) If $d=1$ is a  divisor of $(k,n)$, then due to  \eqref{d-additive}, a $d$-additive sequence must have the form $\alpha(j)=j\frac{\delta_n^k}{k+1}$, and this sequence is $d$-additive, hence there is only one such sequence. 
  
  If $d>1$ is a divisor of $(k,n)$, then using \eqref{d-additive} one shows that restricting any element in the set  $\{\beta \in X_n^k : \beta \ \mbox{is $d$-additive}  \}$ to $\{1,\dots,d-1\}$ is  bijection from this set to the set   (for the case $d=k+1$ we recall that the $d$-additive elements in $X_n^k$ are just the elements with vanishing first coordinate): 
  \begin{gather} \label{help for orbits 333}\left \{ \bd \{1,\dots,d-1\} & \rTo^{\beta} & \left \{0,1,\dots,d \frac{\delta_n^k}{k+1} \right \} \ed : 0\leq \beta_1 \leq \beta_2 \leq \dots \leq \beta_{d-1} \leq d \frac{\delta_n^k}{k+1} \right \}.\end{gather}Now we use Lemma \ref {ci}.
  
  (b) The element $\beta \in X_n^k$ which is $d$-additive and which corresponds to the trivial sequence in the set  \eqref{help for orbits 333} , i.e. $\beta(0)=\beta(1)=\dots = \beta(d-1)=0$ cannot be $d'$ -additive for $d'<d$ (by definition a $d'$-additive sequence does not vanish at $d'$), therefore this $\beta$ has period $d$.
  \epr
  \begin{coro} \label{CGammak-1n} We have (all summands are non-zero):
  	$$ \abs{C^{\Gamma_n}_{\mc T_{k-1}}(\mc T_n)}= \sum_{\{d\geq 1 : d \ \mbox{is a divisor of} \ (k,n) \}} \frac{1}{d} \abs{\left  \{\alpha  \in  X_n^k: \alpha(0)= 0, \ P(\alpha) = d \right \} }. $$
  	
  \end{coro}	
  \bpr This follows from Corollary \ref{number of orbits with d times} and the following equality which we prove below: 
  \begin{gather} \left \{\abs{\mc O}: \mc O \in C^{\Gamma_n}_{\mc T_{k-1}}(\mc T_n)\right \} =\left  \{ d \frac{n+2}{k+1}: d \ \mbox{ is a divisor of $(k,n)$} \right \}. \end{gather}
  
  The inclusion in one direction follows from Proposition \ref{main prop for d-additive}. The inclusion in the other direction follows from Lemma \ref{number of d-additive sequences} (b) and   Corollary  \ref{from number to period} . 
  \epr
  We will find a formula for $ \abs{ \left  \{\alpha  \in  X_n^k: \alpha(0)= 0, \ P(\alpha) = d \right \} }$. To that end we will use Lemma \ref{number of d-additive sequences} (a) and some additional facts. We prove first: 
  \begin{lemma} \label{important lemma for orbits} Let $1\leq k < n+1$ and   let $\{d_i\}_{i\in I}$ be a finite family of divisors of $(k,n)$ in the sense of Definition \ref{def of d-additive}, let $d=g.c.d. \left \{d_i\right \}_{i\in I}$ be the greatest common divisor, then: 
  	
  	{\rm (a)} $d$ is also a divisor of  $(k,n)$;   	
  	{\rm (b)} If $\alpha \in X_n^k$ is $d_i-$additive for each $i\in I$ then $\alpha$ is $d$-additive ;
  	
  	{\rm (c) } The following holds: 
  	$
  	\bigcap_{i\in I}\left  \{\alpha \in X_n^k: \alpha \ \mbox{is $d_i$-additive}\right  \} =\left  \{ \alpha \in X_n^k: \alpha \ \mbox{is $d$-additive}\right \}.
  	$
  	
  \end{lemma}
  \bpr (a) It is clear that $1\leq d \leq d_i\leq k+1 $. Since we have $\frac{k+1}{d_i}, d_i \frac{n+2}{k+1},  \frac{d_i}{d}  \in \ZZ $ for each $i\in \ZZ$ and  $\sum_{i\in I} t_i d_i = d$ for some $\{t_i\}_{i\in \ZZ} \subset \ZZ$, it follows $\frac{k+1}{d}=\frac{k+1}{d_i} \frac{d_i}{d} \in \ZZ$, $d\frac{n+2}{k+1} = \sum_{i\in I} t_i d_i \frac{n+2}{k+1} \in \ZZ$.  
  
  (b) If $d_i=k+1$ for all $i \in I$ the statement holds. Otherwise, since  $\frac{k+1}{d_i}\in \ZZ$ for all $i\in \ZZ$,  we see that $d=g.c.d.\left \{d_i\right \}_{i\in I, d_i\neq k+1} $ and we reduce to the case $d_i < k+1$ for any $i\in I$, in particular $d_i \leq \frac{k+1}{2}$ for any $i \in I$ (see Remark \ref{note for integers}), therefore 
  \begin{gather} \label{help lemma for orbits 333345} \forall i,j \in I \quad d_i + d_j \leq k+1.  \end{gather}  
  We are given that $ \alpha $ is  $d_i$-additive for each $i\in I$, hence (see Definition \ref{def of d-additive}) on the one hand 
  \begin{gather} \label{help lemma for orbits 3333}
  \forall i \in I \qquad \alpha(d_i)= d_i \frac{\delta_n^k}{k+1}
  \end{gather}
  and on the other hand  Lemma \ref{lemma for inc periodic} ensures that for each $i\in \ZZ$  the sequence $\alpha:\{0,1,\dots,k\}\rightarrow \ZZ$ (it has by definition $\alpha(0)=0$) has an unique extension to $\alpha_i : \ZZ \rightarrow \ZZ$ satisfying $\alpha_i(r + t d _i) = \alpha_i(r) + \alpha_i(t d _i)= \alpha_i(r) + t \alpha_i(d _i) $ for any $r,t \in \ZZ$, however due to \eqref{help lemma for orbits 333345} Lemma \ref{equal extensions} ensures that for any $i,j \in \ZZ$ the extensions $\alpha_i, \alpha_j$ coincide, therefore  there exists unique extension of $\alpha$, which we dentoe by the same letter, as follows:
  \begin{gather} \bd \alpha: \ZZ & \rTo^{\alpha} & \ZZ \ed \quad  \forall r,t \in \ZZ \  \forall i \in \ZZ \qquad \alpha(r + t d _i) = \alpha(r) + \alpha(t d _i)= \alpha(r) + t \alpha(d _i). 
  \end{gather} 
  Using that  $\sum_{i\in I} t_i d_i = d$ for some $\{t_i\}_{i\in \ZZ} \subset \ZZ$ it follows that
  $\alpha(r + t d ) = \alpha(r) + \alpha(t d )= \alpha(r) + t \alpha(d)$ for all $r,t \in \ZZ$, and $\alpha (d)=\sum_{i\in I} t_i \alpha(d_i) = \sum_{i\in I} t_i d_i \frac{\delta_n^k}{k+1}= d \frac{\delta_n^k}{k+1} $, therefore  all the properties of being $d$-additive do hold for $\alpha$.
  
  (c) the inclusion in one direction is the content of (b), the other follows from Remark \ref{from d' to d}. 
  \epr
  It is useful to introduce one more notation: 
  \begin{df} \label{PD_d}
  	For any divisor $d$ of $(k,n)$, as defined in Definition \ref{def of d-additive}, we denote a set: 
  	$$PD_d=\left \{p: p \ \mbox{is prime with $\frac{d}{p}\in \ZZ$ and $\frac{d}{p}$ is divisor of $(k,n)$ }  \right \}.$$
  \end{df}
  
  \begin{coro} \label{the seqs of period d}
  	For any divisor $d$ of $(k,n)$ holds:
  	\begin{gather} \nonumber \left  \{\alpha \in X_n^k: \alpha(0)= 0, \ P(\alpha) = d \right  \} = \left \{ \alpha \in X_n^k : \alpha \ \mbox{is $d$-additive} \right \} \setminus \bigcup_{  p \in PD_d   } \left \{ \alpha \in X_n^k : \alpha \ \mbox{is $\frac{d}{p}$-additive} \right \}
  	\end{gather}
  \end{coro}
  \bpr  The inclusion $\subset$ follows from the definition \ref{def of period}. So let us prove the converse $\supset$. Let $\alpha$ be $d$-additive  which is not $\frac{d}{p}$-additive for any prime divisor of $d$, s.t. $\frac{d}{p}$ is a divisor of $(k,n)$. We will show that $\alpha$ is not $d'$-additive for any divisor $d'$ of $ (k,n)$ with $d'<d$ and this shows that $P(\alpha)=d$. Indeed, if $\alpha$ is $d'$-additive for some divisor $d'$ of $(k,n)$ with $d'<d$, then by  Lemma   \ref{important lemma for orbits} (b) we have that $\alpha $ is $d''$-additive, where $d''=g.c.d(d',d)<d$ is a divisor of $(k,n)$, in particular there exists prime $p$ such that $d/p \in \ZZ$ and   $\frac{d/p}{d''} \in \ZZ$.  From this data we deduce: $\frac{k+1}{d/p}=p\frac{k+1}{d} \in \ZZ$ and $\frac{d}{p}\frac{n+2}{k+1}=\frac{d/p}{d''} d''\frac{n+2}{k+1} \in \ZZ$, therefore $\frac{d}{p}$ is a divisor of $(k,n)$, and on the other hand by Remark \ref{from d' to d} and the fact that $\alpha$ is $d''$-additive we deduce that that $\alpha$ is $\frac{d}{p}$-additive, whereas we have chosen $\alpha$ so that it cannot be  $\frac{d}{p}$-additive. The Corollary is proved.
  \epr
  \begin{remark} \label{remark for gcd}
  	Recall that for any  integer $d \in \ZZ$  and any  subset $A\subset \{p:p \ \mbox{is prime and} \ \frac{d}{p} \in \ZZ\}$ the greatest common divisor of $\left \{\frac{d}{p}\right \}_{p\in A}$ is $\frac{d}{\prod_{p\in A} p}$. 
  \end{remark}
  \begin{lemma} \label{lemma for number of d/p} Let  $d$ be a divisor of $(k,n)$.  Recall the definition of $PD_d$ in Definition \ref{PD_d}, then \\ $ \bigcup_{ p \in PD_d   } \left \{ \alpha \in X_n^k : \alpha \ \mbox{is $\frac{d}{p}$-additive} \right \}$ is a subset of $\left \{ \alpha \in X_n^k : \alpha \ \mbox{is $d$-additive} \right \}$ and 
  	\begin{gather}
  	\abs{ \bigcup_{ p \in PD_d   } \left \{ \alpha \in X_n^k : \alpha \ \mbox{is $\frac{d}{p}$-additive} \right \}  }= \sum_{\{B\subset PD_d: B \neq \emptyset \}} (-1)^{\abs{B}+1}\left ( \begin{array}{c} \frac{d}{\prod_{p\in B}p}\frac{n+2}{k+1}-1 \\ \frac{d}{\prod_{p\in B}p}-1 \end{array} \right). 
  	\end{gather} 
  \end{lemma}
  \bpr The inclusion of sets follows from Remark \ref{from d' to d}. Now we apply the formula for number of elements of union of subsets, Lemma \ref{important lemma for orbits}, Remark \ref{remark for gcd}:
  \begin{gather} \nonumber
  \abs{ \bigcup_{ p \in PD_d   } \left \{ \alpha \in X_n^k : \alpha \ \mbox{is $\frac{d}{p}$-additive} \right \} } = \sum_{\{B\subset PD_d: B \neq \emptyset \}} (-1)^{\abs{B}+1} \abs{ \bigcap_{p\in B}  \left \{ \alpha \in X_n^k : \alpha \ \mbox{is $\frac{d}{p}$-additive} \right \}}\\  \nonumber
  =\sum_{\{B\subset PD_d: B \neq \emptyset \}} (-1)^{\abs{B}+1} \abs{   \left \{ \alpha \in X_n^k : \alpha \ \mbox{is $\frac{d}{\prod_{p\in B}p}$-additive} \right \}}.
  \end{gather}
  Now we apply Lemma \ref{number of d-additive sequences} (a) and derive the desired equality.
  \epr 
  \begin{coro} \label{main coro for the formula}
  	Let  $d$ be any divisor of $(k,n)$.  Then (here $\prod_{p\in \emptyset} p = 1$):
  	\begin{gather} \nonumber \abs{\left  \{\alpha \in X_n^k: \alpha(0)= 0, \ P(\alpha) = d \right  \} } =  \sum_{ B\subset PD_d } (-1)^{\abs{B}} \left ( \begin{array}{c} \frac{d}{\prod_{p\in B}p}\frac{n+2}{k+1}-1 \\ \frac{d}{\prod_{p\in B}p}-1 \end{array} \right).
  	\end{gather}
  \end{coro}\bpr To prove this one needs to combine Corollary \ref{the seqs of period d}, Lemma \ref{lemma for number of d/p}, and Lemma \ref{number of d-additive sequences} (a).   \epr
  \subsection{Formula}
  
  Now using Corollaries  \ref{CGammak-1n}, \ref{main coro for the formula} we derive the desired formula for  $\abs{C^\Gamma_{\mc T_{k-1}}(\mc T_n)}$, when $2\leq k < n+1$, however we will show that it holds for $k=1$ as well:
  \begin{prop} \label{final formula} Let $1\leq k < n+1$.   We have ( $PD_d$ is defined in  Definition \ref{def of d-additive})) 
   \begin{gather} \label{help formula for orbits 11234} \abs{C^{\Gamma_n}_{\mc T_{k-1}}(\mc T_n)}= \sum_{\left \{(d,B): \begin{array}{c}  d \ \mbox{is a divisor of} \ (k,n) \\ 
   	B\subset PD_d \end{array}  \right \} } \frac{(-1)^{\abs{B}}}{d}   \left ( \begin{array}{c} \frac{d}{\prod_{p\in B}p}\frac{n+2}{k+1}-1 \\ \frac{d}{\prod_{p\in B}p}-1 \end{array} \right).
   \end{gather}
  \end{prop}	
  \bpr As we explained the formula is proved for $2\leq k < n+1$. If $k=1$, then $k+1$ is prime. On the other hand, if $k+1$ is prime, then using Definition \ref{def of d-additive} we compute:
  \begin{gather} \nonumber 
  \{d\geq 1 : d \ \mbox{is a divisor of} \ (k,n) \}= \left \{ \begin{array}{c c c} \{1, k+1\} &  \mbox{if} \ \frac{n+2}{k+1} \in \ZZ  & \mbox{in this case} \  PD_{k+1}=\{k+1\} \\ \{k+1\}  & \mbox{if} \ \frac{n+2}{k+1} \not \in \ZZ & \mbox{in this case} \  PD_{k+1}=\emptyset .  \end{array} \right.
  \end{gather} 
  Therefore the formula  \eqref{help formula for orbits 11234} for $1\leq k<n+1$, $k+1$ prime  reduces to (note that   always $PD_1=\{\emptyset\}$):
  \begin{gather} \nonumber  \abs{ C^\Gamma_{\mc T_{k-1}}(\mc T_n)}=  \left \{ \begin{array}{c c}  1 + \frac{1}{k+1}\left (\binom{n+2-1}{k+1-1}- \binom{\frac{n+2}{k+1}-1}{0} \right )&  \mbox{if} \ \frac{n+2}{k+1} \in \ZZ  \\ \frac{1}{k+1}\binom{n+2-1}{k+1-1}   & \mbox{if} \ \frac{n+2}{k+1} \not \in \ZZ  \end{array} \right. =\left \{ \begin{array}{c c}  \frac{k}{k+1} + \frac{1}{k+1}\binom{n+1}{k}&  \mbox{if} \ \frac{n+2}{k+1} \in \ZZ  \\ \frac{1}{k+1}\binom{n+1}{k}   & \mbox{if} \ \frac{n+2}{k+1} \not \in \ZZ  \end{array} \right. .
  \end{gather}
  The formula above for $k=1$ is the same as in Corollary \ref{C_A1inA_n+1aut} . Therefore \eqref{help formula for orbits 11234} holds.
  \epr

   \subsection{Remarks and examples} \label{remarks and examples for A_n}
   The notations $\mu(x)$,  $\omega(x)$,  $\Omega(x)$ are from Section \ref{notations}.
   In particular 
   \begin{gather}\label{mu1}  B \subset \{\mbox{Primes}\}, \ \abs{B} < \infty   \ \Rightarrow \ \mu \left ( \prod_{p \in B} p  \right) = (-1)^{\abs{B}} \\ 
   \label{mu2} x \in \ZZ_{\geq 1}, \mu(x)\neq 0 \ \Rightarrow \  x=\prod_{ \{ p: p \ \mbox{is a prime factor of } \ x \} }  p. \end{gather}
   We will give a better version of formula \eqref{help formula for orbits 11234}, Note first that: 
   
    \begin{lemma} \label{divisors of k,n}
    	Let $1\leq k < n+1$ be integers. Let us denote   $D=g.c.d.(k+1,n+2) $, then:
    	$$\left \{d: d \ \mbox{is  divisors of} \ (k,n) \right  \}=\left \{\frac{k+1}{D} x : x \ \mbox{is a  divisors of} \  D \right \}.$$ 
    \end{lemma} \bpr  We have $k+1=D \gamma$,  $n+2 = D z$ for some positive integers $\gamma,z$.
    
    By definition  $d $ is divisor of $(k,n)$ iff $k+1 = d d'$, $n+2= d' d''$ for some positive integers $d',d''$.
    
    If $k+1 = d d'$, $n+2= d' d''$, then  $d'$ is a common divisor of $k+1,n+2$, and therefore there exists divisor  $x$ of $D$ such that $d' x = D$ and from  $k+1 = d d'=\gamma D =\gamma d' x$ follows $d=\gamma x$.
    
    Conversely, let  $d=\gamma x$ for some  divisor  $x$ of $D$, then $D=x y$ for some positive integer $y$ and $k+1= D \gamma =(\gamma x) y $, $n+2=y (x z)$, hence $\gamma x$ is a divisor of $(k,n)$.   \epr
    
    Using Lemma \ref{divisors of k,n} we can give another form of our formula
    
    \begin{coro}
    	\label{final formula third form} Let $1\leq k < n+1$. Let $D=g.c.d(k+1,n+2)$. Then
    	\begin{gather}  \label{formula third form} \abs{C^{\Gamma_n}_{\mc T_{k-1}}(\mc T_n) }= \sum_{\left \{(x,y) :\begin{array}{c} x,y,x/y \in \ZZ_{\geq 1},  \mu(x/y) \neq 0  \\ x, y \ \mbox{are  divisors of} \ D \end{array} \right \}} \frac{ D \mu(x/y)}{(k+1)x}   \left ( \begin{array}{c} y \frac{n+2}{D}-1 \\ y \frac{k+1}{D}-1 \end{array} \right). \end{gather}
    \end{coro}
      \bpr    
   By Lemma \ref{important lemma for orbits} (a),  Remark \ref{remark for gcd}, and the definition of $PD_d$ in Definition \ref{PD_d} follows that  the  assignment from left to right is    well defined:
   \begin{gather}  \label{help formula for orbits 112345}
   \left \{(d,B): \begin{array}{c}  d \ \mbox{is a divisor of} \ (k,n) \\ 
   B\subset PD_d \end{array}  \right \} \leftrightarrow \left \{(d,x) :\begin{array}{c} d,x,d/x \in \ZZ_{\geq 1},  \mu(d/x) \neq 0  \\ d, x \ \mbox{are  divisors of} \ (k,n) \end{array} \right \}\\
   (d,B) \mapsto \left (d, \frac{d}{\prod_{p\in B} p} \right )   \qquad   \left (d, \mbox{the set of  prime factors of} \ d/x\right ) \mapsfrom (d,x)
   \end{gather}
   whereas the assignment from right to left is well defined by the following arguments. Let  $(d,x)$ be an element  in the set on the RHS, in particular $d=x y$ with $x,y \in \ZZ$ and let $p$ be a prime divisor of $d/x=y$, in particular $y=p y'$ with $y'\in \ZZ$ and  $d=p y' x$, therefore $\frac{k+1}{d/p}= p \frac{k+1}{d} \in \ZZ$, $ \frac{d}{p} \frac{n+2}{k+1}= y' x \frac{n+2}{k+1} \in \ZZ$, since $d,x$ are both divisors of $(k,n)$, therefore $d/p$ is a divisor of $(k,n)$ and $p\in PD_d$. Using \ref{mu2} one shows that these assignments are mutually inverse.   
   Therefore we can change the parameters in \eqref{help formula for orbits 11234}  via \eqref{help formula for orbits 112345}, which taking into account  \eqref{mu1}, \eqref{mu2} expresses $\abs{C^{\Gamma_n}_{\mc T_{k-1}}(\mc T_n)}$  as a sum  $\sum_{} \frac{\mu(d/x)}{d}   \left ( \begin{array}{c} x \frac{n+2}{k+1}-1 \\ x-1 \end{array} \right)$, where $(d,x)$ varies in the set $\left \{(d,x) :\begin{array}{c} d,x,d/x \in \ZZ_{\geq 1},  \mu(d/x) \neq 0  \\ d, x \ \mbox{are  divisors of} \ (k,n) \end{array} \right \}$.
  Using  Lemma \ref{divisors of k,n} we obtain formula \eqref{formula third form}.
      \epr
    \begin{coro} \label{n+2 prime coro} Let $1\leq k < n+1$ be integers.
   	
   		{\rm (a)} If  $k+1$, $n+2$ are coprime, then 	$  \abs{ C^{\Gamma_n}_{\mc T_{k-1}}(\mc T_n) }= \frac{1}{k+1}\binom{n+1}{k} =\frac{1}{n+2} \binom{n+2}{k+1}.
   		$

   		{\rm (b)} If  $g.c.d.(n+2,k+1)=p$ is prime, then 
   		$  \abs{ C^{\Gamma_n}_{\mc T_{k-1}}(\mc T_{n})}=  \frac{1}{k+1}\binom{n+1}{k}+ \frac{p-1}{k+1} \binom{(n+2)/p -1}{(k+1)/p-1}. 
   		$
   		
   		{\rm (c)} If	$k+1$ is prime, then 
   		$  \abs{ C^{\Gamma_n}_{\mc T_{k-1}}(\mc T_n)}=  \left \{ \begin{array}{c c}  \frac{k}{k+1} + \frac{1}{k+1}\binom{n+1}{k}&  \mbox{if} \ \frac{n+2}{k+1} \in \ZZ  \\ \frac{1}{k+1}\binom{n+1}{k}   & \mbox{if} \ \frac{n+2}{k+1} \not \in \ZZ.  \end{array} \right.
   		$
   		
   	{\rm	(d) } 	$  \abs{C^{\Gamma_k}_{\mc T_{k-1}}(\mc T_k) }=1.
   		$
   \end{coro}  
   \bpr {\rm (a)} The index set of formula \eqref{formula third form} has a single element, the only summand is $\frac{1}{k+1}\binom{n+1}{k}$.

    {\rm (b)} Now 
    $
    \left \{(x,y) :\begin{array}{c} x,y,x/y \in \ZZ_{\geq 1},  \mu(x/y) \neq 0  \\ x, y \ \mbox{are  divisors of} \ p \end{array} \right \} =  \left \{(p,p), \left (p,1\right ),  \left (1,1 \right ) \right \}.  \nonumber
    $
    Therefore \eqref{formula third form} is:
    \begin{gather} \nonumber 
    \abs{C^{\Gamma_n}_{\mc T_{k-1}}(\mc T_{n})} = \frac{1}{k+1}\binom{n+1}{k} - \frac{1}{k+1}\binom{(n+2)/p-1}{(k+1)/p-1} +  \frac{p}{k+1}\binom{(n+2)/p-1}{(k+1)/p-1}.
    \nonumber
    \end{gather} 
    
    (c)  We proved this in the proof of Proposition \ref{final formula}.
    
    (d)  Now $k=n$ and $g.c.d.(k+1,n+2)=g.c.d.(k+1,k+2)=1$ and we use (a).
   \epr

   Another case, when the formula has a simple form is:
   \begin{coro} \label{n=k+delta}  Let $1\leq k$, $1\leq \Delta$ be integers abd $\Delta+1$ be prime. Then holds {\rm (d)} after  \eqref{formula third form intro}.
   \end{coro}   
   \bpr Now  $D=g.c.d(k+1,n+2)=g.c.d(k+1,k+1+(\Delta+1))=1$, if $(k+1)/(\Delta+1)\not \in \ZZ$, and $D=\Delta+1$, if $(k+1)/(\Delta+1) \in \ZZ$, hence using Corollary  \ref{n+2 prime coro} (a), (b) we get:
   \begin{gather} \nonumber 
   \abs { C^{\Gamma_{k+\Delta}}_{\mc T_{k-1}}(\mc T_{k+\Delta}) } =  \left \{ \begin{array}{c} \frac{1}{k+1}\binom{k+\Delta+1}{k}    \\[4mm]  \frac{1}{k+1}\binom{k+\Delta+1}{k}  + \frac{\Delta}{k+1}\binom{(k+\Delta+2)/(\Delta+1)-1}{(k+1)/(\Delta+1)-1} \end{array} \right.
   =	\left \{ \begin{array}{c}  \frac{1}{\Delta +1} \binom{k+\Delta+1}{\Delta}    \\[4mm]    = \frac{1}{\Delta +1} \binom{k+\Delta+1}{\Delta}  +  \frac{\Delta}{k+1}\frac{k+1}{\Delta+1} \end{array} \right.,
   \nonumber
   \end{gather} where we used the equality $\frac{1}{k+1}\binom{k+\Delta+1}{k}= \frac{1}{\Delta +1} \binom{k+\Delta+1}{\Delta}$   hence the corollary is proved.
   \epr

    \subsection{Interpretation via polygons} \label{Polygons}
    
    Let $1\leq k< n+1$ and let us denote two sets:  \begin{gather} \label{Pn Pnk}
    P_n=\left \{\exp\left (\frac{2 \pi}{n+2} j \right ) \right \}_{j=0}^{n+1}  \qquad \qquad P_n^k = \{U\subset P_n : \abs{U}=k+1\}.
    \end{gather}
    
    \textit{Note that we can view  $P_n$  as the set of vertices of a regular $(n+2)$-gon  embedded into the unit circle, and then $P_n^k$ is the set of $k+1$-subgons of $P_n$. }

    \begin{lemma} 
    	
    	Let $1\leq k< n+1$. 	
    	The following map is a bijection ($X_n^k$ is in \eqref{Xnk}): 	\begin{gather}
    	X_n^k \ni \alpha \mapsto \left \{\exp\left (\frac{2 \pi}{n+2} (\alpha_j+j) \right ) : 0 \leq j \leq k \right \}\in P_n^k.
    	\end{gather} This bijection transports the permutation $S$ of $X_n^k$ from Corollary \ref{Serre on C_0A_kinA_n} to rotation by angle $\frac{2 \pi}{n+2}$.
    \end{lemma}
    \bpr  Denote $Y_n^k=\left \{ \bd \{0,1,\dots,k\} & \rTo^{u} &  \{0,1,\dots,n+1\} \ed : 0\leq u_0 < u_1 < \dots < u_k \leq n+1\right \}$, our map is composition of bijections $X_n^k \rightarrow Y_n^k \rightarrow P_n^k$, which assign  $\alpha \mapsto \{\alpha _j+j\}_{j=0}^k$; $u \mapsto \left \{\exp\left (\frac{2 \pi}{n+2} u_j \right ) \right \}_{j=0}^{k}$, respectively, therefore it is bijection. 
    
    Take now any $\alpha \in X_n^k$. If $\alpha_k < n-k+1$, then using \eqref{Serre 22}  we see that: 
    $$S(\alpha)=\alpha+1 \mapsto  \left \{\exp\left (\frac{2 \pi}{n+2} (\alpha_j+j+1) \right ) \right \}_{j=0}^{k} = \exp\left (\frac{2 \pi}{n+2}  \right ) \left \{\exp\left (\frac{2 \pi}{n+2} (\alpha_j+j) \right ) \right \}_{j=0}^{k}  $$ whereas if $\alpha_k = n-k+1$, then  using \eqref{Serre 23}  we see that:
    \begin{gather}
    S(\alpha)=(0,\alpha_0,\alpha_1\dots. \alpha_{k-1}) \mapsto  \{1\}\cup \left \{\exp\left (\frac{2 \pi}{n+2} (\alpha_{j-1}+j) \right ) \right \}_{j=1}^{k} \nonumber \\ = \{1\}\cup \left \{\exp\left (\frac{2 \pi}{n+2} (\alpha_{j}+j+1) \right ) \right \}_{j=0}^{k-1}  = \exp\left (\frac{2 \pi}{n+2}  \right ) \left \{\exp\left (\frac{2 \pi}{n+2} (\alpha_j+j) \right ) \right \}_{j=0}^{k}.  \nonumber 
    \end{gather}
    Hence the lemma is proved.
    \epr  From this lemma and Corollary \ref{orbits as elements in our set}, it follows
    \begin{coro}  \label{Polygons coro}  Let $1\leq k< n+1$. 	
    	The set  $C^{{\rm Aut}(D^b(A_{n+1}))}_{D^b(A_k)}(D^b(A_{n+1}))$  can be identified with the set of equivalence classes (up to rotation) of $(k+1)$-subgons in a regular $(n+2)$-gon,  where by a ``subgon'' we mean a convex  polygon, whose  vertices are subset of the set of vertices.  
    \end{coro}

    \begin{ex}[ $C^{{\rm Aut}(D^b(A_{6}))}_{D^b(A_3)}(D^b(A_{6}))$] \label{another example}
   Now  $g.c.d. (3+1,6+1)= 1$ and  Corollary \ref{n+2 prime coro} (a) results in $\abs{C^{{\rm Aut}(D^b(A_{6}))}_{D^b(A_3)}(D^b(A_{6}))} = 5$.  Тhe $5$  equivalence classes are shown in  Figure \ref{figure for c36}.
    \begin{figure}	\begin{center}
    		\begin{tikzpicture}[scale=1]

    		\begin{scope}[shift={(-4,0)}]
    		\newdimen\R
    		\R=0.8cm
    		\coordinate (center) at (0,0);
    		\draw (0:\R) node[scale=0.8] {{$\bullet$}} \foreach \x in {0,51.4285,...,360} { -- (\x:\R)node [scale=0.8] {{$\bullet$}} };
    		\draw (51.4285*2:\R) -- (51.4285*-1:\R);
    		\end{scope}

    		\begin{scope}[shift={(-2,0)}]
    		\newdimen\R
    		\R=0.8cm
    		\coordinate (center) at (0,0);
    		\draw (0:\R) node[scale=0.8] {{$\bullet$}} \foreach \x in {0,51.4285,...,360} { -- (\x:\R)node [scale=0.8] {{$\bullet$}} };
    		\draw (51.4285*2:\R) -- (51.4285*-2:\R)--(0:\R);
    		\end{scope}

    		\begin{scope}[shift={(2,0)}]
    		\newdimen\R
    		\R=0.8cm
    		\coordinate (center) at (0,0);
    		\draw (0:\R) node[scale=0.8] {{$\bullet$}} \foreach \x in {0,51.4285,...,360} { -- (\x:\R)node [scale=0.8] {{$\bullet$}} };
    		\draw (51.4285*2:\R) -- (51.4285*-3:\R)--(0:\R);
    		\end{scope}

    		\begin{scope}[shift={(0,0)}]
    		\newdimen\R
    		\R=0.8cm
    		\coordinate (center) at (0,0);
    		\draw (0:\R) node[scale=0.8] {{$\bullet$}} \foreach \x in {0,51.4285,...,360} { -- (\x:\R)node [scale=0.8] {{$\bullet$}} };
    		\draw (51.4285*2:\R) -- (51.4285*-2:\R)  (51.4285*-1:\R)--(51.4285*1:\R);
    		\end{scope}
    		
    		\begin{scope}[shift={(4,0)}]
    		\newdimen\R
    		\R=0.8cm
    		\coordinate (center) at (0,0);
    		\draw (0:\R) node[scale=0.8] {{$\bullet$}} \foreach \x in {0,51.4285,...,360} { -- (\x:\R)node [scale=0.8] {{$\bullet$}} };
    		\draw (51.4285*2:\R) -- (51.4285*-3:\R)  (51.4285*-1:\R)--(51.4285*1:\R)   (51.4285*-3:\R)--(51.4285*-1:\R);
    		\end{scope}

    		\end{tikzpicture}
    	\end{center} \caption{$c_{3,6}=5$} \label{figure for c36}
    	\end{figure}
    \end{ex}
     \subsection{Non-commutative curves in $\mc T_n$} \label{subsection ncc in A_n}
  From Remark \ref{reminder} we know that $C_l(\mc T_n)=\emptyset$ for $l\geq 1$.
 
 The set of genus $0$ non-commutative curves $C_0^{\{\rm Id \}}(\mc T_n)$ is the same as $C^{\{\rm Id \}}_{\mc T_{2-1}}(\mc T_n)$ and from Corollaries \ref{identity group case for A_kinA_n+1},  \ref{n+2 prime coro} (c) it follows: 
  \begin{coro} \label{coro for numC0} Let $n\geq 1$, then $\abs{C^{\{\rm Id\}}_{0}(\mc T_n)}=\binom{n+2}{3}$,   $ 
 	\abs{ C^{\Gamma_n}_{0}(\mc T_n) } =\left \{ \begin{array}{c c}\frac{(n+1)n}{6}  & \mbox{if} \ \frac{n-1}{3} \not \in  \ZZ \\ \frac{(n+1)n+4}{6} & \mbox{if} \ \frac{n-1}{3} \in \ZZ\end{array} \right. .
 	$
 \end{coro}
 It remains to determine $ C^{\{\rm Id\}}_{-1}(\mc T_n),  
  C^{\Gamma_n}_{-1}(\mc T_n) $.

 \begin{lemma} \label{C_-1A_n+1}
 	The set $C^{\{\rm Id\}}_{-1}(\mc T_n)$ (see Definition \ref{C_l}) can be described as follows (using the bijection described in Proposition \ref{bijection 123}): $\abs{C^{\{\rm Id\}}_{-1}(\mc T_n) }=2 \binom{n+2}{4}$ and
 	\begin{gather}\label{C^{{rm Id}}_{-1}}  C^{\{\rm Id\}}_{-1}(\mc T_n)=\left \{\langle  s_{\alpha,\beta}, s_{i,j} \rangle  : 0\leq \alpha \leq \beta <i-1\leq j-1 \leq n-1  \ \mbox{or }  \   0\leq \alpha<i\leq j < \beta \leq n   \right \} \qquad  
 	. \nonumber 
 	 \end{gather}
 \end{lemma}
 \bpr   From Remark \ref{functor from NP tp NP} we see that  $C^{\{\rm Id\}}_{-1}(\mc T_1)$ is emtpy,  so the Lemma holds for $n=1$. Assume that it holds for some $n\geq 1$. 
 
 Let $\mc A \in C^{\{\rm Id\}}_{-1}(\mc T_{n+1})$, due to \eqref{exceptional objects} we can assume that $\mc A=\langle s_{\alpha,\beta}, s_{i,j} \rangle $ for some $0\leq i\leq j \leq n+1$, $0\leq \alpha \leq \beta \leq n+1$, where $ s_{\alpha,\beta}, s_{i,j}$ is an orthogonal exceptional pair. If $j <n+1$ and  $\beta  <n+1$, then since $Rep_k(A_{n+1})$ is a full Serre subcategory in $Rep_k(A_{n+2})$ it follows that $\langle  s_{\alpha,\beta}, s_{i,j} \rangle \in C^{\{\rm Id\}}_{-1}(\mc T_{n})$ and by the induction assumption $\mc A = \langle  s_{\alpha,\beta}, s_{i,j} \rangle  $ for some $\begin{array}{c} 0\leq \alpha \leq \beta <i-1\leq j-1 \leq n-1 \\ \mbox{or }  \   0\leq \alpha<i\leq j < \beta \leq n \end{array} $. So it remains to consider the case ($ \beta = n+1 $ or  $ j=n+1$). Since   $ s_{\alpha,\beta}, s_{i,j}$ is an orthogonal exceptional pair, it is enough to consider the case $j=n+1$. Using Lemma \ref{help lemma 1}  we see that $\beta <n+1$. 
 
  If $\alpha<i$, then Lemma \ref{help lemma 1} (a) implies that $\beta<i-1$ and therefore  $0\leq \alpha \leq \beta <i-1\leq j-1 = n$.  
 If $\alpha \geq i$, then $i\leq \alpha \leq \beta < j$ and Lemma \ref{help lemma 1} (d) implies that $i< \alpha$, hence  $\mc A = \langle   s_{i,j}, s_{\alpha,\beta} \rangle  $ with  $i< \alpha \leq \beta < j$.
 
  From Lemma \ref{help lemma 1} it follows that if $\begin{array}{c} 0\leq \alpha \leq \beta <i-1\leq j-1 \leq n \\ \mbox{or }  \   0\leq \alpha<i\leq j < \beta \leq n+1 \end{array}$, then  $ \langle   s_{\alpha,\beta},  s_{i,j} \rangle \in  C^{\{\rm Id\}}_{-1}(\mc T_{n+1})$.
  
  So far  we showed that the induction assumption implies \eqref{C^{{rm Id}}_{-1}} with $n+1$ instead of $n$ and the first row of the following sequence of equalities\footnote{Note that now if $\langle  s_{\alpha,\beta}, s_{i,j}  \rangle=\langle  s_{\alpha',\beta'}, s_{i',j'}  \rangle$, then $\alpha = \alpha'$, $\beta = \beta'$, $i = i'$, $j =j'$  or  $\alpha = i'$, $\beta = j'$, $i =  \alpha'$, $j =  \beta'$,  (we take into account also Remark \ref{orthogonal pairs})}, for the other rows we use again the induction assumption,  Lemma \ref{C_eoA_n+1}, the Hockey-Stick equality, and Pascal's rule:
  \begin{gather} \nonumber 
  \abs{C^{\{\rm Id\}}_{-1}(\mc T_{n+1}) }= \abs{C^{\{\rm Id\}}_{-1}(\mc T_{n})} + 
  \abs{\left \{\langle  s_{\alpha,\beta}, s_{i,j} \rangle  : \begin{array}{c} 0\leq \alpha \leq \beta <i-1\leq j-1 = n \\ \mbox{or }  \   0\leq \alpha<i\leq j < \beta = n+1 \end{array} \right \} } \end{gather}  \begin{gather} \nonumber 
  = {\scriptstyle \abs{C^{\{\rm Id\}}_{-1}(\mc T_{n})} + 
  \abs{\left \{\langle  s_{\alpha,\beta}, s_{i,j} \rangle  : \begin{array}{c} 0\leq \alpha \leq \beta <i-1\leq j-1 = n \end{array} \right \}}  + \abs{\left \{\langle  s_{\alpha,\beta}, s_{i,j} \rangle  :   0\leq \alpha<i\leq j < \beta = n+1 \right \}} } \end{gather}  \begin{gather}
  \nonumber = {\abs{C^{\{\rm Id\}}_{-1}(\mc T_{n})} + \sum_{i=2}^{n+1} \abs{  C^{\{\rm Id\}}_{\mc T_0, \KK}(\mc T_{i-2})}+ \sum_{\alpha=0}^{n} \abs{  C^{\{\rm Id\}}_{\mc T_0, \KK}(\mc T_{n-\alpha-1})} } \end{gather}  \begin{gather} 
  \nonumber = 2 \binom{n+2}{4} +  \sum_{i=2}^{n+1} \frac{(i-1) i}{2}  + \sum_{\alpha=0}^{n-1}  \frac{(n-\alpha) (n-\alpha+1)}{2} =  2 \binom{n+2}{4} +  \sum_{i=2}^{n+1} \frac{(i-1) i}{2}  + \sum_{j=1}^{n}  \frac{j (j+1)}{2} \end{gather}  \begin{gather} =2\binom{n+2}{4} +  2\sum_{i=2}^{n+1} \binom{i}{2} =
  \nonumber 2 \left (  \binom{n+2}{4} +   \sum_{i=2}^{n+1} \binom{i}{2}\right )  = 2 \left (  \binom{n+2}{4} +   \binom{n+2}{3} \right ) =2   \binom{n+3}{4}.
    \end{gather}
 The lemma is proved. 
 \epr

Next we want to compute $\abs{C^{ \Gamma }_{-1}(\mc T_n)}$ for $\Gamma =  {\rm Aut}_{\KK}(\mc T_n)$.

\begin{lemma} \label{lemma orbits for -1}
	Any orbit in $C_{-1}^{\{\rm Id\}}(\mc T_n)$ of  $\langle S \rangle$ has the form (here  $0\leq \beta < i-1\leq j-1 \leq n-1$):
	\begin{gather}
	\nonumber \langle s_{0,\beta}, s_{i,j} \rangle 
	\mapsto\cdots \mapsto \langle s_{n-j,n-j+\beta}, s_{i+n-j,n} \rangle \\ \nonumber  
	\mapsto  \langle s_{n-j+1,n-j+1+\beta}, s_{0,i+n-j} \rangle 
	\mapsto\dots\mapsto   \langle s_{n-i+1,n-i+1+\beta}, s_{j-i,n} \rangle \\[-2mm] \label{orbit for -1} \\[-2mm] \nonumber 
	\mapsto  \langle s_{0,j-i}, s_{n-i+2,n-i+2+\beta} \rangle
	\mapsto\dots\mapsto   \langle s_{i-2-\beta,j-\beta-2}, s_{n-\beta,n} \rangle \\  \nonumber 
	\mapsto  \langle s_{i-1-\beta,j-1-\beta}, s_{0,n-\beta} \rangle
	\mapsto\dots\mapsto   \langle s_{i-1,j-1}, s_{\beta,n} \rangle   \mapsto \ \mbox{return to} \  \langle s_{0,\beta}, s_{i,j} \rangle.
	\end{gather}
	The orbit reduces to the first and second rows iff \eqref{condition1} holds and then the orbit has $\frac{n}{2}+1$ elements
	\begin{gather} \label{condition1} 0\leq \beta < \frac{n}{2} \ \ \mbox{and}  \ \ i=\frac{n}{2} +1 \ \ \mbox{and}  \ \ j = \frac{n}{2}+1+\beta.  \end{gather}

	Otherwise the elements in \eqref{orbit for -1} are pairwise distinct and in this case the orbit has $n+2$ elements. 
\end{lemma}
\bpr From Lemma  \ref{C_-1A_n+1}  we know that $\langle s_{0,\beta}, s_{i,j} \rangle \in C^{\{\rm Id\}}_{-1}(\mc T_n)$, when $0\leq \beta < i-1\leq j-1 \leq n-1$. Starting from such an element by repeatedly  applying the Serre functor we obtain \eqref{orbit for -1} with the help of Proposition \ref{Serre on exceptional objects}.   

 From  Lemma \ref{C_-1A_n+1} we know that for  $\langle  s_{\alpha,\beta}, s_{i,j} \rangle \in C^{\{\rm Id\}}_{-1}(\mc T_n) $ we have  either $0\leq \alpha \leq \beta <i-1\leq j-1 \leq n-1$  \mbox{or }  $   0\leq \alpha<i\leq j < \beta \leq n $. In the second case applying $n-\beta$ times the Serre functor we obtain an element of the type  $\langle  s_{\alpha',n}, s_{i',j'} \rangle$ with $ 0\leq \alpha'<i'\leq j' < \beta < n $, and now applying $S$ once more  we get $\langle  s_{0,\alpha'}, s_{i'+1,j'+1} \rangle$ with  $0\leq \alpha' < (i'+1)-1\leq (j'+1) \leq n-1$. In the first case, $0\leq \alpha \leq \beta <i-1\leq j-1 \leq n-1$, we apply the inverse of the Serre functor several time to get element of the form   $\langle s_{0,\beta}, s_{i,j} \rangle$, Thus, we proved the first part of the Lemma.

Note that the elements on the first and the third row of \eqref{orbit for -1}  satisfy the first condition  whereas those in the second and the fourth row satisfy the second condition  in \eqref{C^{{rm Id}}_{-1}}. Therefore the orbit reduces to a shorter sequence than the comlplete sequence in \eqref{orbit for -1} iff 
\begin{gather}   \beta = j-i \quad i=n-i+2 \quad j=n-i+2+\beta \nonumber \end{gather}
which is equivalent to $i=\frac{n}{2} +1 $,  $ j = \frac{n}{2}+1+\beta$ and now
$0\leq \beta < i-1$ is equivalent to $\beta <\frac{n}{2}$. Finally note that the first, second, third, fourth row in \eqref{orbit for -1} contain $1+n-j$, $1+j-i$, $1+i-2-\beta$, $1+\beta$ elements, respectively, which sum up to $n+2$.  Hence we proved the lemma.\epr
\begin{coro} \label{two types of orbits for C-1}
	If $n$ is odd, then all orbits of ${\rm Aut}_{\KK}(\mc T_n)$ in $C_{-1}^{\{\rm Id\}}(\mc T_n)$  have $n+2$ elements. Otherwise, there are exactly $\frac{n}{2}$ orbits with  $\frac{n}{2}+1$ elements and all the rest orbits have $n+2$ elements.
\end{coro}
\bpr From Lemma \ref{lemma orbits for -1} we see that if  there is an orbits of length different from $n+2$, then $n$ must be even, and the orbit must have the form 

\begin{gather}
\nonumber \langle s_{0,\beta}, s_{\frac{n}{2} +1,\frac{n}{2} +1+\beta} \rangle 
\mapsto\cdots \mapsto \langle s_{\frac{n}{2} -\beta-1,\frac{n}{2} -1}, s_{n-\beta,n} \rangle
\mapsto  \langle s_{\frac{n}{2} -\beta,\frac{n}{2}}, s_{0,n-\beta} \rangle 
\mapsto\dots\mapsto   \langle s_{\frac{n}{2},\frac{n}{2}+\beta}, s_{\beta ,n} \rangle.
\end{gather}
Furthermore, recalling Remark \ref{orthogonal pairs}  we deduce that  if  $\langle  s_{\alpha,\beta}, s_{i,j}  \rangle=\langle  s_{\alpha',\beta'}, s_{i',j'}  \rangle$, where $\langle  s_{\alpha,\beta}, s_{i,j}  \rangle \in C_{-1}^{\{\rm Id\}}(\mc T_n)$,  then $\alpha = \alpha'$, $\beta = \beta'$, $i = i'$, $j =j'$  or  $\alpha = i'$, $\beta = j'$, $i =  \alpha'$, $j =  \beta'$, and it follows that   for two different $0\leq \beta < \frac{n}{2}$, $0\leq \beta' < \frac{n}{2}$, $\beta \neq  \beta'$ the corresponding orbits are different. 
\epr
We can prove now \eqref{coro for numC-1}.
 From Corollary \ref{two types of orbits for C-1} we know that if $n$ is odd, then all the orbits in $ C^{\{\rm Id\}}_{-1}(\mc T_n) $
 contain the same number of elements $n+2$, on the other hand from Lemma \ref{C_-1A_n+1} we know that $ C^{\{\rm Id\}}_{-1}(\mc T_n)=\frac{(n-1)n (n+1)(n+2)}{12} $, therefore the number of orbits, which is the same as $\abs{ C^{\Gamma_n}_{-1}(\mc T_n)  }$, is $\frac{(n-1)n (n+1)}{12}$.   If $n$ is even, then from Corollary \ref{two types of orbits for C-1}  we know that there are $\frac{n}{2}$ orbits with  $\frac{n}{2}+1$ elements and all the rest orbits have $n+2$ elements. Therfore $\abs{ C^{\Gamma_n}_{-1}(\mc T_n) }-\frac{n}{2}=\left (\frac{(n-1)n (n+1)(n+2)}{12} -\frac{n}{2} \frac{n+2}{2} \right )/(n+2)$
 and it follows that $\abs{ C^{\Gamma_n}_{-1}(\mc T_n)  }=\frac{n(n^2+2)}{12}$.
 
 \begin{remark}
 Using \eqref{coro for numC-1} and Corollary \ref{coro for numC0} one shows that  $\abs{ C^{\Gamma_n}_{-1}(\mc T_n) }>\abs{ C^{\Gamma_n}_{0}(\mc T_n) }$	for $n\geq 4 $. Otherwise  $\abs{ C^{\Gamma_3}_{-1}(\mc T_3) }=\abs{ C^{\Gamma_3}_{0}(\mc T_3) }=2$,  $\abs{ C^{\Gamma_2}_{-1}(\mc T_2) }=1=\abs{ C^{\Gamma_2}_{0}(\mc T_2) }$
 \end{remark}

\subsection{Some valued directed  graphs}  \label{some digraphs}

 For example, in Lemma \ref{C_eoA_n+1}  one sees  that $C_{D^b(A_1), \KK}(D^b(A_2))=C_{D^b(A_1), \KK}(N\PP^0)$   $= \{\langle s_{0,0}\rangle, \langle s_{1,1}\rangle,  \langle s_{0,1}\rangle \}$, whereas from Lemmas \ref{help lemma 2} \ref{Cid_A_kinA_n+1}, \ref{C_-1A_n+1} we know that the semi-orthogonal pairs  in  $C_{D^b(A_1), \KK}(N\PP^0) \times C_{D^b(A_1), \KK}(N\PP^0)$ are $(\langle s_{0,0}\rangle, \langle s_{1,1}\rangle)$, $(\langle s_{1,1}\rangle, \langle s_{0,1}\rangle), (\langle s_{0,1}\rangle, \langle s_{0,0}\rangle)$, therefore we can represent the valued  directed graph $G_{D^b(A_1), \KK}(N\PP^0)$ as in Figure \ref{Figure for digraph of A2}.
 The numbers are everywhere $1$ due to Remark \ref{hom and hom1}, the vanishings \eqref{vanishings general 2} and the fact that there are no orthogonal pairs in $N\PP^0$.
 
 Since $N\PP^{-1}$ is generated by an orthogonal exceptional pair $(E_1,E_2)$, one draws  the  digraph $G_{D^b(A_1), \KK}(N\PP^{-1})$ in Figure \ref{Figure for digraph of NP-1}. 
 
 If $l\geq 1$ in  \cite[Formula (25)]{DK41} is shown the complete lists of exceptional pairs and objects (up to shifts) in $l+1$ and before \cite[Formula (27)]{{DK41}} is explained that for any exceptional pair $(A,B)$ holds $\sum_{i\in \ZZ} \hom(A,B)=l+1$ (only one summand vanishes) and one obtains the  valued digraph drawn in Figure  \ref{Figure for digraph of NPl}.

 In the end we draw the graph $G_{D^b(A_1), \KK}(D^b(A_3))$.  From Lemma \ref{C_eoA_n+1} we know that the derived points in $D^b(A_3)$ are: 
 \begin{gather}\label{Db A1 in DbA3} C^{\{\rm Id\}}_{ D^b(A_1), \KK}(D^b(A_3)) =\{\langle  s_{0,0} \rangle , \langle  s_{1,1} \rangle, \langle  s_{2,2} \rangle, \langle  s_{0,1} \rangle, \langle  s_{0,2}\rangle, \langle  s_{1,2} \rangle  \}. \end{gather}
 Whereas from Lemma \ref{Cid_A_kinA_n+1} the non-commutative curves of genus $0$ in  $D^b(A_3)$ are (for writing the equalities we use Lemma \ref{help lemma 2}): 
 \begin{gather} \label{exceptional pairs in A_3}  C^{\{\rm Id\}}_{0}(D^b(A_3))=\left \{\begin{array}{c}\langle  s_{0,0}, s_{1,1} \rangle = \langle  s_{0,1}, s_{0,0} \rangle=\langle  s_{1,1}, s_{0,1} \rangle \\ \langle  s_{1,1}, s_{2,2} \rangle =  \langle  s_{1,2}, s_{1,1} \rangle=\langle  s_{2,2}, s_{1,2} \rangle \\  \langle  s_{0,1}, s_{2,2} \rangle = \langle  s_{0,2}, s_{0,1} \rangle= \langle  s_{2,2}, s_{0,2} \rangle\\ \langle  s_{0,0}, s_{1,2} \rangle = \langle  s_{0,2}, s_{0,0} \rangle= \langle  s_{1,2}, s_{0,2} \rangle \end{array} \right \} \end{gather} 
  
   Therefore the set of vertices of $G_{D^b(A_1), \KK}(\mc T)$ is   \eqref{Db A1 in DbA3}  and the set of  edges is obtained from the $12$ exceptional pairs in \eqref{exceptional pairs in A_3} and  the   orthogonal exceptional pairs $(\langle  s_{0,0} \rangle , \langle  s_{2,2} \rangle)$,  $(\langle  s_{2,2} \rangle , \langle  s_{0,0} \rangle)$, $(\langle  s_{0,2} \rangle , \langle  s_{1,1} \rangle)$, $(\langle  s_{1,1} \rangle , \langle  s_{0,2} \rangle)$ listed  in Corollary \ref{C_-1A_n+1}.  From the already given data we see that there are $6$ vertices, $12$ one-sided edges, and $2$ double-sided edges. We draw the digraph $G_{D^b(A_1), \KK}(\mc T)$  in $\RR^3$ as shown in Figure \ref{Figure for digraph of A3}, the triangle $\langle s_{1,1} \rangle, \langle s_{2,2} \rangle, \langle s_{0,2} \rangle$ is in a horizontal plane and the arrows from $\langle s_{0,1} \rangle$ to $ \langle s_{0,0} \rangle$ and  from $\langle s_{0,0} \rangle$ to $ \langle s_{1,2} \rangle$ are perpendicular to this plane.
  	The two dimensional simplexes in  $SC_{\mc J, \KK}(\mc T)$ are exactly the triples of vertices in triangles  of arrows which do not form an oriented loop (i.e. not of the form as in Figure \ref{Figure for digraph of A2}), these are triples of arrows as the arrows in the quiver $Q_1$ in \eqref{Q1}. When  a side in a triangle is an  undirected edge (double sided arrow), if for some of the two orientations of the edge   the triple of arrows can be done as in the quiver $Q_1$ in \eqref{Q1}, then this triangle is also a simplex. 
  	Such triangles correspond to exceptional triples in $\mc T$. In  Figure \ref{Figure for digraph of A3}  is presented also a geometric realization of the simplicial complex $SC_{\mc J, \KK}(\mc T)$.

  	\section{Non-commutative  counting in $D^b(D_4)$} \label{the numbers for D_4}
In this section $Q$ is a quiver of type $D_4$ with  orientation of  the arrows to point to the middle vertex, and  we denote: 
$ \mc T = D^b(Rep_\KK(Q)).
$
By  $S$  we will denote the Serre functor of $\mc T$. 
Obviously, the symmetry of the diagram generated by rotation of $2 \pi /3$ determines an auto-equivalence of $Rep_\KK(Q)$ and of $\mc T$, which we denote by $\kappa$.
We note that \cite[Theorem 4.1 and Table I]{Miyachi} imply that: 
\begin{gather} \label{Auto for D4}  {\rm Aut}_{\KK}(\mc T) \cong S_3 \times \ZZ   \qquad    {\rm Aut}_{\KK}(\mc T) =\langle [\kappa], [S], translation functor  \rangle. \end{gather}
Results of this section are   \eqref{table with numbers D4}, \eqref{table with numbers A1 in D4}, \eqref{triples in D4}.

Looking at \cite[Remark 6.1 in Version 2]{DK2} 
 we see that there are $12$ real roots and their dimension vectors are the thin dimension vectors and the  dimension vector with $2$ in the middle vertex and $1$ at the rest. 
Since the assignment to any exceptional representation in $Rep_k(Q)$ its dimension vector is bijection between the set of exceptional objects  (up to isomorphism) and the roots, which we just described, we  can choose representatives of the exceptional representations as follows: 
 \begin{gather}\label{exceptional representations} \{s_i:i=1,2,3\} \ \cup \  \{s_{io}:i=1,2,3\} \ \cup \ \{s_{ij}: 1\leq i < j \leq 3\}  \ \cup \ \{s_{123}, s_o, \delta\}  \end{gather} \begin{gather} 
 \ul{\dim}(s_1) = \begin{diagram}[size=1em]
 &        &     &         &  0    \\
 &        &     & \ldTo   &         \\
 1     & \rTo   & 0 &         &         \\
 &        &     & \luTo   &         \\
 &        &     &         &  0      
 \end{diagram} \quad  \ul{\dim}(s_{1o}) = \begin{diagram}[size=1em]
 &        &     &         &  0    \\
 &        &     & \ldTo   &         \\
 1     & \rTo   & 1 &         &         \\
 &        &     & \luTo   &         \\
 &        &     &         &  0      
 \end{diagram} \quad  \ul{\dim}(s_{12}) = \begin{diagram}[size=1em]
 &        &     &         &  1    \\
 &        &     & \ldTo   &         \\
 1     & \rTo   & 1 &         &         \\
 &        &     & \luTo   &         \\
 &        &     &         &  0      
 \end{diagram} \end{gather} \begin{gather}   \label{kappa}  
 s_2 = \kappa(s_1), s_3 = \kappa(s_2), s_{2o} = \kappa(s_{1o}), s_{3o} = \kappa(s_{2o}),   s_{23} = \kappa(s_{12}), s_{13} = \kappa(s_{23}),
 \end{gather} \begin{gather}  
 \ul{\dim}(s_o) = \begin{diagram}[size=1em]
 &        &     &         &  0    \\
 &        &     & \ldTo   &         \\
 0     & \rTo   & 1 &         &         \\
 &        &     & \luTo   &         \\
 &        &     &         &  0      
 \end{diagram} \quad  \ul{\dim}(s_{123}) = \begin{diagram}[size=1em]
 &        &     &         &  1    \\
 &        &     & \ldTo   &         \\
 1     & \rTo   & 1 &         &         \\
 &        &     & \luTo   &         \\
 &        &     &         &  1      
 \end{diagram} \quad  \ul{\dim}(\delta) = \begin{diagram}[size=1em]
 &        &     &         &  1    \\
 &        &     & \ldTo   &         \\
 1     & \rTo   & 1 &         &         \\
 &        &     & \luTo   &         \\
 &        &     &         &  1      
 \end{diagram}.
 \end{gather}

We will sometimes  notate $s_{ij}$ with $i>j$ but this by definition is the same as $s_{j i}$.

Recalling that the exceptional objects in $\mc T$ are shifts of exceptional representations and using again the bijection \eqref{bijection for eo} in Proposition \ref{bijection 123}  we conclude:
\begin{lemma}  \label{A1 in D4}
	The set $C^{\{\rm Id\}}_{D^b(A_1), \KK}(\mc T)$ has $12$ elements and it can be described as follows: 
	\begin{gather}\nonumber  C^{\{\rm Id\}}_{D^b(A_1), \KK}(\mc T) =\{\langle s_i \rangle\}_{i=1}^3 \ \cup \  \{\langle s_{io} \rangle\}_{i=1}^3 \ \cup \ \{\langle s_{ij} \rangle : 1\leq i < j \leq 3\}  \ \cup \ \{\langle s_{123} \rangle , \langle s_o \rangle , \langle \delta \rangle \} \nonumber. \end{gather}
\end{lemma}
	\begin{lemma} \label{no}
		 Any of the following  pairs  is not exceptional:
		
		$(\delta, s_i)$, $(s_i, \delta)$, $i=1,2,3$
		
		  $(\delta, s_{123})$, $(s_{123}, \delta)$, $(s_o, s_{123})$, $(s_{123}, s_o)$,  $(\delta, s_{o})$, $(s_{o}, \delta)$, $(s_o, s_{123})$,
		  
		   $(s_{123}, s_{io})$, $(s_{io}, s_{123})$, $i=1,2,3$ 
		   
		   $(s_{ij}, s_{o})$, $(s_{o}, s_{ij})$, $1\leq i < j \leq 3$.

		   \end{lemma}
		   \bpr One computes that for any of the listed pairs $(A,B)$ we have  $\scal{\ul{\dim}(B),\ul{\dim}(A)}\neq 0$ and then the lemma follows from Remark \ref{hom and hom1}.\epr 
		   \begin{lemma} \label{exceptional pairs}
		   	Let  $(E_1,E_2)$ be any exceptional pair in $\mc T$,  there is a dichotomy
		   
		   {\rm (a)}  $\hom^p(E_1,E_2)=0$ for all $p\in \ZZ$ which is equivalent to  $(E_1,E_2)\sim X$ or  $(E_1,E_2)\sim \ol{X}$ for some $X\in \{(s_i,s_j), (s_{io}, s_{jo}) :1\leq i <j \leq 3\} \cup \{ (s_{ij},s_{\alpha \beta}), 1\leq  i< j\leq 3, 1\leq  \alpha < \beta \leq 3, (i,j)\neq (\alpha,\beta) \}$, $\ol{(A,B)}=(B,A)$;

		     {\rm (b)}   $\hom^p(E_1,E_2)=1$ for some $p \in \ZZ$ which is equivalent to  $(E_1,E_2)\sim X$  for some $X$ in the following list: 
		     
		     $(s_{io}, \delta)$, $i=1,2,3$ ; 		     
		      $(\delta, s_{ij})$, $1\leq i<j\leq 3$ ; 		     
		      $( s_{ij}, s_{ko})$, $1\leq  i<j \leq 3$, $k\neq i, k \neq j$;
		      
		        $( s_{ij}, s_{123})$, $1\leq  i<j \leq 3$;		        
		         $( s_k, s_{ij})$, $1\leq  i<j \leq 3$, $k\neq i, k \neq j$;		         
		         $(  s_{123}, s_i)$, $i=1,2,3$;
		         
		         $(s_{i}, s_o)$, $(s_{io}, s_{i})$, $(s_o,s_{io})$ $i=1,2,3$;		         
		         $(s_i, s_{jo})$, $(s_{ij}, s_{i})$, $(s_{jo}, s_{ij})$, $1\leq  i,j \leq 3$, $i\neq j$.
		         \end{lemma}
\bpr  One computes that for any of the listed pairs $(A,B)$ we have  $\scal{\ul{\dim}(B),\ul{\dim}(A)}= 0$ and then  from Remark \ref{hom and hom1} it follows that all the pairs are exceptional, furthermore for any pair $(A,B)$ in (a) holds $\scal{\ul{\dim}(A),\ul{\dim}(B)}= 0$ and therefore these pairs are even orthogonal, whereas for any of the pairs in (b) holds $\scal{\ul{\dim}(A),\ul{\dim}(B)}= \pm 1$ and due to  Remark \ref{hom and hom1} we have $\hom^p(A,B)=1$ for $p=0$ or $p=1$ and $\hom^p(A,B)=0$ for all other $p\in \ZZ$. It remains to show that  any exceptional pair $(E_1,E_2)$ is up to shifts equivalent to some of the pairs in (a) or (b). Since $Rep_\KK(D_4)$ is hereditary, we can  assume that $E_1,E_2$ are exceptional representations, hence they are elements in the set \eqref{exceptional representations}. 

If $\delta \in \{E_1,E_2\}$, then due to Lemma \ref{no} $s_i\not \in  \{E_1,E_2\}$, for $i=1,2,3$, $s_o\not \in  \{E_1,E_2\}$,  $s_{123}\not \in  \{E_1,E_2\}$, therefore either  $s_{io}\in \{E_1,E_2\}$ and from the already proved part it follows that $(E_1,E_2)=(s_{io},\delta)$, or   $s_{ij}\in \{E_1,E_2\}$ and  $(E_1,E_2)=(\delta, s_{ij})$. 

So, we can assume that  $\delta \not \in \{E_1,E_2\}$. If $s_{123}\in \{E_1,E_2\}$, then from Lemma \ref{no} it follows that $s_o,s_{io} \not \in \{E_1,E_2\}$, hence either $s_i\in \{E_1,E_2\}$ and $(E_1,E_2)=(s_{123},s_i)$, or $s_{ij} \in  \{E_1,E_2\}$ and $(E_1,E_2)=(s_{ij}, s_{123})$.

We can assume now that  $\delta \not \in \{E_1,E_2\}$ and  $s_{123} \not \in \{E_1,E_2\}$. If $ s_o \in  \{E_1,E_2\}$, then in Lemma \ref{no} we see that  $s_{ij}\not \in  \{E_1,E_2\}$ and therefore  either $s_i\in \{E_1,E_2\}$ and $(E_1,E_2)=(s_{i},s_o)$, or $s_{io} \in  \{E_1,E_2\}$ and $(E_1,E_2)=(s_{o}, s_{io})$.  

So far we reduced to  $\delta \not \in \{E_1,E_2\}$,  $s_{123} \not \in \{E_1,E_2\}$,  $s_{o} \not \in \{E_1,E_2\}$. Assume that $s_{ij} \in  \{E_1,E_2\}$. If $s_{\alpha \beta}\in \{E_1,E_2\}$, then   $(\alpha,\beta) \neq (i,j)$ and the pair is the orthogonal pair $(s_{ij}, s_{\alpha\beta})$. 
 If $s_{k o}\in \{E_1,E_2\}$, then either  $k \neq i, k \neq j$ and $(E_1,E_2)= (s_{ij},s_{ko})$, or $k=i$ and  $(E_1,E_2)= (s_{io}, s_{ij})$, or $k=j$ and  $(E_1,E_2)= (s_{jo}, s_{ij})$. 
 
 So far, we reduced to  $\delta \not \in \{E_1,E_2\}$,  $s_{123} \not \in \{E_1,E_2\}$,  $s_{o} \not \in \{E_1,E_2\}$, $s_{ij}\not \in \{E_1,E_2\}$.  If $s_{io} \in  \{E_1,E_2\}$, then either $s_{jo} \in  \{E_1,E_2\}$ with $j\neq i$ and the pair is the orthogonal pair $(s_{io},s_{jo})$, or $s_i \in \{E_1,E_2\}$ and the pair is $(s_{io},s_i)$, or $s_j \in \{E_1,E_2\}$ with $j\neq i$ and $(E_1,E_2)= (s_j, s_{io})$.
 
 Thus, we can assume  $\delta \not \in \{E_1,E_2\}$,  $s_{123} \not \in \{E_1,E_2\}$,  $s_{o} \not \in \{E_1,E_2\}$, $s_{ij}\not \in \{E_1,E_2\}$, $s_{io} \not \in \{E_1,E_2\}$. And then the pair must be  orthogonal pair of the form $(s_i,s_j)$. 
\epr
\begin{coro} \label{set of -1 curves D4}
 We have $\abs{	C_{-1}^{\{\rm Id \}}(\mc T)}=9$, furthermore ($i,j, \alpha, \beta$ vary in $\{1,2,3\}$): 
 \begin{gather}\label{equation for -1 curves} 	C_{-1}^{\{\rm Id \}}(\mc T)= \{\langle s_i,s_j  \rangle :  i < j \} \cup \{\langle  s_{io},s_{jo} \rangle :  i < j   \}   \cup \{\langle s_{ij},s_{\alpha \beta} \rangle :   i< j,   \alpha < \beta , (i,j)\neq (\alpha,\beta)\}. 
 \end{gather}	
\end{coro}
\bpr  If $\mc A\subset \mc T$, $\mc A \in C_{-1}^{\{\rm Id\}}(\mc T)$, then from bijection \eqref{bijection 1}  there is an orthogonal    exc. pair $(E_1,E_2)$ with  $\langle E_1,E_2 \rangle =\mc A $  hence this exceptional pair  up to shift is some of the pairs in  Lemma \ref{exceptional pairs} (a), and  it follows that $\mc A$ must be an element in the RHS of \eqref{equation for zero curves}. Since all the pairs on the RHS are orthogonal, we see that indeed the sets are equal. From  Remark \ref{orthogonal pairs} it follows that the union is disjoint and that each of the three sets has on the RHS has three elements and the corollary follows.  
\epr
\begin{ex}[The graph $G_{D^b(A_1),\KK}(D^b(D_4))$] \label{example for D4 der points}
 The set of 12 vertices of $G_{D^b(A_1), \KK}(\mc T)$ is  described in Lemma \ref{A1 in D4}.   The set of $54$ edges is obtained from the  exceptional pairs in Lemma \ref{exceptional pairs}(we count each double-sided arrow once). In Figure \ref{F1 for D4}  we draw     full subgraphs of  $G_{D^b(A_1), \KK}(\mc T)$, the set of edges of the complete graph is disjoint union of the edges in these pictures: 
\begin{figure}	
	\begin{tikzpicture}[scale=1] 
	
	\draw[->]  (0,0)--(1,0);	
	\node [left] at (0,0) {$\langle s_{123} \rangle $};
	\node  [right] at (1,0) {$\langle s_k \rangle $};
	\draw[->]  (1.7,0)--(2.7,0);
	\node [right] at (2.7,0) {$\langle s_o \rangle $};
	\draw[<-]  (-0.2,0.2)--(0.3,0.8);
	\node  at (0.5,1) {$\langle s_{ij} \rangle $};
	\draw[->]  (1.1,0.2)--(0.7,0.8);
	\draw[<-]  (1.5,0.2)--(1.8,0.8);
	\node  at (2.1,1) {$\langle s_{ko} \rangle $};
	\draw[->]  (2.8,0.2)--(2.4,0.8);
	\draw[->]  (0.8,1)--(1.7,1);
	\draw[<-]  (1.5,1.8)--(2,1.2);
	\node [above] at (1.3,1.7) {$\langle \delta \rangle $};
	\draw[->]  (1.1,1.8)--(0.6,1.2);	
	\end{tikzpicture} \begin{tikzpicture}  [scale=0.7]
	\draw[<->]  (-1,0)--(+1,0);
	\node [left] at (-1 ,0) {$\langle s_{i}\rangle $};
	\node [right] at (+1 ,0) {$\langle s_{j}\rangle $};
	\draw[<-]   (1.5,0.5)-- (0.5,1.8);
	\node at (0 ,2) {$\langle s_{ij}\rangle $};
	\draw[->]  (-0.5,1.8) --  (-1.5,0.5);
	\node [above] at (0 ,0) {}; 	\node [above] at (0 ,0) {};
	\node [above] at (1.2 ,1) {}; 	\node [above] at (-1.2 ,1) {};
	\end{tikzpicture}		\begin{tikzpicture}  [scale=0.7]
	\draw[<->]  (-1,0)--(+1,0);
	\node [left] at (-1 ,0) {$\langle s_{io}\rangle $};
	\node [right] at (+1 ,0) {$\langle s_{jo}\rangle $};
	\draw[<-]   (1.5,0.5)-- (0.5,1.8);
	\node at (0 ,2) {$\langle s_{k}\rangle $};
	\draw[->]  (-0.5,1.8) --  (-1.5,0.5);
	\node [above] at (0 ,0) {}; 	\node [above] at (0 ,0) {};
	\node [above] at (1.2 ,1) {}; 	\node [above] at (-1.2 ,1) {};
	\end{tikzpicture}		\begin{tikzpicture} [scale=0.7]
	\draw[<->]  (-1,0)--(+1,0);
	\node [left] at (-1 ,0) {$\langle s_{ik}\rangle $};
	\node [right] at (+1 ,0) {$\langle s_{jk}\rangle $};
	\draw[<-]   (1.5,0.5)-- (0.5,1.8);
	\node at (0 ,2) {$\langle s_{k0}\rangle $};
	\draw[->]  (-0.5,1.8) --  (-1.5,0.5);
	\node [above] at (0 ,0) {}; 	\node [above] at (0 ,0) {};
	\node [above] at (1.2 ,1) {}; 	\node [above] at (-1.2 ,1) {};
	\end{tikzpicture}	\caption{Subgraphs of $G_{D^b(A_1),\KK}(D^b(D_4))$} \label{F1 for D4}
	\end{figure}
	where $i,j,k \in \{1,2,3\}$, $i\neq j$, $i\neq k$, $j\neq k$. The complete graph is obtained by gluing all these $12$ pictures in Figure \ref{F1 for D4} along the vertices they share. For the sake of  clarity we draw  three more full subgraphs in Figure \ref{F2 for D4}.  
	\begin{figure}
	\begin{tikzpicture}[scale=0.9] 
	
	\node  at (0.05,0) {$\langle s_{23} \rangle $};
	
	\draw[<->]  (-0.4,1.2)--(0.45,1.2);
	
	\draw[<->]  (-0.1,0.3)--(-0.7,1);	
	\node [above] at (-0.8,1) {$\langle s_{12} \rangle $};

	\node [right] at (0.45,1.2) {$\langle s_{13} \rangle $};
	
	\draw[<->]  (0.6,1)--(0.2,0.3);		
	\node at (0.05,3) {$\langle \delta \rangle $};	
	\draw[->]  (-0.05,2.8)--(-0.8,1.6);	
	\draw[->]  (0.15,2.8)--(0.8,1.6);
	\draw[<-]  (0.05,0.3)--(0.05,2.8);	
	
	\node at (0.05,-1.8) {$\langle s_{123} \rangle $};
	\draw[<-]  (0.05,-1.5)--(0.05,-0.2);

	\draw[->]  (-0.8,1)--(-0.1,-1.5);	
	\draw[->]  (0.7,1)--(0.2,-1.5);		
	\end{tikzpicture} 			
	\begin{tikzpicture}[scale=0.9]

	\node  at (0.05,0) {$\langle s_{1} \rangle $};
	
	\draw[<->]  (-0.4,1.2)--(0.45,1.2);
	
	\draw[<->]  (-0.1,0.3)--(-0.7,1);	
	\node [above] at (-0.8,1) {$\langle s_{2} \rangle $};

	\node [right] at (0.45,1.2) {$\langle s_{3} \rangle $};
	
	\draw[<->]  (0.6,1)--(0.2,0.3);		
	\node at (0.05,3) {$\langle s_{123} \rangle $};	
	\draw[->]  (-0.05,2.8)--(-0.8,1.6);	
	\draw[->]  (0.15,2.8)--(0.8,1.6);
	\draw[<-]  (0.05,0.3)--(0.05,2.8);	
	
	\node at (0.05,-1.8) {$\langle s_o \rangle $};
	\draw[<-]  (0.05,-1.5)--(0.05,-0.2);

	\draw[->]  (-0.8,1)--(-0.1,-1.5);	
	\draw[->]  (0.7,1)--(0.2,-1.5);		
	\end{tikzpicture} 			
	\begin{tikzpicture}[scale=0.9]

	\node  at (0.05,0) {$\langle s_{1o} \rangle $};
	
	\draw[<->]  (-0.4,1.2)--(0.45,1.2);
	
	\draw[<->]  (-0.1,0.3)--(-0.7,1);	
	\node [above] at (-0.8,1) {$\langle s_{2o} \rangle $};

	\node [right] at (0.45,1.2) {$\langle s_{3o} \rangle $};
	
	\draw[<->]  (0.6,1)--(0.2,0.3);		
	\node at (0.05,3) {$\langle s_o \rangle  $};	
	\draw[->]  (-0.05,2.8)--(-0.8,1.6);	
	\draw[->]  (0.15,2.8)--(0.8,1.6);
	\draw[<-]  (0.05,0.3)--(0.05,2.8);	
	
	\node at (0.05,-1.8) {$\langle \delta \rangle$};
	\draw[<-]  (0.05,-1.5)--(0.05,-0.2);

	\draw[->]  (-0.8,1)--(-0.1,-1.5);	
	\draw[->]  (0.7,1)--(0.2,-1.5);		
	\end{tikzpicture} 
	\caption{Subgraphs of $G_{D^b(A_1),\KK}(D^b(D_4))$} \label{F2 for D4}
	\end{figure}	
\end{ex}

\begin{coro}  \label{coro for pairs}
	 {\rm (a) } The following are full exceptional collections: 
	 
	 $(\delta,s_{12},s_{13},s_{23})$, 	$(s_{3o},\delta,s_{13},s_{23})$, $(s_{12},s_{3o},s_{13},s_{23})$  
	 	 
	 $(s_{1o},s_{2o}, s_{3o},\delta)$, 	$(s_{1o},s_{2o}, \delta,s_{12})$, $(s_{1o},s_{2o},s_{12},s_{3o})$
	 
	 in particular (using also \eqref{kappa}) we get for any $k=1,2,3$, $1\leq i < j \leq 3$, $k\neq i$, $k\neq j$: 	
	 \begin{gather} \label{coro for pairs 1} 
	 \langle  \delta ,s_{ij} \rangle = \langle s_{ko}, \delta \rangle = \langle 
	 s_{i j},s_{ko} \rangle = \langle  s_{ik},s_{jk} \rangle^{\perp} =^{\perp} \langle  s_{io},s_{jo} \rangle 
	 \end{gather}

	 {\rm (b) } The following are full exceptional collections: 
	
	$(s_1,s_2,s_3,s_o)$, 	$(s_1,s_2,s_{3o},s_3)$, $(s_1,s_2,s_{o},s_{3o})$ 
	
		$(s_{o},s_{3o},s_{1o}, s_{2o})$, 	$(s_3,s_o,s_{1o}, s_{2o})$, $(s_{3o},s_3,s_{1o}, s_{2o})$

	in particular (using also \eqref{kappa}) we get for any $k=1,2,3$, $1\leq i < j \leq 3$, $k\neq i$, $k\neq j$: 	
	\begin{gather} \label{coro for pairs 2} 
	\langle  s_k,s_o \rangle = \langle s_{ko},s_k \rangle = \langle 
	s_{o},s_{ko} \rangle = ^{\perp}\langle  s_i,s_j \rangle = \langle  s_{io},s_{jo} \rangle^{\perp} 
	\end{gather} 
 {\rm (c) } The following are full exceptional collections: 
	
	$(s_{123},s_1,s_2,s_3)$, 	$(s_{23},s_{123},s_{2},s_3)$, $(s_1,s_{23},s_{2},s_3)$
	
		$(s_{12},s_{13},s_{23},s_{123})$, 	$(s_{12},s_{13},s_{123},s_1)$, 	$(s_{12},s_{13},s_1,s_{23})$ 
		
		in particular (using also \eqref{kappa}) we get for any $k=1,2,3$, $1\leq i < j \leq 3$, $k\neq i$, $k\neq j$:	
	\begin{gather} \label{s_123.s_k} 
	\langle  s_{123},s_k \rangle = \langle s_{i j},s_{123} \rangle = \langle 
	s_{k},s_{i j} \rangle =\langle  s_i,s_j \rangle^{\perp} = ^{\perp}\langle  s_{ki},s_{kj} \rangle
	\end{gather} 
	
	{\rm (d) } For any $k=1,2,3$, $j=1,2,3$, $i=1,2,3$,  $i \neq j$, $k\neq i$, $k\neq j$ the following are full exceptional collections 	
	$(s_{i},s_{jo},s_j,s_{ko})$, 	$(s_{ij},s_{i},s_j,s_{ko})$, $(s_{jo},s_{ij},s_j,s_{ko})$
		in particular:
	\begin{gather} \label{s_i.s_jo} 
	\langle  s_{i},s_{jo} \rangle = \langle s_{i j},s_{i} \rangle = \langle 
	s_{jo},s_{i j} \rangle = 	\langle  s_{j},s_{ko} \rangle^\perp = ^\perp\langle  s_{k},s_{io} \rangle. 	\end{gather} 
\end{coro}

\begin{coro} \label{set of 0 curves D4} We have $\abs{	C_0^{\{\rm Id \}}(\mc T)}=15$, furthermore  : 
	\begin{gather}  \label{equation for zero curves} 	C_0^{\{\rm Id \}}(\mc T)= \{\langle  s_{ko},\delta \rangle,\langle  s_k,s_o \rangle, \langle  s_{123},s_k \rangle : k=1,2,3 \}   \cup \{\langle s_{i},s_{jo} \rangle : 1\leq i , j \leq 3, i\neq j \}. 
	\end{gather}
\end{coro}
\bpr If $\mc A\subset \mc T$, $\mc A \in C_{0}^{\{\rm Id\}}(\mc T)$, then from bijection \eqref{bijection 1}  there is a strong    exc. pair $(E_1,E_2)$ with  $\langle E_1,E_2 \rangle =\mc A $ and $\hom(E_1,E_2)=1$ hence this exceptional pair  up to shift is some of the pairs in  Lemma \ref{exceptional pairs} (b), and then due to  Corollary \ref{coro for pairs} it follows that $\mc A$ must be an element in the RHS of \eqref{equation for zero curves}. Let us now take an element $\mc A$ on the RHS,  by Lemma \ref{exceptional pairs} (b) it is generated by a strong exceptional pair $(E_1,E_2)$ with  $\hom(E_1,E_2)=1$ and using again \eqref{bijection 1} we see that $\mc A \in  C_{0}^{\{\rm Id\}}(\mc T_1)$. 

If $\mc A=\langle  s_{io},\delta \rangle = \langle  s_{jo},\delta \rangle$ , then the orthogonal complement of $\langle \delta \rangle^{\perp}$ in $\mc A$ is $\langle s_{io} \rangle = \langle s_{jo} \rangle$, which implies $i=j$. Similarly one sees that each of  the first tree sets in the union    
 above has 3 elements. If $\langle s_{i},s_{jo} \rangle =\langle s_{\alpha},s_{\beta o} \rangle $, $i\neq j$, $\alpha \neq \beta$ then from Corollary \ref{coro for pairs} (d) we deduce that $(s_{\alpha},s_{\beta o}, s_j,s_{ko})$ is an exceptional collection, where $k\neq i$, $k\neq j$, therefore, since  $s_\alpha, s_{ko}$  is an exceptional pair and from Lemma \ref{exceptional pairs} we know that $\hom^*(s_{ko}, s_k)\neq 0$, we  deduce that $\alpha \neq k$, i. e. $\alpha \in \{i,j\}$. The condition that $(s_{\beta o}, s_j)$ is an exceptional pair together with Lemma \ref{exceptional pairs} imply that $\beta =j$. Now by similar arguments as in the case $\langle  s_{io},\delta \rangle = \langle  s_{jo},\delta \rangle $ it follows that $\alpha = i $. So we proved that the elements in the fourth set of the union above  are $6$. 
 
 From Corollary \ref{coro for pairs} we know that the orthogonal complements of  $ \langle  s_{i},s_{jo} \rangle $ are generated by non-orthogonal pairs, whereas the orthogonal complements of any element in the first three sets of the union above are generated by orthogonal pairs, therefore the fourth set in the union is disjoint with any of the other three sets. 
 
 Finally,  comparing, say the right, orthogonal subcategories  of each of the 9 subcategories in the union of the first three sets  (these orthogonal subcategories are shown in Corollary \ref{coro for pairs}) and using Remark \ref{orthogonal pairs}  we deduce that these  subcategories are pairwise  different  and therefore the subcategories themselves are different. The Corollary is completely proved.   \epr

Next we determine the action of the Serre functor: 

\begin{lemma} \label{Serre of D4} Let $S$ be the Serre functor of $\mc T$. Then we have:
$
	S(\delta) \sim s_o$,  $S(s_o) \sim s_{123}$,  $S(s_{123}) \sim \delta. 
	$
	For any $i,j,k\in \{1,2,3\}$ which are pairwise different we have:
$
S(s_i) \sim s_{jk}$,   $S(s_{jk}) \sim s_{io}$,  $S(s_{io}) \sim s_i. 
$	
\end{lemma}
\bpr We use \eqref{Serre} and appropriate full exceptional   collections. In Corollary \ref{coro for pairs} (a), (c) are shown exceptional collections $(\delta,s_{12},s_{13},s_{23})$, 	$(s_{12},s_{13},s_{23},s_{123})$ and we deduce $S(s_{123})\sim \delta$.  In Corollary \ref{coro for pairs} (b) we have also an  exceptional collection 	$(s_{o},s_{3o},s_{1o}, s_{2o})$, and from Lemma \ref{exceptional pairs} one sees that   $(s_{3o},s_{1o}, s_{2o},\delta)$ ais also full exceptional collection, hence $S(\delta) \sim s_o$. In Corollary \ref{coro for pairs} (b) is shown a full exceptional collection 	$(s_1,s_2,s_3,s_o)$ and  from Lemma \ref{exceptional pairs} we get an exceptional collection $(s_{123}, s_1, s_2, s_3)$ therefore $S(s_0) \sim s_{123}$. 

From Lemma \ref{exceptional pairs} we see that $(s_{23},s_{123}, s_2, s_3)$ is an exceptional collection, and we already have the exceptional collection $(s_{123},s_{2}, s_3, s_1)$, therefore $S(s_1) \sim s_{23}$. From Lemma \ref{exceptional pairs} we see that $(s_{1o},\delta, s_{12}, s_{13})$ is an exceptional collection, and we already have the exceptional collection $(\delta, s_{12},s_{13},s_{23})$, therefore $S(s_{23}) \sim s_{1o}$.  From Lemma \ref{exceptional pairs} we see that $(s_{1}, s_3, s_{2o}, s_{2})$ is an exceptional collection, and we already have used the exceptional collection $(s_3,s_{2o}, s_2, s_{1o})$ in Corollary \ref{coro for pairs} (d), therefore $S(s_{1o}) \sim s_{1}$. Finally, the last statement   follows from the already derived incidences and since $S \circ \kappa \cong \kappa \circ S$ due to \cite[Lemma 3.1]{Miyachi} (recall also \eqref{kappa}). 
\epr
We already know what are the sets $C_{D^b(A_1),\KK}^{\{\rm Id\}}(\mc T)$, $C_0^{\{\rm Id\}}(\mc T)$,  $C_{-1}^{\{\rm Id\}}(\mc T)$. To determine the rest numbers in   tables \eqref{table with numbers D4},  \eqref{table with numbers A1 in D4}   we use again Corollary \ref{group action}. Thus, we need to count the orbits of the action of $\Gamma$ on $C_{D^b(A_1),\KK}^{\{\rm Id\}}(\mc T)$ and on  $C_l^{\{\rm Id\}}(\mc T)$ ( $l=-1,0$)   for the three choices of $\Gamma$. Pictures \ref{Figure for A1 in D4} \ref{Figure for C0}, \ref{Figure for C-1} depict the actions of $[\kappa]$, $[S]$ on $C_{D^b(A_1),\KK}^{\{\rm Id\}}(\mc T)$,  $C_l^{\{\rm Id\}}(\mc T)$ for $l=0,-1$, respectively:  the action of  $[\kappa]$ is depicted via the dashed arrows, whereas the action of  $[S]$  - via non-dashed arrows.

\begin{figure} \centering
	\hspace{-50mm}   \begin{subfigure}[b]{0.2\textwidth}
		\begin{tikzpicture}[scale=0.6] 
		\node   at (0,6.5) { $\langle s_{1}\rangle$\normalsize};
		\draw [|->] (0,6.2)--(0,5.2);
		\node [right] at (0,6) {$S$};
		\node   at (0,5) { $\langle s_{23} \rangle$\normalsize};
		\draw [|->] (0,4.7)--(0,3.8);
		\node [right] at (0,4.5) {$S$};
		\node   at (0,3.5) { $\langle s_{1o} \rangle$\normalsize};
		\draw[|->] (-0.8,3.5) to [out=180,in=195] (-0.6,6.5);
		\node [right] at (-1.5 ,5) {$S$};
		\node   at (0,1.5) { $\langle\delta \rangle$\normalsize};

		\node   at (4,6.5) { $\langle s_{2}\rangle$\normalsize};
		\draw [|->] (4,6.2)--(4,5.2);
		\node [right] at (4,6) {$S$};
		\node   at (4,5) { $\langle s_{13}\rangle$\normalsize};
		\draw [|->] (4,4.7)--(4,3.8);
		\node [right] at (4,4.5) {$S$};
		\node   at (4,3.5) { $\langle s_{2o} \rangle$\normalsize};
		\draw[|->] (3.2,3.8) to [out=135,in=195] (3.4,6.4);
		\node [right] at (2.7 ,4.5) {$S$};
		\node   at (4,1.5) { $\langle s_{o} \rangle$\normalsize};
		
		\node   at (8,6.5) { $\langle s_{3} \rangle$\normalsize};
		\draw [|->] (8,6.2)--(8,5.2);
		\node [right] at (8,6) {$S$};
		\node   at (8,5) { $\langle s_{12} \rangle$\normalsize};
		\draw [|->] (8,4.7)--(8,3.8);
		\node [right] at (8,4.5) {$S$};
		\node   at (8,3.5) { $\langle s_{3o} \rangle$\normalsize};
		\draw[|->] (8.8,3.5) to [out=0,in=315] (8.6,6.5);
		\node [right] at (9.5 ,5) {$S$};
		\node   at (8,1.5) { $\langle s_{123} \rangle$\normalsize};

		\draw [dashed, |->] (0.6,6.5) -- (3.2,6.5);
		\node   at (2,6.2) { $\kappa$\normalsize};
		\draw [dashed, |->] (0.7,5) -- (3.2,5);
		\node   at (2,4.8) { $\kappa$\normalsize};
		\draw [dashed, |->] (0.8,3.5) -- (3.2,3.5);
		\node   at (2,4.8) { $\kappa$\normalsize}; 
		\node   at (2,3.3) { $\kappa$\normalsize};
		
		\draw [dashed, |->] (4.6,6.5) -- (7.2,6.5);
		\node   at (6,6.2) { $\kappa$\normalsize};
		\draw [dashed, |->] (4.6,5) -- (7.2,5);
		\node   at (6,4.8) { $\kappa$\normalsize};
		\draw [dashed, |->] (4.8,3.5) -- (7.2,3.5);
		\node   at (6,4.8) { $\kappa$\normalsize}; 
		\node   at (6,3.3) { $\kappa$\normalsize};
		
		\draw[dashed, |->] (7.4,6.7) to [out=160,in=20] (0.6,6.7);
		\node   at (4,7.2) { $\kappa$\normalsize};
		
		\draw[dashed, |->] (7.3,5.2) to [out=170,in=10] (0.6,5.2);
		\node   at (6,5.5) { $\kappa$\normalsize};
		
		\draw[dashed, |->] (7.2, 3.2) to [out=190,in=350] (0.6,3.2);
		\node   at (4,2.7) { $\kappa$\normalsize};
		
		\draw [|->] (0.7,1.5)-- (3.3,1.5);
		\node   at (2,1.7) { $S$\normalsize};

		\draw [|->] (4.7,1.5)-- (7.3,1.5);
		\node   at (6,1.7) { $S$\normalsize};

		\draw[ |->] (7.4,1.7) to [out=160,in=20] (0.6,1.7);
		\node   at (4,2.2) { $S$\normalsize};
		
		\draw[dashed,  |->]  (-0.2, 1.2) arc [radius=0.5, start angle=110, end angle= 430];
		\node   at (0,0.7) { $\kappa$\normalsize};
		
		\draw[dashed,  |->]  (3.8, 1.2) arc [radius=0.5, start angle=110, end angle= 430];
		\node   at (4,0.7) { $\kappa$\normalsize};
		
		\draw[dashed,  |->]  (7.8, 1.2) arc [radius=0.5, start angle=110, end angle= 430];
		\node   at (8,0.7) { $\kappa$\normalsize};
		
		\end{tikzpicture}
		\caption{on $C_{D^b(A_1),\KK}^{\{\rm Id\}}(\mc T)$}\label{Figure for A1 in D4}
	\end{subfigure}
	\hspace{40mm}
	\begin{subfigure}[b]{0.2\textwidth}
	 \begin{tikzpicture}[scale=0.6]  
	 \node   at (0,6.5) { $\langle s_{1o}, \delta \rangle$\normalsize};
	 \draw [|->] (0,6.2)--(0,5.2);
	 \node [right] at (0,6) {$S$};
	 \node   at (0,5) { $\langle s_{1}, s_o \rangle$\normalsize};
	 \draw [|->] (0,4.7)--(0,3.8);
	 \node [right] at (0,4.5) {$S$};
	 \node   at (0,3.5) { $\langle s_{123}, s_1 \rangle$\normalsize};
	 \draw[|->] (-1.3,3.5) to [out=180,in=195] (-1,6.5);
	 \node [left] at (-1.2 ,5) {$S$};
	 \node   at (0,1.5) { $\langle s_{1}, s_{2o} \rangle$\normalsize};
	 \node   at (0,0) { $\langle s_{1}, s_{3o} \rangle$\normalsize};

	 \node   at (4,6.5) { $\langle s_{2o}, \delta \rangle$\normalsize};
	 \draw [|->] (4,6.2)--(4,5.2);
	 \node [right] at (4,6) {$S$};
	 \node   at (4,5) { $\langle s_{2}, s_o \rangle$\normalsize};
	 \draw [|->] (4,4.7)--(4,3.8);
	 \node [right] at (4,4.5) {$S$};
	 \node   at (4,3.5) { $\langle s_{123}, s_2 \rangle$\normalsize};
	 \draw[|->] (3.2,3.8) to [out=135,in=195] (3.1,6.1);
	 \node [right] at (2.7 ,4.5) {$S$};
	 \node   at (4,1.5) { $\langle s_{2}, s_{3o} \rangle$\normalsize};
	 \node   at (4,0) { $\langle s_{3}, s_{2o} \rangle$\normalsize};
	 
	 \node   at (8,6.5) { $\langle s_{3o}, \delta \rangle$\normalsize};
	 \draw [|->] (8,6.2)--(8,5.2);
	 \node [right] at (8,6) {$S$};
	 \node   at (8,5) { $\langle s_{3}, s_o \rangle$\normalsize};
	 \draw [|->] (8,4.7)--(8,3.8);
	 \node [right] at (8,4.5) {$S$};
	 \node   at (8,3.5) { $\langle s_{123}, s_3 \rangle$\normalsize};
	 \draw[|->] (9.2,3.5) to [out=0,in=315] (9,6.5);
	 \node [left] at (9.7 ,5) {$S$};
	 \node   at (8,1.5) { $\langle s_{3}, s_{1o} \rangle$\normalsize};
	 \node   at (8,0) { $\langle s_{2}, s_{1o} \rangle$\normalsize};

	 \draw [dashed, |->] (0.9,6.5) -- (3,6.5);
	 \node   at (2,6.2) { $\kappa$\normalsize};
	 \draw [dashed, |->] (0.9,5) -- (3,5);
	 \node   at (2,4.8) { $\kappa$\normalsize};
	 \draw [dashed, |->] (1.2,3.5) -- (2.8,3.5);
	 \node   at (2,4.8) { $\kappa$\normalsize}; 
	 \node   at (2,3.3) { $\kappa$\normalsize};
	 
	 \draw [dashed, |->] (4.9,6.5) -- (7,6.5);
	 \node   at (6,6.2) { $\kappa$\normalsize};
	 \draw [dashed, |->] (5,5) -- (7,5);
	 \node   at (6,4.8) { $\kappa$\normalsize};
	 \draw [dashed, |->] (5.2,3.5) -- (6.8,3.5);
	 \node   at (6,4.8) { $\kappa$\normalsize}; 
	 \node   at (6,3.3) { $\kappa$\normalsize};
	 
	 \draw[dashed, |->] (7.1,7) to [out=160,in=20] (0.8,7);
	 \node   at (4,7.2) { $\kappa$\normalsize};
	 
	 \draw[dashed, |->] (7.1,5.5) to [out=170,in=10] (0.8,5.5);
	 \node   at (6,5.5) { $\kappa$\normalsize};
	 
	 \draw[dashed, |->] (7.2, 3.2) to [out=190,in=350] (0.6,3.2);
	 \node   at (4,2.7) { $\kappa$\normalsize};
	 
	 \draw [|->] (1.1,1.6)to [out=10,in=170] (2.9,1.6);
	 \node   at (2,1.9) { $S$\normalsize};
	 
	 \draw [dashed, |->] (1.1,1.4)to [out=350,in=190] (2.9,1.4);
	 \node   at (2,1.4) { $\kappa$\normalsize};
	 
	 \draw [|->] (1.1,0.1)to [out=10,in=170] (2.9,0.1);
	 \node   at (2,0.5) { $S$\normalsize};
	 
	 \draw [dashed, |->] (2.9,-0.1)  to [out=190,in=350] (1.1,-0.1) ;
	 \node   at (2,-0.1) { $\kappa$\normalsize};
	 
	 \draw [|->] (5.1,1.6)to [out=10,in=170] (6.9,1.6);
	 \node   at (6,1.9) { $S$\normalsize};
	 
	 \draw [dashed, |->] (5.1,1.4)to [out=350,in=190] (6.9,1.4);
	 \node   at (6,1.4) { $\kappa$\normalsize};
	 
	 \draw [|->] (5.1,0.1)to [out=10,in=170] (6.9,0.1);
	 \node   at (6,0.5) { $S$\normalsize};
	 
	 \draw [dashed, |->] (6.8,-0.1)  to [out=190,in=350] (5.1,-0.1) ;
	 \node   at (6,-0.1) { $\kappa$\normalsize};
	 
	 \draw[ |->] (7,2) to [out=165,in=15] (1,2);
	 \node   at (4,2.2) { $S$\normalsize};
	 
	 \draw[dashed, |->] (7, 1.1) to [out=190,in=350] (1,1.1);
	 \node   at (4,0.7) { $\kappa$\normalsize};
	 
	 \draw[dashed, |->] (1,-0.4)  to [out=350,in=190]  (7.2, -0.3);
	 \node   at (4,-0.5) { $\kappa$\normalsize};
	 
	 \draw[ |->]  (7.3, -0.6)  to [out=200 ,in=340]  (0.5,-0.6);
	 \node   at (4,-1.5) { $S$\normalsize};
	 
	 \end{tikzpicture}
		\caption{ on $C_0^{\{\rm Id\}}(\mc T)$}\label{Figure for C0}
	\end{subfigure} 
	\caption{Actions}
\end{figure}

 Indeed, the formulas \eqref{kappa}, together with $\kappa( s_{3})\sim s_{1}$, $\kappa( s_{3o})\sim s_{1o}$, and  since $\kappa$ fixes $\langle \delta \rangle $, $\langle s_{123}\rangle $, $\langle s_o \rangle$ amount to action of   $\kappa$ as shown by the dashed arrows. 
 
 For drawing the arrows for $S$ we use Lemma \ref{Serre of D4} and Corollary \ref{coro for pairs}. For example to draw the arrows for $S$ in  Figures  \ref{Figure for A1 in D4} and \ref{Figure for C-1} it is enough to use Lemma \ref{Serre of D4}. In Figure \ref{Figure for C0} the assignment   $S(\langle s_{1o}, \delta \rangle)=\langle s_{1}, s_o \rangle$ is clear from  Lemma \ref{Serre of D4} as well, whereas   $S(\langle s_{1}, s_o \rangle)=\langle s_{23}, s_{123} \rangle$ by Lemma \ref{Serre of D4} and $\langle s_{23}, s_{123} \rangle = \langle  s_{123}, s_1 \rangle$ by Corollary \ref{coro for pairs} (c). Similarly  $S(\langle s_{123}, s_1 \rangle)=\langle \delta, s_{23} \rangle$ by Lemma \ref{Serre of D4} and $\langle \delta, s_{23} \rangle = \langle  s_{1o}, \delta \rangle$ by Corollary \ref{coro for pairs} (a). Similarly  $S(\langle s_{1}, s_{2o} \rangle)=\langle s_{23}, s_{2} \rangle$,  $S(\langle s_{1}, s_{3o} \rangle)=\langle s_{23}, s_{3} \rangle$ by Lemma \ref{Serre of D4} and $\langle s_{23}, s_{2} \rangle = \langle  s_{2}, s_{3 o} \rangle$,  $\langle s_{23}, s_{3} \rangle = \langle  s_{3}, s_{2 o} \rangle$ by Corollary \ref{coro for pairs} (d).

  Since $S$ and $\kappa$ commute we obtain  picture \ref{Figure for C0} and recalling \eqref{Auto for D4}  we  complete  tables \eqref{table with numbers D4}, \eqref{table with numbers A1 in D4}.

\begin{figure}
	
	 \begin{tikzpicture}[scale=0.6]  
	 \node   at (0,6.5) { $\langle s_{12}, s_{23} \rangle$\normalsize};
	 \draw [|->] (0,6.2)--(0,5.2);
	 \node [right] at (0,6) {$S$};
	 \node   at (0,5) { $\langle s_{1o}, s_{3o} \rangle$\normalsize};
	 \draw [|->] (0,4.7)--(0,3.8);
	 \node [right] at (0,4.5) {$S$};
	 \node   at (0,3.5) { $\langle s_{1}, s_3 \rangle$\normalsize};
	 \draw[|->] (-1.3,3.5) to [out=180,in=195] (-1.2,6.5);
	 \node [left] at (-1.2 ,5) {$S$};

	 \node   at (4,6.5) { $\langle s_{23}, s_{13} \rangle$\normalsize};
	 \draw [|->] (4,6.2)--(4,5.2);
	 \node [right] at (4,6) {$S$};
	 \node   at (4,5) { $\langle s_{2o}, s_{1o} \rangle$\normalsize};
	 \draw [|->] (4,4.7)--(4,3.8);
	 \node [right] at (4,4.5) {$S$};
	 \node   at (4,3.5) { $\langle s_{1}, s_2 \rangle$\normalsize};
	 \draw[|->] (3.2,3.8) to [out=135,in=195] (3.1,6.1);
	 \node [right] at (2.7 ,4.5) {$S$};

	 \node   at (8,6.5) { $\langle s_{12}, s_{13} \rangle$\normalsize};
	 \draw [|->] (8,6.2)--(8,5.2);
	 \node [right] at (8,6) {$S$};
	 \node   at (8,5) { $\langle s_{2o}, s_{3o} \rangle$\normalsize};
	 \draw [|->] (8,4.7)--(8,3.8);
	 \node [right] at (8,4.5) {$S$};
	 \node   at (8,3.5) { $\langle s_{2}, s_3 \rangle$\normalsize};
	 \draw[|->] (9.2,3.5) to [out=0,in=315] (9.2,6.5);
	 \node [left] at (9.7 ,5) {$S$};

	 \draw [dashed, |->] (1.2,6.5) -- (2.8,6.5);
	 \node   at (2,6.2) { $\kappa$\normalsize};
	 \draw [dashed, |->] (1.2,5) -- (2.8,5);
	 \node   at (2,4.8) { $\kappa$\normalsize};
	 \draw [dashed, |->] (1.2,3.5) -- (2.8,3.5);
	 \node   at (2,4.8) { $\kappa$\normalsize}; 
	 \node   at (2,3.3) { $\kappa$\normalsize};
	 
	 \draw [dashed, |->] (5.4,6.5) -- (6.7,6.5);
	 \node   at (6,6.2) { $\kappa$\normalsize};
	 \draw [dashed, |->] (5.4,5) -- (6.7,5);
	 \node   at (6,4.8) { $\kappa$\normalsize};
	 \draw [dashed, |->] (5.2,3.5) -- (6.8,3.5);
	 \node   at (6,4.8) { $\kappa$\normalsize}; 
	 \node   at (6,3.3) { $\kappa$\normalsize};
	 
	 \draw[dashed, |->] (7.1,7) to [out=160,in=20] (0.8,7);
	 \node   at (4,7.2) { $\kappa$\normalsize};
	 
	 \draw[dashed, |->] (7.1,5.5) to [out=170,in=10] (0.8,5.5);
	 \node   at (6,5.5) { $\kappa$\normalsize};
	 
	 \draw[dashed, |->] (7.2, 3.2) to [out=190,in=350] (0.6,3.2);
	 \node   at (4,2.7) { $\kappa$\normalsize};

	 \end{tikzpicture}
	\caption{Actions on $C_{-1}^{\{\rm Id\}}(\mc T)$}\label{Figure for C-1}
\end{figure}

For $D^b(D_4)$ we already know these sets for $\mc A$ generated by an exceptional object or by an exceptional pair. Now we determine $C_{\mc A, \KK}^{\{\rm Id\}}(D^b(D_4))$ for $\mc A$ generated by an exceptional triple:
\begin{prop} \label{triples in D4 prop}  If $C_{\mc A, \KK}^{\{\rm Id\}}(D^b(D_4))\neq \emptyset$ and $\mc A$ has a full exceptional triple, then either $\mc A \cong D^b(A_3)$ or  $\mc A \cong D^b(A_1)\oplus  D^b(A_1)\oplus  D^b(A_1)$.  Table \eqref{triples in D4} holds and: 
	\begin{gather} 
	\label{expl for A3 in D4} C^{\{\rm Id\}}_{D^b(A_3), \KK}(D^b(D_4))=\left \{\langle s_i \rangle^{\perp}: i=1,2,3\right \} \cup \left \{\langle s_{io} \rangle^{\perp} : i=1,2,3\right \}\cup\{\langle s_{12} \rangle^{\perp}, \langle s_{13} \rangle^{\perp}, \langle s_{23} \rangle^{\perp} \}  \\
	\label{orthogonal triples in D4}
	C^{\{\rm Id\}}_{D^b(A_1)\oplus D^b(A_1)\oplus D^b(A_1), \KK}(D^b(D_4))=\left \{\langle \delta \rangle^{\perp}, \langle s_o \rangle^{\perp}, \langle s_{123} \rangle^{\perp}\right \}. 
	\end{gather} 
	
\end{prop}
\bpr Let $\mc A$ be generated by an exceptional triple.  If  $\mc B \in C^{\{\rm Id\}}_{\mc A, \KK}(D^b(D_4))$, then  it must be generated by an exceptional triple, say $(A,B,C)$. Due to \cite{WCB1} this triple can be extended to a full exceptional collection   $(A,B,C,D)$ and hence $\mc B= \langle D \rangle^{\perp}$ for some $\langle D \rangle \in C_{D^b(A_1), \KK}(D^b(D_4))$. The set $ C_{D^b(A_1), \KK}(D^b(D_4))$ is determined  in Lemma \ref{A1 in D4} and we have $12$ options for $\langle D \rangle^{\perp}$.

From the exceptional collections  $(s_{1o},s_{2o}, s_{3o},\delta)$, $(s_1,s_2,s_3,s_o)$ and  $(s_{12},s_{13},s_{13},s_{123})$ (see Figures \ref{F1 for D4}, \ref{F2 for D4})  we conclude that the set on RHS in \eqref{orthogonal triples in D4} is a subset of  $ C^{\{\rm Id\}}_{D^b(A_1)\oplus D^b(A_1)\oplus D^b(A_1), \KK}(D^b(D_4))$.

If we show that the RHS in \eqref{expl for A3 in D4} is a subset of  $C^{\{\rm Id\}}_{D^b(A_3), \KK}(D^b(D_4))$, then \eqref{expl for A3 in D4} and  \eqref{orthogonal triples in D4} follow.

Using the actions of $S$ and $\kappa$ shown in Figure \ref{Figure for A1 in D4} we deduce that the subcategories  $\left \{\langle s_i \rangle^{\perp} : i=1,2,3\right \} \cup \left \{\langle s_{io} \rangle^{\perp} : i=1,2,3\right \}\cup\{\langle s_{12} \rangle^{\perp}, \langle s_{13} \rangle^{\perp}, \langle s_{23} \rangle^{\perp} \}$ are pairwise equivalent, therefore it is enogh to  show that  $\langle s_3 \rangle^{\perp}$ is equivalent to $D^b(A_3)$. In Figures \ref{F1 for D4}, \ref{F2 for D4} (or Corollary  \ref{coro for pairs} (c)) we see that $(s_1,s_{23},s_2,s_3)$ is a full exceptional collection, therefore $\langle s_3 \rangle^{\perp} =\langle  s_1,s_{23},s_2 \rangle$. From Figures \ref{F1 for D4}, \ref{F2 for D4} (or using \eqref{equation for -1 curves}, \eqref{s_123.s_k}, \eqref{s_i.s_jo}, \eqref{equation for zero curves}) we deduce that $\langle  s_1,s_{23},s_2 \rangle = \langle A,B,C \rangle$ for some strong exceptional triple with   $\hom(A,B)=\hom(B,C)=1$, $\hom(A,C)=0$, and now  \cite[Corollary 1.9]{Orlov} ensures that $ \langle s_3 \rangle^{\perp} \cong D^b(A_{3})$, hence  \eqref{expl for A3 in D4},  \eqref{orthogonal triples in D4} are proved. The rest follows  from   the actions of $S$ and $\kappa$ describged in  Figure \ref{Figure for A1 in D4}. 
\epr

\section{Examples for intersection relations} \label{subsection for examples of intersections}
Let again $\mc J = \{\mbox{trivial category},  D^b(A_1), D^b(A_2)\}$. 
For two choices of  $\mc T$ we will define  injective  maps
 \begin{gather}  \label{examples for interesection intro}  \bd C_0^{\{\rm Id\}}(\mc T) & \rTo^{\alpha} & \{\mbox{\textit{lines in $\AA^2$}}\} \ed \qquad \qquad  \bd C_{D^b(A_1), \KK}^{\{\rm Id\}}(\mc T) & \rTo^{\beta } &\{\mbox{\textit{points in $\AA^2$}}\}   \ed \end{gather} such that the corresponding 
map $\bd C_{\mc J, \KK}(\mc T) & \rTo^{\gamma} &  \{\mbox{\textit{points and straight lines in  $\AA^2$}}\} \cup \{\emptyset \} \ed $ \footnote{with $\gamma(\mbox{trivial category})=\emptyset $} is an embedding of partially ordered sets, i.e. for any  $x,y \in  C_{\mc J, \KK}(\mc T)$ we have $x \leq y$ iff $\gamma(x) \subset \gamma(y)$.
 Using Lemmas \ref{derived points in a genus 0}, \ref{intersection} one easily shows that such a $\gamma$ has the property that  two different $\mc A, \mc B \in  C_0^{\{\rm Id\}}(\mc T)$ intersect  non-trivially iff $\alpha(\mc A)\cap \alpha (\mc B) \cap  {\rm Im}(\beta) \neq \emptyset $ and   if they intersect  non-trivially $\beta (\mc A \cap \mc B ) = \alpha (\mc A ) \cap \alpha (\mc B )\cap {\rm Im}(\beta)$.

Let us first consider the case  $\mc T = D^b(A_3)$. The derived points, and the genus $0$  nc curves are shown in \eqref{Db A1 in DbA3}, \eqref{exceptional pairs in A_3}. Using Corollary \ref{interesection genus zero} we define   $\alpha$, $\beta$ from \eqref{examples for interesection intro} as shown in Figure  \eqref{Figure for  intersection in A3}.

Till the end of the section we consider the case $\mc T=D^b(D_4)$. The $15$ nc  curves of genus $0$ are described in Corollary \ref{set of 0 curves D4} and the $12$  derived points are described in Lemma \ref{A1 in D4}. To construct  $\alpha, \beta$ from \eqref{examples for interesection intro}  we will use also equalities  \eqref{coro for pairs 1}, \eqref{coro for pairs 2}, \eqref{s_123.s_k}, \eqref{s_i.s_jo}, where one sees what are the three derived points  (Lemma \ref{derived points in a genus 0}) contained in each of the genus $0$ nc curves.  

The curves  $\{\langle s_{i},s_{jo} \rangle : 1\leq i , j \leq 3, i\neq j \}$  can be split  in two three element sets, the curves in each of these sets are pairwise non-intersecting: $\{ \langle s_{1},s_{2o} \rangle, \langle s_{2},s_{3o} \rangle,  \langle s_{3},s_{1o} \rangle \} $,   $\{ \langle s_{1},s_{3o} \rangle, \langle s_{3},s_{2o} \rangle,  \langle s_{2},s_{1o} \rangle \} $, whereas each curve in one of the sets intersects  each of the three curves in the other set in a single derived point, which is seen  from Corollary \ref{coro for pairs} (d)  and Corollary \ref{interesection genus zero}.  Therefore we can present the intersection relations  in   $\{\langle s_{i},s_{jo} \rangle : 1\leq i , j \leq 3, i\neq j \}$ via two triples of parallel lines perpendicular to each other, for example as in Figure \ref{Figure for D4 1}.
\begin{figure}
	\begin{tikzpicture}[scale=0.3] 
	\draw  (-6,2)--(+6,2);
	\draw [fill] (0,2) circle [radius=0.07];
	\node [above right] at (0 ,2) {$\langle s_{23}\rangle $};
	\draw  (-6,4)--(+6,4);
	\draw [fill] (0,4) circle [radius=0.07];
	\node [above right] at (0 ,4) {$\langle s_{2o}\rangle $};		
	\draw  (-6,-2)--(+6,-2);
	\draw [fill] (0,-2) circle [radius=0.07];	
	\node [above left] at (-0.1 ,-2) {$\langle s_{3}\rangle $};
	
	\draw  (0,-3)--(0,+5);
	\node [below] at (0,-3) {$\langle s_{3}, s_{2o} \rangle $};
	
	\draw  (-4,-3)--(-4,+5.2);
	\draw [fill] (-4,-2) circle [radius=0.07];
	\node [below] at (-4,-3) {$\langle s_{1}, s_{3o} \rangle $};
	\node [left] at (-6,-2) {$\langle s_{3}, s_{1o} \rangle $};
	\node [above left] at (-4 ,-2) {$\langle s_{13}\rangle $};
	\draw [fill] (-4,2) circle [radius=0.07];
	\node [left] at (-6,2) {$\langle s_{2}, s_{3o} \rangle $};
	\node [below left] at (-4 ,2) {$\langle s_{3o}\rangle $};
	\draw [fill] (-4,4) circle [radius=0.07];
	
	\node [above right] at (-4 ,4) {$\langle s_{1}\rangle $};
	\node [left] at (-6,4) {$\langle s_{1}, s_{2o} \rangle $};
	
	\draw  (4,-3)--(4,+5);
	\node [below] at (4,-3) {$\langle s_{2}, s_{1o} \rangle $};
	\draw [fill] (4,-2) circle [radius=0.07];
	\node [above right] at (4 ,-2) {$\langle s_{1o}\rangle $};
	\draw [fill] (4,2) circle [radius=0.07];
	\node [below right] at (4 ,2) {$\langle s_{2}\rangle $};
	\draw [fill] (4,4) circle [radius=0.07];
	\node [above right] at (4 ,4) {$\langle s_{12}\rangle $};			
	\end{tikzpicture}
	\caption{Constructing  $\gamma : C_{\mc J,\KK}(D^b(D_4)) \rightarrow \AA^2 $, Step 1}
\label{Figure for D4 1}
\end{figure}
We have embedded already $9$ of the derived points in $C_{D^b(A_1),\KK}^{\{\rm Id\}}(\mc T)$ and $6$ of the genus $0$ nc curves, missing are the three derived points  $\langle \delta \rangle, \langle s_o \rangle, \langle s_{123} \rangle$ and the $9$ nc curves $\{\langle  s_{ko},\delta \rangle : k=1,2,3 \} \cup \{\langle  s_k,s_o \rangle : k=1,2,3 \}   \cup \{\langle  s_{123},s_k \rangle : k=1,2,3 \}$ (see Corollary \ref{set of 0 curves D4} ). Now the curve $\langle s_{3o}, \delta \rangle $ intersects $\langle  s_2, s_{3o}\rangle  $ and   $\langle  s_1, s_{3o}\rangle  $ into $\langle s_{3o} \rangle$, $\langle  s_1, s_{2o}\rangle  $ and   $\langle  s_2, s_{1o}\rangle  $ into $\langle s_{12} \rangle$ and does not intersect the rest of the already drawn curves (we use \eqref{coro for pairs 1}  and Corollary \ref{interesection genus zero}), therefore if we want to continue with the picture \ref{Figure for D4 1}, the nc curve   $\langle s_{3o}, \delta \rangle $  must be mapped  by $\alpha$ from \eqref{examples for interesection intro} to the line spanned by the already embedded points $\langle s_{3o} \rangle$, $\langle s_{12} \rangle$, by similar arguments  $\alpha$ must map the nc curve  $\langle s_{2o}, \delta\rangle $ to the line spanned by  the already embedded points $\langle s_{2o} \rangle$, $\langle s_{13} \rangle$, thus we obtain the two blue lines in the first picture of Figure \ref{Figure for D4 2}.
\begin{figure} 
	\begin{tikzpicture}[scale = 0.5]
	\draw  (-6,3)--(+6,3);
	\draw [fill] (0,3) circle [radius=0.1];
	\node [below right] at (0 ,3) {$\langle s_{23}\rangle $};
	\draw  (-6,4)--(+6,4);
	\draw [fill] (0,4) circle [radius=0.1];
	\node [above right] at (0 ,4) {$\langle s_{2o}\rangle $};		
	\draw  (-6,-4)--(+6,-4);
	\draw [fill] (0,-4) circle [radius=0.1];	
	\node [above left] at (-0.1 ,-4) {$\langle s_{3}\rangle $};
	
	\draw  (0,-5)--(0,+5);
	\node [below] at (0,-5) {$\langle s_{3}, s_{2o} \rangle $};
	
	\draw  (-4,-5)--(-4,+5.2);
	\draw [fill] (-4,-4) circle [radius=0.1];
	\node [below] at (-4,-5) {$\langle s_{1}, s_{3o} \rangle $};
	\node [left] at (-6,-4) {$\langle s_{3}, s_{1o} \rangle $};
	\node [above left] at (-4 ,-4) {$\langle s_{13}\rangle $};
	\draw [fill] (-4,3) circle [radius=0.1];
	\node [left] at (-6,3) {$\langle s_{2}, s_{3o} \rangle $};
	\node [below right] at (-4 ,3) {$\langle s_{3o}\rangle $};
	\draw [fill] (-4,4) circle [radius=0.1];
	
	\node [above right] at (-4 ,4) {$\langle s_{1}\rangle $};
	\node [left] at (-6,4) {$\langle s_{1}, s_{2o} \rangle $};
	
	\draw  (4,-5)--(4,+5);
	\node [below] at (4,-5) {$\langle s_{2}, s_{1o} \rangle $};
	\draw [fill] (4,-4) circle [radius=0.1];
	\node [above right] at (4 ,-4) {$\langle s_{1o}\rangle $};
	\draw [fill] (4,3) circle [radius=0.1];
	\node [below right] at (4 ,3) {$\langle s_{2}\rangle $};
	\draw [fill] (4,4) circle [radius=0.1];
	\node [above right] at (4 ,4) {$\langle s_{12}\rangle $};

	\draw [blue] (-6,2.75)--(-4,3)--(4,4)--(6,4.25);
	\node [below] at (-6,2.75) {$\langle s_{3o}, \delta \rangle $};
	
	\draw [blue] (-4.5,-5)-- (-4,-4)--(0,4)--(0.5,5);	
	\node [left] at (-4.5,-5) {$\langle s_{2o}, \delta \rangle $};

	\end{tikzpicture}  \quad  	\begin{tikzpicture}[scale = 0.5] 
	\draw  (-6,3)--(+6,3);
	\draw [fill] (0,3) circle [radius=0.1];
	\node [below right] at (0 ,3) {$\langle s_{23}\rangle $};
	\draw  (-6,4)--(+6,4);
	\draw [fill] (0,4) circle [radius=0.1];
	\node [above right] at (0 ,4) {$\langle s_{2o}\rangle $};		
	\draw  (-6,-4)--(+6,-4);
	\draw [fill] (0,-4) circle [radius=0.1];	
	\node [above left] at (-0.1 ,-4) {$\langle s_{3}\rangle $};
	
	\draw  (0,-5)--(0,+5);
	\node [below] at (0,-5) {$\langle s_{3}, s_{2o} \rangle $};
	
	\draw  (-4,-5)--(-4,+5.2);
	\draw [fill] (-4,-4) circle [radius=0.07];
	\node [below] at (-4,-5) {$\langle s_{1}, s_{3o} \rangle $};
	\node [left] at (-6,-4) {$\langle s_{3}, s_{1o} \rangle $};
	\node [above left] at (-4 ,-4) {$\langle s_{13}\rangle $};
	\draw [fill] (-4,3) circle [radius=0.1];
	\node [left] at (-6,3) {$\langle s_{2}, s_{3o} \rangle $};
	\node [below right] at (-4 ,3) {$\langle s_{3o}\rangle $};
	\draw [fill] (-4,4) circle [radius=0.1];
	
	\node [above right] at (-4 ,4) {$\langle s_{1}\rangle $};
	\node [left] at (-6,4) {$\langle s_{1}, s_{2o} \rangle $};
	
	\draw  (4,-5)--(4,+5);
	\node [below] at (4,-5) {$\langle s_{2}, s_{1o} \rangle $};
	\draw [fill] (4,-4) circle [radius=0.1];
	\node [above right] at (4 ,-4) {$\langle s_{1o}\rangle $};
	\draw [fill] (4,3) circle [radius=0.1];
	\node [below right] at (4 ,3) {$\langle s_{2}\rangle $};
	\draw [fill] (4,4) circle [radius=0.1];
	\node [above right] at (4 ,4) {$\langle s_{12}\rangle $};

	\draw [blue] (-6,2.75)--(-4,3)--(4,4)--(6,4.25);
	\node [below] at (-6,2.75) {$\langle s_{3o}, \delta \rangle $};
	\draw [blue](-1.1,4.925)-- (4,-4)--(0,3)--(-0.27,3.47)--(4.5,-4.875);	
	\node [above] at (-1.1,4.925) {$\langle s_{1o}, \delta \rangle $};
	\draw [fill] (-0.27,3.47)   circle [radius=0.1];	
	\draw [blue] (-4.5,-5)-- (-4,-4)--(0,4)--(0.5,5);	
	\node [left] at (-4.5,-5) {$\langle s_{2o}, \delta \rangle $};
	\node [above left] at (-0.33,3.4) {$\langle \delta \rangle $};	
	\end{tikzpicture} 
	\caption{Constructing  $\gamma : C_{\mc J,\KK}(D^b(D_4)) \rightarrow \AA^2 $, Step 2}\label{Figure for D4 2}
\end{figure}

 We continue with the  curve   $\langle s_{1o}, \delta\rangle $.  The corresponding line must be spanned by  $\langle s_{1o} \rangle$, $\langle s_{23} \rangle$,  on the other hand the curves $\{ \langle s_{ko}, \delta\rangle \}_{k=1}^3$, have a single intersection derived point $\langle \delta \rangle $, hence if we want to attain $\alpha, \beta$ from the already drawn picture   \ref{Figure for D4 1} the already drawn two blue lines and the line spanned by $\langle s_{1o} \rangle$, $\langle s_{23} \rangle$ must intersect in a single point. To our surprise this is really the case.    Thus we obtain the common point in Figure \ref{Figure for D4 2}. which necessarily must be $\beta(\langle \delta \rangle)$.

Continuing by similar arguments with $\{\langle  s_k,s_o \rangle : k=1,2,3 \}   \cup \{\langle  s_{123},s_k \rangle : k=1,2,3 \}$ we observe two more triples of concurrent lines as seen  in Figure \ref{Figure for D4}. This picture encodes  $\alpha$, $\beta$ from \eqref{examples for interesection intro}.  Note that in this picture ${\rm Im}(\beta)$ contains  exactly the points where more than two lines meet.

\section{Examples of $G_{\mc J,\KK}(\mc T)$ with  $\mc J = \{N\PP^l\}_{l\in \ZZ_{\geq -1}}$} 

\label{last section}

\begin{ex}[$\mc J = \{N\PP^l\}_{l\in \ZZ_{\geq -1}}$,  $\mc T = D^b(D_4)$] From \eqref{vanishings general 2} we see that the set of vertices of the digraph $G_{\mc J, \KK}(\mc T)$ is disjoint union of $C_{-1}(\mc T)$ and $C_0(\mc T)$.  The 15 elements in $C_0(\mc T)$  are listed in  \ref{set of 0 curves D4} (see also Figures \ref{Figure for C0}, \ref{Figure for D4}) whereas the 9 elements in  $C_{-1}(\mc T)$  are  listed in Corollary \ref{set of -1 curves D4} (see also Figure \ref{Figure for C-1}). It remains to determine the pairs $(a,b)$, $a,b \in C_{-1}(\mc T) \cup C_0(\mc T)$ such that $(a,b)$ is a semi-orthogonal pair. To that end we use Figures \ref{F1 for D4}, \ref{F2 for D4} (or Corollary \ref{coro for pairs}). First note that for each directed  edge  $(a,b)$ we have that   $\mc T = \langle a,b \rangle$ is a semi-orthogonal decomposition, and hence $a= b^{\perp}$, $b=^{\perp}a$. On the other hand for any vertex $a$ the left and right orthogonal subcategories $ a^{\perp}$, $^{\perp}a$ are such that $( a^{\perp}, a)$, $(a,^{\perp}a)$ are directed edges. Therefore for any vertex $a$ there is exactly  one directed edge of the form $(a,b)$ and exactly  one directed edge of the form $(c,a)$.  Taking into account that the vertices are finitely many we deduce that the digraph $G_{\mc J, \KK}(\mc T)$  consists of several oriented loops. From  Figures \ref{F1 for D4}, \ref{F2 for D4}, or Corollary \ref{coro for pairs} (d), we get full exceptional collections
	\begin{gather} (s_{1},s_{2o},s_2,s_{3o}), (s_{2},s_{3o},s_3,s_{1o}), (s_{3},s_{1o},s_1,s_{2o})  	 \\
	(s_{1},s_{3o},s_3,s_{2o}), (s_{3},s_{2o},s_2,s_{1o}), (s_{2},s_{1o},s_1,s_{3o})  
	\end{gather}
	Hence  we obtain the first two loops in Figure \ref{Figure 1 for digraph of D4} .  \begin{figure}
\begin{tabular}{c c}\begin{tabular}{c}	\begin{tikzpicture}[scale = 0.5] 
		\draw[->]  (-1,0)--(+1,0);
		\node [left] at (-1 ,0) {$\langle s_{1},s_{2o} \rangle $};
		\node [right] at (+1 ,0) {$\langle s_2,s_{3o}\rangle $};
		\draw[->]   (1.8,0.5)-- (1,1.5);
		\node at (0 ,2) {$\langle s_{3},s_{1o} \rangle $};
		\draw[->]  (-0.7,1.5) --  (-1.7,0.5);
			\end{tikzpicture}  \\
				\begin{tikzpicture}[scale = 0.5] 
				\draw[->]  (-1,0)--(+1,0);
				\node [left] at (-1 ,0) {$\langle s_{1},s_{3o} \rangle $};
				\node [right] at (+1 ,0) {$\langle s_3,s_{2o}\rangle $};
				\draw[->]   (1.8,0.5)-- (1,1.5);
				\node at (0 ,2) {$\langle  s_{2},s_{1o} \rangle  $};
				\draw[->]  (-0.7,1.5) --  (-1.7,0.5);
				\end{tikzpicture} 
		\end{tabular}	\begin{tabular}{c}
				\begin{tikzpicture}[scale = 0.7] 
				\draw[->]  (-1,0)--(0,0);
				\node [left] at (-1 ,0) {$ \langle s_{ko},\delta \rangle $};
				\node [right] at (0 ,0) {$\langle s_{ik},s_{jk}\rangle $};
				\draw[->]   (2.2,0)-- (2.8,0);
				\node [right] at (2.7 ,0) {$\langle s_{123},s_{k} \rangle $};
				\draw[->]  (5,0)--(5.5,0);
			   \node [right] at (5.4,0) {$\langle s_{i},s_{j} \rangle $};
				\draw[->]  (7.3,0)--(7.9,0);
				\node [right] at (7.8,0) {$\langle s_{k},s_{o} \rangle $};
				\draw[->]  (7.8,0.3)--(4.4,2);
				\node [left ] at (4.4,2) {$\langle s_{io},s_{jo} \rangle $};
				\draw[->]  (2.1,2)--(-1,0.3);	
				\end{tikzpicture} \\ $\abs{\{i,j,k\}}=3$ \end{tabular}	\end{tabular} 
	\caption{The graph $G_{{\mc J}, \KK}\left (D^b(D_4)\right )$}\label{Figure 1 for digraph of D4}
	\end{figure}
	Now we use again Figures \ref{F1 for D4}, \ref{F2 for D4}, or Corollary \ref{coro for pairs} (a), (c), (b),  and get the following full exceptional collections
	\begin{gather} (\delta,s_{12},s_{13},s_{23}), (s_{13},s_{23},s_{123},s_{3}), (s_{123},s_{3},s_1,s_{2}), 
	(s_{1},s_{2},s_{3o},s_{3}), (s_{3o},s_{3},s_{1o},s_{2o}), (s_{1o},s_{2o},\delta,s_{12})  \nonumber 
	\end{gather}
	Since $\langle \delta,s_{12} \rangle = \langle s_{3o},\delta \rangle$ (see \eqref{coro for pairs 1} ) and $\langle s_{3o}, s_3 \rangle = \langle s_{3}, s_o \rangle$	(see \eqref{coro for pairs 2}) we obtain a loop in  $G_{\mc J, \KK}(\mc T)$ whose vertices are $ \langle s_{3o},\delta \rangle $, $\langle s_{13},s_{23}\rangle $, $\langle s_{123},s_{3} \rangle $, $\langle s_{1},s_{2} \rangle $, $\langle s_{3},s_{o} \rangle $, $\langle s_{1o},s_{2o} \rangle $.  
	The auto-equivalence $\kappa$ 
	of $\mc T$  ( shown in Figure \ref{Figure for A1 in D4}) acts on the digraph as explained in Corollary  \ref{action on digraphs} and applying  this automorphism of $G_{\mc J, \KK}(\mc T)$ twice on  the obtained loop produces two more loops. Thus we obtain three loops with $6$ vertices  as depicted   in Figure \ref{Figure 1 for digraph of D4}, where the vertices alternate between genus zero and genus $-1$ nc curves.  Thus, the complete  graph   $G_{\mc J, \KK}(\mc T)$ consists of the five loops in Figure \ref{Figure 1 for digraph of D4}. 
	
\end{ex}

\begin{ex}[$\mc J = \{N\PP^l\}_{l\in \ZZ_{\geq -1}}$,  $\mc T = D^b(Q_2)$] From \eqref{vanishings general} we see that the set of vertices of  $G_{\mc J, \KK}(\mc T)$ is disjoint union of $C_{-1}(\mc T)$, $C_0(\mc T)$, $C_1(\mc T)$, and these are determined in Corollaries \ref{coro for Q2 genus 1}, \ref{coro for Q2 genus 0}, \ref{coro for Q2 genus -1}. To determine the pairs $(a,b)$, $a,b \in C_{-1}(\mc T) \cup C_0(\mc T) \cup C_0(\mc T)$ such that $\Hom(b,a)=0$  we use Figures \ref{cube for Q2}, \ref{four pyra 1 for Q2}, \ref{four pyra 2 for Q2}  (or Lemma  \ref{exceptional pairs in T2}). It is useful to keep in mind that $( F_+,F_-)$,  $ ( G_+,G_-)$,  $(F_+,G_+),$ $( F_-,G_-)$, $(d^m, c^m)$, $(a^m,b^{m+1})$ are orthogonal pairs. As in the previous example for any vertex $a$ there is exactly  one directed edge of the form $(a,b)$ and exactly  one  edge of the form $(c,a)$.  It follows that   $G_{\mc J, \KK}(\mc T)$  consists of oriented loops and infinite sequences. From  Figures \ref{cube for Q2}, \ref{four pyra 1 for Q2}, \ref{four pyra 2 for Q2}  we get full exceptional collections:
	\begin{gather} (c^m,c^{m+1}, G_-,F_-), (G_-,F_-,d^m,d^{m+1}), (d^m,d^{m+1},G_+,F_+), (G_+,F_+,c^m,c^{m+1})  	 \\
	(a^m,a^{m+1}, F_+,F_-), ( F_+,F_-,b^m,b^{m+1}), (b^m,b^{m+1},G_+,G_-), (G_+,G_-,a^m,a^{m+1}).  
	\end{gather}
	Hence  we obtain the two loops in Figure \ref{Figure 1 for digraph of Q2 curves}  as part of $G_{\mc J, \KK}(\mc T)$, where the vertices alternate between genus $1$  and genus $-1$ non-commutative curves.  \begin{figure}
		\begin{tikzpicture}[scale = 0.7] 
		\draw[->]  (-1,0)--(0,0);
		\node [left] at (-1 ,0) {$ \langle c^m,c^{m+1} \rangle $};
			\node [right] at (0 ,0) {$\langle G_-,F_-\rangle $};
		\draw[->]   (2.4,0)-- (2.8,0);
		\node [right] at (2.7 ,0) {$\langle  d^m,d^{m+1} \rangle $};

			\node [right ] at (0,2) {$\langle F_+,G_+ \rangle $};
			\draw[->]  (0.5,1.7)--(-1.5,0.5);	
				\draw[->]  (3,0.5)--(1.3,1.7);	
			\end{tikzpicture}  	\begin{tikzpicture}[scale = 0.7] 
			\draw[->]  (-1,0)--(0,0);
			\node [left] at (-1 ,0) {$ \langle a^m,a^{m+1} \rangle $};
			\node [right] at (0 ,0) {$\langle F_+,F_-\rangle $};
			\draw[->]   (2.4,0)-- (2.8,0);
			\node [right] at (2.7 ,0) {$\langle  b^m,b^{m+1} \rangle $};

			\node [right ] at (0,2) {$\langle G_+,G_- \rangle $};
			\draw[->]  (0.5,1.7)--(-1.5,0.5);	
			\draw[->]  (3,0.5)--(1.3,1.7);	
			\end{tikzpicture}  
	\caption{Full subgraphs in $G_{{\mc J}, \KK}\left (D^b(Q_2)\right )$}\label{Figure 1 for digraph of Q2 curves}
	\end{figure}
	 Looking at the cube in  Figures \ref{cube for Q2} one checks that $(d^m, c^m, a^m,b^{m+1})$, $( a^m,b^{m+1}, d^{m+1}, c^{m+1})$ are exceptional collections for each $m$, which amounts to the subgraph in Figure \ref{Figure 2 for digraph of Q2 curves} where the vertices are genus $-1$ curves. 
	 \begin{figure}
	  	\begin{tikzpicture} 
	  	\node [left] at (-7.5,0) {$\cdots$};
	  	\draw[->]  (-7.5,0)--(-5.5,0);	
	  	\node [above] at (-6,0) {};	
	  	\draw[->]  (-3.8,0)--(-1.8,0);
	  	\node [above] at (-2.8,0) {};	
	  	\node [left] at (-3.8,0) {$\langle d^m, c^m \rangle $};
	  	\draw[->]  (0.2,0)--(2.2,0);
	  	\node [right] at (-1.8 ,0) {$\langle a^m, b^{m+1} \rangle  $};
	  	\node [right] at (2.2,0) {$\langle d^{m+1}, c^{m+1} \rangle $};
	  	\draw[->]  (4.6,0)--(6.6,0);
	  	\node [above] at (3.2,0) {};	
	  	\node [right] at (6.6,0) {$\cdots $};	
	  	\node [above] at (0,0) {};	
	  	
	  	\end{tikzpicture}
	 		\caption{Full subgraph in $G_{{\mc J}, \KK}\left (D^b(Q_2)\right )$}\label{Figure 2 for digraph of Q2 curves}
	 \end{figure} It remains to find the full subgraph  whose vertices are the  genus $0$ curves listed in Corollary \ref{coro for Q2 genus 0}
	 Using  again Figures \ref{cube for Q2}, \ref{four pyra 1 for Q2}, \ref{four pyra 2 for Q2}   we obtain for each $m\in \ZZ$ the exceptional collections:
	 
	 	$ (c^m, G_-, a^{m},F_-), ( a^{m},F_-,d^{m+1},F_+), (d^{m+1},F_+,b^{m+2},G_+), (b^{m+2},G_+,c^{m+2},G_-)  	 
		 	$
	 which amounts to the first two   digraphs in Figure \ref{Figure 3 for digraph of Q2 curves}. 
	 The auto-equivalence $\zeta$ 
	 of $\mc T$ from Proposition \ref{autoeq sigma prop}  produces  the other two subgraphs in Figure \ref{Figure 3 for digraph of Q2 curves} from the obtained subgraphs.
	 
	 \begin{figure}
	 	\begin{tikzpicture} 
	 	\node [left] at (-3.1 ,0) {$ \cdots $};
	 	\draw[->]  (-3.1,0)--(-2.1,0);
	 	\draw[->]  (-0.7,0)--(0.3,0);
	 	\node [left] at (-1 ,0) {$  aF^m_+  $};
	 	 	\node [right] at (0.2 ,0) {$ cF^{m+1}_- $};
	 	\draw[->]   (1.7,0)-- (2.7,0);
		 	\node [right] at (2.9 ,0) {$ bG^{m+2}_- $};
	 	\draw[->]  (4.4,0)--(5.4,0);
		 	\node [right] at (5.4,0) {$ dG^{m+2}_+ $};
	 	\draw[->]  (6.8,0)--(7.8,0);
	 	\node [right] at (7.8,0) {$aF^{m+2}_+$};
	 	\draw[->]  (9.3,0)--(10.3,0);
	 	\node [right] at (10.3,0) {$ \cdots $};	
	 	\end{tikzpicture}  
	 	
	 	\begin{tikzpicture} 
	 	\node [left] at (-3.1 ,0) {$ \cdots $};
	 	\draw[->]  (-3.1,0)--(-2.1,0);
	 	\draw[->]  (-0.7,0)--(0.3,0);
	 	\node [left] at (-0.7,0) {$  aF^{m-1}_+  $};
	 	 	\node [right] at (0.2 ,0) {$ cF^{m}_- $};
	 	\draw[->]   (1.7,0)-- (2.7,0);
	 	 	\node [right] at (2.9 ,0) {$ bG^{m+1}_- $};
	 	\draw[->]  (4.4,0)--(5.4,0);
	 	 	\node [right] at (5.4,0) {$ dG^{m+1}_+ $};
	 	\draw[->]  (6.8,0)--(7.8,0);
		 	\node [right] at (7.8,0) {$aF^{m+1}_+$};
	 	\draw[->]  (9.3,0)--(10.3,0);
	 	\node [right] at (10.3,0) {$ \cdots $};	
	 	\end{tikzpicture}

	 	\begin{tikzpicture} 
	 	\node [left] at (-3.1 ,0) {$ \cdots $};
	 	\draw[->]  (-3.1,0)--(-2.1,0);
	 	\draw[->]  (-0.7,0)--(0.3,0);
	 	\node [left] at (-1 ,0) {$  dF^m_+  $};
	 	\node [right] at (0.2 ,0) {$ bG^{m+1}_+ $};
	 	\draw[->]   (1.7,0)-- (2.7,0);
	  	\node [right] at (2.9 ,0) {$ cG^{m+1}_- $};
	 	\draw[->]  (4.4,0)--(5.4,0);
	 	\node [right] at (5.4,0) {$ aF^{m+1}_- $};
	 	\draw[->]  (6.8,0)--(7.8,0);
	 	 	\node [right] at (7.8,0) {$dF^{m+2}_+$};
	 	\draw[->]  (9.3,0)--(10.3,0);
	 	\node [right] at (10.3,0) {$ \cdots $};	
	 	\end{tikzpicture}  
	 	
	 	\begin{tikzpicture} 
	 	\node [left] at (-3.1 ,0) {$ \cdots $};
	 	\draw[->]  (-3.1,0)--(-2.1,0);
	 	\draw[->]  (-0.7,0)--(0.3,0);
	 		\node [left] at (-0.7,0) {$  dF^{m-1}_+  $};
	 	 	\node [right] at (0.2 ,0) {$ bG^{m}_+ $};
	 	\draw[->]   (1.7,0)-- (2.7,0);
	 	 	\node [right] at (2.9 ,0) {$ cG^{m}_- $};
	 	\draw[->]  (4.4,0)--(5.4,0);
	 	 	\node [right] at (5.4,0) {$ aF^{m}_- $};
	 	\draw[->]  (6.8,0)--(7.8,0);
	 	 	\node [right] at (7.8,0) {$dF^{m+1}_+$};
	 	\draw[->]  (9.3,0)--(10.3,0);
	 	\node [right] at (10.3,0) {$ \cdots $};	
	 	\end{tikzpicture}  
	 	\caption{Full subgraphs in $G_{{\mc J}, \KK}\left (D^b(Q_2)\right )$}\label{Figure 3 for digraph of Q2 curves}
	 \end{figure}

	   Thus, in Figures \ref{Figure 1 for digraph of Q2 curves}, \ref{Figure 2 for digraph of Q2 curves}, \ref{Figure 3 for digraph of Q2 curves}  are shown all the components of the  digraph   $G_{\mc J, \KK}(\mc T)$. 
	
\end{ex}
\let\oldaddcontentsline\addcontentsline
\renewcommand{\addcontentsline}[3]{}

\let\addcontentsline\oldaddcontentsline

\end{document}